%% file: main.tex
\documentclass[11pt]{caltech_thesis}
\usepackage[hyphens]{url}
\usepackage{lipsum}
\usepackage{graphicx}
\usepackage{listings}
\usepackage[dvipsnames]{xcolor}

\usepackage{todonotes}

\usepackage[utf8]{inputenc}
\usepackage[T1]{fontenc}
\usepackage{newtxtext,newtxmath}
\usepackage{microtype}
\usepackage{comment}

\usepackage[
    backend=biber,natbib,
    style=alphabetic,
    backref=true
]{biblatex}

\definecolor{linkcol}{rgb}{0.8,0.4,0.1}
\definecolor{citecol}{rgb}{0.1,0.1,0.75}
\usepackage[unicode=true,
 bookmarks=true,bookmarksnumbered=false,bookmarksopen=false,
 breaklinks=false,pdfborder={0 0 1},colorlinks=true,
 linkcolor=linkcol,
 citecolor=citecol,linktoc=all]
 {hyperref}

\newenvironment{constraint}[1][]{\par\medskip\noindent \textbf{Constraint #1.}}{\medskip}

\newenvironment{optprob}[1][]{\par\medskip\noindent \textbf{Optimization Problem #1.}}{\par\bigskip}

\input{thmmacros}

\begin{document}

\title{Combinatorial and Algebraic Properties of Nonnegative Matrices}
\author{Jenish C. Mehta}

\degreeaward{PhD in Computer Science}                 
\university{California Institute of Technology}    
\address{Pasadena, California}                     
\unilogo{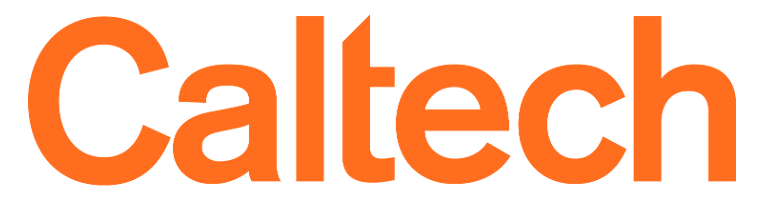}                                 
\copyyear{2022}  
\defenddate{Dec 8, 2021}          


\rightsstatement{All rights reserved
}

\maketitle[logo]

\begin{abstract}

We study the combinatorial and algebraic properties of Nonnegative
Matrices. Our results are divided into three different categories.

1. We show the first quantitative generalization of the 100 year-old
Perron-Frobenius theorem, a fundamental theorem which has been used
within diverse areas of mathematics. The Perron-Frobenius theorem
states that every irreducible nonnegative matrix $R$ has a largest
positive eigenvalue $r$, and every other eigenvalue $\lambda$ of $R$ is
such that $\text{Re}\lambda<r$ and $|\lambda|\leq r$. We capture
the notion of irreducibility through the widely studied notion of
edge expansion $\phi$ of $R$ which intuitively measures how well-connected
the underlying digraph of $R$ is, and show a quantitative relation
between the spectral gap $\Delta=1-\text{Re}\lambda/r$  (where $\lambda\not=r$
is the eigenvalue of $R$ with the largest real part) and the edge expansion $\phi$ as
follows. 
\[
\dfrac{1}{15}\cdot\dfrac{\Delta(R)}{n}\leq\phi(R)\leq\sqrt{2\cdot\Delta(R)}.
\]
This also provides a more general result than the Cheeger-Buser inequalities
since it applies to any nonnegative matrix.

2. We study constructions of specific nonsymmetric matrices (or nonreversible
Markov Chains) that have small  edge expansion  but large spectral gap, taking
us in a direction more novel and unexplored than studying symmetric matrices with
constant  edge expansion  that have been extensively studied. We first analyze
some known but less studied Markov Chains, and then provide a novel
construction of a nonreversible chain for which
\[
\phi(R)\leq\dfrac{\Delta(R)}{\sqrt{n}},
\]
obtaining a bound exponentially better than known bounds. We also
present a candidate construction of matrices for which 
\[
\phi(R)\leq2\dfrac{\Delta(R)}{n},
\]
which is the most beautiful contribution of this thesis. We believe
these matrices have properties remarkable enough to deserve study
in their own right.

3. We connect the edge expansion and spectral gap to other combinatorial
properties of nonsymmetric matrices. The most well-studied property
is mixing time, and we provide elementary proofs of known results and/or new results relating  
mixing time, edge expansion, spectral gap, and  spectral gap of the additive symmetrization. Further, we provide a unified view of
the notion of capacity and normalized capacity, and show the monotonicity
of capacity of nonreversible chains amongst other results for nonsymmetric
matrices. We finally discuss and prove interesting lemmas about different
notions of expansion and show the first results for tensor walks or
nonnegative tensors.
\end{abstract}

\begin{acknowledgements} 	 

It took me considerable time to decide on doing a PhD, and I am deeply
thankful to all the people in my life that were a source of support
in different forms and aspects after I decided to do so.

Applying to different places from India without any actual interaction
with persons at those places felt like a blackbox process with uncertain
outcomes, and I thank all the people who gave me an opportunity
to do a PhD -- Lance Fortnow (Georgia Tech), Laszlo Babai (U Chicago),
David Zuckerman (UT Austin), Anup Rao (U Washington, 2011), and Thomas
Vidick and Leonard Schulman (Caltech).

I would like to thank Thomas Vidick for being a near perfect advisor
-- energetic, hardworking, a true quantum expert, extremely sharp
to be able to slice and see through ideas, super smart, and with a
superhuman ability to do multiple tasks at the same time. He infused
life into the theory group at Caltech and went out of his way to help
students. He ensured that I always had all the resources in every
form so that nothing, however tiny, was an impediment to my research.
We had many interesting discussions about commuting Hamiltonians in
my first years at Caltech, and even after I started working on other
problems in circuit complexity, he would always check in regularly
to see that everything was fine. He also never stopped me from indulging
in extremely difficult problems although it was almost certain to
lead to very little output.

I would like to thank Chris Umans for some enjoyable conversations
on different problems in complexity theory at different points, and in a certain sense, my
way of thinking about a problem for the first few weeks is closest
to his. His thought process is extremely clear and precise, he goes into the intricacies of problems immediately, and
goes fiercely after difficult problems for years with the vigor of
a grad student, which was always motivating.

I thank Shuki Bruck for serving on my thesis committee, for being
generous with his time, and for always asking very insightful questions
during talks.

For the last, since I can write in more detail having worked with
him for the past few years, I would like to deeply thank Leonard Schulman,
for everything. He has a penchant for interesting problems irrespective
of where they come from, the acumen to see the fundamental essence
in them, and the breadth in diverse areas of mathematics to be able
to solve them. In fact, during my first visit to Caltech, we only talked about a high school project of mine on differently weighted
balls in Newton's cradle, and nothing about mainstream theory.
He has extreme kindness (not evident from his terse emails), has cultured
and refined tastes, has his own unique form of slightly dry wit that
adds charm to every conversation, and has deep wisdom that becomes
visible when he talks philosophically about any topic. He has many
powerful and original intuitions about almost every kind of random
object, and the uncommon ability to question the very
obvious and (almost always) very lucidly grasp the complicated. He asked a
question about Cheeger-type inequalities for directed graphs at one
point which I found extremely interesting, and working on it with
him for the last few years ultimately led to this thesis. He always
made me very comfortable and made every discussion so interesting
and enjoyable that I looked forward to every meeting with him, and
he always ensured that we kept learning newer things about the
problem each time. It was extremely helpful that our meetings were
almost always flexible, and could last for two or three or four or
more hours.  I
shall always remain thankful to him for making the experience of doing
research with him wonderful.


Above all, I would like to thank Thomas Vidick and Leonard Schulman
for giving me the freedom to develop my own 
sensibilities,  intuitions, and perceptions
of theory, and letting me choose problems based on my own taste and
temperament. This was invaluable and will stay with me forever.

Having mostly self-studied complexity theory before I came to Caltech,
there were many gaps in both my knowledge and understanding of concepts
in different areas, and my first few years at Caltech were like an
extended undergrad, and I would like to thank many different people
for some extraordinary courses -- Joel Tropp (Linear Algebra), Chris
Umans (Complexity Theory, arguably the best complexity course on the
planet), Alexei Kitaev (Quantum Computation), Leonard Schulman (Randomness
and Algorithms), Thomas Vidick (Hamiltonian Complexity), Gil Cohen
(Randomness Extractors), Adam Sheffer (Additive Combinatorics; Polynomial Method and Incidences), and Erik Winfree (Biomolecular Computation)
among others. I also thank Zeyu Guo for getting me interested in Additive Combinatorics. 

I spent the first few years working on circuit complexity (trying
to show exponential lower bounds on the size of constant depth $AC^{0}\cdot[\oplus]$
circuits computing the Inner Product function), and although I was
not able to obtain many meaningful results, it led to some deep insights
and understandings about the behavior of circuits and the limits of
different tools, which I hope will be helpful in future work on the
problem. During this time, I would like to thank Johan Hastad, Ben
Rossman, and Srikanth Srinivasan for some inspiring and/or fruitful
conversations. I would also like to thank Jaikumar Radhakrishnan and
Prahlad Harsha for hosting me twice at TIFR during two summers, and
Meena Mahajan and Samir Datta for hosting me for two summers during
my undergrad, and Samir Datta for formulating my first problem in
theory.

Caltech is a beautiful place, and I have spent many pleasant hours
at the turtle pond and at different places ruminating over all kinds
of thoughts, and I thank it for giving me this privilege to do so.
I also thank all the staff -- specifically Maria and Linda -- for
ensuring that nothing organizational became a concern, Linda for many warm
and nice conversations, and the Caltech International Students
Program for making all governmental matters simple. I would
also like to thank the entire Caltech West Coast Swing group for teaching
a stone like me how to move.

I would like to thank many friends, both new and old, and although an exhaustive list is infeasible, I will mention a few close and/or easily recollectable names, from grad school
-- Andrea (closest grad friend), Jing, Zofii, Tal, Palma, Swagato,
Florian, Fabien, Sudhi, Rohit, Ben Lee, Anand, Saeed and Pranav;
 from undergrad -- Anant (always takes me on adventures), Praneeth, and everyone still in touch; from middle/high school -- Miren (like a brother by now), Dhruman,
Saagar, and Preet; and towards the end of my PhD, Corinne, for being the wonderful being that she is. 

Finally, I would like to thank my family for their unbridled support.
I would like to thank my father, whose raw, untarnished, undifferentiated,
authentic, and almost sublime intellect always made me feel like a mere trained mind that was moving around symbols
(or symbolically expressible thoughts) to produce desired strings,
and he was thus a constant inspiration during my time at Caltech,
specifically when problems (circuit lower bounds) were not getting
solved. The manner in which he has supported everything I do goes
beyond words in both count and significance, and it would not have
been possible to start on the PhD without his rock-solid presence
and support. Unfortunately, he passed away in 2021 due to Covid, just a few months
before I defended. We shared a deep and inexpressible connection,
and he called and talked to me every single day (although he would
have a thousand things to do) through which he unknowingly, unintentionally,
and indirectly provided the catalyst for me to keep going on. I would
like to thank my mother for her deep care, her unconditional and sacrificial
love, her ethereal and calming presence, and the way she has always smiled and stood behind me. She has 
an empathetic understanding of the world and  she interacts with all forms
of reality very delicately and genuinely, and consequently has deeply
affected my own aesthetic sense. I would like to thank my sister Rutvi for
her joyful exuberance and pure innocence, and our mutual craziness
has helped me protect a little bit of my own. I thank my cousin Rishabh for our camaraderie and
for always being with me in everything, and Dhwani, Chintan, Dhruvi, and the
entire extended family.

I would like to conclude by thanking everyone whom I might have missed,
and all the direct and indirect correlations that I'm unable to recall
at present.

\end{acknowledgements}

\vspace*{5cm}

\begin{center}
For my parents,
\par\end{center}
\begin{center}
Chetan and Bina Mehta
\par\end{center}
\begin{center}
(and whoever they might dedicate it to,
\par\end{center}
\begin{center}
ad infinitum)
\par\end{center}

\pagebreak{}

\tableofcontents

\mainmatter

\chapter{Introduction}

\vspace{1.5cm}

\begin{quote}
That the sun will not rise tomorrow is no less intelligible a proposition,
and implies no more contradiction than the affirmation, that it will
rise.
\begin{flushright}
\textasciitilde{} David Hume, \emph{An Enquiry Concerning Human Understanding}\\
\par\end{flushright}

\end{quote}

\vspace{3cm}

\section{On chance and its evolution}

A primary idea in the history of human thought is that of chance,
and is further made indispensable due to the functionality accorded
by its \emph{measurability}. This measurability of chance is distinct
from its \emph{meaning} -- the latter being its different interpretations,
and amongst the two notable ones are the Bayesian perspective --
purely for rational decision making -- and the Frequentist perspective
-- a much stronger and less justified statement about the nature
of sequences of events -- and there are multiple other interpretations
none of which are satisfactory or comprehensive enough to provide
a rigorous and justified universal meaning of chance. From a sufficiently
critical perspective, we find a certain form of conditioned non-determinism
to be a more fundamental and universal concept for understanding observed
events than chance. However, although deep and arguably more primary
and important than anything that we will write further on, our concern
in this thesis will not be the meaning of chance (and not even its
measurability), but primarily its \emph{evolution}.

To proceed further, we assume we have a possibility of $n$ states,
where $n$ may be finite or infinite ($\aleph_{0}$, $\aleph_{1}$,
or more) such that there is an associated chance $p$ over the collection
of states. Note that the different states could correspond to states
of any system or universe in general. Let the chance $p$ evolve to
$p'$ at a different point in time (we have an implicit notion of
\emph{distinct} events which produces an implicit notion of time or
steps for us). To understand the evolution of chance more mathematically,
we shall assume that $p'=Ap$, and the operator $A$ which takes $p$
to $p'$ will be our primary object of study.

\section{Markov Chains and the Perron-Frobenius theorem}

It is possibly to consider a large number of properties of $A$ (and
$p$), but to understand the operator $A$ in a mathematically meaningful
manner, we make various simplifying assumptions about $A$ and $p$.
Our first assumption is that we will consider only a finite state
space in this thesis, always referred to by $n$. We will assume that
$A$ is \emph{independent of $p$} and \emph{finite-dimensional}.
Thus $A$ could be represented as a finite collection of constants.
Throughout the thesis (except the final section on tensors), we further
assume that $A$ is \emph{linear}, and thus can be expressed as an
$n\times n$ matrix. At this point, the study becomes restricted to
general linear operators, and we make the final assumption that $A$
is \emph{entry-wise nonnegative}, to say something stronger about
$A$ than general linear operators.

Such an $A$ is referred to as a Markov chain, first studied by Andrey
Markov in 1906, motivated by a disagreement with Pavel Nekrasov who
had claimed that independence was necessary for the weak law of large
numbers to hold \cite{markov1906rasprostranenie,seneta1996markov}.
Markov chains have an illustrious history and have been thoroughly
studied over the past 100 years, and have been applied to diverse
areas within mathematics.

Arguably the most fundamental theorem about Markov chains is the Perron-Frobenius
theorem. For simplicity, assume that $A$ is entry-wise positive.
Then Perron's theorem \cite{perron1907theorie} says that there is
a positive number $r$ which is an eigenvalue of $A$, all other eigenvalues
$\lambda\not=r$ of $A$ are such that $|\lambda|<r$, and the corresponding
left and right eigenvectors corresponding to $r$ are entry-wise positive.
This theorem was extended to all irreducible $A$ by Frobenius \cite{frobenius1912matrizen},
where $A$ is irreducible if the underlying graph with non-negative edge weights
$A_{i,j}$ is strongly-connected, i.e. there is a directed path of non-zero
weights between every pair of vertices. In the irreducible case, there
is one specific difference -- for all other eigenvalues $\lambda\not=r$
of $A$, it is the case that $|\lambda|\leq r$ and  $\text{Re}\lambda<r$. The resulting
theorem is referred to as the Perron-Frobenius theorem.

The Perron-Frobenius theorem is fundamental to many areas within mathematics,
and has been extensively used in a large number of applications (see
\cite{maccluer2000many}, \cite{berman1994nonnegative} for instance),
such as Markov Chains \cite{LevinPW09} (traditionally, Markov chains
have been nonnegative row stochastic matrices in order to preserve the
probability simplex); theory of Dynamical Systems \cite{dynsys-katok1997introduction};
a large number of results in economics such as Okishio's theorem \cite{okishio-bowles1981technical},
Hawkins Simon condition \cite{hawkins1949note}, Leontiev Input/Output
Economic Model \cite{leontief1986input}, and Walrasian Stability
of Competitive Markets \cite{econ-takayama1985mathematical}; Leslie
population age distribution model in Demography \cite{leslie1945use,leslie1948some};
DeGroot learning process in social networks \cite{degroot1974reaching,degroot-french1956formal,degroot-harary1959criterion};
PageRank and internet search engines \cite{page1999pagerank,pagerank-langville2011google};
Thurston\textquoteright s classification of surface diffeomorphisms
in Low-Dimensional Topology \cite{thurston1988geometry}; Kermack--McKendrick
threshold in Epidemiology \cite{epid-kermack1991contributions,epid-kermack1932contributions,epid-kermack1933contributions};
Statistical Mechanics (specifically partition functions) \cite{statmech-tolman1979principles};
the Stein--Rosenberg theorem and Seidel versus the Jacobi iterative
methods for solving linear equations in Matrix Iterative Analysis
\cite{iteranaly-varga1962iterative}; and see the comprehensive \cite{berman1994nonnegative}
for more.

However, already at this juncture, there is an interesting question that will be the main focus of the first part of this thesis. The
Perron-Frobenius theorem provides a \emph{qualitative} result about
the irreducibility and eigenvalues of the matrix, i.e., it says that
if $A$ is irreducible, then every eigenvalue $\lambda$ of $A$ except
$r$ is such that $\text{Re}\lambda<r$. Our main question is:
what is the \emph{quantitative} version of this result? Assume that
$\lambda$ is some eigenvalue except $r$ such that it has the maximum
real part. Then specifically, we can ask, 
\noindent \begin{center}
\textbf{Question 1.} Is there a quantitative relation between some
measure of irreducibility of $A$ and the gap $r-\text{Re}\lambda$?
\\
\par\end{center}

This question will be our first primary consideration. In fact, this
question has been completely answered in the case that $A$ is symmetric.
The result holds even for reversible $A$ where $A$ is reversible
if $D_{u}AD_{v}=D_{v}A^{T}D_{u}$ where $u$ and $v$ are the left and right eigenvectors of $R$ for eigenvalue $r$, but syntactically (and intuitively)
both reversibility and symmetricity are exactly the same (see Lemma \ref{lem:transRtoA}). To understand
the symmetric version of this question, we will require some more
definitions.

\section{The Cheeger-Buser Inequalities}

We start by defining the spectral gap and edge expansion of symmetric
matrices $A$. For the sake of simplicity of exposition, assume that
$A$ is irreducible and the largest eigenvalue of $A$ is 1, and the
corresponding positive eigenvector is $w$. The first notion is that
of the \emph{spectral gap}, which is again simple to define for symmetric
matrices. Since $A$ is symmetric, let the real eigenvalues of $A$
be $1\geq\lambda_{2}(A)\geq...\geq\lambda_{n}(A)$. Define the spectral
gap of $A$ as 
\[
\Delta(A)=1-\lambda_{2}(A).
\]
The second notion is the edge expansion (or expansion) of $A$ which helps to 
quantify the notion of irreducibility of $A$, and is defined as
\[
\phi(A)=\min_{S\subset[n]:\sum_{i\in S}w_{i}^{2}\leq\frac{1}{2}\sum_{i}w_{i}^{2}}\dfrac{\sum_{i\in S,j\in\overline{S}}A_{i,j}w_{i}w_{j}}{\sum_{i\in S,j}A_{i,j}w_{i}w_{j}}.
\]
It might appear notationally cumbersome, but it is intuitively simple
to understand by looking at the equivalent definition for a symmetric
doubly stochastic matrix $A$, in which case $w$ is just the all
ones vector. In that case, we would have 
\[
\phi(A)=\min_{S:|S|\leq n/2}\dfrac{\sum_{i\in S,j\in\overline{S}}A_{i,j}}{|S|}.
\]
Thus, given a set $S$, the  edge expansion  of the set is exactly the total
average weight of edges that leaves the set $S$, i.e. goes from $S$
to $\overline{S}$. The same intuition holds for general
$A$, except that the edges are re-weighted by the eigenvector $w$
(the reasons for which will become clear later). This notion of  edge expansion 
is a simple and sufficiently general quantitative notion of the irreducibility
of $A$. Note that if $\phi(A)=0$, then there is some set $S$ such
that the total mass of edges leaving the set is 0, implying that there
is no way to reach a vertex in $\overline{S}$ from the set $S$,
implying that the matrix is disconnected and thus has 0 irreducibility
or is reducible. Similarly, if $\phi(A)\approx0$, it means that there
are two \emph{almost} disconnected sets in $A$ and thus $A$ is badly
connected or has high reducibility. Similarly, if $\phi(A)\geq c$ (where $c$ any constant independent of the size of $A$),
it implies that for every set $S$, there is a constant fraction of
the mass of edges leaving $S$. This means that a random walk starting at a uniformly random vertex in any set $S$ is very likely
to \emph{leave} $S$, and it will not be restrained to $S$, implying that $A$ is very
well-connected, or has high irreducibility. Thus the  edge expansion  $\phi$
captures the notion of irreducibility neatly. Further, we remark that
given $A$, it is NP-hard to decide if the  edge expansion  of $A$ is greater
than some given input number \cite{garey1974some}.

Given these two definitions, we are ready to answer Question 1 for
symmetric matrices by relating $\phi(A)$ with $\Delta(A)$. The upper
bound on $\phi(A)$ in terms of $\Delta(A)$ was obtained by Cheeger
\cite{cheeger70}, albeit over Reimannian manifolds, and it was translated
to symmetric matrices by \cite{AlonM85,SinclairJ89,dodziuk84,Nilli91}.
The lower bound on $\phi(A)$ in terms of $\Delta(A)$ was obtained
by Buser \cite{Buser82} for Reimannian manifolds, and the translation
to symmetric matrices is a direct consequence of the variational characterization
of eigenvalues for symmetric matrices. Together, the following bounds
are referred to as the Cheeger-Buser inequalities:
\begin{equation}
\frac{1}{2}\Delta(A)\leq\phi(A)\leq\sqrt{2\cdot\Delta(A)}.\label{eq:intro-cheeg}
\end{equation}
The bounds are tight, the lower bound exactly holding for the hypercube,
and the upper bound for the undirected cycle up to constants. Note
that this helps us understand the connectivity or irreducibility of
$A$ in terms of the spectral gap, exactly as we had hoped, and provides
a quantitative version of the Perron-Frobenius theorem for symmetric
(or reversible) matrices. The Cheeger-Buser inequalities have been
used in many different results in clustering and expander graph construction
amongst others \cite{tolliver2006graph,arora2009expander,hoory2006expander,spielman1996spectral},
and higher-order Cheeger inequalities have also been shown relatively
recently \cite{lee2014multiway}.

The general problem for any nonnegative matrix (not necessarily symmetric/reversible) had been open, and that will indeed lead to our
main result.

\section{Quantitative generalization of the Perron-Frobenius theorem}

We first redefine the edge expansion and spectral gap of a general irreducible
nonnegative matrix $R$, that is not necessarily symmetric. We will
use $R$ for matrices whose left and right eigenvectors for the largest
eigenvalue are not the same and $A$ otherwise. Again for the sake
of simplicity of exposition, assume that $R$ has largest eigenvalue
1, and let the eigenvalues of $R$ (roots of the equation $\text{det}(\lambda I-R)=0$)
be arranged so that $1\geq\text{Re}\lambda_{2}(R)\geq...\geq\text{Re}\lambda_{n}(R).$
Then define the spectral gap of $R$ as 
\[
\Delta(R)=1-\text{Re}\lambda_{2}(R).
\]
Unlike the eigenvalues which are basis-independent, the  edge expansion 
depends on the explicit entries of the matrix (for it to be meaningful in terms of probability),
and we will assume thus that $R$ is nonnegative, and also irreducible (for purely technical reasons),
and the positive (from the Perron-Frobenius theorem) left and right
eigenvectors for eigenvalue 1 of $R$ are $u$ and $v$ respectively.
Then define the edge expansion of $R$ as 
\[
\phi(R)=\min_{S\subset[n]:\sum_{i\in S}u_{i}v_{i}\leq\frac{1}{2}\sum_{i}u_{i}v_{i}}\dfrac{\sum_{i\in S,j\in\overline{S}}R_{i,j}u_{i}v_{j}}{\sum_{i\in S,j\in [n]}R_{i,j}u_{i}v_{j}}.
\]
We extend these definitions to reducible $R$, by taking limit infimums
within a ball of irreducible matrices around $R$, and the formal
definition of $\phi(R)$ is given in \ref{def:(Edge-expansion-doubstoc}
and \ref{def:gen-edge-expansion}, and that of $\Delta(R)$ is given
in \ref{def:(Spectral-gap-irred} and \ref{def:(gen-Spectral-gap-and-edgeexp}.
With these definitions, we can present the main theorem of this thesis here.

\begin{theorem}
\label{thm:intro-phi-delta} Let $R$ be an $n\times n$ nonnegative
matrix, with edge expansion $\phi(R)$ defined as \ref{def:gen-edge-expansion}
and \ref{def:(gen-Spectral-gap-and-edgeexp}, and the spectral gap $\Delta(R)$
defined as \ref{def:(Spectral-gap-irred} and \ref{def:(gen-Spectral-gap-and-edgeexp}.
Then 
\[
\dfrac{1}{15}\cdot\dfrac{\Delta(R)}{n}\leq\phi(R)\leq\sqrt{2\cdot\Delta(R)}.
\]
\end{theorem}

We remark that this exactly helps to answer Question 1, and provides
a generalization of the Perron-Frobenius theorem. It is more general
than the Cheeger-Buser inequalities since it holds for all nonnegative
matrices but is not a generalization of it since the lower bound becomes
weaker for the case of symmetric matrices.

We make some comments about Theorem \ref{thm:intro-phi-delta}. The
first is that the upper bound is exactly same as in the case of symmetric
matrices. The upper bound was shown for doubly stochastic matrices
$A$ by Fiedler, and it is straightforward to extend it to all $R$.
In fact, given the Cheeger inequality, it is simple to use the inequality
for showing the upper bound on $\phi(R)$ in Theorem \ref{thm:intro-phi-delta}.
However, unlike the symmetric case, obtaining a lower bound on $\phi(R)$
or the equivalent of the Buser inequality is much more difficult.
The primary source of the difficulty is that there is no Courant-Fisher
variational characterization of eigenvalues for nonsymmetric matrices.
In fact, this characterization for symmetric matrices is at the heart
of most theorems in spectral theory about symmetric matrices, and
the latter are completely understood by it. However, different ideas
are required to deal with the nonsymmetric case.

The lower bound on $\phi(R)$ in Theorem \ref{thm:intro-phi-delta}
is proven using a sequence of lemmas each of which is independently usable in different contexts. The main idea is to show that
$\phi(R^{k})\leq k\cdot\phi(R)$ (for all integers $k$), and use the Schur decomposition
of $R$ to to relate $\phi(R^{k})$ to the spectral gap of $R$. A
detailed exposition on this is given in Section \ref{subsec:Lower-bound-on-phi-actual}
before the theorem is proved. The main difference in Theorem \ref{thm:intro-phi-delta}
from the Cheeger-Buser inequalities is the factor of $n$, which makes
the bound extremely weak. However, the strength of the lower bound
on $\phi(R)$ in the theorem is in the factors that are \emph{not} present.

The first straightforward thing to note is that the dependence on $n$ is linear and not exponential. This
is crucial, since our initial attempt to prove this result was by
using tools from perturbation theory -- since $\phi(R)$ can be understood
as a perturbation to a matrix of disconnected components. It is possible
to understand the change in the coefficients of the characteristic
polynomial from the change in the matrix, and further relate it to
the change in the spectral gap. However, this brings in factors exponentially
depending on $n$. Another direction is to use results similar to
Baur-Fike \cite{bauer1960norms} (see Lemma \ref{lem:(Baur-Fike-}),
but it is very limited since it holds only for diagonalizable $R$,
and depends on the eigenvector condition number. The eigenvector condition
number is $\mathcal{K}(R)=\|V\|_{2}\|V^{-1}\|_{2}$ where $V$ is an
invertible matrix that transforms $R$ to its Jordan decomposition
$R=VJV^{-1}$. Since the eigenvector condition number can be enormous
even for small matrices, the dependence on $\mathcal{K}(R)$ makes
the Baur-Fike result much weaker than what we seek. There is an extension
of Baur-Fike to non-diagonalizable matrices by Saad \cite{Saad11},
but it has an exponential dependence on $n$ and a dependence on $\mathcal{K}(R)$.
Note, however, that our result in Theorem \ref{thm:intro-phi-delta}
does not depend on $\mathcal{K}(R)$. We also reprove Baur-Fike using
Schur decomposition, which removes $\mathcal{K}(R)$, but instead
has a dependence on $\exp(n)$ and $\sigma(R)=\|R\|_{2}$, where $\sigma(R)$ is the
largest singular value of $R$ (See Lemma \ref{lem:(Baur-Fike-extension)-Let}).
There are two other factors that could appear in our proof of Theorem \ref{thm:intro-phi-delta}, and were
indeed present in the initial version of our proofs but were sequentially
removed in subsequent iterations -- the largest singular value $\sigma(R)$
of $R$, and the eigenvalue condition number $\kappa(R)=\min_{i}u_{i}\cdot v_{i}$
where $u$ and $v$ are the positive eigenvectors corresponding to
the largest eigenvalue 1 of $R$, and normalized so that $\sum_{i}u_{i}v_{i}=1$.
The removal of $\sigma$ is possible by a neat trick (Lemma \ref{lem:transRtoA}),
but the removal of $\kappa(R)$ is nontrivial. In fact, in the published
version of our result \cite{mehts2019edge}, we have a weaker bound
of 
\[
\dfrac{1}{15}\cdot\dfrac{\Delta(R)}{n+\ln\left(\frac{1}{\kappa(R)}\right)}\leq\phi(R).
\]
In the attempt to construct matrices showing the necessity of $\kappa(R)$
in the lower bound, it turned out that for all possible "extreme cases" of
matrices, $\kappa(R)$ was indeed unnecessary. This made us reconsider
all the lemmas in our proof, and indeed, it was possible to rewrite
them to remove the dependence on $\kappa(R)$ entirely and arrive
at Theorem \ref{thm:intro-phi-delta}.

We would also like to comment on the regime where the lower bound
on $\phi(R)$ in Theorem \ref{thm:intro-phi-delta} is useful. Due
to the factor of $n$, the bound is most meaningful when $\phi(R)\ll1/n$.
As such, the utility of the bound is perhaps in treating it as a perturbation
result, where tiny perturbations -- say $O(1/n^{2})$ or $O(1/n\log n)$
or $O(1/n\log^{2}(n))$ -- closely affect the spectral gap, and the
gap between the eigenvalues and the positive eigenvalue in the perturbed
matrix can be bounded by the resulting perturbation expressed in $\phi$ (see Section \ref{sec:limits-perturb}).

Thus, the lack of dependence on $\exp(n)$ (since we do not use perturbation
theory), $\mathcal{K}(R)$ (since we do not use the Jordan decomposition
but the Schur decomposition), $\sigma(R)$ (since we transform $R$
to an $A$ that preserves the eigenvalues and edge expansion and has equal
maximum singular and eigenvalues), and $\kappa(R)$ (since we ensure
all the intermediate lemmas are tight and crucially that the right eigenvector $v$ for eigenvalue 1 is in the kernel of $I-R$) is where the primary importance
of the lower bound on $\phi(R)$ lies, since it helps to answer the
first question and exactly quantify the Perron-Frobenius theorem.
It means $\Delta(R)$ is proportional to $\phi(R)$ up to the scaling
factor of $n$. However, since there is no dependence on $n$ in the
symmetric case (inequalities \ref{eq:intro-cheeg}), it raises the
second main question -- is the dependence on $n$ in Theorem \ref{thm:intro-phi-delta}
necessary, or is it more an artefact of our proof method? More specifically,
\begin{center}
\textbf{Question 2.} Are there matrices $R$ for which $\phi(R)\approx\dfrac{\Delta(R)}{n}$?\\
\par\end{center}

\section{Constructions of nonreversible matrices}

The attempt to Question 2 takes us in novel and sparsely explored areas.
To answer it in the strictest sense, we will now seek \emph{doubly stochastic
matrices} that achieve the bound required in Question 2. The main reason
is that they have the uniform vector (which has no variance) as both the left and right eigenvector
for eigenvalue 1, and it would seem possible on the outset that
at least for doubly stochastic matrices, the lower bound on $\phi$
in Theorem \ref{thm:intro-phi-delta} is possibly much weaker. To understand this, we will start by studying the edge expansion of \emph{specific}
irreversible or non-symmetric matrices carefully.

The edge expansion of symmetric matrices has been amply studied, and the
matrices with constant edge expansion are called (combinatorial) expanders.
These matrices have remarkable properties, making them a fundamental
building block within combinatorics, and matrices with constant  edge expansion 
and few non zero entries (edges) are complicated (in construction
or proof) and have a rich history, since they have ``magical'' properties
and help to construct seemingly impossible objects or algorithms (see \cite{hoory2006expander}).
Question 2 however is irrelevant in the symmetric case since the  edge expansion 
and spectral gap are tightly related, and we seek a gap of a factor
of $n$ between the two quantities.

The first steps in this direction were already taken in the work
of Maria Klawe \cite{klawe1981non,klawe1984limitations} who showed
that certain affine-linear constructions of matrices have inverse
polylog$(n)$  edge expansion (see Construction \ref{constr:klawe-vaz}). It was observed by Umesh Vazirani \cite{UVaz-pers},
that it is possible to orient the edges in the construction of Klawe
to obtain doubly stochastic matrices with constant spectral gap but
expansion that was $1/\text{polylog}(n)$! This is quite remarkable,
since it shows that there are doubly stochastic matrices with 
\[
\phi(A_{KV})\approx\dfrac{\Delta(A_{KV})}{(\log n)^{c}}
\]
for some constant $c$. In fact, there are many other affine-linear
constructions, all achieving similar bounds. A construction to note
in particular is that of de Bruijn matrices (described and analyzed
in detail in Section \ref{subsec:debruijn}). For these matrices, it
turns out that 
\[
\phi(A_{dB})\approx\dfrac{\Delta(A_{KV})}{\log n}.
\]
These are beautiful matrices, and their properties are listed in Lemma \ref{lem:debruijn}. These constructions indicate that a dependence
on some function of $n$ cannot be avoided. However, Question 2 can
be rephrased as follows.

\noindent \begin{center}
\textbf{Question 3.} Are there doubly stochastic matrices $A$ with 
\[
\dfrac{\phi(A)}{\Delta(A)}\in o\left(\dfrac{1}{\log n}\right)
\]
 or is it the case that for all doubly stochastic matrices, for some constant
$c$, 
\[
c\cdot\dfrac{\Delta(A)}{\log(n)}\leq\phi(A)?
\]
\\
\par\end{center}

This is the primary question for consideration, since if it is true
that the edge expansion is at least the spectral gap by $\log(n)$, it
would mean that our lower bound on $\phi(R)$ in Theorem \ref{thm:intro-phi-delta}
is exponentially far from the truth (at least for doubly stochastic matrices), and newer techniques would be
required for a stronger bound. Further, it would also mean that the
spectral gap provides a good estimate of the edge expansion which in itself
is difficult to compute since the $\log(n)$ factors are essentially
negligible in most applications. However, we show the following striking
construction.

\begin{theorem}
\label{thm:intro-rootn}There is a family of doubly stochastic matrices
$A_{n}$, called Rootn matrices, such that for every $n$,
\[
\phi(A_{n})\leq\dfrac{\Delta(A_{n})}{\sqrt{n}}.
\]
\end{theorem}

This construction of Rootn matrices is presented in Section \ref{subsec:Rootn-matrices-=002013}
and arrived at, presented, and discussed in detail in Sections \ref{subsec:Constraints-on-the},
\ref{subsec:Special-Cases-of}, \ref{subsec:The-general-case}, \ref{subsec:Rootn-matrices-=002013}.
This construction is exponentially better than known constructions,
and shows that the factor of $n$ in Theorem \ref{thm:intro-phi-delta}
is indeed closer to the truth. The construction for Rootn matrices
is arrived at by carefully fixing the eigenvalues in the triangular
matrix in the schur decomposition to create a constant spectral gap,
and choosing a unitary such that the resulting matrix is doubly stochastic
and has minimum  edge expansion.

As such, the final question that remains is the gap between $\sqrt{n}$
and $n$. To resolve this problem, we take a direction different from
the one that was used to construct Rootn matrices. Observing Rootn
matrices, we learn a possible structure of the matrix that helps to fix the  edge expansion,
and then set the entries in the matrix very carefully to maximize
the spectral gap, as discussed in Section \ref{subsec:Observations-from-the}.
This leads us to the most beautiful contribution of this thesis --
Chet Matrices -- that help to show the following. 

\begin{theorem}
\label{thm:intro-chet}There is a family of matrices $C_{n}$ called Chet matrices, for
which $\sum_{i}C_{n}(i,j)=1$ and $\sum_{j}C_{n}(i,j)=1$ for all
$i,j$, and if  $\phi(C_{n})$
is defined as for doubly stochastic matrices,
\[
\phi(C_n)=\min_{S,|S|\leq n/2}\frac{\sum\limits _{i\in S,j\in\overline{S}}C_{n}(i,j)}{|S|},
\]
then for all $n$, 
\[
\phi(C_{n})\leq 2\cdot\dfrac{\Delta(C_{n})}{n}.
\]
\end{theorem}

The construction is presented in Section \ref{subsec:Chet-Matrices-=002013}.
However, we are unable to show that Chet Matrices are nonnegative for all $n$.
It seems to be the case for $n$ up till 500 that we tested numerically (see Appendix \ref{subsec:Chet-Matrix-for-500}),
but we do not have a proof for all $n$. In spite of this, these matrices
have many remarkable properties listed in Section \ref{subsec:Chet-Observations-and-Properties},
and proving their nonnegativity is one of the main open problems of
this thesis. To show their nonnegativity, we phrase a sequence of
Trace conjectures in Section \ref{subsec:Conjectures-chet}, which
if true would (almost) imply the nonnegativity of Chet Matrices. We
believe these conjectures are interesting in their own right.

The results so far give a complete answer to the first question we
had asked initially, since we have tight bounds and (almost) matching
constructions. We now proceed to other combinatorial properties and
relate them to the edge expansion and spectral gap.

\section{Relations with other combinatorial properties}

\subsection{Mixing Time}

The most widely studied property of Markov chains that is indispensable
to algorithms is mixing time. To simplify exposition, consider
an irreducible doubly stochastic matrix $A$ such that it is $\frac{1}{2}$-lazy,
i.e., for every $i$, we have that $A_{i,i}\geq\frac{1}{2}$, which
is required purely for a technical reason (ergodicity), and any other constant
would also be sufficient. This condition ensures that except the eigenvalue
1, all other eigenvalues $\lambda$ of $A$ are such that $|\lambda|<1$.
Thus, if we write the Jordan decomposition of $A=VJV^{-1}$, it is
clear that for large enough $t$, $A^{t}\approx\frac{1}{n}\mathbf{1}\cdot\mathbf{1}^{T}$
(where $\mathbf{1}$ is the all ones vector), since all other eigenvalues
in $A$ will become close to zero and only the space corresponding
to eigenvalue 1 remains. Note that this also means, that if we start
with some probability distribution $p$, for large enough $t$, $A^{t}p\approx_{\epsilon}\frac{1}{n}\mathbf{1}$ (in say, $\ell_1$ norm).
The smallest such $t$ which works for \emph{any }starting $p$ is
referred to as the mixing time of $A$ or $\tau_{\epsilon}$. Note that similarly, we can
formally define the mixing time for any irreducible and $\frac{1}{2}$-lazy
nonnegative matrix $R$ with largest eigenvalue 1 and corresponding
left and right eigenvectors $u$ and $v$, such that $\langle u,v\rangle=1$.
This will be the smallest $t$ such that $R^{t}\approx_{\epsilon}v\cdot u^{T}$.

There is extensive literature on mixing time, different methods to
bound it and different tools that can be used for specific types of
chains. Our aim will mostly be to obtain general bounds on  mixing
time and relate it to other combinatorial and algebraic quantities.

The mixing time of symmetric (or reversible) $A$ is comprehensively studied and sufficiently well understood 
(see \cite{aldous2002reversible,LevinPW09,MonT06}). 
Due to the spectral decomposition of $A$, the mixing time $\tau_{\epsilon}$ for a reversible chain is
bounded as 
\[
c_{1}\cdot\dfrac{1}{\Delta(A)}\leq\tau_{\epsilon}(A)\leq c_{2}\cdot\dfrac{\ln(\frac{n}{\epsilon\cdot\kappa(A)})}{\Delta(A)}
\]
for some constants $c_{1}$ and $c_{2}$. Thus, the mixing time is
approximately the inverse of the spectral gap, up to the factor of $\log(n)$
and the eigenvalue condition number $\kappa(A)$. The factor $\kappa(A)$
will in fact appear in all bounds related to the mixing time, and
it is easy to see that it cannot be removed by taking even a $3\times3$
matrix which is very close to being reducible.

Regarding general chains that are not necessarily symmetric, again
many different results are known for specific chains, and the most
general is the result of Mihail \cite{mihail1989conductance}, which
directly relates the mixing time to edge expansion (Mihail shows it
for row stochastic $R$, but it directly extends to all $R$ by syntactic
changes).
\[
\tau_{\epsilon}(R)\leq c\cdot\dfrac{\ln\left(\frac{n}{\epsilon\cdot\kappa(R)}\right)}{\phi^{2}(R)}.
\]
This shows that the mixing time is inversely proportional to the  edge expansion 
of $R$. Our first result is essentially a one-line proof of the above
result, that can be derived from one of our main lemmas relating the
mixing time of a matrix $A$ (derived from $R$) which has identical
left and right eigenvector $w$ for eigenvalue 1 to the second largest \emph{singular
value} of $A$. We also show a lower bound on $\tau_{\epsilon}(R)$
in terms of $\phi(R)$, which is also known (see \cite{LevinPW09}),
but we write the bound for general $\epsilon$.

We also relate the mixing time to the spectral gap in the irreversible
case, and get the following theorem.

\begin{theorem}
\label{thm:intro-mixt}Let $\tau_{\epsilon}(R)$ be the mixing time
of an irreducible $\frac{1}{2}$-lazy nonnegative matrix $R$. Then
\[
\tau_{\epsilon}(R)\leq c\cdot\dfrac{n+\ln\left(\frac{1}{\epsilon\cdot\kappa(R)}\right)}{\Delta(R)},
\]
and our constructions of Rootn and Chet Matrices in Sections \ref{subsec:Rootn-matrices-=002013}
and \ref{subsec:Chet-Matrices-=002013} show that the factor of $n$
is again necessary.
\end{theorem}

The interesting thing in the result \ref{thm:intro-mixt} is that
$n$ is additive and not multiplicative in the numerator. The next
thing we study is the relation between the mixing time of $A$ (irreducible,
nonnegative, $\frac{1}{2}$-lazy with left and right eigenvector $w$
for eigenvalue 1) and the mixing time of $\tilde{A}=\frac{1}{2}(A+A^{T})$,
and show two-sided and tight bounds between them. We also relate the
mixing time of $A$ to that of $\exp(t\cdot(I-A))$.

We mostly present simple/elementary proofs of many known bounds on Mixing time, and present some new results in Section \ref{subsec:Mixing-Time}. Our main contribution is to show how all the (optimal) bounds are achievable using only a few key tools
and lemmas, and it also helps us obtain many bounds in a form that we have not come across. 

\subsection{Capacity and Normalized Capacity}

The next quantity we explore has been studied extensively for symmetric
matrices in different communities under different names at different
points of time, which we refer to as \emph{capacity}. The term and
the related questions arise from the study of harmonic functions and the Dirichlet
problem. Again to present a simplified definition, assume we have
an irreducible nonnegative matrix $A$ with largest eigenvalue $1$
and the corresponding left and right eigenvector $w$. Let the Laplacian
of $A$ be $L=I-A$, let $U$ be some subset of the vertices, and
let $a\in\mathbb{R}^{|U|}$ be some real vector over $U$. We want
to find a vector $q\in\mathbb{R}^{n}$ such that $q_{i}=a_{i}$ for
$i\in U$, i.e. $q$ is same as $a$ on $U$, and $(Lq)_{i}=0$ for
$i\in\overline{U}$. Since $A$ is irreducible, it is not difficult
to see that there is a unique vector $q$ satisfying the equations.
Given such a $q$, the capacity is defined as follows.
\[
\text{cap}_{A}(U,\overline{a})=\langle q,Lq\rangle.
\]
The vector $\overline{a}$ is such that $\overline{a}_{i}=a_{i}/w_{i}$.
This gives us the capacity of the vector $\overline{a}$ on the set
$U$ for the matrix $A$ (It is not the capacity of $a$ due to a
technical condition). Our definition is completely general, and special
cases of capacity for symmetric matrices have been studied (see \cite{doyle1984random})
by interpreting the graph as an electrical network, with the edge
weights $A_{i,j}=A_{j,i}$ denoting the conductance of the edge $\{i,j\}$.
Further, if $a\in\{0,1\}^{|U|},$ the $U$ can be written as $U=S\cup T$
with $a$ being 1 on $S$ and $0$ on $T$, and the resulting capacity
can be referred to as $\text{cap}(S,T)$. For symmetric matrices,
the vector $q$ is exactly the voltages at each vertex when the vertices
in $S$ are put at voltage 1 and the vertices in $T$ are put at voltage
0. Further, if the entire graph is modified such that there is one
vertex $s$ for $S$ and $t$ for $T$, and one edge between them,
then the effective resistance between $s$ and $t$ or the sets $S$
and $T$ in the original graph is exactly 
\[
\dfrac{1}{\text{cap}_{S,T}(A)}.
\]
In fact, it is possible to create many different results and algorithms
through the usage of capacity for symmetric matrices, and there have
been many results in this direction, see for instance \cite{doyle1984random,spielman2011graph,lyons1983simple,chandra1996electrical,tetali1991random}.

Our aim is to study capacity in its full generality for nonsymmetric
matrices. We show many basic results some of which are folklore for
symmetric matrices, but our main result is the following, discussed
and shown in Section \ref{subsec:Capacity-and-Schur}.

\begin{theorem}
\label{thm:intro-mon-cap}Let $A_{\alpha}=\alpha A+(1-\alpha)A^{T}$,
then for every $U$ and $a$, if   $\ 0\leq\beta\leq\alpha\leq\frac{1}{2}$,
then 
\[
\text{cap}_{U,a}(A_{\alpha})\leq\text{cap}_{U,a}(A_{\beta}).
\]
\end{theorem}

We remark a few things again about this result. The first is that
the result is simple if we compare a nonsymmetric matrix to a symmetric
matrix, i.e. it is simple to show that for $\tilde{A}=\frac{1}{2}(A+A^{T})$,
\[
\text{cap}_{U,\overline{a}}(\tilde{A})\leq\text{cap}_{U,\overline{a}}(A).
\]
This is in fact a direct consequence of Dirichlet's theorem, which
itself is a consequence of the Cauchy-Schwarz inequality. The remarkable
thing in Theorem \ref{thm:intro-mon-cap} is that it helps to compare the
capacity of two \emph{nonsymmetric} matrices, and thus the general
tool for symmetric matrices -- specifically the Cauchy-Schwarz inequality in this case
-- cannot be used. The proof is relatively long and intricate, but
we find it striking that capacity monotonically increases as $\alpha$
decreases and goes from $\frac{1}{2}$ to 0. Other related results
can also be found in Section \ref{subsec:Capacity-and-Schur}.

\subsection{A different notion of expansion and Tensor walks}

In the last section, we take the first steps towards a new definition
of expansion inspired by the recent studies of expansion in high-dimensional
expanders \cite{bafna2020high}, and we show a neat lemma similar
to the standard notion of edge expansion for it in Section \ref{subsec:A-different-notion} which we believe  could be useful in the study of high dimensional
objects. The second and final thing we study are tensor walks. Due
to the meteoric rise in interest in both machine learning and quantum
information, tensors have become a fundamental tool in both the areas.
However, the results related to them are extremely specific to the
application, and the task of showing general results is in its infancy.
This is mostly due to the behavior of tensors that is different from matrices in almost every aspect. To highlight this, consider
a $3$-tensor $T$ indexed as $T_{i,j,k}$, with $i,j,k\in[n]$. Given
probability vectors $p_{t}\in\mathbb{R}_{\geq0}^{n}$, we can define
one step of the tensor walk as 
\[
p_{t}(i)=\sum_{j,k}T_{i,j,k}p_{t-1}(j)p_{t-2}(k).
\]
To ensure $p_{t}$ is also a probability distribution, we can enforce
that for any $j,k$, $\sum_{i}T_{i,j,k}=1$. Note that one step of
the tensor walk is no longer a linear operation, but it is a very meaningful
walk, and there are an enormous number of areas where such walks naturally
appear. The same definition can be extended to any $k$-tensor. The
surprising thing is that even a 3-tensor could have an \emph{exponential}
(in $n$) number of eigenvalues (defined appropriately). There are
even multiple variants of the Perron-Frobenius theorem for nonnegative
tensors (see \cite{chang2008perron,friedland2013perron} for instance), however, even many basic questions remain unanswered. The
first question is, do positive tensors have a unique fixed point?
In fact, surprisingly, this turns out to be false, and an example
was shown in \cite{changzhang2013counterexample}. Our primary result
here is to delineate a condition for which it is sufficient for the
walk to have a unique fixed point. We say a $k$-tensor is stochastic over index $i$, if for every fixing of the other indices $l_1,l_2,...,l_{k-1}$, 
\[
\sum_{i=1}^n T_{l_1,l_2,i,l_3,...,l_{k-1}}=1.
\]

We say $T$ is 2-line stochastic if it is stochastic over the output index $i$ and any input index. Our main theorem is the following.

\begin{theorem}
\label{thm:intro-tens}Let $T$ be a $k$-tensor in $n$ dimensions
with positive entries, and let $T$ be $2$-line stochastic. Then
$T$ has a unique positive fixed point.
\end{theorem}

This result can be found in Section \ref{subsec:Tensors-and-Beyond}.
It raises an extremely important question: how fast does the walk
converge to the fixed point? We leave this question for future considerations.

The thesis is organized as follows. We first give the Preliminary
definitions in Section \ref{sec:Preliminaries}, proceed to prove our   
Theorem \ref{thm:intro-phi-delta}  in Section \ref{sec:Generalizations-of-the},
show theorems \ref{thm:intro-rootn} in Section \ref{subsec:Rootn-matrices-=002013}
and \ref{thm:intro-chet} in Section \ref{subsec:Chet-Matrices-=002013},
and theorems \ref{thm:intro-mixt}, \ref{thm:intro-mon-cap} and \ref{thm:intro-tens}
in Sections \ref{subsec:Mixing-Time}, \ref{subsec:Capacity-and-Schur}
and \ref{subsec:Tensors-and-Beyond} respectively. We conclude with
some summarizing final thoughts in Section \ref{sec:finth}.

\chapter{\label{sec:Preliminaries} Preliminaries -- Spectral Gap, Edge Expansion,
and the Perron-Frobenius Theorem}

\vspace{1.5cm}
\begin{quote}
As for the rest of my readers, they will accept such portions as
apply to them. I trust that none will stretch the seams in putting
on the coat, for it may do good service to him whom it fits.
\begin{flushright}
\textasciitilde{} Henry David Thoreau, \emph{Walden}\\
\par\end{flushright}
\end{quote}
\vspace{3cm}

We will consider nonnegative matrices throughout this thesis using
$R$ or $A$, also referred to as chains or graphs. For any nonnegative
matrix $R\in\mathbb{R}_{\geq0}^{n\times n}$, there is an implicit
underlying graph on $n$ vertices with the edge $(i,j)$ having weight
$R_{j,i}$, and we assume there is no edge if $R_{j,i}=0$. Note that
the subscripts are reversed since we assume right multiplication by
a vector, although this will not matter anywhere in this thesis except
for internal consistency in lemmas and proofs.

Given a matrix $R$, we say that $R$ is \emph{strongly connected}
or \emph{irreducible}, if there is a path from $s$ to $t$ for every
pair $(s,t)$ of vertices in the underlying digraph on edges with
positive weight, i.e. for every $(s,t)$ there exists $k>0$ such
that $R^{k}(s,t)>0$. We say $R$ is \emph{weakly connected}, if there
is a pair of vertices $(s,t)$ such that there is a path from $s$
to $t$ but no path from $t$ to $s$ in the underlying digraph (on
edges with positive weight).

We start by stating the Perron-Frobenius theorem. The theorem was
shown for positive matrices $R$ by Perron in \cite{perron1907theorie},
and for irreducible matrices $R$ by Frobenius in \cite{frobenius1912matrizen}.
In the last hundred years, the theorem has been used extensively in
many different areas of mathematics as discussed in the Introduction,
and has become a fundamental tool within spectral theory and dynamical
systems.

\begin{theorem}
\emph{\label{thm:(Perron-Frobenius)}(Perron-Frobenius theorem \cite{perron1907theorie,frobenius1912matrizen})}
Let $R\in\mathbb{R}^{n\times n}$ be a nonnegative matrix. Then the
following hold for $R$. 
\begin{enumerate}
\item $R$ has some nonnegative eigenvalue $r$, such that all other eigenvalues
have magnitude at most $r$, and $R$ has nonnegative left and right
eigenvectors $u$ and $v$ for $r$. 
\item If $R$ has some \emph{positive} left and right eigenvectors $u'$
and $v'$ for some eigenvalue $\lambda$, then $\lambda=r$. 
\item If $R$ is irreducible, then the eigenvalue $r$ is positive and simple
(unique), the corresponding left and right eigenvectors $u$ and $v$ are positive and unique, and all other eigenvalues
$\lambda$ are such that $|\lambda|\leq r$, and $\text{Re}\lambda<r$.
\end{enumerate}
\end{theorem}

Many nice proofs of the theorem can be found in different lecture
notes online, and proving it directly is also not difficult. We will state all our results for irreducible matrices, and
they will extend for general matrices by using limit infimums. By
the Perron-Frobenius (Theorem \ref{thm:(Perron-Frobenius)}, part
3), an irreducible nonnegative matrix $R$ will have a simple positive
eigenvalue $r$ such that all eigenvalues have magnitude strictly
less than $r$, and it will be called the \emph{trivial }or \emph{stochastic}
or PF eigenvalue of $R$, and all other eigenvalues of $R$ will be
called \emph{nontrivial. } The left and right eigenvectors corresponding
to $r$ will be called the \emph{trivial or PF left eigenvector (generally
referred to as $u$) }and \emph{trivial or PF right eigenvector} (generally
referred to as $v$). This leads us to the following definition.
\medskip

\begin{definition}
\emph{\label{def:(Spectral-gap-irred}(Spectral gap of irreducible
nonnegative matrices)} Let $R$ be an $n\times n$ irreducible nonnegative
matrix. Let the eigenvalues $\lambda_{1}$ to $\lambda_{n}$ of $R$
be arranged so that $\lambda_{1}>\text{Re}\lambda_{2}\geq\text{Re}\lambda_{3}\geq...\geq\text{Re}\lambda_{n}$,
where due to the Perron-Frobenius Theorem \ref{thm:(Perron-Frobenius)},
$\lambda_{1}$ is real and positive and $\text{Re}\lambda_{2}<\lambda_{1}$.
Define the spectral gap of $R$ as
\[
\triangle(R)=1-\dfrac{\text{Re}\lambda_{2}(R)}{\lambda_{1}(R)}.
\]
\end{definition}

Note that $0\leq\Delta(R)\leq2$ from the Perron-Frobenius theorem,
since $|\lambda_{i}(R)|<\lambda_{1}(R)$. We remind the reader that
the eigenvalues of $R$ are simply the roots of the characteristic
polynomial, i.e. of the equation $\text{det}(\lambda I-R)=0$. We
will also consider singular values of nonnegative matrices $A$ with
identical positive left and right eigenvector $w$ for PF eigenvalue
1, and denote them as $1=\sigma_{1}(A)\geq\sigma_{2}(A)\geq\cdots\geq\sigma_{n}(A)\geq0$
(see Lemma \ref{lem:transRtoA} for proof of $\sigma_{1}(A)=1$).
We denote $(i,j)$'th entry of $M\in\mathbb{C}^{n\times n}$ by $M_{i,j}$,
or if $M$ already has a subscript, then as $M_{t}(i,j)$ as will be
clear from context. We denote the conjugate-transpose of $M$ as $M^{*}$
and the transpose of $M$ as $M^{T}$. Any $M\in\mathbb{C}^{n\times n}$
has a \emph{Schur decomposition} (see, e.g., \cite{lax07}) $M=UTU^{*}$
where $T$ is an upper triangular matrix whose diagonal entries are
the eigenvalues of $M$, and $U$ is a unitary matrix ($UU^{*}=U^{*}U=I$).
When we write ``vector'' we mean by default a column vector. For
a vector $v$, we again write or $v_{i}$ to denote its $i$'th entry
or $v_{t}(i)$ in case a subscript is already present. For any two
vectors $x,y\in\mathbb{C}^{n}$, we use the standard \emph{inner product
}$\langle x,y\rangle=\sum_{i=1}^{n}x_{i}^{*}\cdot y_{i}$ defining
the norm $\|x\|_{2}=\sqrt{\langle x,x\rangle}$. We write $u\perp v$
to indicate that $\langle u,v\rangle=0$. Note that $\langle x,My\rangle=\langle M^{*}x,y\rangle$.
We denote the operator norm of $M$ by $\|M\|_{2}=\max_{u:\|u\|_{2}=1}\|Mu\|_{2}$,
and recall that the operator norm is at most the Frobenius norm, i.e.,
$\|M\|_{2}\leq\|M\|_{F}:=\sqrt{\sum_{i,j}|M_{i,j}|^{2}}.$ We write
$D_{u}$ for the diagonal matrix whose diagonal contains the vector
$u$. We will use the phrase ``positive'' to mean entry-wise positive,
and the phrase ``expansion'' to mean edge expansion. We note the
Birkhoff-von Neumann theorem, which states that every doubly stochastic
matrix $A$ can be written as $A=\sum\alpha_{i}P_{i}$, where $P_{i}$
are permutation matrices and the $\alpha_{i}\geq0$ and $\sum\alpha_{i}=1$.

Recall the Courant-Fischer variational characterization of eigenvalues
for symmetric real matrices, applied to the second eigenvalue: 
\[
\max_{u\perp v_{1}}\frac{\langle u,Mu\rangle}{\langle u,u\rangle}=\lambda_{2}(M),
\]
where $v_{1}$ is the eigenvector for the largest eigenvalue of $M$.
We will use the symbol $J$ for the all $1$'s matrix divided by $n$,
i.e., $J=\frac{1}{n}\mathbf{1}\cdot\mathbf{1}^{T}$. We will denote
the standard basis vectors as $e_{1},...,e_{n}$. We denote the all
$1$'s vector by $\boldsymbol{1}$, and we say that any subset $S\subseteq[n]$
is a \emph{cut}, denote its complement by $\overline{S}$, and denote
the \emph{characteristic vector of a cut} as $\boldsymbol{1}_{S}$,
where $\boldsymbol{1}_{S}(i)=1$ if $i\in S$ and 0 otherwise. We
first present the definition of edge expansion for doubly stochastic matrices. 
\medskip

\begin{definition}
\emph{\label{def:(Edge-expansion-doubstoc}(Edge expansion of doubly
stochastic matrices)} For a doubly stochastic matrix $A$, the \emph{edge
expansion of the cut $S$} is defined as 
\[
\phi_{S}(A):=\dfrac{\langle\mathbf{1}_{S},A\mathbf{1}_{\overline{S}}\rangle}{\langle\mathbf{1}_{S},A\mathbf{1}\rangle}
\]
and the \emph{edge expansion of $A$} is defined as 
\[
\phi(A)=\min_{S,|S|\leq n/2}\phi_{S}(A)=\min_{S,|S|\leq n/2}\frac{\sum\limits _{i\in S,j\in\overline{S}}A_{i,j}}{|S|}.
\]
\end{definition}

We wish to extend these notions to general nonnegative matrices $R$.
Since eigenvalues and singular values of real matrices remain unchanged
whether we consider $R$ or $R^{T}$, the same should hold of a meaningful
definition of edge expansion. However, note that Definition \ref{def:(Edge-expansion-doubstoc}
has this independence only if the matrix is Eulerian, i.e., $R\mathbf{1}=R^{T}\mathbf{1}$.
Thus, to define edge expansion for general matrices, we transform
$R$ using its left and right eigenvectors $u$ and $v$ to obtain
$D_{u}RD_{v}$, which is indeed Eulerian, since 
\[
D_{u}RD_{v}\mathbf{1}=D_{u}Rv=rD_{u}v=rD_{u}D_{v}\mathbf{1}=rD_{v}D_{u}\mathbf{1}=rD_{v}u=D_{v}R^{T}u=D_{v}R^{T}D_{u}\mathbf{1}.
\]
Since $D_{u}RD_{v}$ is Eulerian, we can define the edge expansion
of $R$ similar to that for doubly stochastic matrices: 
\medskip

\begin{definition}
\label{def:gen-edge-expansion} \emph{(Edge expansion of irreducible
nonnegative matrices)} Let $R\in\mathbb{R}^{n\times n}$ be an irreducible
nonnegative matrix with positive left and right eigenvectors $u$
and $v$ for the PF eigenvalue $r$. The \emph{edge expansion of the
cut $S$} is defined as 
\begin{align}
\phi_{S}(R) & :=\dfrac{\langle\mathbf{1}_{S},D_{u}RD_{v}\mathbf{1}_{\overline{S}}\rangle}{\langle\mathbf{1}_{S},D_{u}RD_{v}\mathbf{1}\rangle}\label{phiSR}
\end{align}
and the \emph{edge expansion of $R$} is defined as 
\[
\phi(R)=\min_{S:\sum_{i\in S}u_{i}\cdot v_{i}\leq\frac{1}{2}\sum_{i}u_{i}v_{i}}\phi_{S}(R)=\min_{S:\sum_{i\in S}u_{i}\cdot v_{i}\leq\frac{1}{2}\sum_{i}u_{i}v_{i}}\dfrac{\sum\limits _{i\in S,j\in\overline{S}}R_{i,j}\cdot u_{i}\cdot v_{j}}{\sum\limits _{i\in S,j}R_{i,j}\cdot u_{i}\cdot v_{i}}.
\]
\end{definition}

The edge expansion $\phi(R)$ can be intuitively understood as the
\emph{measure of irreducibility of $R$} (though the correspondence
is not exact in the strictest sense). If $\phi(R)$ is high, say a
constant or even $1/\text{polylog}(n)$, then $R$ is well-connected
(strongly). As $\phi(R)$ becomes smaller, it implies there are two
(strongly) disconnected components within $R$, or that $R$ has a
sink-like component.

Our next aim is to define $\Delta(R)$ and $\phi(R)$ for nonnegative
matrices $R$ that are not irreducible. We do this by considering the
ball of irreducible matrices in a small ball around $R$, and taking
the limit of the infimum to get a well-defined quantity. 
\medskip

\begin{definition}
\emph{\label{def:(gen-Spectral-gap-and-edgeexp}(Spectral gap and Edge Expansion
of any nonnegative matrix)} Let $E$ be the set of all irreducible
nonnegative matrices. For any nonnegative matrix $R$, let 
\[
R_{\epsilon}=\{H:H\in E,\ \|R-H\|_{F}\leq\epsilon\}.
\]
Then we define the spectral gap of $R$ as 
\[
\Delta(R)=\lim_{\epsilon\rightarrow0}\inf_{H\in R_{\epsilon}}\Delta(H),
\]
and the edge expansion of $R$ as 
\[
\phi(R)=\lim_{\epsilon\rightarrow0}\inf_{H\in R_{\epsilon}}\phi(H).
\]
\end{definition}

The first simple thing we can say is the following.

\begin{lemma}
\label{lem:if_irred_then_phi_morethan_0}Let $R$ be a nonnegative
matrix. If $R$ is irreducible, then $\phi(R)>0$.
\end{lemma}

\begin{proof}
Since $R$ is irreducible, the left and right PF eigenvectors $u$
and $v$ are positive. For the sake of contradiction, if $\phi(R)=0$,
there is some $S$ with $\phi_{S}(R)=0$, implying that $\langle\mathbf{1}_{S},D_{u}RD_{v}\mathbf{1}_{\overline{S}}\rangle=0$,
and since $u$ and $v$ are positive, it implies that $R_{i,j}=0$
for $i\in S,j\in\overline{S}$, and similarly $R_{i,j}=0$ for $i\in\overline{S},j\in S$
since $D_{u}RD_{v}$ is Eulerian. This implies there is no path from
$i\in S$ to $j\in\overline{S}$, implying $R$ is not irreducible,
a contradiction.
\end{proof}

According to Perron-Frobenius (Theorem \ref{thm:(Perron-Frobenius)}, part 3), if $R$ is irreducible
then $\text{Re}\lambda_{2}(R)<\lambda_{1}(R)$. However, since irreducibility
implies positivity of edge expansion from Lemma \ref{lem:if_irred_then_phi_morethan_0},
it means that the limiting case of Theorem \ref{thm:intro-phi-delta}
showing that $\phi(R)>0\Leftrightarrow\Delta(R)>0$ implies that we
get a much stronger and tighter statement than the Perron-Frobenius
theorem.

\chapter{\label{sec:Generalizations-of-the}Generalizations of the Perron-Frobenius
Theorem and the Cheeger-Buser Inequalities}

\vspace{1.5cm}
\begin{quote}
Human reason has the peculiar fate in one species of its cognitions
that it is burdened with questions which it cannot dismiss, since
they are given to it as problems by the nature of reason itself, but
which it also cannot answer, since they transcend every capacity of
human reason.
\begin{flushright}
\textasciitilde{} Immanuel Kant, \emph{Critique of Pure Reason}
\par\end{flushright}

\end{quote}
\vspace{3cm}

The aim of this section is to prove our main Theorem \ref{thm:intro-phi-delta}.
As stated in the Introduction, there are two well known theorems in
Spectral Theory. The first is the Perron-Frobenius Theorem, which
provides us, among other results, a \emph{qualitative} statement connecting
the edge expansion of matrices and their second eigenvalue. More specifically,
it tells us that for any irreducible nonnegative matrix $R$ -- for
which it is always the case that $\phi(R)>0$ -- the spectral gap
$\Delta(R)>0$. The second are the Cheeger-Buser inequalities as follows.

\begin{theorem}
\label{thm:cheeger-buser} \emph{(Cheeger-Buser Inequality \cite{cheeger70,Buser82,AlonM85,SinclairJ89,dodziuk84,Nilli91})}
Let $R$ be a \emph{reversible} nonnegative matrix. Then 
\[
\frac{1}{2}\cdot\Delta(R)\leq\phi(R)\leq\sqrt{2\cdot\Delta(R)}.
\]
\end{theorem}

This gap in our understanding, between the qualitative result for
the nonreversible case and a quantitative result for the reversible
case, is the primary focus of this section. Our main aim in theorem
\ref{thm:intro-phi-delta} is to show an inequality similar to the
form above but for \emph{any} nonnegative matrix (not necessarily
reversible).

We will prove Theorem \ref{thm:intro-phi-delta}, for the case
of irreducible nonnegative matrices first, and extend it to all matrices (not necessary irreducible) in Section \ref{sec:extension-to-gen-mats}. We restate the main theorem that
we want to prove in this chapter.

\begin{theorem}
\label{thm:main_irred}Let $R$ be an irreducible nonnegative matrix.
Then 
\[
\dfrac{1}{15}\cdot\dfrac{\Delta(R)}{n}\leq\phi(R)\leq\sqrt{2\cdot\Delta(R)}.
\]
\end{theorem}

The upper bound on $\phi(R)$ akin to Cheeger's inequality follows
by a straightforward extension of Fiedler's proof for doubly stochastic
matrices \cite{fiedler1995estimate}, and the lower bound will require
substantial work.

However, before we go into the details of the proof of
Theorem \ref{thm:intro-phi-delta}, we will first prove the Cheeger-Buser
Inequalities both for completeness and concreteness, since our Definition \ref{def:gen-edge-expansion}
of edge expansion is different from some sources. We
will start by proving a crucial lemma that will be used in all subsequent
sections.

\begin{lemma}
\label{lem:transRtoA} Let $R$ be an irreducible nonnegative matrix
with positive (left and right) eigenvectors $u$ and $v$ for the
PF eigenvalue 1, normalized so that $\langle u,v\rangle=1$. Define
$A=D_{u}^{\frac{1}{2}}D_{v}^{-\frac{1}{2}}RD_{u}^{-\frac{1}{2}}D_{v}^{\frac{1}{2}}$.
Then the following hold for $A$: 
\begin{enumerate}
\item $\phi(A)=\phi(R)$. 
\item For every $i$, $\lambda_{i}(A)=\lambda_{i}(R)$.
\item $\|A\|_{2}=1$.
\item If $R$ is reversible, i.e. $D_{u}RD_{v}=D_{v}R^{T}D_{u}$, then $A$
is symmetric.
\end{enumerate}
\end{lemma}

\begin{proof}
Let the matrix $A$ be as defined, and let $w=D_{u}^{\frac{1}{2}}D_{v}^{\frac{1}{2}}\mathbf{1}$.
Then it is easily checked that $Aw=w$ and $A^{T}w=w$. Further, 
\[
\langle w,w\rangle=\langle D_{u}^{\frac{1}{2}}D_{v}^{\frac{1}{2}}\mathbf{1},D_{u}^{\frac{1}{2}}D_{v}^{\frac{1}{2}}\mathbf{1}\rangle=\langle D_{u}\mathbf{1},D_{v}\mathbf{1}\rangle=\langle u,v\rangle=1
\]
where we used the fact that the matrices $D_{u}^{\frac{1}{2}}$ and
$D_{v}^{\frac{1}{2}}$ are diagonal, and so they commute, and are
unchanged by taking transposes. Let $S$ be any set. The condition
$\sum_{i\in S}u_{i}\cdot v_{i}\leq\frac{1}{2}$ translates to $\sum_{i\in S}w_{i}^{2}\leq\frac{1}{2}$
since $u_{i}\cdot v_{i}=w_{i}^{2}$. Thus, for any set $S$ for which
$\sum_{i\in S}u_{i}\cdot v_{i}=\sum_{i\in S}w_{i}^{2}\leq\frac{1}{2}$,
\begin{align*}
\phi_{S}(R) & =\frac{\langle\mathbf{1}_{S},D_{u}RD_{v}\mathbf{1}_{\overline{S}}\rangle}{\langle\mathbf{1}_{S},D_{u}D_{v}\mathbf{1}_{S}\rangle}\\
 & =\frac{\langle\mathbf{1}_{S},D_{u}D_{u}^{-\frac{1}{2}}D_{v}^{\frac{1}{2}}AD_{u}^{\frac{1}{2}}D_{v}^{-\frac{1}{2}}D_{v}\mathbf{1}_{\overline{S}}\rangle}{\langle\mathbf{1}_{S},D_{u}^{\frac{1}{2}}D_{v}^{\frac{1}{2}}D_{u}^{\frac{1}{2}}D_{v}^{\frac{1}{2}}\mathbf{1}_{S}\rangle}\\
 & =\frac{\langle\mathbf{1}_{S},D_{u}^{\frac{1}{2}}D_{v}^{\frac{1}{2}}AD_{u}^{\frac{1}{2}}D_{v}^{\frac{1}{2}}\mathbf{1}_{\overline{S}}\rangle}{\langle\mathbf{1}_{S},D_{u}^{\frac{1}{2}}D_{v}^{\frac{1}{2}}D_{u}^{\frac{1}{2}}D_{v}^{\frac{1}{2}}\mathbf{1}_{S}\rangle}\\
 & =\frac{\langle\mathbf{1}_{S},D_{w}AD_{w}\mathbf{1}_{\overline{S}}\rangle}{\langle\mathbf{1}_{S},D_{w}^{2}\mathbf{1}_{S}\rangle}\\
 & =\phi_{S}(A)
\end{align*}
and (1) holds. Further, since $A$ is a similarity transform of $R$,
all eigenvalues are preserved and (2) holds. For (3), consider the
matrix $H=A^{T}A$. Since $w$ is the positive left and right eigenvector
for $A$, i.e. $Aw=w$ and $A^{T}w=w$, we have $Hw=w$. But since
$A$ was nonnegative, so is $H$, and since it has a positive eigenvector
$w$ for eigenvalue 1, by Perron-Frobenius (Theorem \ref{thm:(Perron-Frobenius)},
part 2), $H$ has PF eigenvalue 1. But $\lambda_{i}(H)=\sigma_{i}^{2}(A)$,
where $\sigma_{i}(A)$ is the $i$'th largest singular value of $A$.
Thus, we get $\sigma_{1}^{2}(A)=\lambda_{1}(H)=1$, and thus $\|A\|_{2}=1$.
For (4), 
\begin{align*}
D_{u}RD_{v} & =D_{v}R^{T}D_{u}\Leftrightarrow\\
D_{u}D_{u}^{-\frac{1}{2}}D_{v}^{\frac{1}{2}}AD_{u}^{\frac{1}{2}}D_{v}^{-\frac{1}{2}}D_{v} & =D_{v}D_{u}^{\frac{1}{2}}D_{v}^{-\frac{1}{2}}A^{T}D_{u}^{-\frac{1}{2}}D_{v}^{\frac{1}{2}}\Leftrightarrow\\
D_{u}^{\frac{1}{2}}D_{v}^{\frac{1}{2}}AD_{u}^{\frac{1}{2}}D_{v}^{\frac{1}{2}} & =D_{u}^{\frac{1}{2}}D_{v}^{\frac{1}{2}}A^{T}D_{u}^{\frac{1}{2}}D_{v}^{\frac{1}{2}}\Leftrightarrow\\
A & =A^{T}
\end{align*}
as required.
\end{proof}
Given Lemma \ref{lem:transRtoA}, we can restrict to showing proofs
for $A$ instead of $R$, since both the spectral gap and  edge expansion 
remain unaffected.

We will now prove the Cheeger-Buser Inequalities. There
are many well-known proofs in literature, but the definition of $\phi(R)$
for general matrices has often been unsatisfying (although it works
technically), and we show complete proofs for our Definition \ref{def:gen-edge-expansion}
of $\phi(R)$ which is the same definition as in \cite{mihail1989conductance}.
This exposition will also be helpful in contrasting with similar lemmas
for the nonreversible case.

\section{\label{subsec:Buser-Inequality-=002013}Buser Inequality -- the lower bound on $\phi(R)$}

We start by showing the lower bound on $\phi(R)$ for reversible $R$.
To achieve this, note that for the case of
reversible $R$, the definition of $\phi(R)$ for reducible $R$ is
exactly same as the irreducible case (by using the variational characterization for the second eigenvalue as stated in the Preliminaries), and thus we show the theorem for
irreducible reversible $R$. For any such $R$, note from the Perron-Frobenius
Theorem \ref{thm:(Perron-Frobenius)} that $\lambda_{1}(R)>0$, and
by rescaling the matrix with a constant -- which does not change
$\Delta(R)$ or $\phi(R)$ (Lemma \ref{lem:transRtoA}) -- we assume
that $\lambda_{1}(R)=1$. Further, we will work with $A=D_{u}^{\frac{1}{2}}D_{v}^{-\frac{1}{2}}RD_{u}^{-\frac{1}{2}}D_{v}^{\frac{1}{2}}$
which is symmetric since $R$ is reversible from Lemma \ref{lem:transRtoA},
is irreducible, and has the same spectral gap and  edge expansion as $R$. Thus,
we will show a lower bound on $\phi(A)$ where we assume $Aw=w$ and
$A^{T}w=w$ for a positive $w$ (Theorem \ref{thm:(Perron-Frobenius)}).
The lower bound easily follows from the spectral decomposition of
A. Note that the following lemma does not require $A$ to be reversible.

\begin{lemma}
\label{lem:phi_low_bound_gamma}Let $A$ be any irreducible nonnegative
matrix with largest eigenvalue 1 and corresponding positive left and right eigenvector
$w$. If $\langle v,Av\rangle\leq\gamma\langle v,v\rangle$ with $0\leq\gamma\leq1$
for every vector $v$ with $v\perp w$, then 
\[
\dfrac{1-\gamma}{2}\leq\phi(A).
\]
\end{lemma}

\begin{proof}
Let $A=ww^{*}+B$ with $Bw=0$, then $\langle v,Av\rangle=\langle v,Bv\rangle$,
and for any set $S$, let $D_{w}1_{S}=c_{0}w+v$ where $v\perp w$.
Then 
\[
\langle w,D_{w}1_{S}\rangle=c_{0}
\]
and 
\[
\langle D_{w}1_{S},D_{w}1_{S}\rangle=|c_{0}|^{2}+\langle v,v\rangle
\]
and note that 
\[
\langle w,D_{w}1_{S}\rangle=\langle D_{w}1_{S},D_{w}1_{S}\rangle=\langle D_{w}1,D_{w}1_{S}\rangle.
\]
Then for any set $S$ with $\langle D_{w}1_{S},D_{w}1_{S}\rangle\ensuremath{\leq1/2}$,
we have 
\begin{align*}
\phi_{S}(A) & =\dfrac{\langle D_{w}1_{\overline{S}},AD_{w}1_{S}\rangle}{\langle D_{w}1,D_{w}1_{S}\rangle}\\
 & =1-\dfrac{\langle D_{w}1_{S},AD_{w}1_{S}\rangle}{\langle D_{w}1,D_{w}1_{S}\rangle}\\
 & =1-\dfrac{\langle c_{0}w+v,A(c_{0}w+v)\rangle}{\langle D_{w}1,D_{w}1_{S}\rangle}\\
 & =1-\dfrac{\langle c_{0}w+v,c_{0}w+Av\rangle}{\langle D_{w}1,D_{w}1_{S}\rangle}\\
 & =1-\dfrac{|c_{0}|^{2}+\langle v,Av\rangle}{\langle D_{w}1,D_{w}1_{S}\rangle}\\
 & \geq1-\dfrac{|c_{0}|^{2}+\gamma\langle v,v\rangle}{\langle D_{w}1,D_{w}1_{S}\rangle}\\
 & =1-\dfrac{\langle D_{w}1_{S},D_{w}1_{S}\rangle^{2}+\gamma(\langle D_{w}1_{S},D_{w}1_{S}\rangle-\langle D_{w}1_{S},D_{w}1_{S}\rangle^{2})}{\langle D_{w}1,D_{w}1_{S}\rangle}\\
 & =1-(\gamma+(1-\gamma)\langle D_{w}1_{S},D_{w}1_{S}\rangle)\\
 & \ \ \ \ [\text{since \ensuremath{\langle D_{w}1_{S}},\ensuremath{D_{w}1_{S}\rangle}\ensuremath{\ensuremath{\leq}1/2}}]\\
 & \geq\dfrac{1-\gamma}{2}
\end{align*}
as required.
\end{proof}
The Buser inequality easily follows from the following.

\begin{lemma}
\label{lem:buser2}For any $v\perp w$, $\langle v,Av\rangle\leq\lambda_{2}(A)\langle v,v\rangle$
\end{lemma}

\begin{proof}
Let $u_{1}=w,u_{2},u_{3},...,u_{n}$ be an orthogonal set of eigenvectors
for $A$ (See Preliminaries \ref{sec:Preliminaries}), and let $v=\sum_{i}c_{i}u_{i}$
with $c_{1}=0$ since $v\perp w$, and 
\[
\langle v,Av\rangle=\langle\sum_{i}c_{i}u_{i},A\sum_{i}c_{i}u_{i}\rangle=\langle\sum_{i}c_{i}u_{i},\sum_{i}c_{i}\lambda_{i}u_{i}\rangle=\sum_{i\geq2}c_{i}^{2}\lambda_{i}\leq\lambda_{2}\sum_{i}c_{i}^{2}=\lambda_{2}\langle v,v\rangle
\]
as required.
\end{proof}
Combining Lemma \ref{lem:phi_low_bound_gamma} and Lemma \ref{lem:buser2}
gives the lower bound on $\phi(A)$ in Theorem \ref{thm:cheeger-buser}.

\section{Cheeger Inequality -- the upper bound on $\phi(R)$}

We now proceed with the upper bound on $\phi(A)$ which is relatively
complicated. Our proof method is similar to that in \cite{mihail1989conductance},
\cite{chung2007four} and going earlier to \cite{AlonM85}.

\begin{lemma}
\label{lem:cheeger}$\phi(A)\leq\sqrt{2\cdot\Delta(A)}$.
\end{lemma}

\begin{proof}
Let $v$ be the eigenvector for eigenvalue $\lambda_{2}(A)$ with
$v\perp w$ and assume $\langle w,w\rangle=1$. Then we have that
\[
1-\lambda_{2}=\dfrac{\langle v,(I-A)v\rangle}{\langle v,v\rangle}.
\]
Let $u=D_{w}^{-1}v-c\cdot\mathbf{1}$, and note that $\langle D_{w}u,D_{w}u\rangle=\langle v,v\rangle+c^{2}$
since $v\perp w$. Our aim is to choose $c$ in order to divide the
positive and negative entries of $u$ into groups $S$ and $\overline{S}$,
such that $\sum_{i\in S}w_{i}^{2}\leq1/2$ and $\sum_{i\in\overline{S}}w_{i}^{2}\leq1/2$.
To achieve this, without loss of generality, assume that the entries
of $D_{w}^{-1}v$ are arranged in decreasing order (as indices go
from 1 to $n$), and let $r$ be the smallest index such that $\sum_{i=1}^{r}w_{i}^{2}>1/2$.
Thus $\sum_{i=1}^{r-1}w_{i}^{2}\leq1/2$ and $\sum_{i=r+1}^{n}w_{i}^{2}\leq1/2$
since $\sum_{i}w_{i}^{2}=1$. We then choose $c=v_{r}/w_{r}$, and
this gives us $u_{i}\geq0$ for $i<r$, $u_{r}=0$, and $u_{i}\leq0$
for $i>r$. Thus, letting $x$ be the vector with $x_{i}=u_{i}$ for $i<r$
and 0 otherwise, letting $y$ be the vector with $y_{i}=u_{i}$ for $i>r$ and $0$ otherwise, we
get $u=x-y$, with the property that $\sum_{i}w_{i}^{2}\leq1/2$ when
$i$ runs over nonzero entries of $x$ or $y$, and $\langle x,y\rangle=0$.
Thus, we have that
\begin{align}
1-\lambda_{2} & =\dfrac{\langle v,(I-A)v\rangle}{\langle v,v\rangle}\nonumber \\
 & =\dfrac{\langle D_{w}u,(I-A)D_{w}u\rangle}{\langle D_{w}u,D_{w}u\rangle+c^{2}}\nonumber \\
 & \ \ \ [\text{since \ensuremath{(I-A)w=0} and \ensuremath{(I-A)^{T}w=0}}]\nonumber \\
 & \geq\dfrac{\langle D_{w}x,(I-A)D_{w}x\rangle+\langle D_{w}y,(I-A)D_{w}y\rangle-2\langle D_{w}x,(I-A)D_{w}y\rangle}{\langle D_{w}x,D_{w}x\rangle+\langle D_{w}y,D_{w}y\rangle}\nonumber \\
 & \geq\dfrac{\langle D_{w}x,(I-A)D_{w}x\rangle+\langle D_{w}y,(I-A)D_{w}y\rangle}{\langle D_{w}x,D_{w}x\rangle+\langle D_{w}y,D_{w}y\rangle}\nonumber \\
 & \ \ \ 
 [\text{since } \langle D_{w}x,D_{w}y \rangle - 
 \langle D_{w}x,A\cdot D_{w}y \rangle
 = 0-\langle D_{w}x,A\cdot D_{w}y \rangle \leq 0, \nonumber \\
 & \ \ \ \ \text{as all entries of } x,y,w,A \text{ are nonnegative}
 ]\nonumber \\
 & \geq\dfrac{\langle D_{w}x,(I-A)D_{w}x\rangle}{\langle D_{w}x,D_{w}x\rangle}\label{eq:cheeg-1}
\end{align}
where we assume in the last line that 
\[
\dfrac{\langle D_{w}x,(I-A)D_{w}x\rangle}{\langle D_{w}x,D_{w}x\rangle}\leq\dfrac{\langle D_{w}y,(I-A)D_{w}y\rangle}{\langle D_{w}y,D_{w}y\rangle}.
\]
Note that $Aw=w$ and $A^{T}w=w$, so for any $i$, $w_{i}^{2}=\sum_{j}a_{i,j}w_{i}w_{j}$,
for any $j$, $w_{j}^{2}=\sum_{i}a_{i,j}w_{i}w_{j}$. Let $a'_{i,j}=a_{i,j}w_{i}w_{j}$.
So
\begin{align*}
\langle D_{w}x,D_{w}x\rangle & =\sum_{i}w_{i}^{2}x_{i}^{2}=\dfrac{1}{2}\cdot \left(\sum_{i}w_{i}^{2}x_{i}^{2}+\sum_{j}w_{j}^{2}x_{j}^{2} \right)\\
 & =\dfrac{1}{2}\cdot\left(\sum_{i,j}a_{i,j}w_{j}w_{i}x_{i}^{2}+\sum_{j,i}a_{i,j}w_{i}w_{j}x_{j}^{2}\right)\\
 & =\dfrac{1}{2}\cdot\sum_{i,j}a_{i,j}w_{i}w_{j}(x_{i}^{2}+x_{j}^{2})\\
 & \geq\dfrac{1}{4}\cdot\sum_{i,j}a'_{i,j}(x_{i}+x_{j})^{2}
\end{align*}
and 
\[
\langle D_{w}x,(I-A)D_{w}x\rangle=\frac{1}{2}\sum_{i,j}a'_{i,j}(x_{i}-x_{j})^{2}
\]
and multiplying LHS of the two equations, we get 
\begin{align*}
\langle D_{w}x,D_{w}x\rangle\langle D_{w}x,(I-A)D_{w}x\rangle & \geq\frac{1}{8}\sum_{i,j}a'_{i,j}(x_{i}+x_{j})^{2}\sum_{i,j}a'_{i,j}(x_{i}-x_{j})^{2}\\
 & \geq\frac{1}{8}\left(\sum_{i,j}\sqrt{a'_{i,j}}\left|x_{i}-x_{j}\right|\cdot\sqrt{a'_{i,j}}\left|x_{i}+x_{j}\right|\right)^{2}\\
 & \ \ \ \text{[Cauchy-Shwarz]}\\
 & =\frac{1}{8}\left(\sum_{i,j}a'_{i,j}\left|x_{i}^{2}-x_{j}^{2}\right|\right)^{2}
\end{align*}
Since $x_{i}\geq x_{j}$ for $i<j$, then for $m'_{i,j}=a'_{i,j}+a'_{j,i}$
\begin{align*}
\sum_{i,j}a'_{i,j}\left|x_{i}^{2}-x_{j}^{2}\right| & =\sum_{i<j}m'_{i,j}(x_{i}^{2}-x_{j}^{2})\\
 & =\sum_{i=1}^{n}\sum_{j=i+1}^{n}\sum_{k=i}^{j-1}m'_{i,j}(x_{k}^{2}-x_{k+1}^{2})\\
 & =\sum_{i=1}^{n}\sum_{k=i}^{n-1}\sum_{j=k+1}^{n}m'_{i,j}(x_{k}^{2}-x_{k+1}^{2})\\
 & =\sum_{k=1}^{n-1}\sum_{i=1}^{k}\sum_{j=k+1}^{n}m'_{i,j}(x_{k}^{2}-x_{k+1}^{2})\\
 & =\sum_{k=1}^{r-1}(x_{k}^{2}-x_{k+1}^{2})\sum_{i=1}^{k}\sum_{j=k+1}^{n}m'_{i,j}\\
 & \ \ \ \text{[since for \ensuremath{k\geq r}, \ensuremath{x_{k}=0}]}\\
 & =2\sum_{k=1}^{r-1}(x_{k}^{2}-x_{k+1}^{2})\cdot\phi_{k}\cdot\mu_{k}\\
 & \ \ \ \text{[factor of 2 since \ensuremath{m_{i,j}} is already sum of two entries,}\\
 & \ \ \ \text{and \ensuremath{\mu_{k}=\sum_{i=1}^{k}w_{i}^{2}\leq1/2} for any \ensuremath{k\leq r}]}\\
 & \geq2\alpha\cdot\sum_{k=1}^{r-1}(x_{k}^{2}-x_{k+1}^{2})\cdot\mu_{k}\\
 & \text{\ \ \ [letting \ensuremath{\alpha}=\ensuremath{\min_{k}\phi_{k}}]}\\
 & =2\alpha\cdot\sum_{k=1}^{r-1}w_{k}^{2}x_{k}^{2}\\
 & =2\alpha\cdot\langle D_{w}x,D_{w}x\rangle
\end{align*}
which finally gives after combining everything, 
\begin{align*}
1-\lambda_{2} & =\dfrac{\langle v,(I-A)v\rangle}{\langle v,v\rangle}\\
 & \geq\dfrac{\langle D_{w}x,(I-A)D_{w}x\rangle}{\langle D_{w}x,D_{w}x\rangle}\\
 & \ \ \ \text{[from equation \ref{eq:cheeg-1}]}\\
 & \geq\frac{1}{8}\cdot\dfrac{\left(2\alpha\cdot\langle D_{w}x,D_{w}x\rangle\right)^{2}}{\langle D_{w}x,D_{w}x\rangle^{2}}\\
 & =\frac{1}{2}\cdot\alpha^{2}
\end{align*}
as required.
\end{proof}

\section{\label{subsec:Tightness-of-cheeger-bounds}Tightness of the Cheeger-Buser Inequality}

It follows from straightforward calculations that for the Hypercube $H_{n}$
on $n$ vertices where $n$ is a power of 2, we get 
\[
\phi(H_{n})=\dfrac{1}{2}\cdot\Delta(H_{n})
\]
showing that the Buser inequality is \emph{exactly} tight, and for
the undirected cycle $C_{n}$ on $n$ vertices, we get that 
\[
\phi(C_{n})\geq\sqrt{\Delta(C_{n})}
\]
showing that the Cheeger inequality is tight up to the constant. We proceed to prove the upper bound on $\phi(R)$ in Theorem \ref{thm:main_irred}.

\section{Fiedler's Proof -- Upper bound on $\phi(R)$ in terms of $\Delta(R)$}

Given Cheeger's inequality (Lemma \ref{lem:cheeger}), the upper bound
on $\phi(R)$ is relatively straightforward. We show the bound for
$A$ defined in Lemma \ref{lem:transRtoA}.

\begin{lemma}
\label{lem:fied} (Extension of Fiedler \cite{fiedler1995estimate}) Let $R$ be an irreducible nonnegative
matrix with positive (left and right) eigenvectors $u$ and $v$ for
eigenvalue 1, and let $A=D_{u}^{\frac{1}{2}}D_{v}^{-\frac{1}{2}}RD_{u}^{-\frac{1}{2}}D_{v}^{\frac{1}{2}}$.
Then 
\[
\phi(A)\leq\sqrt{2\cdot\Delta(A)}.
\]
\end{lemma}

\begin{proof}
Given $R$ as stated, and letting
\[
A=D_{u}^{\frac{1}{2}}D_{v}^{-\frac{1}{2}}RD_{u}^{-\frac{1}{2}}D_{v}^{\frac{1}{2}},
\]
note that 
\[
w=D_{u}^{\frac{1}{2}}D_{v}^{\frac{1}{2}}\boldsymbol{1},
\]
where we use positive square roots for the entries in the diagonal
matrices. Since $A$ has $w$ as both the left and right eigenvector
for eigenvalue 1, so does $M=\frac{A+A^{T}}{2}$.

As explained before Definition \ref{def:gen-edge-expansion} for edge
expansion of general nonnegative matrices, $D_{w}AD_{w}$ is Eulerian,
since $D_{w}AD_{w}\mathbf{1}=D_{w}Aw=D_{w}w=D_{w}^{2}\mathbf{1}=D_{w}w=D_{w}A^{T}w=D_{w}A^{T}D_{w}\mathbf{1}$.
Thus, for any $S$, 
\[
\langle\mathbf{1}_{S},D_{w}AD_{w}\mathbf{1}_{\overline{S}}\rangle=\langle\mathbf{1}_{\overline{S}},D_{w}AD_{w}\mathbf{1}_{S}\rangle=\langle\mathbf{1}_{S},D_{w}A^{T}D_{w}\mathbf{1}_{\overline{S}}\rangle,
\]
and thus for any set $S$ for which $\sum_{i\in S}w_{i}^{2}\leq\frac{1}{2}\sum_{i}w_{i}^{2}$,
\[
\phi_{S}(A)=\frac{\langle\mathbf{1}_{S},D_{w}AD_{w}\mathbf{1}_{\overline{S}}\rangle}{\langle\mathbf{1}_{S},D_{w}AD_{w}\mathbf{1}\rangle}=\frac{1}{2}\cdot\frac{\langle\mathbf{1}_{S},D_{w}AD_{w}\mathbf{1}_{\overline{S}}\rangle+\langle\mathbf{1}_{S},D_{w}A^{T}D_{w}\mathbf{1}_{\overline{S}}\rangle}{\langle\mathbf{1}_{S},D_{w}AD_{w}\mathbf{1}\rangle}=\phi_{S}(M)
\]
and thus 
\begin{equation}
\phi(A)=\phi(M).\label{eq:up2-1}
\end{equation}
For any matrix $H$, let 
\[
R_{H}(x)=\frac{\langle x,Hx\rangle}{\langle x,x\rangle}.
\]
For every $x\in\mathbb{C}^{n}$, 
\begin{equation}
\text{Re}R_{A}(x)=R_{M}(x),\label{eq:up1-1}
\end{equation}
since $A$ and $\langle x,x\rangle$ are nonnegative and we can write
\[
2\cdot R_{M}(x) 
= \frac{\langle x,Ax\rangle}{\langle x,x\rangle}+\frac{\langle x,A^{*}x\rangle}{\langle x,x\rangle}
=\frac{\langle x,Ax\rangle}{\langle x,x\rangle}+\frac{\langle Ax,x\rangle}{\langle x,x\rangle}
=\frac{\langle x,Ax\rangle}{\langle x,x\rangle}+\frac{\langle x,Ax\rangle^{*}}{\langle x,x\rangle}
=2\cdot\text{Re}R_{A}(x).
\]
Also, 
\begin{equation}
\text{Re}\lambda_{2}(A)\leq\lambda_{2}(M).\label{eq:up3-1}
\end{equation}
To see this, first note that $\lambda_{2}(A)<1$ since $A$ is irreducible.
Thus, since $\lambda_{2}(A)\not=1$, then let $v$ be the eigenvector
corresponding to $\lambda_{2}(A)$. Then since $Aw=w$, we have for $\lambda_2 = \lambda_2(A)$,
\[
Av=\lambda_{2}v\Rightarrow\langle w,Av\rangle=\langle w,\lambda_{2}v\rangle\Leftrightarrow\langle A^{T}w,v\rangle=\lambda_{2}\langle w,v\rangle\Leftrightarrow(1-\lambda_{2})\langle w,v\rangle=0
\]
which implies that $v\perp w$. Thus, we have that 
\[
\text{Re}\lambda_{2}(A)=\text{Re}\frac{\langle v,Av\rangle}{\langle v,v\rangle}=\frac{\langle v,Mv\rangle}{\langle v,v\rangle}\leq\max_{u\perp w}\frac{\langle u,Mu\rangle}{\langle u,u\rangle}=\lambda_{2}(M)
\]
where the second equality uses equation \ref{eq:up1-1}, and the last
equality follows from the variational characterization of eigenvalues
stated in the Preliminaries \ref{sec:Preliminaries}. Thus,
using equation \ref{eq:up3-1}, equation \ref{eq:up2-1} and Cheeger's
inequality for $M$ (Theorem \ref{thm:cheeger-buser}), we get 
\[
\phi(A)=\phi(M)\leq\sqrt{2\cdot(1-\lambda_{2}(M))}\leq\sqrt{2\cdot(1-\text{Re}\lambda_{2}(A))}=\sqrt{2\cdot\Delta(A)}
\]
as required.
\end{proof}

\section{Role of singular values}

In this section, we try and obtain a lower bound on $\phi$ by mimicking
the proof of Buser's inequality in Theorem \ref{thm:cheeger-buser}.
The following lemma immediately follows.

\begin{lemma}
\label{lem:phi_sing}Let $A$ be an irreducible nonnegative matrix
with largest eigenvalue 1 and $w$ the corresponding positive left
and right eigenvector. Then 
\[
\dfrac{1-\sigma_{2}(A)}{2}\leq\phi(A),
\]
where $\sigma_{2}(A)$ is the second largest singular value of $A$.
\end{lemma}

\begin{proof}
The proof immediately follows by nothing that for $v\perp w$, $\langle v,Av\rangle\leq\sigma_{2}(A)\langle v,v\rangle$
by basic linear algebra, and using lemmas \ref{lem:phi_low_bound_gamma}
gives the result.
\end{proof}
We note that Lemma \ref{lem:phi_sing} is not very meaningful in the
following sense -- let $A$ be a directed cycle, then $\sigma_{2}(A)=1$
but $\phi(A)\approx1/n$, showing the bound is not continuous. A larger
gap can be obtained by de Bruijn matrices -- discussed in Section
\ref{subsec:debruijn} -- for which $\sigma_{2}(A)=1$ but $\phi(A)\approx1/\log n$.

Another thought might be to first ensure that $A$ is lazy, or write
$A'=\frac{1}{2}A+\frac{1}{2}I$, in which case $\sigma_{2}(A')$ might
become more meaningful for $\phi(A')=\frac{1}{2}\phi(A)$. However,
noting that $A=w\cdot w^{T}+B$ and $A^{T}=w\cdot w^{T}+B^{T}$, we
have that the second largest eigenvalue of $\tilde{A}=\frac{1}{2}(A+A^{T})$
\[
\frac{1}{2}(A+A^{T})=w\cdot w^{T}+\frac{1}{2}(B+B^{T})
\]
is exactly the largest eigenvalue of $\tilde{B}=\frac{1}{2}(B+B^{T})$.
Further, since $A$ is such that $\|A\|_{2}=1$ from Lemma \ref{lem:transRtoA}
since it has the same left and right eigenvector $w$ for eigenvalue
1, it implies that $\sigma_{2}(A)=\|B\|_{2}$. Thus, let $x$ with
$\|x\|_{2}=1$ be the eigenvector for the largest eigenvalue of $\frac{1}{2}(B+B^{T})$,
then 
\[
\lambda_{2}\left(\frac{1}{2}(A+A^{T})\right)=\langle x,\frac{1}{2}(B+B^{T})x\rangle\leq\frac{1}{2}\|x\|_{2}\|Bx\|_{2}+\frac{1}{2}\|x\|_{2}\|B^{T}x\|_{2}=\|B\|_{2}=\sigma_{2}(A).
\]
Thus, using easy cheeger, we would have 
\[
\frac{1}{2}-\frac{1}{2}\sigma_{2}(A)\leq\frac{1}{2}-\frac{1}{4}\lambda_{2}(A+A^{T})\leq\phi(A),
\]
and singular values are again not meaningful since the lower bound
provided by them is weaker.

\section{Limits of perturbation theory} \label{sec:limits-perturb}

In this section, we discuss one perturbation bound to illustrate the
difficulty in proving a lower bound on $\phi(R)$ in terms of the
spectral gap using tools from perturbation theory, and also illuminate
the contrast with Theorem \ref{thm:intro-phi-delta}. The Baur-Fike
bounds state the following.

\begin{lemma}
(Baur-Fike \cite{bauer1960norms}) \label{lem:(Baur-Fike-}Let $A$
be any diagonalizable matrix with such that $A=VGV^{-1}$ is the Jordan
form for $A$ with $\mathcal{K}_{p}(A)=\|V\|_{p}\|V^{-1}\|_{p}$.
Let $E$ be a perturbation of $A$, such that $A'=A+E$. Then for
every eigenvalue $\mu$ of $A'$, there is an eigenvalue $\lambda$
of $A$ such that 
\[
|\mu-\lambda|\leq\mathcal{K}_{p}(A)\cdot\|E\|_{p}.
\]
\end{lemma}

The proof is simple but we do not rewrite it here. To contrast
with Theorem \ref{thm:intro-phi-delta}, let $S$ be the set that
achieves $\phi$ for $A$ (assume it has largest eigenvalue 1), and
let $A'$ be two disconnected components where the mass going out
and coming into the set $S$ has been removed and put back within
the sets to ensure that the largest eigenvalue of $A'$ is 1. Note
that $\lambda_{2}(A')=\lambda_{2}(A')=1=\lambda_{1}(A)$, and thus,
from the above lemma, we get that 
\[
\Delta(A)\leq\mathcal{K}_{p}(A)\cdot\|E\|_{p}.
\]
Note that even if we could have a bound where $\|E\|_{p}\approx\phi(A)$,
the dependence of $\mathcal{K}(A)$ cannot be avoided, and further,
it is still limited to diagonalizable matrices. An extension to general matrices was obtained by Saad \cite{Saad11}. 

\begin{lemma}
\label{lem:(Baur-Fike-extension)-Saad} (Baur-Fike extension by Saad \cite{Saad11}) Let $A$
be any matrix (not necessarily diagonalizable), and let $\mathcal{K}_{2}(A)=\|V\|_{2}\|V^{-1}\|_{2}$
for the Jordan decomposition $A=VGV^{-1}$. Let $E$ be a perturbation
of $A$, such that $A'=A+E$. Then for every eigenvalue $\mu$ of
$A'$, there is an eigenvalue $\lambda$ of $A$ such that
\[
1\leq\frac{1}{|\lambda-\mu|}\sum_{i=0}^{l-1}\left(\frac{1}{|\lambda-\mu|}\right)^{i}\mathcal{K}_{2}(A)\|E\|_{2}
\]
where $l$ is the size of the largest Jordan block.
\end{lemma}

We can reprove the lemma using Schur decomposition instead of the
Jordan decomposition, and we get the following.

\begin{lemma}
\label{lem:(Baur-Fike-extension)-Let} Let $A$
be any matrix (not necessarily diagonalizable) with largest singular
value $\sigma=\|A\|_{2}$. Let $E$ be a perturbation
of $A$, such that $A'=A+E$. Then for every eigenvalue $\mu$ of
$A'$, there is an eigenvalue $\lambda$ of $A$ such that
\[
1\leq\frac{1}{|\lambda-\mu|}\left(1+\frac{\sigma}{|\lambda-\mu|}\right)^{n-1}\|E\|_{2}.
\]
\end{lemma}

\begin{proof}
Let $A=UTU^{*}$ be the Schur decomposition of $A$. Let $T$ be an
upper triangular matrix with $|T_{i,i}|\geq\alpha>0$ implying that
$T$ is invertible, and $\|T\|\leq\sigma$. Let maximum entry at distance
$k$ from diagonal in $T^{-1}$ be $\gamma_{k}$. Note that for the
diagonal entries of $T^{-1}$, we have that $|T^{-1}(k,k)|\leq\frac{1}{\alpha}.$
Thus $\gamma_{0}\leq\frac{1}{\alpha}$. Assume for $k\geq1$, 
\[
\gamma_{k}\leq\frac{\sigma}{\alpha^{2}}\left(1+\frac{\sigma}{\alpha}\right)^{k-1}.
\]
For $k=1$, it is easy to verify that for any row $r$, 
\[
T^{-1}(r,r)T(r,r+1)+T^{-1}(r,r+1)T(r+1,r+1)=0
\]
and thus 
\[
|T^{-1}(r,r+1)|\leq\frac{|T^{-1}(r,r)|\cdot|T(r,r+1)|}{|T(r+1,r+1)|}\leq\frac{\frac{1}{\alpha}\cdot\sigma}{\alpha}=\frac{\sigma}{\alpha^{2}}
\]
as required. Assume the above equation holds for all $\gamma_{l}$
for $l\leq k$, we will show the equation for $\gamma_{k+1}.$ Consider
any entry $x$ in $T^{-1}$ at distance $x=k+1$ from the diagonal.
Then since $T^{-1}T=I$, we get that for some fixed row $r$, and
column $r+x$, 
\begin{align*}
\sum_{j}T^{-1}(r,j)T(j,r+x) & =0\\
\sum_{j=r}^{r+x}T^{-1}(r,j)T(j,r+x) & =0\\
\sum_{j=r}^{r+x-1}T^{-1}(r,j)T(j,r+x)+T^{-1}(r,r+x)T(r+x,r+x) & =0
\end{align*}
and thus 
\begin{align*}
|T^{-1}(r,r+x)| & \leq\frac{1}{T(r+x,r+x)}\sum_{j=r}^{r+x-1}|T^{-1}(r,j)|\cdot|T(j,r+x)|\\
 & \leq\frac{1}{\alpha}\left(\frac{\sigma}{\alpha}+\sum_{j=r+1}^{r+x-1}\gamma_{j-r}\cdot\sigma\right)\\
 & =\frac{\sigma}{\alpha}\left(\frac{1}{\alpha}+\sum_{j-r=1}^{j-r=x-1}\gamma_{j-r}\right)\\
 & =\frac{\sigma}{\alpha}\left(\frac{1}{\alpha}+\sum_{l=1}^{l=k}\gamma_{l}\right)\\
 & =\frac{\sigma}{\alpha}\left(\frac{1}{\alpha}+\sum_{l=1}^{l=k}\frac{\sigma}{\alpha^{2}}\left(1+\frac{\sigma}{\alpha}\right)^{l-1}\right)\\
 & =\frac{\sigma}{\alpha}\left(\frac{1}{\alpha}+\frac{\sigma}{\alpha^{2}}\frac{(1+\frac{\sigma}{\alpha})^{k}-1}{(1+\frac{\sigma}{\alpha})-1}\right)\\
 & =\frac{\sigma}{\alpha^{2}}\left(1+\frac{\sigma}{\alpha}\right)^{k}\\
 & =\gamma_{k+1}.
\end{align*}
Thus we have
\begin{equation}
\|T^{-1}\|_{2}\leq n\sum_{i=0}^{n-1}\gamma_{i}\leq\frac{1}{\alpha}\left(1+\frac{\sigma}{\alpha}\right)^{n-1}\label{eq:in-baur-fike-ext}
\end{equation}
and reproving Baur-Fike using Schur decomposition, we get that for
some eigenvalue $\mu$ of $A+E$ and corresponding eigenvector $v$,
we have
\[
(A+E)v=\mu v,
\]
\[
1\leq\|(A-\mu I)^{-1}\|_{2}\|E\|_{2},
\]
and thus, writing $A=UTU^{*}$, we get from equation \ref{eq:in-baur-fike-ext}
that the eigenvalue $\lambda$ of $A$ closed to $\mu$ is such that,
\[
1\leq\frac{1}{|\lambda-\mu|}\left(1+\frac{\|A\|_{2}}{|\lambda-\mu|}\right)^{n-1}\|E\|_{2}.
\]
Similarly, expanding $A=VGV^{-1}$ with the Jordan decomposition,
we get by expanding each block using 
\[
(|\lambda-\mu|I+E)^{-1}=\frac{1}{|\lambda-\mu|}\sum_{i=0}^{l-1}\left(\frac{1}{|\lambda-\mu|}\right)^{i},
\]
that 
\[
1\leq\frac{1}{|\lambda-\mu|}\sum_{i=0}^{l-1}\left(\frac{1}{|\lambda-\mu|}\right)^{i}\mathcal{K}_{2}(A)\|E\|_{2}.
\]
\end{proof}

Thus, for general matrices, there is a loss of a factor of $\exp(n)$, and a further loss of a factor of $\mathcal{K}(A)$ if Lemma \ref{lem:(Baur-Fike-extension)-Let} is used or $\sigma(A)$ if Lemma \ref{lem:(Baur-Fike-extension)-Let} is used. These give extremely weak lower bounds for $\phi$ even if we choose
$E$ appropriately. Thus we need different ideas.

\section{\label{subsec:Lower-bound-on-phi-actual}Lower bound on $\phi(R)$
in terms of $\Delta(R)$}

As a consequence of the discussion in the previous two sections, we
want to obtain a more meaningful bound than that provided by singular
values and also one that is independent of the condition numbers or
exponential factors or singular values, so that combined with Lemma
\ref{lem:fied}, it provides us with a complete generalization of
the Perron-Frobenius theorem. The lower bound will be the result of
a sequence of lemmas that we state next. Before that, we want to give
an overview of the proof. For ease of understanding, we will discuss
it for the case of double stochastic matrices -- which have the uniform
distribution as the fixed point -- and they turn out to be a sufficiently
rich special case that contains most of the properties of the general case.

Consider an irreducible doubly stochastic matrix $A$. Note that all
our spectral reasoning so far relied crucially on the spectral decomposition
of some matrix in all of our Lemmas \ref{lem:cheeger}, \ref{thm:cheeger-buser},
\ref{lem:phi_sing}. However, we cannot use it for showing a result
which has form as in Theorem \ref{thm:intro-phi-delta} for the lower
bound on $\phi(A)$ since $A$ cannot be diagonalized by a unitary,
and attempting to mimic the proof of Buser's inequality (Section \ref{subsec:Buser-Inequality-=002013})
only gives a result similar to Lemma \ref{lem:phi_sing}. As discussed
in the previous section, other theorems in literature could
be used in limited contexts -- Baur-Fike for instance in cases where
$A$ is diagonalizable, however extending them using the Schur or
Jordan decomposition to general non-diagonalizable matrices or using
other tools from perturbation theory gives lower bounds with factors
exponential in $n$, and/or a dependence on the eigenvalue/eigenvector
condition numbers of $A$, and/or a dependence on the largest singular
value of $A$, as shown in Lemma \ref{lem:(Baur-Fike-extension)-Let}
by using the Schur decomposition instead of the Jordan decomposition.
The singular values and the eigenvalue condition numbers of $A$ are
less meaningful if $A$ is doubly stochastic, but they become dominant
when we consider general matrices as required in Theorem \ref{thm:intro-phi-delta}.

The key idea that unlocks the entire proof is to try and understand
$\phi(A^{2})$ and $\phi(A^{k})$. It turns out that $\phi(A^{k})\leq k\cdot\phi(A)$.
Note that this helps us lower bound the edge expansion of $A$ by understanding
the edge expansion of $A^{k}$. Write $A=J+UTU^{*}$ where $J=\frac{1}{n}\mathbf{1}\cdot\mathbf{1}^{T}$,
and $U$ is the unitary and $T$ the upper triangular matrix in the
Schur decomposition of $A-J$. Further, it is simple to see that $A^{k}=J+UT^{k}U^{*}$.
We need to find a sufficiently small $k$ such that the expansion
of $A^{k}$ is easy to bound. To do this, we will try to ensure that
for some $k$, $\|T^{k}\|\approx0$, and thus $\phi(A^{k})\approx\phi(J)=\text{constant}$.

Note that since the eigenspace corresponding to the largest eigenvalue
1 is removed in $J$, the diagonal of $T$ contains all the eigenvalues
of $A$ with eigenvalue 1 replaced by 0. Let the eigenvalue with largest
magnitude be $\lambda_{m}$ and assume that $|\lambda_{m}|<1$. If
we raise $T^{1/1-|\lambda_{m}|}$, then it is clear that every entry
in its diagonal will be about $e^{-1}$, and the diagonal will exponentially
decrease by taking further powers. Thus the diagonal entries of $T^{\log n/(1-|\lambda_{m}|)}$
will all be inverse polynomially small, and assume at this point that
they are approximately 0. Thus, the matrix will now behave in a manner
similar to a nilpotent matrix, and it will follow that $\|T^{n\log n/(1-|\lambda_{m}|)}\|\approx0$.
Thus for $k\approx\dfrac{n\cdot\log n}{1-|\lambda_{m}|}$, we will
have 
\[
\phi(A)\geq\dfrac{1}{k}\phi(A^{k})\approx\dfrac{1}{k}\phi(J)\approx\dfrac{1-|\lambda_{m}|}{n\cdot\log n},
\]
qualitatively giving us the kind of bound that we need. There are
many technical caveats that produce a dependence on the largest singular
value of $A$ and the eigenvalue condition number of the largest eigenvalue
of $A$, but a careful analysis (with sufficient tricks) helps us
avoid both the factors, and also helps us remove the factor of $\log n$
and transform the dependence from $|\lambda_{m}|$ to $\text{Re}\lambda_{2}$.

We proceed now to prove the lower bound in Theorem \ref{thm:main_irred}.

The first lemma states that $\phi$ is sub-multiplicative in the following
sense.

\begin{lemma}
\emph{(Submultiplicativity of $\phi_{S}$)}\label{lem:(Submultiplicativity-of-phiS)}
Let $R$ and $B$ be nonnegative matrices that have the same left
and right eigenvectors $u$ and $v$ for eigenvalue 1. Then for every
cut $S$, we have that 
\[
\phi_{S}(RB)\leq\phi_{S}(R)+\phi_{S}(B).
\]
\end{lemma}

\begin{proof}
Let
\[
\gamma_{S}(R)=\langle\boldsymbol{1}_{\overline{S}},D_{u}RD_{v}\boldsymbol{1}_{S}\rangle+\langle\boldsymbol{1}_{S},D_{u}RD_{v}\boldsymbol{1}_{\overline{S}}\rangle
\]
and similarly $\gamma_{S}(B)$. Without loss of generality, we assume
that the largest eigenvalue of $R$ and $B$ is 1, since the  edge expansion 
does not change by scaling by Lemma  \ref{lem:transRtoA}. Then for every cut
$S$, we will show that 
\[
\gamma_{S}(RB)\leq\gamma_{S}(R+B)=\gamma_{S}(R)+\gamma_{S}(B).
\]
Fix any cut $S$. Assume $R=\left[\begin{array}{cc}
P & Q\\
H & V
\end{array}\right]$ and $B=\left[\begin{array}{cc}
X & Y\\
Z & W
\end{array}\right]$ naturally divided based on cut $S$. For any vector $u$, let $D_{u}$
be the diagonal matrix with $u$ on the diagonal. Since $Rv=v,R^{T}u=u,Bv=v,B^{T}u=u$,
we have
\begin{equation}
P^{T}D_{u_{S}}\boldsymbol{1}+H^{T}D_{u_{\overline{S}}}\boldsymbol{1}=D_{u_{S}}\boldsymbol{1},\label{eq:gen-pow-1}
\end{equation}
\begin{equation}
ZD_{v_{S}}\boldsymbol{1}+WD_{v_{\overline{S}}}\boldsymbol{1}=D_{v_{\overline{S}}}\boldsymbol{1},\label{eq:gen-pow-2}
\end{equation}
\begin{equation}
XD_{v_{S}}\boldsymbol{1}+YD_{v_{\overline{S}}}\boldsymbol{1}=D_{v_{S}}\boldsymbol{1},\label{eq:gen-pow-3}
\end{equation}
\begin{equation}
Q^{T}D_{u_{S}}\boldsymbol{1}+V^{T}D_{u_{\overline{S}}}\boldsymbol{1}=D_{u_{\overline{S}}}\boldsymbol{1},\label{eq:gen-pow-4}
\end{equation}
where $u$ is divided into $u_{S}$ and $u_{\overline{S}}$ and $v$
into $v_{S}$ and $v_{\overline{S}}$ naturally based on the cut $S$.
Further, in the equations above and in what follows, the vector $\mathbf{1}$
is the all 1's vector with dimension either $|S|$ or $|\overline{S}|$
which should be clear from the context of the equations, and we avoid
using different vectors to keep the notation simpler. Then we have
from the definition of $\gamma_{S}$,
\begin{align*}
\gamma_{S}(R) & =\langle\boldsymbol{1}_{\overline{S}},D_{u}RD_{v}\boldsymbol{1}_{S}\rangle+\langle\boldsymbol{1}_{S},D_{u}RD_{v}\boldsymbol{1}_{\overline{S}}\rangle\\
 & =\langle\boldsymbol{1},D_{u_{S}}QD_{v_{\overline{S}}}\boldsymbol{1}\rangle+\langle\boldsymbol{1},D_{u_{\overline{S}}}HD_{v_{S}}\boldsymbol{1}\rangle
\end{align*}
and similarly 
\[
\gamma_{S}(B)=\langle\boldsymbol{1},D_{u_{S}}YD_{v_{\overline{S}}}\boldsymbol{1}\rangle+\langle\boldsymbol{1},D_{u_{\overline{S}}}ZD_{v_{S}}\boldsymbol{1}\rangle.
\]
The matrix $RB$ also has $u$ and $v$ as the left and right eigenvectors
for eigenvalue 1 respectively, and thus,
\begin{align*}
\gamma_{S}(RB) & =\langle\boldsymbol{1},D_{u_{S}}PYD_{v_{\overline{S}}}\boldsymbol{1}\rangle+\langle\boldsymbol{1},D_{u_{S}}QWD_{v_{\overline{S}}}\boldsymbol{1}\rangle+\langle\boldsymbol{1},D_{u_{\overline{S}}}HXD_{v_{S}}\boldsymbol{1}\rangle+\langle\boldsymbol{1},D_{u_{\overline{S}}}VZD_{v_{S}}\boldsymbol{1}\rangle\\
 & =\langle P^{T}D_{u_{S}}\boldsymbol{1},YD_{v_{\overline{S}}}\boldsymbol{1}\rangle+\langle\boldsymbol{1},D_{u_{S}}QWD_{v_{\overline{S}}}\boldsymbol{1}\rangle+\langle\boldsymbol{1},D_{u_{\overline{S}}}HXD_{v_{S}}\boldsymbol{1}\rangle+\langle V^{T}D_{u_{\overline{S}}}\boldsymbol{1},ZD_{v_{S}}\boldsymbol{1}\rangle\\
 & =\langle D_{u_{S}}\boldsymbol{1}-H^{T}D_{u_{\overline{S}}}\boldsymbol{1},YD_{v_{\overline{S}}}\boldsymbol{1}\rangle+\langle\boldsymbol{1},D_{u_{S}}Q(D_{v_{\overline{S}}}\boldsymbol{1}-ZD_{v_{S}}\boldsymbol{1})\rangle\\
 & \ \ \ \ \ +\langle\boldsymbol{1},D_{u_{\overline{S}}}H(D_{v_{S}}\boldsymbol{1}-YD_{v_{\overline{S}}}\boldsymbol{1}\rangle+\langle D_{u_{\overline{S}}}\boldsymbol{1}-Q^{T}D_{u_{S}}\boldsymbol{1},ZD_{v_{S}}\boldsymbol{1}\rangle\\
 & \,\,\,\,\,\,\,\,\,\,\,\,\,\,\,\,\,\,\,\,\,\,\text{[from equations \ref{eq:gen-pow-1}, \ref{eq:gen-pow-2}, \ref{eq:gen-pow-3}, \ref{eq:gen-pow-4} above]}\\
 & =\langle D_{u_{S}}\boldsymbol{1},YD_{v_{\overline{S}}}\boldsymbol{1}\rangle+\langle\boldsymbol{1},D_{u_{S}}QD_{v_{\overline{S}}}\boldsymbol{1}\rangle-\langle H^{T}D_{u_{\overline{S}}}\boldsymbol{1},YD_{v_{\overline{S}}}\boldsymbol{1}\rangle-\langle\boldsymbol{1},D_{u_{S}}QZD_{v_{S}}\boldsymbol{1}\rangle\\
 & \ \ \ \ \ +\langle\boldsymbol{1},D_{u_{\overline{S}}}HD_{v_{S}}\boldsymbol{1}\rangle+\langle D_{u_{\overline{S}}}\boldsymbol{1},ZD_{v_{S}}\boldsymbol{1}\rangle-\langle\boldsymbol{1},D_{u_{\overline{S}}}HYD_{v_{\overline{S}}}\boldsymbol{1}\rangle-\langle Q^{T}D_{u_{S}}\boldsymbol{1},ZD_{v_{S}}\boldsymbol{1}\rangle\\
 & \leq\langle D_{u_{S}}\boldsymbol{1},YD_{v_{\overline{S}}}\boldsymbol{1}\rangle+\langle\boldsymbol{1},D_{u_{S}}QD_{v_{\overline{S}}}\boldsymbol{1}\rangle+\langle\boldsymbol{1},D_{u_{\overline{S}}}HD_{v_{S}}\boldsymbol{1}\rangle+\langle D_{u_{\overline{S}}}\boldsymbol{1},ZD_{v_{S}}\boldsymbol{1}\rangle\\
 & \,\,\,\,\,\,\,\,\,\,\,\,\,\,\,\,\,\,\,\,\,\,\text{[since every entry of the matrices is nonnegative]}\\
 & =\gamma_{S}(R)+\gamma_{S}(B)
\end{align*}
as required. Note that since $R$ and $B$ have the same left and
right eigenvectors, so does $RB$, and since $\langle\boldsymbol{1}_{\overline{S}},D_{u}RD_{v}\boldsymbol{1}_{S}\rangle=\langle\boldsymbol{1}_{S},D_{u}RD_{v}\boldsymbol{1}_{\overline{S}}\rangle$,
we exactly have that 
\[
\phi_{S}(R)=\frac{1}{2}\dfrac{\gamma_{S}(R)}{\sum_{i\in S}u_{i}v_{i}}
\]
which gives the lemma.
\end{proof}
The following lemma follows directly as a corollary of Lemma \ref{lem:(Submultiplicativity-of-phiS)},
and is essentially the heart of the entire proof.

\begin{lemma}
\emph{\label{lem:phi_pow_bound} (Submultiplicativity of $\phi$)} Let $R\in\mathbb{R}^{n\times n}$
be an irreducible nonnegative matrix with left and right eigenvectors
$u$ and $v$ for the PF eigenvalue 1. Then 
\[
\phi(R^{k})\leq k\cdot\phi(R).
\]
\end{lemma}

\begin{proof}
Noting that $R^{k}$ has $u$ and $v$ as left and right eigenvectors
for any $k$, we inductively get using Lemma \ref{lem:(Submultiplicativity-of-phiS)}
that 
\[
\phi_{S}(R^{k})\leq\phi_{S}(R)+\phi_{S}(R^{k-1})\leq\phi_{S}(R)+(k-1)\cdot\phi_{S}(R)=k\cdot\phi_{S}(R),
\]
and we get by letting $S$ be the set that minimizes $\phi(R)$, that
\[
\phi(R^{k})\leq k\cdot\phi(R).
\]
\end{proof}
For the case of \emph{symmetric} doubly stochastic matrices $R$,
Lemma \ref{lem:phi_pow_bound} follows from a theorem of Blakley and
Roy \cite{BlakleyR65} (it does not fall into the framework of an
extension of that result to the nonsymmetric case \cite{Pate12}).
Lemma \ref{lem:phi_pow_bound} helps to lower bound $\phi(R)$ by
taking powers of $R$, which is useful since we can take sufficient
powers in order to make the matrix simple enough that its edge expansion
is easily calculated.

The next two lemmas follow by technical calculations.

\begin{lemma}
\label{lem:T_norm_bound} (Bounded norm of powers) Let $T\in\mathbb{C}^{n\times n}$
be an upper triangular matrix with $\|T\|{}_{2}=\sigma$ and for every
$i$, $|T_{i,i}|\leq\alpha<1$. Then 
\[
\|T^{k}\|_{2}\leq n\cdot\sigma^{n}\cdot{k+n \choose n}\cdot\alpha^{k-n}.
\]
\end{lemma}

\begin{proof}
Let $g_{r}(k)$ denote the maximum of the absolute value of entries
at distance $r$ from the diagonal in $T^{k}$, where the diagonal
is at distance 0 from the diagonal, the off-diagonal is at distance
1 from the diagonal and so on. More formally, 
\[
g_{r}(k)=\max_{i}|T^{k}(i,i+r)|.
\]
We will inductively show that for $\alpha\leq1$, and $r\geq1$, 
\begin{equation}
g_{r}(k)\leq{k+r \choose r}\cdot\alpha^{k-r}\cdot\sigma^{r},\label{eq:lo2-1-1}
\end{equation}
where $\sigma=\|T\|_{2}$. First note that for $r=0$, since $T$
is upper triangular, the diagonal of $T^{k}$ is $\alpha^{k}$, and
thus the hypothesis holds for $r=0$ and all $k\geq1$. Further, for
$k=1$, if $r=0$, then $g_{0}(1)\leq\alpha$ and if $r\geq1$, then
$g_{r}(1)\leq\|T\|_{2}\leq\sigma$ and the inductive hypothesis holds
also in this case, since $r\geq k$ and $\alpha^{k-r}\geq1$. For
the inductive step, assume that for all $r\geq1$ and all $j\leq k-1$,
$g_{r}(j)\leq{j+r \choose r}\cdot\alpha^{j-r}\cdot\sigma^{r}$. We
will show the calculation for $g_{r}(k)$.

Since $|a+b|\leq|a|+|b|$, 
\begin{align*}
g_{r}(k) & \leq\sum_{i=0}^{r}g_{r-i}(1)\cdot g_{i}(k-1)\\
 & =g_{0}(1)\cdot g_{r}(k-1)+\sum_{i=0}^{r-1}g_{r-i}(1)\cdot g_{i}(k-1).
\end{align*}
The first term can be written as, 
\begin{align}
g_{0}(1)\cdot g_{r}(k-1) & =\alpha\cdot{k-1+r \choose r}\cdot\alpha^{k-1-r}\cdot\sigma^{r}\nonumber \\
 & \ \ \ \ \ \ \ \ \text{[using that \ensuremath{g_{0}(1)\leq\alpha} and the inductive hypothesis for the second term]}\nonumber \\
 & \leq\alpha^{k-r}\cdot\sigma^{r}\cdot{k+r \choose r}\cdot\dfrac{{k-1+r \choose r}}{{k+r \choose r}}\nonumber \\
 & \leq\frac{k}{k+r}\cdot{k+r \choose r}\cdot\alpha^{k-r}\cdot\sigma^{r}\label{eq:inproofeq-1}
\end{align}
and the second term as 
\begin{align*}
g_{r}(k) & \leq\sum_{i=0}^{r-1}g_{r-i}(1)\cdot g_{i}(k-1)\\
 & \leq\sigma\cdot\sum_{i=0}^{r-1}{k-1+i \choose i}\cdot\alpha^{k-1-i}\cdot\sigma^{i}\\
 & \ \ \ \ \ \text{[using \ensuremath{g_{r-i}(1)\leq\sigma} and the inductive hypothesis for the second term]}\\
 & \le\sigma^{r}\cdot\alpha^{k-r}\cdot{k+r \choose r}\cdot\sum_{i=0}^{r-1}\cdot\dfrac{{k-1+i \choose i}}{{k+r \choose r}}\cdot\alpha^{r-1-i}\\
 & \ \ \ \ \ \text{[using \ensuremath{\sigma^{i}\leq\sigma^{r-1}} since \ensuremath{\sigma\geq1} and \ensuremath{i\leq r-1}]}\\
 & \leq\sigma^{r}\cdot\alpha^{k-r}\cdot{k+r \choose r}\cdot\sum_{i=0}^{r-1}\cdot\dfrac{{k-1+i \choose i}}{{k+r \choose r}}\\
 & \ \ \ \ \ \text{[using \ensuremath{\alpha\leq1} and \ensuremath{i\leq r-1}]}
\end{align*}
We will now show that the quantity inside the summation is at most
$\frac{r}{k+r}$. Inductively, for $r=1$, the statement is true,
and assume that for any other $r$, 
\[
\sum_{i=0}^{r-1}\frac{{k-1+i \choose i}}{{k+r \choose r}}\leq\frac{r}{k+r}.
\]
Then we have 
\begin{align*}
\sum_{i=0}^{r}\frac{{k-1+i \choose i}}{{k+r+1 \choose r+1}} & =\frac{{k+r \choose r}}{{k+r+1 \choose r+1}}\cdot\sum_{i=0}^{r-1}\frac{{k-1+i \choose i}}{{k+r \choose r}}+\frac{{k-1+r \choose k-1}}{{k+r+1 \choose k}}\\
 & \leq\frac{r+1}{k+r+1}\cdot\frac{r}{k+r}+\frac{(r+1)\cdot k}{(k+r+1)\cdot(k+r)}\\
 & =\frac{r+1}{k+r+1}
\end{align*}
Thus we get that the second term is at most $\sigma^{r}\cdot\alpha^{k-r}\cdot{k+r \choose r}\cdot\frac{r}{k+r}$,
and combining it with the first term (equation \ref{eq:inproofeq-1}),
it completes the inductive hypothesis.

Noting that the operator norm is at most the Frobenius norm, and since
$g_{r}(k)$ is increasing in $r$ and the maximum value of $r$ is
$n$, we get using equation \ref{eq:lo2-1-1}, 
\begin{align*}
\|T^{k}\|_{2} & \leq\sqrt{\sum_{i,j}|T^{k}(i,j)|^{2}}\\
 & \leq n\cdot\sigma^{n}\cdot\alpha^{k-n}\cdot{k+n \choose n}
\end{align*}
as required.
\end{proof}
Using Lemma \ref{lem:T_norm_bound}, we can show the following lemma
for the special case of upper triangular matrices with operator norm
at most 1.

\begin{lemma}
\label{lem:T_norm_1_bound} Let $T\in\mathbb{C}^{n\times n}$ be an
upper triangular matrix with $\|T\|_{2}\leq1$ and $|T_{i,i}|\leq\alpha<1$
for every $i$. Assume $T^{k}$ is well defined for every real $k$.
Then $\|T^{k}\|\leq\epsilon$ for 
\[
k\geq\dfrac{4n+2\ln(\frac{n}{\epsilon})}{\ln\left(\dfrac{1}{\alpha}\right)}.
\]
\end{lemma}

\begin{proof}
Let $X=T^{k_{1}}$, where $k_{1}=\frac{c_{1}}{\ln(\frac{1}{\alpha})}$
for $\alpha<1$. Then $\|X\|_{2}\leq\|T\|_{2}^{k_{1}}\leq1$, and
for every $i$, 
\[
|X_{i,i}|\leq|T_{i,i}|^{k_{1}}\leq|\lambda_{m}|^{c_{1}/\ln(\frac{1}{\alpha})}=e^{-c_{1}}.
\]
Using Lemma \ref{lem:T_norm_bound} for $X$ with $\sigma=1$ and
$\beta=e^{-c_{1}}$, we get that for $k_{2}=c_{2}\cdot n$, 
\begin{align*}
\|X^{k_{2}}\| & \leq n\cdot{k_{2}+n \choose n}\cdot e^{-c_{1}(k_{2}-n)}\\
 & \leq n\cdot e^{n}\cdot(c_{2}+1)^{n}\cdot e^{-c_{1}(c_{2}-1)n}\\
 & \ \ \ \ \ \ \ \text{[using \ensuremath{{a \choose b}\leq\left(\frac{ea}{b}\right)^{b}}]}\\
 & =\exp\left(n\cdot\left(\frac{\ln n}{n}+\ln(c_{2}+1)+1+c_{1}-c_{1}c_{2}\right)\right)
\end{align*}
and to have this quantity less than $\epsilon$, we require 
\begin{align}
 & \exp\left(n\cdot\left(\frac{\ln n}{n}+\ln(c_{2}+1)+1+c_{1}-c_{1}c_{2}\right)\right)\leq\epsilon\nonumber \\
\Leftrightarrow & \left(\frac{1}{\epsilon}\right)^{\frac{1}{n}}\leq\exp\left(-1\cdot\left(\frac{\ln n}{n}+\ln(c_{2}+1)+1+c_{1}-c_{1}c_{2}\right)\right)\nonumber \\
\Leftrightarrow & \frac{1}{n}\ln\frac{n}{\epsilon}+1+c_{1}+\ln(c_{2}+1)\leq c_{1}c_{2}\label{eq:proofmodbound-1}
\end{align}
and we set 
\[
c_{1}=1+\frac{1}{2.51}\cdot\frac{1}{n}\ln\left(\frac{n}{\epsilon}\right)
\]
and $c_{2}=3.51$ which always satisfies inequality \ref{eq:proofmodbound-1}.
As a consequence, for 
\begin{align*}
k & =k_{1}\cdot k_{2}\\
 & =\frac{c_{1}\cdot c_{2}\cdot n}{\ln(\frac{1}{\alpha})}\\
 & =\dfrac{3.51\cdot n+1.385\cdot\ln\left(\frac{n}{\epsilon}\right)}{\ln(\frac{1}{\alpha})}
\end{align*}
we get, 
\[
\|T^{k}\|_{2}\leq\|X^{k_{2}}\|\leq\epsilon
\]
as required.
\end{proof}
Given lemmas \ref{lem:phi_pow_bound} and \ref{lem:T_norm_1_bound},
we can lower bound $\phi(R)$ in terms of $\ln(\frac{1}{|\lambda_{m}|})$
(where $\lambda_{m}$ is the nontrivial eigenvalue that is maximum
in magnitude). Our aim is to lower bound $\phi(R)$ by $\phi(R^{k})$,
but since the norm of $R^{k}$ increases by powering, we cannot use
the lemmas directly, since we do not want a dependence on $\sigma(R)$
in the final bound. To handle this, we transform $R$ to $A$, such
that $\phi(R)=\phi(A)$, the eigenvalues of $R$ and $A$ are the
same, but $\sigma(A)=\|A\|_{2}=1$ irrespective of the norm of $R$,
using Lemma \ref{lem:transRtoA}. 

Given Lemma \ref{lem:transRtoA}, we lower bound $\phi(A)$ using
$\phi(A^{k})$ in terms of $1-|\lambda_{m}(A)|$, to obtain the corresponding
bounds for $R$.

\begin{lemma}
\label{lem:mod_lambda_bound} Let $R$ be an irreducible nonnegative
matrix with positive (left and right) eigenvectors $u$ and $v$ for
the PF eigenvalue 1, normalized so that $\langle u,v\rangle=1$. Let $\lambda_{m}$ be the nontrivial eigenvalue of $R$ that is maximum
in magnitude, and let $R^{k}$ be well-defined for every real $k$.
Then 
\[
\frac{1}{15}\cdot\frac{\ln\left(\dfrac{1}{|\lambda_{m}|}\right)}{n}\leq\phi(R).
\]
\end{lemma}

\begin{proof}
From the given $R$ with positive left and right eigenvectors $u$
and $v$ for eigenvalue 1 as stated, let $A=D_{u}^{\frac{1}{2}}D_{v}^{-\frac{1}{2}}RD_{u}^{-\frac{1}{2}}D_{v}^{\frac{1}{2}}$
and $w=D_{u}^{\frac{1}{2}}D_{v}^{\frac{1}{2}}\mathbf{1}$ as in Lemma
\ref{lem:transRtoA}. Note that $w$ is positive, and 
\[
\langle w,w\rangle=\langle D_{u}^{\frac{1}{2}}D_{v}^{\frac{1}{2}}\mathbf{1},D_{u}^{\frac{1}{2}}D_{v}^{\frac{1}{2}}\mathbf{1}\rangle=\langle u,v\rangle=1.
\]
Further, $Aw=w$ and $A^{T}w=w$, and since $A^{k}=D_{u}^{\frac{1}{2}}D_{v}^{-\frac{1}{2}}R^{k}D_{u}^{-\frac{1}{2}}D_{v}^{\frac{1}{2}}$
it is well-defined for every real $k$.

Let $A=w\cdot w^{T}+B$. Since $(w\cdot w^{T})^{2}=w\cdot w^{T}$
and $Bw\cdot w^{T}=w\cdot w^{T}B=0$, we get 
\[
A^{k}=w\cdot w^{T}+B^{k}.
\]
Let $B=UTU^{*}$ be the Schur decomposition of $B$, where the diagonal
of $T$ contains all but the stochastic eigenvalue of $A$, which
is replaced by 0, since $w$ is both the left and right eigenvector
for eigenvalue 1 of $A$, and that space is removed in $w\cdot w^{T}$.
Further, the maximum diagonal entry of $T$ is at most $|\lambda_{m}|$
where $\lambda_{m}$ is the nontrivial eigenvalue of $A$ (or $R$)
that is maximum in magnitude. Note that $|\lambda_{m}|<1$ since $A$
(and $R$) is irreducible. Since $w\cdot w^{T}B=Bw\cdot w^{T}=0$
and $\|A\|_{2}\leq1$ from Lemma \ref{lem:transRtoA}, we have that
$\|B\|_{2}\leq1$.

Thus, using Lemma \ref{lem:T_norm_1_bound} (in fact, the last lines
in the proof of Lemma \ref{lem:T_norm_1_bound}), for 
\[
k\geq\dfrac{3.51\cdot n+1.385\cdot\ln\left(\frac{n}{\epsilon}\right)}{\ln\left(\dfrac{1}{|\lambda_{m}|}\right)},
\]
we get that 
\[
\|B^{k}\|_{2}=\|T^{k}\|_{2}\leq\epsilon,
\]
and for $e_{i}$ being the vector with 1 at position $i$ and zeros
elsewhere, we get using Cauchy-Schwarz 
\[
|B^{k}(i,j)|=\left|\langle e_{i},B^{k}e_{j}\rangle\right|\leq\|e_{i}\|_{2}\|B\|_{2}\|e_{j}\|_{2}\leq\epsilon.
\]
Further, note that $B^{k}w=0$ and $w^{T}B^{k}=0$. This means that
for any $i$ and $S\subseteq[n]$, 
\begin{align*}
\sum_{j\in S}B^{k}(i,j)w_{j} & =-\sum_{j\in\overline{S}}B^{k}(i,j)w_{j}
\end{align*}
Thus, for any set $S$ for which $\sum_{i\in S}w_{i}^{2}\leq\frac{1}{2}$,
we get 
\begin{align*}
\phi_{S}(A^{k}) & =\frac{\langle\mathbf{1}_{S},D_{w}A^{k}D_{w}\mathbf{1}_{\overline{S}}\rangle}{\langle\mathbf{1}_{S},D_{w}D_{w}\mathbf{1}\rangle}\\
 & =\frac{\sum_{i\in S,j\in\overline{S}}A^{k}(i,j)\cdot w_{i}\cdot w_{j}}{\sum_{i\in S}w_{i}^{2}}\\
 & =\frac{\sum_{i\in S,j\in\overline{S}}\left(w\cdot w^{T}(i,j)+B^{k}(i,j)\right)\cdot w_{i}\cdot w_{j}}{\sum_{i\in S}w_{i}^{2}}\\
 & =\dfrac{\sum_{i\in S}w_{i}^{2}\sum_{j\in\overline{S}}w_{j}^{2}-\sum_{i\in S,j\in S}B^{k}(i,j)\cdot w_{i}\cdot w_{j}}{\sum_{i\in S}w_{i}^{2}}\\
 & \geq\sum_{j\in\overline{S}}w_{j}^{2}-\epsilon\dfrac{(\sum w_{i}^{2})^{2}}{\sum_{i\in S}w_{i}^{2}}\\
 & \geq\dfrac{1}{2}(1-\epsilon)
\end{align*}
since $\sum_{i\in S}w_{i}^{2}\leq\frac{1}{2}.$ Note that this holds
for \emph{every} set $S$. Thus, we get that for 
\[
k\geq\dfrac{3.51n+1.385\cdot\ln\left(\frac{n}{\epsilon}\right)}{\ln\left(\dfrac{1}{|\lambda_{m}|}\right)}
\]
the edge expansion 
\[
\phi(A^{k})\geq\frac{1}{2}(1-\epsilon),
\]
and thus using Lemma \ref{lem:phi_pow_bound}, we get that 
\begin{align*}
\phi(A) & \geq\frac{1}{k}\cdot\phi(A^{k})\geq\frac{1-\epsilon}{2}\cdot\frac{\ln\left(\dfrac{1}{|\lambda_{m}|}\right)}{3.51\cdot n+1.385\cdot\ln\left(n/\epsilon\right)}
\end{align*}
and setting $\epsilon=\frac{1}{2\cdot e}$, and using $\ln(e\cdot x)\leq x,$we
get that
\begin{align*}
\phi(A) & \geq\frac{1}{15}\cdot\frac{\ln\left(\dfrac{1}{|\lambda_{m}|}\right)}{n},
\end{align*}
 and it carries to $R$ due to Lemma \ref{lem:transRtoA}.
\end{proof}
We almost have the proof that we seek, but we need to take care of
two issues. We need a dependence on $\text{Re}\lambda_{2}$ instead
of $|\lambda_{m}|$, and we need to ensure that we can take real powers
of the underlying matrix on which the lemmas are applied. To solve
both the issues, we will now focus on the exponential version of $R$,
or $E_{R}=\exp(R-I)$. We note the following properties of this matrix.

\begin{lemma}
\label{lem:expt(I-A)bound}Let $E_{R}$ be the exponential version
of $R$, i.e. $E_{R}=\exp(R-I)$. Then the following hold for $E_{R}$.
\begin{enumerate}
\item $E_{A}=D_{u}^{\frac{1}{2}}D_{v}^{-\frac{1}{2}}E_{R}D_{u}^{-\frac{1}{2}}D_{v}^{\frac{1}{2}}$
for $A=D_{u}^{\frac{1}{2}}D_{v}^{-\frac{1}{2}}RD_{u}^{-\frac{1}{2}}D_{v}^{\frac{1}{2}}$
\item The largest eigenvalue of $E_{R}$ is 1, the left and right eigenvectors
of $E_{R}$ and $R$ are same, and $E_{R}$ is irreducible.
\item $(E_{R})^{k}$ is well-defined for all real $k$
\item $\phi(E_{R}^{t})\leq t\cdot\phi(R)$
\end{enumerate}
\end{lemma}

\begin{proof}
The first three properties are straightforward and follow from the
definitions. For the last property, we show for $A$ and it will be
implied for $R$ by Lemma \ref{lem:transRtoA}. For any set $S$,
we can write 
\begin{align*}
\phi_{S}(E_{A}^{t}) & =\phi_{S}(\exp(t(A-I)))\\
 & =e^{-t}\phi_{S}(\sum_{i=0}\dfrac{1}{i!}t^{i}A^{i})\\
 & =e^{-t}\left(\phi_{S}(I)+\sum_{i=1}^{\infty}\dfrac{t^{i}\cdot\phi_{S}(A^{i})}{i!}\right)\\
 & \leq e^{-t}\left(0+\sum_{i=1}^{\infty}i\dfrac{t^{i}\cdot\phi_{S}(A)}{i!}\right)\\
 & \ \ \ [\text{using Lemma \ref{lem:phi_pow_bound}}]\\
 & =e^{-1}\left(t\cdot\phi_{S}(A)\sum_{i=1}^{\infty}\dfrac{t^{i-1}}{(i-1)!}\right)\\
 & =t\cdot\phi_{S}(A)
\end{align*}
and thus 
\[
\phi(E_{A})\leq t\cdot\phi(A)
\]
as required, and it extends to $R$ by Lemma \ref{lem:transRtoA}.
\end{proof}
Given Lemma \ref{lem:mod_lambda_bound}, we use the exponential version
of $R$ to obtain the lower bound on edge expansion in Theorem \ref{thm:main_irred},
restated below.

\begin{lemma}
\label{lem:real-lambda-bound} Let $R$ be an irreducible nonnegative
matrix with positive (left and right) eigenvectors $u$ and $v$ for
the PF eigenvalue 1, normalized so that $\langle u,v\rangle=1$. Then
\[
\frac{1}{15}\cdot\dfrac{1-\text{Re}\lambda_{2}(R)}{n}\leq\phi(R).
\]
\end{lemma}

\begin{proof}
We show the theorem for $A$ and it extends to $R$ by Lemma \ref{lem:transRtoA}.
Note that for every eigenvalue $\lambda=a+i\cdot b$ of $A$, the
corresponding eigenvalue of $E_{A}$ is 
\begin{align*}
\exp(\lambda-1) & =\exp(a-1)\exp(i\cdot b)\\
|\exp(\lambda-1)| & =\exp(a-1).
\end{align*}
Thus, to maximize the magnitude $|\exp(\lambda-1)|$ of the eigenvalue
of $E_{A}$, we need to maximize $\exp(a-1)$, and since $a<1$ as
$A$ is irreducible, the value $|\exp(\lambda-1)|$ is maximized when
$a$ is closest to 1. In other words, letting $\lambda_{m}(E_{A})$
be the eigenvalue of $E_{A}$ that is largest in magnitude, we have
that 
\[
|\lambda_{m}(E_{A})|=\exp(\text{Re}\lambda_{2}(A)-1).
\]

Applying Lemma \ref{lem:mod_lambda_bound} on $E_{A}$, we get that
\[
\phi(E_{A})\geq\frac{1}{15}\cdot\frac{\ln\left(\dfrac{1}{|\lambda_{m}(E_{A})|}\right)}{n}=\frac{1}{15}\cdot\frac{1-\text{Re}\lambda_{2}(A)}{n},
\]
and from Lemma \ref{lem:expt(I-A)bound}, we finally get 
\[
\phi(A)\geq\frac{1}{15}\cdot\frac{1-\text{Re}\lambda_{2}(A)}{n}
\]
as required, and it extends to $R$ through Lemma \ref{lem:transRtoA}.
\end{proof}
This completes the proof of our first main theorem, the lower bound
on $\phi$ for irreducible nonnegative matrices. A simple extension
to all matrices follows after extending the definition of the spectral
gap and  edge expansion.

\section{\label{sec:extension-to-gen-mats}Extension of Theorem \ref{thm:main_irred} to all nonnegative matrices}

The following simple lemma follows from the definition of $\lim\inf$.

\begin{lemma}
\label{lem:liminf}Let $p:E\rightarrow\mathbb{R}$ and $q:E\rightarrow\mathbb{R}$
be two functions such that for all $H\in E$, $p(H)\leq q(H)$. Then
for all $R$, 
\[
\lim_{\epsilon\rightarrow0}\inf_{H\in R_{\epsilon}}p(H)\leq\lim_{\epsilon\rightarrow0}\inf_{H\in R_{\epsilon}}q(H).
\]
\end{lemma}

\begin{proof}
Note that if $R\in E$, it holds trivially. Assume $R\not\in E$.
Fix $\epsilon$. Let $a_{\epsilon}=\inf_{H\in R_{\epsilon}}p(H)$
and $b_{\epsilon}=\inf_{H\in R_{\epsilon}}q(H)$. First we claim that
$a_{\epsilon}\leq b_{\epsilon}$. For the sake of contradiction, assume
\begin{equation}
a_{\epsilon}>b_{\epsilon}.\label{eq:assm1}
\end{equation}
Let $\delta_{\epsilon}=(a_{\epsilon}-b_{\epsilon})/2>0$, then there
exists $T\in E$ such that $q(T)\leq b_{\epsilon}+\delta_{\epsilon}$
(else $\inf_{H\in R_{\epsilon}}q(H)\geq b_{\epsilon}+\delta_{\epsilon}>b_{\epsilon}$).
Further, 
\[
a_{\epsilon}\leq p(T)\leq q(T)\leq b_{\epsilon}+\delta_{\epsilon}
\]
or $a_{\epsilon}\leq b_{\epsilon}$, a contradiction to \ref{eq:assm1}.
Thus, $a_{\epsilon}\leq b_{\epsilon}$.

Now let $c_{\epsilon}=a_{\epsilon}-b_{\epsilon}$, and we know that
for all $\epsilon>0$, $c_{\epsilon}\leq0$. It follows that 
\begin{equation}
\lim_{\epsilon\rightarrow0}c_{\epsilon}\leq0\label{eq:assm2}
\end{equation}
Again for the sake of contradiction assume the contrary, i.e $\lim_{\epsilon\rightarrow0}c_{\epsilon}>0$
and $\delta=\frac{1}{2}\cdot\lim_{\epsilon\rightarrow0}c_{\epsilon}>0$.
By the definition of limit, there exists $\epsilon>0$ such that $c_{\epsilon}>\delta$,
or $c_{\epsilon}>0$, a contradiction to \ref{eq:assm2}. Thus the
lemma follows.
\end{proof}
\begin{claim}
The main Theorem \ref{thm:main_irred} holds for all nonnegative matrices
$R$.
\end{claim}

\begin{proof}
Follows from Lemma \ref{lem:liminf} after setting $p$ and $q$ appropriately.
\end{proof}
This concludes the proof of our main Theorem \ref{thm:intro-phi-delta}.

\chapter{Constructions of Nonreversible chains}

\vspace{1.5cm}
\begin{quote}
We have now perceived, that all the explanations commonly given
of nature are mere modes of imagining, and do not indicate the true
nature of anything, but only the constitution of the imagination.
I do not attribute to nature either beauty or deformity, order or
confusion. Only in relation to our imagination can things be called
beautiful or ugly, well-ordered or confused.
\begin{flushright}
\textasciitilde{} Baruch Spinoza, \emph{Ethics}
\par\end{flushright}

\end{quote}
\vspace{3cm}

The starting point of this section is to understand the optimality
of Theorem \ref{thm:intro-phi-delta}. The main difference between
Theorem \ref{thm:intro-phi-delta} and the Cheeger-Buser inequality (Theorem 
\ref{thm:cheeger-buser}) is the loss of a factor of $n$ in the nonreversible
case, and we want to understand whether this loss is indeed necessary
or whether it is a relic of the limitations of our proof techniques.
We will also seek constructions of \emph{doubly stochastic} matrices
with these properties, to test the optimality of our theorem even
in the case of a uniform principal eigenvector for the nonnegative
matrix, where one might expect a polylog$(n)$ loss instead of a loss
of $n$ in Theorem \ref{thm:intro-phi-delta}. We remark that we explain
in detail in Sections \ref{subsec:Constraints-on-the}, \ref{subsec:Special-Cases-of},
\ref{subsec:The-general-case} and \ref{subsec:Observations-from-the}
how we arrive at our constructions in Sections \ref{subsec:Rootn-matrices-=002013}
and \ref{subsec:Chet-Matrices-=002013} so that they do not seem mysterious.
Although brevity is the hallmark of wit, we think comprehensiveness
is the hallmark of understanding, and we'll adhere to the latter for
the most part, and as a consequence detail the thought process to
arrive at our constructions, which might be instructive in a search
for other similar constructions.

Before we proceed towards constructions, we need to understand the
\emph{non-expansion} of graphs.

\section{\label{subsec:Non-expansion}Non-expansion}

The quantity edge expansion $\phi$ as defined in \ref{def:gen-edge-expansion}
has been extensively studied in the last 70 years within combinatorics
and spectral theory, albeit mostly for reversible matrices, and undirected
graphs (or symmetric matrices) with a few edges (or few nonzero entries in the matrix) and constant expansion -- also called
(combinatorial) expanders -- have many remarkable, almost magical
properties (see \cite{hoory2006expander}), making them a fundamental
combinatorial object from which many other optimal (up to lower order terms)
pseudorandom objects can be constructed, such as error-correcting
codes and pseudorandom generators amongst others, and further they
serve as a building block in a large number of constructions within
mathematics.

What we seek now is a direction opposite to that of expansion of graphs
-- the \emph{non-expansion} of graphs. This is uninteresting for
undirected graphs, since non-expansion of symmetric matrices simply
implies a small spectral gap due to the Cheeger-Buser inequality \ref{thm:cheeger-buser},
and it is elementary to construct such graphs. However, for the case
of directed graphs (or irreversible matrices), the property of non-expansion
is non-trivial, since at the outset, it is possible for a nonreversible
matrix to have constant spectral gap but edge expansion that diminishes
with the matrix size.

All the constructions that we present in this section are of non-expanding
irreversible matrices, and although they go furthest from the remarkable
and highly sought-after property of constant expansion, they are truly
beautiful in a different and unique manner and deserve study in their
own right. Our aim of this section is two-fold -- to analyze and
give an exposition of a few interesting known constructions, and present
new constructions that are \emph{exponentially} better than known constructions with regards
to the lower bound in Theorem \ref{thm:intro-phi-delta}, and with
regards to the non-expansion of graphs as defined next. In fact, even
the sub-optimal (in the sense of the lower bound in Theorem \ref{thm:intro-phi-delta})
constructions will have many interesting and aesthetically pleasing
properties. To systematize exposition, we define the following quantity.
\medskip

\begin{definition}
\label{def:nonexpansion} \emph{(Non-expansion of matrices)} Let $E_{n}$ be the set of all $n\times n$
doubly stochastic matrices. For any $A\in E_{n}$, define the non-expansion
of $A$ as

\[
\Gamma(A)=\dfrac{\phi(A)}{1-\text{Re}\lambda_{2}(A)}
\]
and overriding notation, define 
\[
\Gamma(n)=\inf_{A\in E_{n}}\Gamma(A).
\]
\end{definition}

\begin{lemma}
\label{lem:gamma-basic}For nonreversible matrices, $\dfrac{1}{15n}\leq\Gamma(n)\leq\dfrac{1}{2}$
and for reversible matrices, $\Gamma^{\text{rev}}(n)=\dfrac{1}{2}$
\end{lemma}

\begin{proof}
The lower bounds follow from Theorem \ref{thm:intro-phi-delta} and
the Cheeger-Buser inequality \ref{thm:cheeger-buser}, and the upper
bound follows by taking the hypercube on $n$ vertices in $E_{n}$,
as discussed in Section \ref{subsec:Tightness-of-cheeger-bounds}.
\end{proof}

\section{Beyond $1/2$}

The first immediate question is to find a doubly stochastic $A$ such
that $\Gamma(A)<\frac{1}{2}$, since we know that for any symmetric
$A$, $\Gamma(A)\geq\frac{1}{2}$ from \ref{lem:gamma-basic}. In
fact, it is not difficult to find such matrices numerically. A straightforward
manner of achieving this is to start with $J$, and for some random
$(i,j)$ and $(k,l)$ and random $\delta$, set $A_{i,j}\pm\delta$,
$A_{k,l}\pm\delta$, $A_{i,l}\mp\delta$, and $A_{k,j}\mp\delta$
(where the range of $\delta$ is chosen to ensure the entries of $A$
do not become negative). With this simple algorithm, it is possible
to find many examples that beat $1/2$, and the following example
was constructed after observing some of these matrices.

\begin{lemma}
\emph{(Beyond 1/2)} \label{lem:gamma-onethird}Let the matrix $A$ be as
follows (where the $\cdot$ are zeros):
\[
A=\begin{bmatrix}\cdot & \frac{1}{3} & \frac{1}{3} & \frac{1}{3}\\
\cdot & \frac{1}{3} & \frac{1}{3} & \frac{1}{3}\\
\cdot & \frac{1}{3} & \frac{1}{3} & \frac{1}{3}\\
1 & \cdot & \cdot & \cdot
\end{bmatrix}.
\]
Then $\Gamma(A)=\dfrac{1}{3}<\dfrac{1}{2}$, improving upon Lemma  \ref{lem:gamma-basic}.
\end{lemma}

\begin{proof}
Note that $A$ is doubly stochastic, and for $S=\{2,3\}$, $\phi(A)=\phi_{S}(A)=\dfrac{1}{3}$,
and all nontrivial eigenvalues of $A$ are 0. Thus, the lemma follows.
\end{proof}
\medskip

\begin{remark}
\emph{(Symmetry about the opposite diagonal)} The matrix $A$ in Lemma \ref{lem:gamma-onethird}
represents an interesting chain, whose adjacency matrix is \emph{symmetric
about the opposite diagonal}. Although this property seems unrelated
to edge expansion or eigenvalues, it will curiously appear again in our
final example.
\end{remark}

\section{Affine-linear constructions}

There are many affine-linear constructions known in literature, and
one such construction is a result of Maria Klawe \cite{klawe1984limitations}.
In fact, the purpose of Klawe's paper was indeed to show non-expansion
of certain constructions of $d$-regular graphs, albeit undirected,
but it was observed by Umesh Vazirani \cite{UVaz-pers} in the 80's
that there is a natural way to orient the edges to create directed
graphs with $d/2$ in-and-out degrees with similar expansion properties,
and we learned of this construction from him. There is an entire family
of constructions with different parameters, but all are equivalent
for us from the perspective of minimizing $\Gamma(n)$. We state one
specific (and neat) construction below.
\medskip

\begin{construction} \label{constr:klawe-vaz} \emph{(Klawe-Vazirani Matrices \cite{klawe1984limitations,UVaz-pers})} Let $n>2$ be a prime, and
create the graph on $n$ vertices with in-degree and out-degree 2
by connecting every vertex $v\in\mathbb{Z}/n$ to two vertices, $1+v$
and $2v$, each with edge-weight 1/2.\end{construction}
\medskip

\begin{example} 
An example for $n=7$ is shown below for right multiplication
by a vector (take the transpose of the matrix for left multiplication):

\[
A_{KV}=\begin{bmatrix}\frac{1}{2} & \cdot & \cdot & \cdot & \cdot & \cdot & \frac{1}{2}\\
\frac{1}{2} & \cdot & \cdot & \cdot & \frac{1}{2} & \cdot & \cdot\\
\cdot & 1 & \cdot & \cdot & \cdot & \cdot & \cdot\\
\cdot & \cdot & \frac{1}{2} & \cdot & \cdot & \frac{1}{2} & \cdot\\
\cdot & \cdot & \frac{1}{2} & \frac{1}{2} & \cdot & \cdot & \cdot\\
\cdot & \cdot & \cdot & \cdot & \frac{1}{2} & \cdot & \frac{1}{2}\\
\cdot & \cdot & \cdot & \frac{1}{2} & \cdot & \frac{1}{2} & \cdot
\end{bmatrix}.
\]

\end{example}

We note the following properties of these graphs.

\begin{lemma}
\label{lem:klawe_vaz}The matrices $A_{KV}$ in Construction \ref{constr:klawe-vaz}
have the following properties:
\begin{enumerate}
\item The matrices have one eigenvalue 1, one eigenvalue 0, and $n-2$ eigenvalues
$\lambda$ such that $\lambda=\dfrac{1}{2}\exp\left(\dfrac{2\pi I}{n-1}\cdot k\right)$
for $k=1,..,n-2$.
\item \cite{klawe1984limitations} $\phi(A)\leq c\cdot\left(\frac{\log\log n}{\log n}\right)^{1/5}$
where $c$ is a constant.
\item $\Gamma(n)\in O\left(\frac{\log\log n}{\log n}\right)^{1/5}$ improving
upon Lemma \ref{lem:gamma-onethird}.
\end{enumerate}
\end{lemma}

\begin{proof}
(1) Let the matrix be as follows: $A=\dfrac{1}{2}(A'+A'')$ where
$A'$ is the directed cycle on $n$ vertices and $A''$ represent
the cycle (since $n>2$ is a prime) that goes from vertex $i$ to
$(2\cdot i)\mod n$. Considering the matrices corresponding to right
multiplication by a vector with vertices $\{0,...,n-1\}$, we have
$A'_{i,i-1}=1$ and otherwise $A'_{i,j}=0$, and $A''_{i,j}=1$ if
$i=2\cdot j$, and else $A''_{i,j}=0$. Let $U$ denote the Fourier
transform over the field $\mathbb{Z}_{n}$, with $U_{i,j}=\dfrac{1}{\sqrt{n}}\omega^{i\cdot j}$
where $\omega=\exp(2\pi I/n)$ is the $n$'th root of unity and $I=\sqrt{-1}$.
Let $B=U^{*}AU$, $B'=U^{*}A'U$, and $B''=U^{*}A''U$. Note that
$B'$ is a diagonal matrix with $B'_{i,i}=w^{i}$ and 0 otherwise.
Similarly, 
\[
B''_{i,j}=\langle U_{i},A''U_{j}\rangle=\dfrac{1}{n}\sum_{k=0}^{n-1}\omega^{-i\cdot k}\omega^{j\cdot(k/2)}=\dfrac{1}{n}\sum_{k=0}^{n-1}\omega^{(j/2-i)\cdot k}
\]
and thus $B''_{i,2\cdot i}=1$ and 0 otherwise. Note $B_{0,0}=1$
and $B_{0,j}=0$ and $B_{i,0}=0$ for all $i,j$, which is the trivial
block corresponding to eigenvalue 1, and the rest of the eigenvalues
of $B$ are in the $n-1\times n-1$ block. For any such eigenvalue
$\lambda$ and corresponding eigenvector $v$ of $B$, we have for
$1\leq i\leq n-1$ that 
\[
\dfrac{1}{2}(\omega^{i}v_{i}+v_{2\cdot i})=\lambda v_{i}
\]
or 
\[
(2\lambda-\omega^{i})v_{i}=v_{2\cdot i}.
\]
Note that $(2\lambda-\omega^{i})\not=0$, else $v$ will be the all
zeros vector. Thus, we get that the equation
\begin{equation}
\prod_{i=1}^{n-1}(2\lambda-\omega^{i})=1,\label{eq:klawe-1}
\end{equation}
whose roots will be the non-trivial eigenvalues of $B$ and $A$.
Note that $\lambda=0$ is a root, since $\prod_{i=1}^{n-1}\omega^{i}=1$
and $(-1)^{n-1}=1$ since $n$ is a prime. We claim that the remaining
$n-2$ roots are all such that $|\lambda|=\dfrac{1}{2}$ with $\lambda\not=\dfrac{1}{2}$,
i.e. $\lambda=\dfrac{1}{2}\exp(2\pi k\cdot I/(n-1))$ for $k\in[n-1]$.
To see this, assume $2\lambda\not=1,$and multiplying and dividing
equation \ref{eq:klawe-1} by $(2\lambda-1)$, we get that 
\[
1=\dfrac{\prod_{i=0}^{n-1}(2\lambda-\omega^{i})}{2\lambda-1}=\dfrac{(2\lambda)^{n}-1}{2\lambda-1}
\]
or 
\[
2\lambda=(2\lambda)^{n}
\]
proving the claim.

(2) Let $s(A)=\dfrac{A+A^{T}}{2}$. It is shown in \cite{klawe1984limitations}
{[}Theorem 2.1{]} that the vertex expansion $\mu(s(A))$ is bounded
as follows:
\[
\mu(s(A))\leq c_{1}\cdot\left(\frac{\log\log n}{\log n}\right)^{1/5}.
\]
Note that since the degree of $s(A)$ is 4, we get that 
\[
\phi(A)=\phi(s(A))\leq4\cdot c_{1}\left(\frac{\log\log n}{\log n}\right)^{1/5}
\]
which proves the claim.

(3) Since $\text{Re}\lambda_{2}\leq\dfrac{1}{2}$, we get the claimed
bound on $\Gamma(n)$ from (1) and (2).
\end{proof}

In fact, there are many different affine-linear constructions that
are known. Recall that Alon and Boppana \cite{alon86,Nilli91} showed
that for any infinite family of $d$-regular undirected graphs, the
adjacency matrices, normalized to be doubly stochastic and with eigenvalues
$1=\lambda_{1}\geq\lambda_{2}\geq\ldots\geq-1$, have $\lambda_{2}\geq\frac{2\sqrt{d-1}}{d}-o(1)$.
Feng and Li \cite{Li92} showed that undirectedness is essential to
this bound: they provide a construction of cyclically-directed $r$-partite
($r\geq2$) $d$-regular digraphs (with $n=kr$ vertices for $k>d$,
$\gcd(k,d)=1$), whose normalized adjacency matrices have (apart from
$r$ ``trivial'' eigenvalues), only eigenvalues of norm $\leq1/d$.
The construction is of an affine-linear nature quite similar to the
matrices in Klawe-Vazirani, and to our knowledge does not give an
upper bound on $\Gamma$ any stronger than those.

\section{\label{subsec:debruijn} The de Bruijn construction}

One of the most beautiful constructions is the de Bruijn graphs or
(implicitly) the de Bruijn sequences. Consider the following problem:
Let $s$ be a string of bits such that every $x\in\{0,1\}^{3}$ appears
as a continuous substring within $s$ exactly once. The earliest reference
to such a string comes from a Sanskrit prosody in the work of Pingala
\cite{brown1869sanskrit}, and one such string $s$ is $s=$yam\={a}t\={a}r\={a}jabh\={a}nasalag\={a}m,
in which each three-syllable pattern occurs starting at its name:
'yam\={a}t\={a}' has a short--long--long pattern, 'm\={a}t\={a}r\={a}'
has a long--long--long pattern, and so on, until 'salag\={a}m' which
has a short--short--long pattern. In general, a de Bruijn sequence
$\text{dB}(n,k)$ is a string containing every $n$-letter sequence
from an alphabet of size $k$ exactly once as a contiguous subsequence,
and the prosody contains a string in $\text{dB}(3,2)$. These were
described by de Bruijn \cite{de1946combinatorial} and I J Good \cite{good1946normal}
independently, and previously by Camille Flye Sainte-Marie for an
alphabet of size 2 \cite{de1975acknowledgement}. These sequences
have also been called ``shortest random-like sequences'' by Karl
Popper \cite{popper2005logic}.
\medskip

\begin{construction} \label{constr:debruijn} \emph{(de Bruijn Matrices \cite{de1946combinatorial})} Define de Bruijn graphs on vertices $\{0,1\}^{k}$ and let every vertex
$|a\rangle=|a_{1}a_{2}\cdots a_{k}\rangle$, have two outgoing edges
to vertices $|a_{2}\cdots a_{k}0\rangle$ and $|a_{2}\cdots a_{k}1\rangle$
with weight $1/2$ each. In words, starting with a $k$-bit vertex
$v$, uniformly go to one of the following two vertices: shift $v$
one bit to the left, and append 0 or 1 uniformly.
\end{construction}

\medskip
\begin{example}
For $k=3$ ($n=8$), the adjacency matrix (for right-multiplication
by a vector) looks like the following:

\[
A_{dB}=\begin{bmatrix}\frac{1}{2} & \cdot & \cdot & \cdot & \frac{1}{2} & \cdot & \cdot & \cdot\\
\frac{1}{2} & \cdot & \cdot & \cdot & \frac{1}{2} & \cdot & \cdot & \cdot\\
\cdot & \frac{1}{2} & \cdot & \cdot & \cdot & \frac{1}{2} & \cdot & \cdot\\
\cdot & \frac{1}{2} & \cdot & \cdot & \cdot & \frac{1}{2} & \cdot & \cdot\\
\cdot & \cdot & \frac{1}{2} & \cdot & \cdot & \cdot & \frac{1}{2} & \cdot\\
\cdot & \cdot & \frac{1}{2} & \cdot & \cdot & \cdot & \frac{1}{2} & \cdot\\
\cdot & \cdot & \cdot & \frac{1}{2} & \cdot & \cdot & \cdot & \frac{1}{2}\\
\cdot & \cdot & \cdot & \frac{1}{2} & \cdot & \cdot & \cdot & \frac{1}{2}
\end{bmatrix}.
\]

\end{example}

In fact, it is possible to construct many different matrices with
different shifts and different base fields instead of $\mathbb{F}_{2}$,
but all the constructions have similar properties.

The matrices in Construction \ref{constr:debruijn} have many remarkable
and beautiful properties. The first observation is that the walk mixes
exactly in $k$ steps, or $A^{k}=J$, since after $k$ steps, the
vertex that is reached is completely uniform independent of the starting
vertex. In fact, the most remarkable thing about these matrices is
that all the (nontrivial) eigenvalues are 0. This implies that each
of the Jordan blocks are of size at most $k$ and nilpotent. Further,
half the singular values are 1, and the other half are all 0. And
the remarkable thing is that the edge expansion of these matrices is $O(1/k)=O(1/\log n)$,
which shows that $\Gamma(n)\in O(1/\log n)$. The edge expansion of these  matrices was studied in \cite{delorme1998spectrum}, but it does not give the actual non-expanding set. In the following lemma, we prove all the relevant properties of de Bruijn matrtices.

\begin{lemma}
\label{lem:debruijn}Let A be the de Bruijn matrix on $n=2^{k}$ vertices.
Then the following hold for $A$:
\begin{enumerate}
\item $\dfrac{1}{2\log n}\leq\phi(A)\leq\dfrac{8}{\log n}$.
\item $A$ is doubly stochastic, and all the nontrivial eigenvalues of $A$
are zero.
\item $\Gamma(n)\leq\dfrac{8}{\log n}\in O\left(\dfrac{1}{\log n}\right)$,
improving upon Lemma \ref{lem:klawe_vaz}.
\item The Jordan and Schur forms are the same, and $A$ has the trivial
block of size 1 for eigenvalue 1, and it has exactly $2^{k-1-r}$
Jordan blocks of size $r\in[k-1]$ and one block of size $k$ for
the 0 eigenvalues.
\item $A$ has $n/2$ singular values that are 1, and $n/2$ singular values
that are 0.
\end{enumerate}
\end{lemma}

\begin{proof}
(1) Let $S$ be the set of all $k$-bit strings that have $r=\lceil k/2\rceil$
contiguous ones. Let $S=\bigcup_{i=r}^{k}T_{i}$ where $T_{i}$ is
the set of all strings that have ones at all the $r$ positions ending
at $i$, and there is no $j>i$ with the same property, i.e. there
is no $j>i$ such that the substring {[}$j-r+1,j${]} is all ones.
Then note that the $T_{i}$'s are disjoint. Further, the definition
of $T_{i}$ implies that the $r$ positions ending at $i$ are all
ones, and the $i+1$'st position (except for $T_{k}$) is 0. In fact,
these two constraints are \emph{sufficient} to define the $T_{i}$'s
and ensure that they do not overlap, since to overlap, there must
be a string with $r$ contiguous ones ending at $i$ and $r$ contiguous
ones ending at $j$, but if the $i+1$ position in the string is $0$,
then there are $k-r-1<k/2\leq r$ positions left for $r$ ones, which
is an impossibility.

Thus, for $T_{i},$ all strings have $r+1$ positions fixed -- the
indices $i-r+1$ to $i$ contain 1, and the position $i+1$ contains
0. Note that the positions $\leq i-r$ could be any bits, and also
all positions $>i+1$, since the length of the suffix is $k-(i+1)<r$
for all $i\geq r$. Thus, we get that $|T_{i}|=2^{k-r-1}$ for $r\leq i\leq k-1$,
and $|T_{k}|=2^{k-r}$. Since $\frac{k}{2}\leq r\leq\frac{k}{2}+1$,
we have 
\begin{align*}
|S| & =(k-r)\cdot2^{k-r-1}+2^{k-r}\\
 & \leq\dfrac{k}{2}2^{k/2-1}+2^{k/2}\\
 & \leq\dfrac{5}{4}\cdot k\cdot2^{k/2}\\
 & \ll2^{k-1}\\
 & =\dfrac{n}{2}
\end{align*}
and thus, $S$ is a valid set to consider for determining $\phi$.
Moreover, 
\begin{align*}
|S| & =(k-r)\cdot2^{k-r-1}+2^{k-r}\\
 & \geq\left(\dfrac{k}{2}-1\right)2^{k/2-2}+2^{k/2-1}\\
 & \geq\dfrac{k}{8}2^{k/2}
\end{align*}
The key property of all the sets is the following: for any $i>r$,
the vertices in the sets $T_{i}$ have all outgoing edges inside the
set $S$, i.e. there are no edges from the vertices in the set $T_{i}$
to $\overline{S}$, since the strings will always have a contiguous
substring of $r$ ones after a one-bit left shift. The only set that
could have edges to $\overline{S}$ is $T_{r}$, and it could have
at most $2\cdot|T_{r}|$ outgoing edges. Thus, we get that
\[
\phi_{S}(A)\leq\dfrac{2\cdot|T_{r}|}{|S|}\leq\dfrac{2\cdot2^{k/2-1}}{\dfrac{1}{8}\cdot k\cdot2^{k/2}}\leq\dfrac{8}{k}
\]
as required.

Further, since $A^{k}=J$ from (2), we have from Lemma \ref{lem:phi_pow_bound}
that $\phi(A^{k})\leq k\cdot\phi(A)$ or 
\[
\phi(A)\geq\dfrac{1}{2k}.
\]

(2) It is straightforward to see from the definition that $A\mathbf{1}=\mathbf{1}$
and $A^{T}\mathbf{1}=\mathbf{1}$ and thus $A$ is doubly stochastic.
Let $a=a_{1}a_{2}\cdots a_{k}$ be a $k$-bit string, let $|a\rangle$
denote the corresponding standard basis vector in an $n$-dimensional
space, and for any $S\subseteq[k]$, let $a_{S}=\sum_{i\in S}a_{i}$
and
\[
v_{S}=\sum_{a}(-1)^{a_{S}}|a\rangle.
\]
Let $V=\{v_{S}:S\subseteq[k]\}$ and $W=\{v_{S}:1\in S,S\subseteq[k]\}$.
Then we show that for any $v_{S}\in W$, 
\[
Av_{S}=0.
\]
Note that by definition, 
\[
A|a_{1}...a_{k}\rangle=\dfrac{1}{2}|a_{2}...a_{k}0\rangle+\dfrac{1}{2}|a_{2}...a_{k}1\rangle
\]
By direct calculation, we have letting $T=S\backslash\{1\}$, 
\begin{align*}
Av_{S} & =\sum_{a}(-1)^{a_{S}}A|a_{1}\cdots a_{k}\rangle\\
 & =\sum_{a}(-1)^{a_{1}}(-1)^{a_{T}}A|a_{1}\cdots a_{k}\rangle\\
 & =\dfrac{1}{2}\sum_{a}(-1)^{a_{T}}(-1)^{a_{1}}\left(|a_{2}\cdots a_{k}0\rangle+|a_{2}\cdots a_{k}1\rangle\right)\\
 & =\dfrac{1}{2}\sum_{a_{2},\cdots a_{k}=0}^{1}(-1)^{a_{T}}\left(|a_{2}\cdots a_{k}0\rangle+|a_{2}\cdots a_{k}1\rangle\right)\sum_{a_{1}=0}^{1}(-1)^{a_{1}}\\
 & =0
\end{align*}
as claimed. Further, for $S\not=T$, let $U=(S\backslash(S\cap T))\cup(T\backslash(S\cap T))\not=\emptyset$,
then it is simple to see that
\begin{align*}
\langle v_{S},v_{T}\rangle & =\sum_{a}\sum_{b}(-1)^{a_{S}}(-1)^{b_{T}}\langle a|b\rangle=\sum_{a}(-1)^{a_{S}}(-1)^{a_{T}}=\sum_{a}(-1)^{a_{U}}=0.
\end{align*}
Thus the $v_{S}$ are all orthogonal, and $|W|=2^{n-1}$, and since
$Av_{S}=0$ for $v_{S}\in W$, it implies that the kernel $W$ of
$A$ has dimension $n-1$. Since $A$ has 1 as the trivial eigenvalue
for $v_{\emptyset}$, it means that all the nontrivial eigenvalues
of $A$ are 0.

(3) Combining (1) and (2) gives the bound on gamma since $\lambda_{2}(A)=0$.

(4) We will use notation from (2). The set $V$ as defined contains
orthogonal vectors as shown in (2), and if we express $A$ in the
basis of the vectors in $V$, we get a matrix of Jordan blocks, and
as such both the Jordan and Schur forms of $A$ will exactly be the
same (since the vectors of $V$ form a unitary matrix). Our aim is
to now understand the effect of $A$ on vectors not in $W$. This
will help us understand the chain of generalized eigenvectors. Consider
$S\subseteq[k]$ where $S=\{i_{1},...,i_{r}\}$ with $i_{1}<i_{2}<...<i_{r}$,
and let $T=\{i_{1}-1,i_{2}-1,...,i_{r}-1\}$. Let $a=a_{1}...a_{k}$,
$b=b_{1}...b_{k}$ with $b_{i}=a_{i}+1$ and $b_{k}=0$, and $c=c_{1}...c_{k}$
with $c_{i}=a_{i}+1$ with $c_{k}=1$. Then note that if $i_{1}=1$,
then from part (2), $Av_{S}=0$, and if $i_{1}>1$, then 
\begin{align*}
Av_{S} & =A\sum_{a}(-1)^{a_{S}}|a_{1}...a_{r}\rangle\\
 & =\dfrac{1}{2}\sum_{a}(-1)^{a_{S}}|a_{2}...a_{r}0\rangle+\dfrac{1}{2}\sum_{a}(-1)^{a_{S}}|a_{2}...a_{r}1\rangle\\
 & =\sum_{a_{2},...,a_{k}=0}^{1}(-1)^{a_{S}}|a_{2}...a_{r}0\rangle+\sum_{a_{2},...,a_{k}=0}^{1}(-1)^{a_{S}}|a_{2}...a_{r}1\rangle\\
 & \ \ [\text{since \ensuremath{1\not\in S}}]\\
 & =\sum_{b_{1},...,b_{k-1}=0}^{1}(-1)^{b_{T}}|b\rangle+\sum_{c_{1},...,c_{k-1}=0}^{1}(-1)^{c_{T}}|c\rangle\\
 & =\sum_{a}(-1)^{a_{T}}|a\rangle\\
 & =v_{T}.
\end{align*}
Thus, for any $S\subseteq[k]$ denote $S+j=\{i\in[k]:i-j\in S\}$,
and let $1\in S$ and $\max S=r$, then we have for $1\leq j\leq k-r$
\[
Av_{S+j}=v_{S+j-1}
\]
giving us a chain of $k+1-r$ vectors that end with a vector in the
kernel. Thus, this forms one specific Jordan block. In general, we
have that any vector $v_{S}$ will belong to a Jordan block of size
$k-(\max S-\min S)$. To count the number of distinct Jordan blocks,
consider $S$ such that $1\in S$. Then the vector $v_{S}$ in the
kernel will be the last vector in the chain that starts with $v_{S+k-\max S}$,
in a block of size $k+1-\max S$. Thus, for each $r\in[k-1]$, the
number of Jordan blocks of size $r$ is the number of vectors $v_{S}$
with $1\in S$ and $k+1-r\in S$, which is exactly $2^{k-1-r}$. Further,
the only block of size $k$ is obtained for the chain ending in $v_{\{1\}}$,
and there is the trivial block of size 1 corresponding to $v_{\emptyset}$.
As a sanity check, if we sum the sizes of all blocks, we get 
\[
\sum_{r=1}^{k-1}r\cdot2^{k-1-r}+k+1=2^{k}\left(1-2^{-(k-1)}-(k-1)2^{-k}\right)+k+1=2^{k}=n
\]
as expected.

(5) Note that with the unitary $U$ formed by the vectors from $V=\{v_{S}:S\subseteq[k]\}$
(after suitably normalizing), we get that $A=UTU^{*}$ where $T$
is a collection of Jordan blocks as shown in (4) above. Since a Jordan
block of size $r$ can be converted to a diagonal matrix of size $r$
with $r-1$ ones and one zero entry by multiplying with a permutation
matrix, we get that $T=DP$ where $P$ is a permutation matrix that
converts each Jordan block to a diagonal matrix. Thus we have the
singular value decomposition of $A=UDPU^{*}=UDQ$ since $Q$ is a
unitary. To count the number of zeros on the diagonal in $D$, note
that each nontrivial Jordan block contributes exactly one 0, and thus
the total number of zeros is equal to the number of nontrivial Jordan
blocks, which is exactly 
\[
\sum_{r=1}^{k-1}2^{k-1-r}+1=2^{k-1}=\dfrac{n}{2}
\]
implying that there are $n-n/2$ ones on the diagonal. Thus $A$ has
$n/2$ singular values that are 1 and $n/2$ singular values that
are 0.
\end{proof}

The beautiful thing about this construction is that it is extremely
simple to describe, and still has the remarkable properties in the
lemma above. Also, this construction is the benchmark for other constructions,
and although it does not achieve a low enough value of $\Gamma(n)$
that would be sufficient for Theorem \ref{thm:intro-phi-delta}, it
will be the starting point for the construction in Section \ref{subsec:Constraints-on-the}.
\medskip

Our aim now will be to beat the upper bound on $\Gamma$ in Lemma
\ref{lem:debruijn}, and understand if our lower bound on $\phi$
in Theorem \ref{thm:intro-phi-delta} is tight or whether it is exponentially
worse than the truth. Note that if it is true that 
\[
\Gamma(n)\in\Omega\left(\dfrac{1}{\text{polylog}(n)}\right),
\]
it would mean that our techniques for the proof of Theorem \ref{thm:intro-phi-delta}
are extremely weak, and different techniques will be required to get
a tighter bound. From a utilitarian perspective, it would imply that
the spectral gap is a good estimate for $\phi$ (up to $\text{polylog}(n)$
terms which are essentially negligible in succinctly defined chains
where the input length is $O(\log n)$) even in the nonreversible
case. Towards this end, we will try to find constructions of matrices
that try to surpass the upper bound on $\Gamma$ in Lemma \ref{lem:debruijn}.

\section{\label{subsec:Constraints-on-the}Constraints on the Search Space}

At this point, we are in search of a doubly stochastic matrix $A$
that helps to improve the bound in Lemma \ref{lem:debruijn}. Since
the search space (all doubly stochastic matrices) is difficult to
understand in terms of non-expansion, we will systematically try to
impose meaningful constraints on it to arrive at the type and form
of matrices that we want. Towards this end, our first question, that
will turn out to be sufficient to be the last, is to understand the
following:
\medskip

\noindent \textbf{Main question:}\emph{ How small can the  edge expansion 
of a doubly stochastic matrix be if all its nontrivial eigenvalues
are 0?}
\medskip

If we can show that $\Gamma(n)\in\Omega(1/\log n)$ for all $n\times n$
matrices that have all nontrivial eigenvalues 0, then it will imply that de Bruijn
matrices are optimal in the sense of non-expansion as described in
Definition \ref{def:nonexpansion} in Section \ref{subsec:Non-expansion},
that is, they have the least edge expansion amongst all matrices that have
all (nontrivial) eigenvalues 0. This will be our first constraint,
and further, this restriction is also sufficiently general in the
sense described next.

\begin{constraint}[1]
\label{Constraint-1:-Restrict}
Restrict
all nontrivial eigenvalues to 0.
\end{constraint}

The rationale behind choosing all eigenvalues
0 is as follows -- as seen in the first steps in the proof of Lemma
\ref{lem:T_norm_1_bound}, if for a doubly stochastic matrix $A$,
every nontrivial eigenvalue has magnitude at most $1-c$ for some
constant $c$, then powering just $O(\log n)$ times will make the
diagonal entries inverse polynomially small in magnitude, and thus
it would seem that the matrix should have behavior similar to matrices
with all eigenvalues 0. Thus, if $A$ had all eigenvalues with magnitude
less than $1-c$, we can simply consider the matrix $A^{k}$ with
$k\in O(\log n)$ as our starting matrix, and we know that its expansion will be at most $O(\log n)$ times the expansion of $A$ (from Lemma
\ref{lem:phi_pow_bound}), and thus if $A$ had small expansion --
about $O(1/\sqrt{n})$ -- $A^{k}$ will have similar expansion (albeit
off by a factor of $\log n$). Given any doubly stochastic $A$ with
$\text{Re}\lambda_{2}<1$, it is simple to ensure that $|\lambda_{i}|\leq1-c$
for some specific $c$ by lazification, i.e. by considering $\exp(t\cdot(A-I))$
or $\frac{1}{2}(A+I)$, which is essentially the first step in the
proof of Lemma \ref{lem:real-lambda-bound} or the lower bound in
Theorem \ref{thm:intro-phi-delta}.

Restricting to matrices with all 0 (nontrivial) eigenvalues, our primary
concern throughout will be the matrix $B=A-J$. If all nontrivial
eigenvalues of $A$ are 0, it means all eigenvalues of $B$ are 0,
and we will further write the Schur form of $B$ as $B=UTU^{*}$ where
$U$ is a unitary and $T$ is upper triangular and all its diagonal
entries are 0.

\begin{lemma}
If $A=J+B$ has all (nontrivial) eigenvalues 0, then $B^{n-1}=0$.
\end{lemma}

\begin{proof}
The proof is immediate since $T$ is nilpotent.
\end{proof}
Note that the entries in the rows and columns of the matrix $A$ will
always sum to 1, since we have removed the all ones eigenvector in
$J$ and we ensure that the first column of $U$ is the vector $1/\sqrt{n}\cdot\mathbf{1}$,
which would imply all other columns of $U$ are orthogonal to the
vector $\mathbf{1}$ (since $U$ is a unitary), which would further
give us $B\mathbf{1}=0$ and $B^{T}\mathbf{1}=0$. Also, there is
another useful consequence of constraint \hyperref[Constraint-1:-Restrict]{1}
-- since all the eigenvalues of $A$ are real, we can restrict to
real unitary matrices $U$ (by looking at the process of obtaining
a Schur decomposition, since $A$ is real and has real eigenvalues).
This actually gives us our second implied condition: 

\begin{constraint}[2] \label{(Implied)-Constraint-2:unitary-real}
The
unitary $U$ is real, and has $\frac{1}{\sqrt{n}}\mathbf{1}$ as the
first column.
\end{constraint}

Given this observation, the next step is to look at the proof of Lemma
\ref{lem:T_norm_1_bound}. In Lemma \ref{lem:T_norm_1_bound}, the
main aim was to upper bound $k$ for which $T^{k}\approx0$ and it
was seen that a value of $k$ about $n/1-\text{Re}\lambda_{2}$ was
sufficient. Our main aim is to now find a construction of matrices
such that this lemma is tight. In other words, we now seek a construction
of $T$ (and $B$) such that for sufficiently small $k$, the norm
$\|T^{k}\|_{2}$ is sufficiently far from zero. Consider the de Bruijn
matrices and the matrix $B_{dB}=A_{dB}-J$. We know that $B_{dB}^{k}=0$
for $k=\log n$. Further from Lemma \ref{lem:debruijn}, since $\phi(A_{dB})\leq8/k$,
it also implies that for any $j\ll k$, $\|B^{j}\|\gg0$, since otherwise
from the proof of Lemma \ref{lem:debruijn}, we would get that $\phi(A)>1/j$,
which would be a contradiction. 

It is not difficult to construct a $T$ such that $T^{k}\gg0$ for
$k\in o(n)$. For instance, it is simple to obtain this if we let
the Schur form be one large Jordan block, that is, $T$ is simply
the matrix with 1's above the diagonal.
\[
T=\begin{bmatrix}\mathbf{0} & 0 & 0 & 0 & 0 & 0\\
0 & \mathbf{0} & 1 & 0 & 0 & 0\\
0 & 0 & \mathbf{0} & 1 & 0 & 0\\
0 & 0 & 0 & \mathbf{\ddots} & \ddots & 0\\
0 & 0 & 0 & 0 & \mathbf{0} & 1\\
0 & 0 & 0 & 0 & 0 & \mathbf{0}
\end{bmatrix},\ \ T^{2}=\begin{bmatrix}\mathbf{0} & 0 & 0 & 0 & 0 & 0\\
0 & \mathbf{0} & 0 & 1 & 0 & 0\\
0 & 0 & \mathbf{0} & 0 & \ddots & 0\\
0 & 0 & 0 & \mathbf{\ddots} & \ddots & 1\\
0 & 0 & 0 & 0 & \mathbf{0} & 0\\
0 & 0 & 0 & 0 & 0 & \mathbf{0}
\end{bmatrix},\ \ T^{3}=\begin{bmatrix}\mathbf{0} & 0 & 0 & 0 & 0 & 0\\
0 & \mathbf{0} & 0 & 0 & 1 & 0\\
0 & 0 & \mathbf{0} & 0 & \ddots & 1\\
0 & 0 & 0 & \ddots & \ddots & 0\\
0 & 0 & 0 & 0 & \mathbf{0} & 0\\
0 & 0 & 0 & 0 & 0 & \mathbf{0}
\end{bmatrix}
\]

For this matrix, it is clear that $\|T^{k}\|=1$ for any $k<n-1$,
since the diagonal consisting of ones keeps shifting away from the
main diagonal with higher powers of $T$. However, the immediate problem
in choosing a matrix $T$ as mentioned above is that there might not be any
unitary $U$ such that $J+UTU^{*}=A$ is a nonnegative matrix. In
fact, this is one of the primary problems in constructing these
matrices -- ensuring that they are positive. Note that the entries
in the rows and columns of the matrix will always sum to 1 due to
condition 2. At this point, it is difficult to find a unitary $U$
that transforms the matrix $T$ above to a doubly stochastic matrix,
and we will need to modify $T$ in some manner. Note that essentially,
the only requirement for $T$ is that it is upper triangular with
diagonal entries 0, but having a completely general $T$ is extremely
difficult to handle. Thus, our next relaxation/constraint on $T$
will be as follows:

\begin{constraint}[3]
For $A=J+UTU^{*}$, let $T$ have
some number $r$ above the diagonal instead of 1, and let the rest
of the entries in $T$ be 0.
\end{constraint}

From the constraint above, we get that 
\[
T=\begin{bmatrix}\mathbf{0} & 0 & 0 & 0 & 0 & 0\\
0 & \mathbf{0} & r & 0 & 0 & 0\\
0 & 0 & \mathbf{0} & r & 0 & 0\\
0 & 0 & 0 & \mathbf{\ddots} & \ddots & 0\\
0 & 0 & 0 & 0 & \mathbf{0} & r\\
0 & 0 & 0 & 0 & 0 & \mathbf{0}
\end{bmatrix},\ \ T^{2}=\begin{bmatrix}\mathbf{0} & 0 & 0 & 0 & 0 & 0\\
0 & \mathbf{0} & 0 & r^{2} & 0 & 0\\
0 & 0 & \mathbf{0} & 0 & \ddots & 0\\
0 & 0 & 0 & \mathbf{\ddots} & \ddots & r^{2}\\
0 & 0 & 0 & 0 & \mathbf{0} & 0\\
0 & 0 & 0 & 0 & 0 & \mathbf{0}
\end{bmatrix},\ \ T^{3}=\begin{bmatrix}\mathbf{0} & 0 & 0 & 0 & 0 & 0\\
0 & \mathbf{0} & 0 & 0 & r^{3} & 0\\
0 & 0 & \mathbf{0} & 0 & \ddots & r^{3}\\
0 & 0 & 0 & \ddots & \ddots & 0\\
0 & 0 & 0 & 0 & \mathbf{0} & 0\\
0 & 0 & 0 & 0 & 0 & \mathbf{0}
\end{bmatrix}
\]

Note that for $r\approx1-1/n$, since $\|T^{k}\|=r^{k}$, for any
$k<n/2$, we would have that $\|T^{k}\|=r^{k}\gtrsim e^{-1/2}$ which
is much larger than 0. Thus if we can transform $T$ to a doubly stochastic
matrix using any valid unitary $U$ with $r\approx1-1/n$, we will
get the type of matrix that we are looking for. At this point, instead
of fixing $r$ to being about $O(1/n)$ far from 1, our aim will be
to have a matrix $T$ for which there is \emph{some} unitary that
transforms it to a doubly stochastic matrix and has $r$ as large
as possible.

Also, the rationale for $T$ having non-zero entries only on the off-diagonal
is as follows: In any upper triangular matrix $T$ in which the diagonal
has zeros, the off-diagonal entries \emph{affect} the entries in powers
of $T$ the most, since entries far from the diagonal will become
ineffective after a few powers of $T$. Thus, choosing $T$ with the
non-zeros pattern of a Jordan block will not be far from optimal.

Note that our structure of $T$ gives the following lemma which helps
to illustrate the main issue with trying to have constructions that
beat the upper bound on $\Gamma$ in Lemma \ref{lem:debruijn}.

\begin{lemma}
Let $A=J+UTU^{*}$ be a doubly stochastic matrix with $T$ being the
all zeros matrix but with entry $r$ in the off-diagonal entries from
row $2$ to $n-1$. Then 
\[
r^{2}=\dfrac{\sum A_{i,j}^{2}-1}{n-2}
\]
\end{lemma}

\begin{proof}
Looking at the trace of $AA^{T}$, we get 
\[
\sum_{i,j}A_{i,j}^{2}=\text{Tr}(AA^{T})=1+r^{2}(n-2)
\]
which gives the lemma.
\end{proof}
Note that from this lemma, it might seem that maximizing $r$ is a
simple task -- maximizing the sum of squares of entries of $A$,
and it is easy to construct such an $A$, for instance $A=(1-p)I+pJ$
for an appropriate $p$. This illustrates two issues. The lemma does
not use our constraints on $A$ -- nonnegativity, and the fact that
all nontrivial eigenvalues 0. This is the primary issue in all attempts
of construction of the intended matrices. 
\medskip

\begin{remark}
\emph{(Problems with simulation)} We would also like to
re-state, and mentioned in the introduction, that it is not possible
to use simulations and find different $A$'s that are doubly stochastic
and have small expansion but large spectral gap, and it is even harder
to simulate $T$ and $U$ to have the type of properties that we care
about. The main issue is the extreme sensitivity of these matrices
to small perturbations, and in fact this actually is our very aim
-- to find matrices that \emph{are} heavily affected (in terms of
change in eigenvalues) by a small perturbations (which corresponds
to a small change in the expansion). Since simulations with finite
number of bits act as perturbations themselves, it would be highly improbable
to arrive at such matrices numerically.
\end{remark}

\section{\label{subsec:Special-Cases-of}Special Cases of specific matrices
satisfying constraints}

To get started on our aim of finding $T$ with a large value of $r$
and some unitary that transforms it to a doubly stochastic matrix,
we are going to start with a very simple case. We were averse to looking
at specific small examples initially, since it seemed that for any
finite $n$ (sufficiently small enough, say $n\leq12$), the values
of  edge expansion  and the spectral gaps would be off by large constant
factors, and would essentially be uninformative since they are scalar
values. However, it turned out that these matrices indicated something
that we did not initially expect. 
\medskip

Consider the case of $n=3$. Let our $T$ and $U$ be as follows:
\[
T=\begin{bmatrix}0 & 0 & 0\\
0 & 0 & r\\
0 & 0 & 0
\end{bmatrix},\ \ U=\begin{bmatrix}\dfrac{1}{\sqrt{3}} & a_{1} & b_{1}\\
\dfrac{1}{\sqrt{3}} & a_{2} & b_{2}\\
\dfrac{1}{\sqrt{3}} & a_{3} & b_{3}
\end{bmatrix}
\]
where the $a_{i}$'s and $b_{i}$'s are real, as stated in constraint
\hyperref[(Implied)-Constraint-2:unitary-real]{2}. We can now treat these
as 3 points in the x-y plane, and get the following equations using
the fact that $U$ is a unitary:
\begin{align*}
a_{1}+a_{2}+a_{3} & =0\\
b_{1}+b_{2}+b_{3} & =0\\
a_{1}^{2}+b_{1}^{2} & =\dfrac{2}{3}\\
a_{2}^{2}+b_{2}^{2} & =\dfrac{2}{3}\\
a_{3}^{2}+b_{3}^{2} & =\dfrac{2}{3}
\end{align*}

The last three equations tell us that each of the three points are
on a circle of radius $t=\sqrt{2/3}$, and the first two equations
tells us that the vector $(a_{3},b_{3})=-(a_{1},b_{1})-(a_{2},b_{2})$.
Since all the three points are on the circle, it implies that the
third point lies on the bisector of the angle between the vectors
from the origin to the first and second points. However, since the
equations are completely symmetric, it implies that this must be satisfied
by all three points, and thus it shows that the three points lie at
the the corners of an equilateral triangle. Thus, we get that the
points are at an angle of $2\pi/3$ from each other, and this gives
\begin{align*}
(a_{1},b_{1}) & =(a,b)\\
(a_{2},b_{2}) & =(-\dfrac{1}{2}a-\dfrac{\sqrt{3}}{2}b,-\dfrac{1}{2}b+\dfrac{\sqrt{3}}{2}a)\\
(a_{3},b_{3}) & =(-\dfrac{1}{2}a+\dfrac{\sqrt{3}}{2}b,-\dfrac{1}{2}b-\dfrac{\sqrt{3}}{2}a)
\end{align*}

Thus we have our unitary $U$, and we need to choose $(a,b)$ (with
$a^{2}+b^{2}=\frac{2}{3}$) in order to choose as large a value of
$r$ as possible keeping the matrix nonnegative. Note that due to
the structure of $T$ and $U$, we get that the entries $B_{i,j}=r\cdot a_{i}\cdot b_{j}$,
and we need to ensure $B_{i,j}+\frac{1}{3}\geq0$. Thus our problem
becomes solving the following:

\textbf{maximize} $r$

\textbf{s.t.} $B_{i,j}+\dfrac{1}{3}\geq0$ for all $i,j$, and $a^{2}+b^{2}=\dfrac{2}{3}$.\\

The above problem has a simple solution, obtained by
\[
a=b=\dfrac{1}{\sqrt{3}}
\]
 giving 
\[
r=\dfrac{\sqrt{3}+1}{\sqrt{3}+2}=1-\dfrac{1}{\sqrt{3}+2}.
\]

Note that at this point, we have a doubly stochastic matrix $A$ with
all eigenvalues 0, with a specific value of $r$ stated above,
but all these numbers could be off by large factors and are essentially
meaningless.

What is indeed meaningful, is to look at the unitary $U$:
\[
U=\begin{bmatrix}\dfrac{1}{\sqrt{3}} & \dfrac{1}{\sqrt{3}} & \dfrac{1}{\sqrt{3}}\\
\dfrac{1}{\sqrt{3}} & -\dfrac{1}{2}(1+\dfrac{1}{\sqrt{3}}) & \dfrac{1}{2}(1-\dfrac{1}{\sqrt{3}})\\
\dfrac{1}{\sqrt{3}} & -\dfrac{1}{2}(1-\dfrac{1}{\sqrt{3}}) & -\dfrac{1}{2}(1+\dfrac{1}{\sqrt{3}})
\end{bmatrix}
\]

The fascinating thing is that in addition to having the first column
as a multiple of $\mathbf{1}$, \emph{even the first row in $U$ is
a multiple of the all ones vector.} And since $U$ indicates a rotation
in the space of $T$ and is a vector unlike the numbers $\phi$, $r$,
$\lambda_{2}$ that are scalars, it indicates that the maximum value
of $r$ is obtained in the direction of the all 1's vector. This gives
us our next crucial constraint:

\begin{constraint}[4]
The maximum value of $r$ is obtained
when $U$ has the first row containing the all ones vector. 
\end{constraint}

\section{\label{subsec:The-general-case}The general case and relevant reasoning}

We now wish to understand how large the value of $r$ can be, for
any unitary with the first row and column vectors that are multiples
of the all ones vector. Indeed, we can show the following:

\begin{lemma}
\label{lem:phi-max-1byrootn}Let $A=J+UTU^{*}$ be a doubly stochastic
matrix where $T$ contains zeros on the diagonal and $r$ on the off-diagonal
entries (from row 2 to $n$) and $U$ be some real unitary with the
first row and column being a multiple of the all ones vector. Then
\[
r \leq 1-\dfrac{\sqrt{n}+2}{n-4} \in1-\Omega\left(\dfrac{1}{\sqrt{n}}\right).
\]
\end{lemma}

\begin{proof}
Let the $k$'th column of $U$ be $u_{k}$, where $u_{i,k}$ represents
the $i$'th entry of column $k$ or $(i,k)$'th entry of $U$, then
based on our choice of $T$ and the conditions on $U$, we get 
\[
A=J+r\cdot\sum_{k=2}^{n-1}u_{k}u_{k+1}^{*}
\]
and 
\[
A_{i,j}=\dfrac{1}{n}+r\sum_{k=2}^{n-1}u_{i,k}u_{j,k+1}
\]
since the unitaries are real as observed in constraint \hyperref[(Implied)-Constraint-2:unitary-real]{2}.
Note that every entry in the first row and column of $U$ is $\dfrac{1}{\sqrt{n}}$,
and thus $u_{i,1}=\dfrac{1}{\sqrt{n}}$ for all $i$, $u_{1,k}=\dfrac{1}{\sqrt{n}}$
for all $k$, and since the first row and column of $U$ are multiples
of the all ones vector, the sum of entries in any other row or column
of $U$ is 0. Thus we get that 
\begin{align*}
A_{1,1} & =\dfrac{1}{n}+r\cdot{\dfrac{n-2}{n}}\\
1-A_{1,1} & =1-\dfrac{1}{n}-r+\dfrac{2r}{n}
\end{align*}
and for $j\geq2$, 
\begin{align*}
A_{1,j} & =\dfrac{1}{n}+r\sum_{k=2}^{n-1}u_{1,k}u_{j,k+1}\\
 & =\dfrac{1}{n}+\dfrac{r}{\sqrt{n}}\sum_{k=2}^{n-1}u_{j,k+1}\\
 & =\dfrac{1}{n}+\dfrac{r}{\sqrt{n}}(-u_{j,1}-u_{j,2})\\
 & \ \ \ \text{[since \ensuremath{\sum_{k=1}^{n}u_{j,k}=0}]}\\
 & =\dfrac{1-r}{n}-\dfrac{r}{\sqrt{n}}\cdot u_{j,2}
\end{align*}
which gives for $j\geq2$
\[
u_{j,2}=\dfrac{\sqrt{n}}{r}\cdot\left(\dfrac{1-r}{n}-A_{1,j}\right).
\]
Now consider $u_{2}$ (or the second column of $U$), then we have
\begin{align}
1 & =\sum_{j=1}^{n}u_{j,2}^{2}\nonumber \\
 & =\dfrac{1}{n}+\dfrac{n}{r^{2}}\sum_{j=2}^{n}\left(\dfrac{1-r}{n}-A_{1,j}\right)^{2}\nonumber \\
\dfrac{n-1}{n^{2}}\cdot r^{2} & =\dfrac{n-1}{n^{2}}\cdot(1-r)^{2}+\sum_{j=2}^{n}A_{1,j}^{2}-2\cdot\dfrac{1-r}{n}\cdot\sum_{j=2}^{n}A_{1,j}\label{eq:phi-max-1byrootn}
\end{align}
Note that since the matrix is double stochastic, $\sum_{j=2}^{n}A_{1,j}=1-A_{1,1}$
and each $A_{1,j}\geq0$, and we get
\[
\sum_{j=2}^{n}A_{1,j}^{2}\leq\left(\sum_{j=2}^{n}A_{1,j}\right)^{2}=(1-A_{1,1})^{2}.
\]

Replacing these in equation \ref{eq:phi-max-1byrootn}, we get 
\begin{align*}
\dfrac{n-1}{n^{2}}\cdot r^{2} & \leq\dfrac{n-1}{n^{2}}\cdot(1-r)^{2}+(1-A_{1,1})^{2}-2\cdot\dfrac{1-r}{n}\cdot(1-A_{1,1})\\
\dfrac{n-1}{n^{2}}\cdot r^{2} & \leq\dfrac{n-1}{n^{2}}\cdot(1-r)^{2}+(1-A_{1,1})\left(1-A_{1,1}-2\cdot\dfrac{1-r}{n}\right)\\
(n-1)(2r-1) & \leq(n(1-r)+2r-1)(n(1-r)+4r-3)\\
0 & \leq(n^{2}-6n+8)r^{2}+(-2n^{2}+8n-8)r+n^{2}-3n+2
\end{align*}
and since $r\leq1$ (else the matrix will have norm larger than
1 and will not be doubly stochastic, see Lemma \ref{lem:transRtoA}),
the valid range of $r$ solving the above quadratic is obtained by
\[
r\leq1-\dfrac{\sqrt{n}+2}{n-4}
\]
which implies that 
\[
r\in1-\Omega\left(\dfrac{1}{\sqrt{n}}\right)
\]
and completes the proof.
\end{proof}

Rephrasing Lemma \ref{lem:phi-max-1byrootn} again, it says that after fixing the first row of the unitary matrix
to the all ones vector, given the form of $T$ that we have fixed,
the best value of $r$ that we can hope to achieve is about $1-\frac{1}{\sqrt{n}}$,
and this will henceforth be our aim. Note that if possible, this will
give us an \emph{exponential} improvement over the deBruin construction
for $\Gamma(n)$ in Lemma \ref{lem:debruijn}. Thus, to get a maximum value of $r$, our optimization problem becomes
the following:

\begin{optprob}[1] \label{opt:first}

\textbf{maximize} $r$

\textbf{such that:}

$\dfrac{1}{n}+r\sum_{k=2}^{n-1}u_{i,k}u_{j,k+1}\geq0$ for all $2\leq i,j\leq n$,

$\dfrac{1}{n}+\dfrac{r}{\sqrt{n}}\sum_{k=2}^{n-1}u_{j,k+1}\geq0$
for all $2\leq j\leq n$

$\dfrac{1}{n}+\dfrac{r}{\sqrt{n}}\sum_{k=2}^{n-1}u_{i,k}\geq0$ for
all $2\leq i\leq n$

$u_{i,k}=\langle e_{i},Ue_{k}\rangle$

$UU^{T}=U^{T}U=I$.

$u_{1,j}=\dfrac{1}{\sqrt{n}},$$u_{i,1}=\dfrac{1}{\sqrt{n}}$, for
all $1\leq i,j\leq n$. 
\end{optprob}

Although we can attempt and solve the above problem, it is unwieldy
and the conditions on the $(n-1)^{2}$ variables in $U$ make it almost
intractable. To solve this issue and understand how the solutions
to this problem would look like, we will simplify it by fixing other
entries of the unitary $U$. This fixing will be mildly creative and
empirical, based on our observations in the $n=3$ and other cases
with small $n$. Note that intuitively, based on the constraints on the unitary matrix $U$, it is a rotation
in space of two distinct ``Schur'' blocks, where a Schur block
is same as a Jordan block except that the off-diagonal entries are
$r$ and not 1. The first is the trivial block containing only the
eigenvalue 1, and the second block of size $n-1$  contains
the remaining zero eigenvalues and has $r$ above the diagonal and zeros elsewhere. Since only two types of \emph{actions} are performed
on the Schur blocks, we are going to set the entries (at positions
greater than 1) in any column of the unitary $U$ to consist only
of \emph{two} distinct values. With this restriction, since $U$ has
$1/\sqrt{n}$ in the first row and column and since $U$ has to be
unitary, there is exactly one unitary matrix that has only two distinct
values in each column, and it turns out to be a \emph{symmetric} unitary,
with $U_{i,i}=\alpha$, and $U_{i,j}=\beta$ for $2\leq i,j\leq n$.

Having fixed the unitary matrix $U$, and the upper triangular matrix
$T$, our aim is to choose an $r$ as large as possible such that
the resulting matrix $A$ is nonnegative. Thus, optimization problem
\hyperref[opt:first]{1} becomes:

\begin{optprob}[2] \label{opt:sec}

\textbf{maximize} $r$

\textbf{such that:}

$\dfrac{1}{n}+r\sum_{k=2}^{n-1}u_{i,k}u_{j,k+1}\geq0$ for all $2\leq i,j\leq n$

$\dfrac{1}{n}+\dfrac{r}{\sqrt{n}}\sum_{k=2}^{n-1}u_{j,k+1}\geq0$
for all $2\leq j\leq n$

$\dfrac{1}{n}+\dfrac{r}{\sqrt{n}}\sum_{k=2}^{n-1}u_{i,k}\geq0$ for
all $2\leq i\leq n$
\end{optprob}

Subject to the constraints provided by these inequalities, we aim
to \emph{minimize} $\phi$ or \emph{maximize} $r$. Due to our restrictions
on $U$, the cut $S=\{1\}$ in the resulting matrix $A$ is special,
and we solve the above problem by \emph{minimizing} the edge expansion
$1-A_{1,1}$ of this cut. With the set of possible values that $r$
can take, we note that a set of extreme points of the resulting optimization
problem of minimizing $1-A_{1,1}$ or maximizing $A_{1,1}$ are obtained
if we \emph{force} the values of all the entries $A_{1,i}$ for $3\leq i\leq n$
to 0. We then maximize $r$ for the resulting matrix (indeed, there
are exactly two possible doubly stochastic matrices at this point),
and the result is the following construction.

\section{\label{subsec:Rootn-matrices-=002013}Rootn matrices -- fixing the
spectral gap and minimizing  edge expansion }
\medskip

\begin{construction}\label{constr:Rootn-matrices-=002013} Let $m=\sqrt{n}$,
\begin{align*}
a_{n} & =\dfrac{m^{2}+m-1}{m\cdot(m+2)}, & b_{n} & =\dfrac{m+1}{m\cdot(m+2)}, & c_{n} & =\dfrac{1}{m\cdot(m+1)},\\
\\
d_{n} & =\dfrac{m^{3}+2m^{2}+m+1}{m\cdot(m+1)\cdot(m+2)}, & e_{n} & =\dfrac{1}{m\cdot(m+1)\cdot(m+2)}, & f_{n} & =\dfrac{2m+3}{m\cdot(m+1)\cdot(m+2)},
\end{align*}

\noindent and define the $n\times n$ matrix
\[
A_{n}=\left[\begin{array}{ccccccccc}
a_{n} & b_{n} & 0 & 0 & 0 & 0 & 0 & \cdots & 0\\
0 & c_{n} & d_{n} & e_{n} & e_{n} & e_{n} & e_{n} & \cdots & e_{n}\\
0 & c_{n} & e_{n} & d_{n} & e_{n} & e_{n} & e_{n} & \cdots & e_{n}\\
0 & c_{n} & e_{n} & e_{n} & d_{n} & e_{n} & e_{n} & \cdots & e_{n}\\
0 & c_{n} & e_{n} & e_{n} & e_{n} & d_{n} & e_{n} & \cdots & e_{n}\\
0 & c_{n} & e_{n} & e_{n} & e_{n} & e_{n} & d_{n} & \cdots & e_{n}\\
\vdots & \vdots & \vdots & \vdots & \vdots & \vdots & \vdots & \ddots & \vdots\\
0 & c_{n} & e_{n} & e_{n} & e_{n} & e_{n} & e_{n} & \ldots & d_{n}\\
b_{n} & f_{n} & c_{n} & c_{n} & c_{n} & c_{n} & c_{n} & \cdots & c_{n}
\end{array}\right]
\]

\end{construction}

\begin{example} Construction \ref{constr:Rootn-matrices-=002013} for
$n=9$ looks like the following:

\[
A_{9}=\dfrac{1}{60}\begin{bmatrix}44 & 16 & 0 & 0 & 0 & 0 & 0 & 0 & 0\\
0 & 5 & 49 & 1 & 1 & 1 & 1 & 1 & 1\\
0 & 5 & 1 & 49 & 1 & 1 & 1 & 1 & 1\\
0 & 5 & 1 & 1 & 49 & 1 & 1 & 1 & 1\\
0 & 5 & 1 & 1 & 1 & 49 & 1 & 1 & 1\\
0 & 5 & 1 & 1 & 1 & 1 & 49 & 1 & 1\\
0 & 5 & 1 & 1 & 1 & 1 & 1 & 49 & 1\\
0 & 5 & 1 & 1 & 1 & 1 & 1 & 1 & 49\\
16 & 9 & 5 & 5 & 5 & 5 & 5 & 5 & 5
\end{bmatrix}.
\]
Note $\lambda_{2}(A_{9})=0$, and $\phi(A_{9})=\dfrac{16}{60}<\dfrac{1}{\sqrt{9}}.$
\end{example}

\begin{theorem}
\label{thm:rootn}The following hold for the matrices $A_{n}$ in
Construction \ref{constr:Rootn-matrices-=002013}:
\begin{enumerate}
\item $A_{n}$ is doubly stochastic.
\item Every nontrivial eigenvalue of $A_{n}$ is 0.
\item The edge expansion is bounded as 
\[
\dfrac{1}{6\sqrt{n}}\leq\phi(A_{n})\leq\dfrac{1}{\sqrt{n}}.
\]
\item As a consequence of 1,2,3, 
\[
\phi(A_{n})\leq\dfrac{1-\text{Re}\lambda_{2}(A_{n})}{\sqrt{n}}
\]
and thus 
\[
\Gamma(n)\leq\dfrac{1}{\sqrt{n}},
\]
exponentially improving upon the bound in Lemma \ref{lem:debruijn}.
\end{enumerate}
\end{theorem}

\begin{proof}
The following calculations are easy to check, to see that $A_{n}$
is a doubly stochastic matrix: 
\begin{enumerate}
\item $a_{n}\geq0$, $b_{n}\geq0$, $c_{n}\geq0$, $d_{n}\geq0$, $e_{n}\geq0$,
$f_{n}\geq0$. 
\item $a_{n}+b_{n}=1$. 
\item $c_{n}+d_{n}+(n-3)e_{n}=1$. 
\item $b_{n}+f_{n}+(n-2)c_{n}=1$. 
\end{enumerate}
This completes the proof of (1).

$A_{n}$ is triangularized as $T_{n}$ by the unitary $U_{n}$, i.e.
\[
A_{n}=U_{n}T_{n}U_{n}^{*},
\]
with $T_{n}$ and $U_{n}$ defined as follows. Recall that $m=\sqrt{n}$.
Let 
\[
r_{n}=1-\dfrac{1}{m+2},
\]

\begin{align*}
\alpha_{n} & =\dfrac{-n^{2}+2n-\sqrt{n}}{n\cdot(n-1)}=-1+\frac{1}{m\cdot(m+1)},
\end{align*}
\begin{align*}
\beta_{n} & =\dfrac{n-\sqrt{n}}{n\cdot(n-1)}=\frac{1}{m\cdot(m+1)},
\end{align*}
\[
T_{n}=\begin{bmatrix}1 & 0 & 0 & 0 & 0 & 0 & 0 & 0 & 0\\
0 & 0 & r_{n} & 0 & 0 & 0 & 0 & \cdots & 0\\
0 & 0 & 0 & r_{n} & 0 & 0 & 0 & \cdots & 0\\
0 & 0 & 0 & 0 & r_{n} & 0 & 0 & \cdots & 0\\
0 & 0 & 0 & 0 & 0 & r_{n} & 0 & \cdots & 0\\
0 & 0 & 0 & 0 & 0 & 0 & r_{n} & \cdots & 0\\
\vdots & \vdots & \vdots & \vdots & \vdots & \vdots & \ddots & \ddots & \vdots\\
0 & 0 & 0 & 0 & 0 & 0 & \cdots & 0 & r_{n}\\
0 & 0 & 0 & 0 & 0 & 0 & 0 & 0 & 0
\end{bmatrix},
\]
and 
\[
U_{n}=\begin{bmatrix}\frac{1}{\sqrt{n}} & \frac{1}{\sqrt{n}} & \frac{1}{\sqrt{n}} & \frac{1}{\sqrt{n}} & \frac{1}{\sqrt{n}} & \frac{1}{\sqrt{n}} & \cdots & \frac{1}{\sqrt{n}} & \frac{1}{\sqrt{n}}\\
\frac{1}{\sqrt{n}} & \alpha_{n} & \beta_{n} & \beta_{n} & \beta_{n} & \beta_{n} & \cdots & \beta_{n} & \beta_{n}\\
\frac{1}{\sqrt{n}} & \beta_{n} & \alpha_{n} & \beta_{n} & \beta_{n} & \beta_{n} & \cdots & \beta_{n} & \beta_{n}\\
\frac{1}{\sqrt{n}} & \beta_{n} & \beta_{n} & \alpha_{n} & \beta_{n} & \beta_{n} & \cdots & \beta_{n} & \beta_{n}\\
\frac{1}{\sqrt{n}} & \beta_{n} & \beta_{n} & \beta_{n} & \alpha_{n} & \beta_{n} & \cdots & \beta_{n} & \beta_{n}\\
\frac{1}{\sqrt{n}} & \beta_{n} & \beta_{n} & \beta_{n} & \beta_{n} & \alpha_{n} & \cdots & \beta_{n} & \beta_{n}\\
\vdots & \vdots & \vdots & \vdots & \vdots & \vdots & \ddots & \vdots & \vdots\\
\frac{1}{\sqrt{n}} & \beta_{n} & \beta_{n} & \beta_{n} & \beta_{n} & \beta_{n} & \cdots & \alpha_{n} & \beta_{n}\\
\frac{1}{\sqrt{n}} & \beta_{n} & \beta_{n} & \beta_{n} & \beta_{n} & \beta_{n} & \cdots & \beta_{n} & \alpha_{n}
\end{bmatrix}.
\]
To show that $U_{n}$ is a unitary, the following calculations can
be easily checked: 
\begin{enumerate}
\item $\frac{1}{n}+\alpha_{n}^{2}+(n-2)\cdot\beta_{n}^{2}=1$. 
\item $\frac{1}{\sqrt{n}}+\alpha_{n}+(n-2)\cdot\beta_{n}=0$. 
\item $\frac{1}{n}+2\cdot\alpha_{n}\cdot\beta_{n}+(n-3)\cdot\beta_{n}^{2}=0$. 
\end{enumerate}
Also, to see that $A_{n}=U_{n}T_{n}U_{n}^{*}$, the following calculations
are again easy to check: 
\begin{enumerate}
\item 
\[
A_{n}(1,1)=a_{n}=\langle u_{1},Tu_{1}\rangle=\frac{1}{n}+\frac{1}{n}\cdot(n-2)\cdot r_{n}.
\]
\item 
\[
A_{n}(1,2)=b_{n}=\langle u_{1},Tu_{2}\rangle=\frac{1}{n}+\frac{1}{\sqrt{n}}\cdot(n-2)\cdot r_{n}\cdot\beta_{n}.
\]
\item 
\[
A_{n}(n,1)=b_{n}=\langle u_{n},Tu_{1}\rangle=\frac{1}{n}+\frac{1}{\sqrt{n}}\cdot(n-2)\cdot r_{n}\cdot\beta_{n}.
\]
\item For $3\leq j\leq n$, 
\[
A_{n}(1,j)=0=\langle u_{1},Tu_{j}\rangle=\frac{1}{n}+\frac{1}{\sqrt{n}}\cdot\alpha_{n}\cdot r_{n}+(n-3)\cdot\frac{1}{\sqrt{n}}\cdot\beta_{n}\cdot r_{n}.
\]
\item For $2\leq i\leq n-1$, 
\[
A_{n}(i,1)=0=\langle u_{i},Tu_{1}\rangle=\frac{1}{n}+\frac{1}{\sqrt{n}}\cdot\alpha_{n}\cdot r_{n}+\frac{1}{\sqrt{n}}\cdot(n-3)\cdot\beta_{n}\cdot r_{n}.
\]
\item For $2\leq i\leq n-1$, 
\[
A_{n}(i,2)=c_{n}=\langle u_{i},Tu_{2}\rangle=\frac{1}{n}+\alpha_{n}\cdot\beta_{n}\cdot r_{n}+(n-3)\cdot\beta_{n}^{2}\cdot r_{n}.
\]
\item For $3\leq j\leq n$, 
\[
A_{n}(n,j)=c_{n}=\langle u_{n},Tu_{j}\rangle=\frac{1}{n}+\alpha_{n}\cdot\beta_{n}\cdot r_{n}+(n-3)\cdot\beta_{n}^{2}\cdot r_{n}.
\]
\item For $2\leq i\leq n-1$, 
\[
A_{n}(i,i+1)=d_{n}=\langle u_{i},Tu_{i+1}\rangle=\frac{1}{n}+\alpha_{n}^{2}\cdot r_{n}+(n-3)\cdot\beta_{n}^{2}\cdot r_{n}.
\]
\item For $2\leq i\leq n-2$, $3\leq j\leq n$, $i+1\not=j$, 
\[
A_{n}(i,j)=e_{n}=\langle u_{i},Tu_{j}\rangle=\frac{1}{n}+2\cdot\alpha_{n}\cdot\beta_{n}\cdot r_{n}+(n-4)\cdot\beta_{n}^{2}\cdot r_{n}.
\]
\item 
\[
A_{n}(n,2)=f_{n}=\langle u_{n},Tu_{2}\rangle=\frac{1}{n}+(n-2)\cdot r_{n}\cdot\beta_{n}^{2}.
\]
\end{enumerate}
We thus get a Schur decomposition for $A_{n}$, and since the diagonal
of $T_{n}$ contains only zeros except the trivial eigenvalue 1, we
get that all nontrivial eigenvalues of $A_{n}$ are zero. This completes
the proof of (2).

If we let the set $S=\{1\}$, then we get that 
\[
\phi(A_{n})\leq\phi_{S}(A_{n})=b_{n}<\frac{1}{\sqrt{n}}.
\]
Further, since $T_{n}$ can be written as $\Pi_{n}D_{n}$, where $D_{n}(1,1)=1$,
$D_{n}(i,i)=r_{n}$ for $i=2$ to $n-1$, and $D_{n}(n,n)=0$ for
some permutation $\Pi_{n}$, we get that $A_{n}=(U_{n}\Pi_{n})D_{n}U_{n}^{*}$
which gives a singular value decomposition for $A_{n}$ since $U_{n}\Pi_{n}$
and $U_{n}^{*}$ are unitaries. Thus, $A_{n}$ has exactly one singular
value that is 1, $n-2$ singular values that are $r_{n}$, and one
singular value that is 0. Thus, from Lemma \ref{lem:phi_sing}, we
get that 
\[
\phi(A)\geq\frac{1-r_{n}}{2}=\frac{1}{2\cdot\left(\sqrt{n}+2\right)}\geq\frac{1}{6\sqrt{n}}
\]
and this completes the proof of (3).
\end{proof}
\begin{remark}
We remark that for the matrices $A_{n}$ constructed in Theorem \ref{thm:rootn},
it holds that 
\[
\phi(A_{n})\leq\frac{1-|\lambda_{i}(A)|}{\sqrt{n}}
\]
for any $i\not=1$, giving a stronger guarantee than that required
for Theorem \ref{thm:intro-phi-delta}. \\
\end{remark}

We reiterate that it would be unlikely to arrive at such a construction
by algorithmic simulation, since the eigenvalues of the matrices $A_{n}$
are extremely sensitive. Although $\lambda_{2}(A_{n})=0$, if we shift
only $O(1/\sqrt{n})$ of the mass in the matrix $A_{n}$ to create
a matrix $A'_{n}$, by replacing $a_{n}$ with $a'_{n}=a_{n}+b_{n}$,
$b_{n}$ with $b'_{n}=0$, $f_{n}$ with $f'_{n}=f_{n}+b_{n}$ and
keeping $c_{n},d_{n},e_{n}$ the same, then $\lambda_{2}(A'_{n})=1$.
Thus, since perturbations of $O(1/\sqrt{n})$ (which is tiny for large
$n$) cause the second eigenvalue to jump from 0 to 1 (and the spectral
gap from 1 to 0), it would not be possible to make tiny changes to
random matrices to arrive at a construction satisfying the required
properties in Theorem \ref{thm:rootn}.

\section{\label{subsec:Observations-from-the}Observations from the construction
of Rootn Matrices}

At this point, we know that our lower bound on $\phi$ in theorem
\ref{thm:intro-phi-delta} is close to optimal, and a loss of $n^{\alpha}$
is necessary for $\frac{1}{2}\leq\alpha\leq1$. Our aim now is to
find a construction that achieves a bound $\Gamma(n)\in O(1/n)$.
However, this is going to take us in a direction different from the
one we have taken so far, but to start out in that direction, we need
some observations from the construction of Rootn Matrices. Our first
observation comes from observing the underlying markov chain. Observe
that the Rootn matrices in Construction \ref{constr:Rootn-matrices-=002013}
look approximately like the following, where $\epsilon\in\Theta(1/\sqrt{n})$.

\[
A_{n}=\left[\begin{array}{ccccccccc}
1-\epsilon & \epsilon & 0 & 0 & 0 & 0 & 0 & \cdots & 0\\
0 & \epsilon^{2} & 1-\epsilon & \epsilon^{3} & \epsilon^{3} & \epsilon^{3} & \epsilon^{3} & \cdots & \epsilon^{3}\\
0 & \epsilon^{2} & \epsilon^{3} & 1-\epsilon & \epsilon^{3} & \epsilon^{3} & \epsilon^{3} & \cdots & \epsilon^{3}\\
0 & \epsilon^{2} & \epsilon^{3} & \epsilon^{3} & 1-\epsilon & \epsilon^{3} & \epsilon^{3} & \cdots & \epsilon^{3}\\
0 & \epsilon^{2} & \epsilon^{3} & \epsilon^{3} & \epsilon^{3} & 1-\epsilon & \epsilon^{3} & \cdots & \epsilon^{3}\\
0 & \epsilon^{2} & \epsilon^{3} & \epsilon^{3} & \epsilon^{3} & \epsilon^{3} & 1-\epsilon & \cdots & \epsilon^{3}\\
\vdots & \vdots & \vdots & \vdots & \vdots & \vdots & \vdots & \ddots & \vdots\\
0 & \epsilon^{2} & \epsilon^{3} & \epsilon^{3} & \epsilon^{3} & \epsilon^{3} & \epsilon^{3} & \ldots & 1-\epsilon\\
\epsilon & \epsilon^{2} & \epsilon^{2} & \epsilon^{2} & \epsilon^{2} & \epsilon^{2} & \epsilon^{2} & \cdots & \epsilon^{2}
\end{array}\right].
\]

Note that for any vertex $k\in[3,n]$, since $A_{i,j}$ is the probability
of going from vertex $j$ to vertex $i$ (for right multiplication
by a vector), the chain from vertex $k$ almost always goes backward
from $k$ to $k-1$ with probability $\approx1-\epsilon$, and otherwise
goes approximately uniformly to any other vertex except vertex 1.
Thus we have the following.
\medskip

\begin{remark}\label{Observation:-The-Rootn} The
Rootn matrix in Construction \ref{constr:Rootn-matrices-=002013}
is \emph{almost} a cycle.
\end{remark}
\medskip

This observation in fact leads further to the following concrete lemma.

\begin{lemma}
\label{lem:all-sets-phi-same-rootn}For the Rootn Matrix $A_{n}$,
for \emph{any} set $S=\{i,i+1,\dots,j\}$, 
\[
\phi_{S}(A_{n})\in O\left(\dfrac{1}{\sqrt{n}}\right).
\]
\end{lemma}

\begin{proof}
This can immediately be observed by directly summing the entries in
the matrix $A$ in Construction \ref{constr:Rootn-matrices-=002013}.
\end{proof}
The main observation from Lemma \ref{lem:all-sets-phi-same-rootn}
is that although we pick a set of size $n/2$, say $\{2,...,n/2+1\}$,
there is contribution from \emph{exactly one entry} $d_{n}$ of about
$1$ that is in the complement set, and the rest of the $O(n^{2})$
numbers contribute about $O(\sqrt{n})$ mass that leaves the set,
and thus the average mass leaving the set is approximately $n^{2}\epsilon^{3}/(n/2)\in O(1/\sqrt{n})$.
Thus, if we had zeros above the diagonal instead of the entries $e_{n}$,
for a set of size $n/2$, we would have contribution from exactly
\emph{one} entry $d_{n}$ of about $1$, and the resulting  edge expansion 
of the matrix would become $O(1/n)$ as we would want. Note that the
entries $e_{n}$ are less meaningful in terms of the expansion, since
removing a factor of $n\cdot e_{n}\cdot J$ from the matrix will only
keep the entries $d_{n}$ in the off-diagonal block.

This takes us in a novel direction, since by observing construction
\ref{constr:Rootn-matrices-=002013} and noting \ref{Observation:-The-Rootn}
and Lemma \ref{lem:all-sets-phi-same-rootn}, we see that we can \emph{fix
the  edge expansion  of the matrix} and somehow \emph{set the entries of
the matrix to have desired eigenvalues}, going in a direction opposite
to that which we had for construction  \ref{constr:Rootn-matrices-=002013} of 
Rootn Matrices. One simple manner to achieve this is to have the following
matrix (discussed previously for another requirement) -- 
\[
A_{n}=\begin{bmatrix}\mathbf{0} & 1 & 0 & 0 & 0 & 0\\
0 & \mathbf{0} & 1 & 0 & 0 & 0\\
0 & 0 & \mathbf{0} & 1 & 0 & 0\\
0 & 0 & 0 & \ddots & \ddots & 0\\
0 & 0 & 0 & 0 & \mathbf{0} & 1\\
0 & 0 & 0 & 0 & 0 & \mathbf{0}
\end{bmatrix}
\]
which has  edge expansion  $O(1/n)$ and all its eigenvalues are 0. The matrix
however is not doubly stochastic, and making it so by putting $A_{n,1}=1$
makes it a cycle and which brings the spectral gap close to $\phi^{2}$,
putting it strictly within the Cheeger regime (see Section \ref{subsec:Tightness-of-cheeger-bounds}).
However, inspecting the heaviest permutation within the Rootn Matrix
(Construction \ref{constr:Rootn-matrices-=002013}) through observation
\ref{Observation:-The-Rootn}, we can start by fixing the \emph{structure}
of the matrix to be following:
\begin{equation}
\label{structure-chet-1}
A_{n}=\begin{bmatrix}
\mathbf{1-r} & r & 0 & 0 & 0 & 0\\
\cdot & \mathbf{\cdot} & r & 0 & 0 & 0\\
\cdot & \cdot & \cdot & r & 0 & 0\\
\cdot & \cdot & \cdot & \ddots & \ddots & 0\\
\cdot & \cdot & \cdot & \cdot & \mathbf{\cdot} & r\\
\cdot & \cdot & \cdot & \cdot & \cdot & \mathbf{1-r}
\end{bmatrix}
\end{equation}

Note that the matrix is almost exactly similar to the Rootn matrix
in Construction \ref{constr:Rootn-matrices-=002013} as remarked in observation
\ref{Observation:-The-Rootn}, where the entries represented by $\cdot$
are yet to be filled. Since the matrix has to be doubly stochastic,
$r\leq1$, and irrespective of the value of $r$, 
\[
\phi(A_{n})\in O(1/n),
\]
since by considering the set $S=\{1,...,n/2\}$, 
\[
\phi_{S}(A)=\dfrac{r}{n/2}\leq \dfrac{2}{n}\in O(1/n).
\]
The key thing to note is that this is almost exactly similar to the
structure of the triangular matrix in the Schur form of Rootn matrices, and in essence,
is the same structure as for Rootn matrices up to rotation by a unitary. All we need to do is
to fill the entries of $A_{n}$ so that it remains doubly stochastic,
and has the second eigenvalue much far from 1, which will give us
$\Gamma(n)\in O(1/n)$

However, having fixed this structure, even if we choose, say, a $5\times5$
matrix, it seems difficult to set variables to obtain some desired
eigenvalues, since the characteristic polynomial has degree $\geq 5$,
and it will not have an analytic solution, which although not necessary,
its lack makes the analysis of the matrix nearly impossible. To overcome this, the main
idea is as follows: instead of trying to directly control the eigenvalues
of the matrix, we try and \emph{control the eigenvalues indirectly}
by controlling the coefficients of the characteristic polynomial.
For any given matrix $A_{n}$, let the coefficient of $\lambda^{n-k}$
in the characteristic polynomial be $a_{k}$, then 
\[
a_{k}=\dfrac{(-1)^{k}}{k!}\begin{vmatrix}\text{Tr}A & 1 & 0 & 0 & 0 & 0\\
\text{Tr}A^{2} & \text{Tr}A & 2 & 0 & 0 & 0\\
\text{Tr}A^{3} & \text{Tr}A^{2} & \text{Tr}A & \ddots & 0 & 0\\
\ddots & \text{Tr}A^{3} & \text{Tr}A^{2} & \text{Tr}A & k-2 & 0\\
\text{Tr}A^{k-1} & \ddots & \text{Tr}A^{3} & \text{Tr}A^{2} & \text{Tr}A & k-1\\
\text{Tr}A^{k} & \text{Tr}A^{k-1} & \ddots & \text{Tr}A^{3} & \text{Tr}A^{2} & \text{Tr}A
\end{vmatrix}
\]
where $|\cdot|$ is the determinant. Thus, we realize that controlling
the traces of the powers of the matrix $A_{n}$ helps to control the
coefficients of the characteristic polynomial, which in turn controls
the eigenvalues. This is not surprising in retrospect, since $\text{Tr}A^{k}=\sum\lambda_{i}^{k}$,
but just this fact is not very useful unless the values of $\sum\lambda_{i}^{k}$
are simple. However, for us, with constraint \hyperref[Constraint-1:-Restrict]{1},
to have all eigenvalues 0, we only need to set $\text{Tr}A^{k}=1$
for all $k\geq1$, and this implies that the matrix has one (trivial)
eigenvalue 1, and all other eigenvalues of the matrix will be 0. Thus,
given the structure \ref{structure-chet-1}, we are going to set the
remaining values of $A$ to ensure that $\text{Tr}A^{k}=1$.

To set the values concretely, note that the underlying graph of our
matrix is the following, letting $A_{i,j}$ be the weight
of the edge from $j\rightarrow i$ (for right multiplication by a
vector): the only ``back'' edges in the graph are from vertex $i$
to $i-1$ (for $i\geq2$) of weight $r$, and all the other edges
are ``forward'' edges that go from $i$ to $j$ ($j\geq i$). Consider
the combinatorial meaning of $\text{Tr}A^{k}$. It sums the weight
of every length $k$ path (walk) from any vertex to itself (where the weight of a path is the weight of the product of the edges). Fix any
vertex $i$, and consider any path of length $k$ from $i$ to itself.
How far can this path go forward? Note that if the path went from
$i$ to $i+k$ in one step, since it can only go back one vertex in
one step, the path cannot go back to vertex $i$ in $k$ steps. Thus,
the maximum that a path of length $k$ starting at any vertex can
go forward in one step is $k-1$. Thus, in $\text{Tr}A^{k}$, there
is \emph{exactly} one path of length $k$ from vertex $i$ to itself
that goes back $k-1$ steps, and it has weight $r^{k-1}a_{i,i-k+1}$.
Similarly, the weight of any path of length $k$ that goes back $l$
steps (with $l\leq k-1$) will be a term of the form $r^{l}\cdot c$
where $c$ is the product of other forward edges in the graph along
the path. As a consequence, we get that for any $k\geq1$, 
\begin{equation}
\text{Tr}A^{k}=\sum_{i=0}^{k-1}r^{i}w_{i}(k),\label{eq:disc-chet-traceeq}
\end{equation}
or that $\text{Tr}A^{k}$ is a polynomial in $r$ of degree $k-1$
where $w_{i}(k)$ are functions of all the other nonzero entries of
$A$. As stated above, the coefficient of $r^{k-1}$ in $\text{Tr}A^{k}$
is easy to find, and is exactly 
\[
[r^{k-1}]\text{Tr}A^{k}=w_{k-1}(k)=\sum_{j=k-1}^{n}a_{j,j-k+1}
\]
and in particular, contains only those entries that are at a distance
of $k-1$ from the diagonal. This key fact is extremely important,
and we state it explicitly as a fact that we are going to exploit:

\begin{fact}\label{fact:chet-1} For the structure
of $A$ as discussed (see \ref{structure-chet-1}), the coefficient of $r^{k-1}$ in $\text{Tr}A^{k}$
is linear in the entries of the matrix.
\end{fact}

Observation \ref{fact:chet-1} suggests to us an \emph{inductive}
manner of assigning forward edge weights. Note that if we have set
all edge weights that are at a distance of $k-2$ from the diagonal
in $A$ using only equations $\text{Tr}A^{l}=1$ for $1\leq l\leq k-1$,
 then equation \ref{eq:disc-chet-traceeq} becomes 
\[
\text{Tr}A^{k}=W+r^{k-1}\cdot\sum_{j=k-1}^{n}a_{j,j-k+1}
\]
and setting it to 1 gives 
\begin{equation}
\sum_{j=k-1}^{n}a_{j,j-k+1}=\dfrac{1-W}{r^{k-1}},\label{eq:disc-chet-traceeq-2}
\end{equation}
where $W$ is the sum of weights of all paths that go at most $k-2$
times back.

Since $A^{n}=J$, we have $n-1$ independent equations, $\text{Tr}A^{k}=1$
for $1\leq k\leq n-1$. Given these $n-1$ equations,\emph{ if our
matrix consists only of $n-1$ different variables, we will be able
to set each variable using one equation}, and since we can use equation
\ref{eq:disc-chet-traceeq-2} (that comes from the equation $\text{Tr}A^{k}=1$)
to set the variables at distance $k-1$ from the diagonal, it gives
us our $n-2$ different variables. Thus, \emph{we set all variables
at distance $k$ from the origin to have value $c_{k}$}. This is
the most important and final conclusion.

We further need to ensure that our matrix remains doubly stochastic,
the first column and last row will have special values, but otherwise
the values are set as stated. We thus obtain the following matrix:
\[
A_{n}=\left[\begin{array}{ccccccccc}
b_{0} & r & 0 & 0 & 0 & 0 & 0 & \cdots & 0\\
b_{1} & c_{0} & r & 0 & 0 & 0 & 0 & \cdots & 0\\
b_{2} & c_{1} & c_{0} & r & 0 & 0 & 0 & \cdots & 0\\
b_{3} & c_{2} & c_{1} & c_{0} & r & 0 & 0 & \cdots & 0\\
b_{4} & c_{3} & c_{2} & c_{1} & c_{0} & r & 0 & \cdots & 0\\
\vdots & \vdots & c_{3} & c_{2} & c_{1} & c_{0} & r & \cdots & 0\\
b_{n-3} & c_{n-4} & \ddots & \ddots & \ddots & \ddots & \ddots & \ddots & \vdots\\
b_{n-2} & c_{n-3} & c_{n-4} & \cdots & c_{3} & c_{2} & c_{1} & c_{0} & r\\
b_{n-1} & b_{n-2} & b_{n-3} & \cdots & b_{4} & b_{3} & b_{2} & b_{1} & b_{0}
\end{array}\right]
\]
where for $0\leq i\leq n-2$
\[
b_{i}=1-r-\sum_{j=0}^{i-1}c_{i}
\]
and 
\[
b_{n-1}=1-\sum_{i=0}^{n-2}b_{i}.
\]

With this structure, we can set $c_{0}$ using $\text{Tr}A=1$, and
$c_{1}$ using $\text{Tr}A^{2}=1$, and similarly $c_{k-1}$ using $\text{Tr}A^{k}=1$,
and note that $\text{Tr}A^{k}$ contains $b_{k-1}$ but $b_{k-1}$
does not depend on $c_{k-1}$. Thus, using equation \ref{eq:disc-chet-traceeq-2},
we get each of the values of $c_{k}$ as functions (polynomials, albeit
with negative exponents) of $r$. And finally, the magic value of
$r$ that we use to evaluate each of the $c_{k}$'s, and that which
makes the entire construction work, is to set $r$ to be the $(n-1)$'th
root of $1/n$: 
\[
r=\left(\dfrac{1}{n}\right)^{\dfrac{1}{n-1}}
\]

The value of $r$ can be obtained in multiple ways, for instance by
using the equation $\text{Tr}A^{n-1}=1$ which will be a polynomial
only in $r$ since all of the $a_{i}$'s have already been fixed,
and choosing the positive real root of $r$ in the equation. Alternately
it can be obtained by taking $A^{n-1}=J$ and looking at the $a_{1,n}$.
We thus arrive at the following construction.

\section{\label{subsec:Chet-Matrices-=002013}Chet Matrices -- fixing the
 edge expansion  and maximizing the spectral gap}

\begin{construction}
\label{constr-chet}
For any $n$, let $C_{n}$ have the following structure.
\medskip

\[
C_{n}=\left[\begin{array}{ccccccccc}
b_{0} & r & 0 & 0 & 0 & 0 & 0 & \cdots & 0\\
b_{1} & c_{0} & r & 0 & 0 & 0 & 0 & \cdots & 0\\
b_{2} & c_{1} & c_{0} & r & 0 & 0 & 0 & \cdots & 0\\
b_{3} & c_{2} & c_{1} & c_{0} & r & 0 & 0 & \cdots & 0\\
b_{4} & c_{3} & c_{2} & c_{1} & c_{0} & r & 0 & \cdots & 0\\
\vdots & \vdots & c_{3} & c_{2} & c_{1} & c_{0} & r & \cdots & 0\\
b_{n-3} & c_{n-4} & \ddots & \ddots & \ddots & \ddots & \ddots & \ddots & \vdots\\
b_{n-2} & c_{n-3} & c_{n-4} & \cdots & c_{3} & c_{2} & c_{1} & c_{0} & r\\
b_{n-1} & b_{n-2} & b_{n-3} & \cdots & b_{4} & b_{3} & b_{2} & b_{1} & b_{0}
\end{array}\right]
\]
where for $0\leq i\leq n-2$,
\[
b_{i}=1-r-\sum_{j=0}^{i-1}c_{i}
\]
and 
\[
b_{n-1}=1-\sum_{i=0}^{n-2}b_{i}.
\]
Inductively for $k=0$ to $n-3$, set $c_k$ as a polynomial in $r$ using the equation $\text{Tr}(C_n)^{k+1}=1$. More specifically, let $C_{n}|_{k}$ be the matrix with structure of $C_{n}$, in which
all $c_{l}$ with $l\geq k$ are set to 0, and all $b_{l}$ with $l\geq k+1$
are set to 0. Inductively for $k=0$ to $n-3$, set the polynomial
\[
c_{k}=\dfrac{1-\text{Tr}(C_{n}|_{k})^{k+1}}{(n-k-2)\cdot(k+1)\cdot r^{k}}.
\]
Set $r$ to be the $(n-1)$'th root of $1/n$, i.e.,
\[
r=\left(\dfrac{1}{n}\right)^{\dfrac{1}{n-1}}.
\]
\end{construction}

We now show certain specific cases, so that we can actually compute
and see how the matrices $A_{n}$ look. 
\medskip

\begin{example}

For $n=5$, we get 
\begin{align*}
c_{0} & =\dfrac{2}{3}r-\dfrac{1}{3}\\
c_{1} & =\dfrac{5}{6}r-\dfrac{1}{3}\cdot\dfrac{1}{r}\\
c_{2} & =\dfrac{55}{27}r+\dfrac{2}{9}\cdot\dfrac{1}{r}-\dfrac{20}{27}\cdot\dfrac{1}{r^{2}}\\
r & =5^{-1/4}
\end{align*}
and the matrix approximated to 5 decimal places looks like:
\[
C_{5}\approx\begin{bmatrix} 0.33126 & 0.66874 \\
 0.21870 & 0.11249 & 0.66874 \\
 0.15993 & 0.05886 & 0.11249 & 0.66874 \\
 0.12170 & 0.03820 & 0.05886 & 0.11249 & 0.66874 \\
 0.16810 & 0.12170 & 0.15993 & 0.21870 & 0.33126 
\end{bmatrix}
\]
\end{example}

Similarly, we can write another example.
\medskip

\begin{example}
For $n=8$, we get
\begin{align*}
c_{0} & =\dfrac{1}{6}\cdot\left(2r-1\right)\\
c_{1} & =\dfrac{1}{5r}\cdot\left(\dfrac{4}{3}r^{2}-\dfrac{7}{12}\right)\\
c_{2} & =\dfrac{1}{4r^{2}}\cdot\left(\dfrac{152}{135}r^{3}+\dfrac{7}{45}r-\dfrac{14}{27}\right)\\
c_{3} & =\dfrac{1}{3r^{3}}\cdot\left(\dfrac{778}{675}r^{4}+\dfrac{49}{900}r^{2}+\dfrac{7}{27}r-\dfrac{1981}{3600}\right)\\
c_{4} & =\dfrac{1}{2r^{4}}\cdot\left(\dfrac{112}{75}r^{5}+\dfrac{133}{4050}r^{3}+\dfrac{7}{135}r^{2}+\dfrac{1981}{5400}r-\dfrac{1043}{1620}\right)\\
c_{5} & =\dfrac{1}{r^{5}}\cdot\left(\dfrac{27056}{10125}r^{6}+\dfrac{2863}{121500}r^{4}+\dfrac{28}{1215}r^{3}+\dfrac{1981}{81000}r^{2}+\dfrac{1043}{2430}r-\dfrac{193417}{243000}\right)\\
r & =8^{-1/7}
\end{align*}

and the matrix approximated to 5 decimals looks like:
\[
C_{8}\approx\begin{bmatrix}0.25700 & 0.74300\\
0.17601 & 0.08100 & 0.74300\\
0.13489 & 0.04112 & 0.08100 & 0.74300\\
0.10823 & 0.02666 & 0.04112 & 0.08100 & 0.74300\\
0.08901 & 0.01922 & 0.02666 & 0.04112 & 0.08100 & 0.74300\\
0.07431 & 0.01470 & 0.01922 & 0.02666 & 0.04112 & 0.08100 & 0.74300\\
0.06261 & 0.01170 & 0.01470 & 0.01922 & 0.02666 & 0.04112 & 0.08100 & 0.74300\\
0.09800 & 0.06261 & 0.07431 & 0.08901 & 0.10823 & 0.13489 & 0.17601 & 0.25700
\end{bmatrix}
\]

\end{example}


We give the equations for $n=16$ in Appendix \ref{subsec:Chet-Matrix-for-16}

Given construction  \ref{constr-chet}, we
have the following theorem.

\begin{theorem}
\label{lem:chet}Let $C_{n}$ be the matrix generated by algorithm
\ref{constr-chet}. Then the following hold for $C_{n}$.
\begin{enumerate}
\item (Conjecture) Every entry of $C_{n}$ is nonnegative.
\item All nontrivial eigenvalues of $C_{n}$ are 0.
\item $\phi(C_{n})\leq2/n\in O(1/n)$
\item As a consequence of 2 and 3, we get 
\[
\Gamma(n)\in O\left(\dfrac{1}{n}\right)
\]
improving the bound in Lemma \ref{thm:rootn} and combined with Theorem
\ref{thm:intro-phi-delta}, it implies that 
\[
\Gamma(n)\in\Theta\left(\dfrac{1}{n}\right)
\]
and that the lower bound on $\phi(R)$ in Theorem \ref{thm:intro-phi-delta}
is tight up to constants.
\end{enumerate}
\end{theorem}

\begin{proof}
Part (1) of lemma is a conjecture and we assume it for the remaining
parts. For (3), since $C_{n}$ is nonnegative due to (1), $\phi(C_{n})$
is well defined, and considering $S=\{1,...,n/2\}$ we have $\phi_{S}(C_{n})\leq\dfrac{r}{n/2}\leq\dfrac{2}{n}\in O(1/n)$.
For (2), since the sum of every row and column of $A_{n}$ is 1, it
has one eigenvalue 1 corresponding to the all ones eigenvector, and
further since $\text{Tr}C_{n}^{k}=\sum_{i=1}^{n}\lambda_{i}^{k}(C_{n})=1+\sum_{i=2}^{n}\lambda_{i}^{k}(C_{n})=1$
for all $k\geq1$, it implies all other eigenvalues of $C_{n}$ except
1 are 0. 
\end{proof}

There is only one uncertainty about this construction, and that is
to show that all $c_{i}$ and $b_{i}$ in Construction \ref{constr-chet}
are nonnegative. We do not know how to show this at present, and this
will be the main conjecture of this thesis, stated and discussed in
Section \ref{subsec:Conjectures-chet}. However, we try to list some
of the properties of Chet Matrices first.

\section{\label{subsec:Chet-Observations-and-Properties}Observations and
Properties of Chet Matrices}

We find Chet matrices beautiful enough to warrant their independent
study, and we shall explore many of their properties.

\noindent \textbf{1.} \emph{(Hessenberg-Toeplitz Structure)} Chet
matrices have a Hessenberg-Toeplitz structure. Hessenberg matrices
$H_{k,l}$ are matrices for which all the entries that are below the
diagonal at a distance greater than $k$ and entries above the diagonal
at a distance greater than $l$ are zero. Further, the matrix has
a Topelitz-like structure where entries are the same along every diagonal,
except the first column and the last row.
\par\medskip

\noindent \textbf{2.} \emph{(Permanent)} Let $C_{n}'$ be the matrix
which has all the entries of $C_{n}$, but the $n-1$ entries above
the diagonal that are $r$ are replaced by $-r$. Then 
\[
\text{permanent}(C_{n})=\text{determinant}(C_{n}')
\]
and the proof follows from straightforward induction.
\par\medskip

\noindent \textbf{3. }\emph{(Approximate values)} The approximate
values of the entries of $C_{n}$, by numerical and analytical observations,
are as follows:

\[
r\approx1-\dfrac{\log n}{n-1}.
\]
For the entries of $C_{n}$, 
\[
c_{i}\in O\left(\dfrac{1}{i\cdot n}\right).
\]
See Appendix \ref{subsec:Chet-Matrix-for-500} for a table of the
above for $n=500$. Thus, note that 
\[
\sum_{i=0}^{O(\log n)}c_{i}\in O\left(\dfrac{\log\log n}{n}\right),
\]
 and thus for $i\in O(\log n)$, 
\[
b_{i}=1-r-\sum_{j=0}^{i}c_{i}\approx\dfrac{\log n-\log\log n}{n}>0.
\]

A weaker statement is possible to show directly.

\begin{lemma}
For any Chet Matrix $C_{n}$, if $c_{i}\geq0$ and $b_{i}\geq0$ for
$i\in O(n/\ln n)$, then for $i\in O(n/\ln n)$,
\[
c_{i}\in O\left(\dfrac{1}{i\cdot(n-i)}\right).
\]
\end{lemma}

\begin{proof}
The proof follows immediately by noting that since $c_{i}$ is set
by using $\text{Tr}(C_{n}^{i+1})=1$, we get that 
\begin{align*}
1 & =\text{Tr}C_{n}^{i+1}=S+r^{i}\cdot c_{i}\cdot i\cdot(n-i)
\end{align*}
where $S$ is the weighted sum of all other paths of length $i+1$
that do not use $c_{i}$. Since $S\geq0$ by the assumption of $b_{i},c_{i}\geq0$,
and 
\[
r^{i}=n^{-\frac{i}{n-1}}\geq\exp\left(-\dfrac{\ln n}{n-1}\cdot i\right)\in\Omega(1)
\]
since $i\in O(n/\ln n)$, the claim immediately follows.
\end{proof}
\par\medskip

\noindent \textbf{4.} \emph{(Analytical values)} The first three entries,
computed analytically for any $n$, are as follows. Let $r=n^{-1/(n-1)}$.

\[
c_{0}(n)=\dfrac{2\cdot r-1}{n-2}
\]

\[
c_{1}(n)=\dfrac{(2r^{2}-1)\cdot n+1}{2(n-3)(n-2)\cdot r}
\]

\[
c_{2}(n)=\dfrac{(2r^{3}-1)\cdot n^{3}+(-8r^{3}+3r+4)\cdot n^{2}+(12r^{3}-15r-3)\cdot n+12r}{3(n-4)(n-3)(n-2)^{2}r^{2}}
\]
It can be shown that these are always nonnegative.
\par\medskip

\noindent \textbf{5.} \emph{(Distinct entries)} An interesting thing
about this construction, is that there are $\Theta(n)$ distinct entries
in the matrix, while Rootn matrices in Construction \ref{constr:Rootn-matrices-=002013}
has 6 distinct entries in the matrix, and de Bruijn has only
2. Note that starting with the structure of the matrix in construction
\ref{constr:Rootn-matrices-=002013}, and again using the equations
$\text{Tr}(A^{k})=1$ for $k=1$ to $6$ (since there are 6 distinct
variables), we obtain the same values of the variables as obtained
from the Schur transformation of the matrix, avoiding the Schur decomposition
entirely for Rootn matrices.
\par\medskip

\noindent \textbf{6. }\emph{(Route to the construction, and infinite
similar constructions)} The path taken
 to arrive at Construction \ref{constr-chet} of Chet Matrices is quite remarkable in retrospect.
We started with wanting to understand the edge expansion of matrices with
all eigenvalues 0 since it is not far from optimal (see Section \ref{subsec:Constraints-on-the}),
and ended up carefully choosing $T$ and $U$ in the schur decomposition
of doubly stochastic matrices to arrive at Rootn Matrices (construction
\ref{constr:Rootn-matrices-=002013}). The construction informed us
of a possible structure of matrices to consider, and with that structure
and subsequent reasoning we arrived at having $\text{Tr}A^{k}=1$
as our primary condition and diagonally setting the values in the
matrix to arrive at Chet Matrices. At this point, we can start with
any matrix with small edge expansion and use the same algorithm \ref{constr-chet}
to set values at distance $k$ from the diagonal by using $\text{Tr}A^{k}=1$.
Moreover, we can construct infinite such matrices with the structure
of Chet Matrices but with different values of $r$ and the $b_{i}$'s
and $c_{i}$'s, simply by choosing some set of eigenvalues -- say
$\frac{1}{2},\frac{1}{4},...,\frac{1}{2^{n-1}}$ -- and using corresponding
trace inequalities.
\par\medskip

\noindent \textbf{7.} \emph{(Matrices from approximate equations)}
Note that since it is difficult grasp how the coefficients of the
polynomials of the entries $c_{i}$ in the matrix are behaving, it
is tempting to try and approximate the entries of the Chet matrices
in Construction \ref{constr-chet} and see if the
approximation works. However, it turns out that the eigenvalues of
the matrix, as expected, are extremely sensitive, and the spectral
gap quickly diminishes to zero with approximation. Another idea, is
to instead approximate the equations themselves instead of the matrix
entries, by letting $\text{Tr}A^{k}\approx1$. However, it is simple
to show that this cannot work.

\begin{lemma}
There are doubly stochastic matrices $A$ with $\text{Tr}A^{k}\approx1$
but $\text{Re}\lambda_{2}(A)\approx1$.
\end{lemma}

\begin{proof}
Consider the adjacency matrix $C$ of the graph which consists of
1 isolated vertex, and a directed cycle of length $n-1$. Note that
$C$ is doubly stochastic and disconnected. The eigenvalues of $C$
are $1,\omega_{i}$ for $i=0$ to $n-2$, where $\omega_{i}=e^{\frac{2\pi\cdot i}{n-1}}$.
Let $A=\alpha J+(1-\alpha)C$. Note that the eigenvalues of $A$ are
$1,(1-\alpha)\omega_{i}$, for $0\leq i\leq n-2$, since the eigenvector
for the first eigenvalue 1 is the all ones vector, and $J$ is a projection
on it. The other eigenvalue that was 1, corresponding to $\omega_{0}$,
is in the space orthogonal to the all ones vector and just shrinks
by a factor of $(1-\alpha)$. Let $\alpha=\frac{2\log n}{n-1}$,
then the second eigenvalue of $A$ is $\lambda_{2}=1-\alpha=1-\frac{2\log n}{n-1}\approx1$.

We have that for any $k\leq n-2$, $\text{Tr}A^{k}=1+(1-\alpha)^{k}\sum_{i=0}^{n-2}\omega_{i}^{k}=1$
since $\sum_{i=0}^{n-2}\omega_{i}^{k}=0$ for $k\leq n-2$. For $k=n-1$,
we have 
\begin{align*}
\text{Tr}A^{n-1} & =1+(1-\alpha)^{n-1}\sum_{i=0}^{n-2}\omega_{i}^{n-1}\\
 & =1+(n-1)(1-\alpha)^{n-1}\\
 & \leq1+(n-1)\exp(-2\log n)\\
 & \leq1+\dfrac{1}{n}\\
 & \approx1
\end{align*}
which completes the proof. Note that such equations might help to
construct matrices with $\Gamma(n)\approx\dfrac{1}{\log n}$, but
will not help to obtain $\Gamma(n)\in o\left(\dfrac{1}{\log n}\right)$.
\end{proof}
\par\medskip

\noindent \textbf{8.} \emph{(Jordan Form)} We now try and understand
the Jordan form of Chet Matrices. The Jordan forms of Hessenberg matrices
are well understood (see \cite{zemke2006hessenberg} for instance)
and are easy to derive using the resolvent $(t\cdot I-A)^{-1}$, and
we only state the final forms forms here. Assume that $A=VQV^{-1}$
is the Jordan decomposition of $A$, where $Q$ contains two blocks,
1 block of length 1 for eigenvalue 1 and another of length $n-1$
for eigenvalue 0. Specifically, $Q$ contains all zeros, except $Q_{1,1}=1$,
and $Q_{i,i+1}=1$ for $i=2$ to $n-1$. The columns of $V$ are the
most interesting parts, and can be described as follows. Let $A_{[1,j]}$
denote the submatrix of $A$ consisting of all rows and columns with
indices $\{1,2,...,j\}$, with $A_{[1,0]}=1$. Let the columns of
$V$ be $V_{i}$ for $1\leq i\leq n$. Note that since $A$ is doubly
stochastic, the first column $V_{1}$ is just a multiple of the all
ones vector. Consider the second column. The entries of the second
column, for $i=1$ to $n$, are given by 
\[
V_{2}(i)=\dfrac{\det(t\cdot I-A_{[1,i-1]})}{r^{i-1}}\Bigg|_{t=0}
\]
(Note $I$ has dimension that of $A_{[1,i-1]}$, thus for $i=1$ the
numerator is simply $-1$). Similarly, the third column, $V_{3}$,
can be described as follows, 
\[
V_{3}(i)=\dfrac{1}{r^{i-1}}\cdot\dfrac{d}{dt}(\det(t\cdot I-A_{[1,i-1]}))\Bigg|_{t=0}
\]
and in general, for $2\leq j\leq n$,
\[
V_{j}(i)=\dfrac{1}{(j-2)!\cdot r^{i-1}}\cdot\dfrac{d^{(j-2)}}{dt^{(j-2)}}(\det(t\cdot I-A_{[1,i-1]}))\Bigg|_{t=0}.
\]
Thus, note that even the Jordan form has Hessenberg structure! This
follows simply because $\det(t\cdot I-A_{[1,i-1]})$ has degree $i-1$
in $t$, and when it is differentiated more than $i-1$ times, the
term becomes zero. We find this property beautiful, but do not know
how to employ it to show nonnegativity of variables at present.
\par\medskip

\noindent \textbf{9.} \emph{(Generating functions)} We now express
all the entries of the matrix in terms of generating functions, in
order to link it to the combinatorics of the underlying chain, and
obtain closed form expressions for some of the quantities. Since the
matrix has all eigenvalues 0 except the trivial eigenvalue, let $\chi_{t}(C_{n})$
be the characteristic polynomial of $C_{n}$, then $\chi_{t}(C_{n})=t^{n-1}\cdot(t-1)$,
and thus except the coefficient of $t^{n}$ ($=1$) and $t^{n-1}$
($=-1$), all the coefficients are 0. However, we know that $\chi_{t}(C_{n})=\det(t\cdot I-C_{n})$,
and since $t\cdot I-C_{n}$ also has the Hessenberg structure (see
point (1)), letting $C'_{n}$ denote $C_{n}$ with the entries $r$
replaced by $-r$, we get from point (2) that $\text{perm}(t\cdot I-C'_{n})=\det(t\cdot I-C{}_{n})$.
Thus, we have that 
\begin{equation}
\text{perm}(t\cdot I-C'_{n})=t^{n-1}\cdot(t-1).\label{eq:gen-func-chet-charpoly}
\end{equation}
The important thing about the permanent is that the underlying combinatorics
is much more amenable to generating functions than the determinant.
The permanent is the sum of weights of all the cycle covers of the
matrix which is much simpler to deal with than the number of weighted
cycles of length $k$ (for all $k$) as is the case with the determinant.

We are going to let $Z_{n}=t\cdot I-C'_{n}$, and count the cycle
covers in $Z$. Let $Z_{1,1}=Z_{n,n}=q_{0}=t-b_{0}$, $Z_{i,i}=a_{0}=t-c_{0}$
for $2\leq i\leq n-1$, and we let $q_{i}=-b_{i}$ and $a_{i}=-c_{i}$
for $i\geq2$ for the corresponding entries in $Z$. Note that $Z$
has $r$ above the diagonal, since $C'_{n}$ had $-r$. To count the
cycle covers in $Z_{n}$, note that due to the structure of our matrix,
starting at vertex $i$, if a cycle goes back $k$ steps to vertex
$i-k$, then one of the cycles must be the cycle $\{i,i-k,i-k+1,i-k+2,\ldots,i-1,i\}$
of length $k+1$, since this is the only manner in which the cycle
can be a part the cover. Note that the weight of this cycle is $c\cdot r^{k}$
where $c$ denotes the weight of the edge from vertex $i$ to vertex
$i-k$. Thus, note that any particular cycle disconnects the graph
into two parts, and thus, we can define the following two polynomials:
\begin{align*}
R(x) & =\sum_{i\geq0}a_{i}r^{i}x^{i+1}\\
G(x) & =\sum_{i\geq0}q_{i}r^{i}x^{i+1}
\end{align*}
where assume that $a_{l}=0$ for $l\geq n-1$ and $q_{l}=0$ for $l\geq n$.
Note that $R$ and $G$ do not have constant terms, and since the
$b_{i}$'s were chosen to ensure the first column and last row of
$C_{n}$ sum to 1, we have that 
\begin{align}
c_{i} & =b_{i}-b_{i+1}\nonumber \\
a_{i} & =q_{i}-q_{i+1}\nonumber \\
R & =G-\dfrac{G-q_{0}x}{rx}\nonumber \\
G-R & =\dfrac{G'}{rx}\label{eq:gen-func-chet-0}
\end{align}
where $G'=G-q_{0}x$. This equation will be useful later.

The coefficients of $x^{k+1}$ in the two polynomials for $k\geq1$
are 
\begin{align*}
[x^{k}]R(x) & =a_{k}r^{k}\\{}
[x^{k}]G(x) & =q_{k}r^{k}.
\end{align*}

The first coefficient is exactly the length of a cycle of length $k+1$
starting at any vertex except the first and last vertex, and the second
coefficient is the length of a cycle of length $k+1$ starting at
the first or last vertex. Thus, since the permanent of $Z$ is the
sum of the weights of all cycle covers of length $n$ in our graph,
we get that 
\begin{align*}
H(x) & =G(x)^{2}\sum_{i\geq0}R(x)^{i}+G(x)\\
 & =\dfrac{G(x)^{2}}{1-R(x)}+G(x)
\end{align*}
where the first term counts the weights of all cycle covers by fixing
the cover, and considering the cycle in which the first vertex appears,
the cycle in which the last vertex appears, and then the number of
different ways in which to divide the remaining vertices into $i$
cycles, and the second term counts the one cycle cover in which all
vertices appear in a single cycle. Note that the coefficients of $x^{k}$
in $H(x)$ are polynomials in $t$. Thus we have
\begin{align*}
t^{n-1}\cdot(t-1)=\text{perm}(Z) & =[x^{n}]H(x)\\
 & =[x^{n}]\dfrac{G(x)^{2}}{1-R(x)}-b_{n}r^{n-1}
\end{align*}

Thus, to compare coefficients on both sides of the above equation,
we have the following $n+1$ equations:
\begin{align}
1 & =\dfrac{1}{n!}\cdot\dfrac{d^{n}}{dt^{n}}[x^{n}]\dfrac{G(x)^{2}}{1-R(x)}\Bigg|_{t=0}\nonumber \\
-1 & =\dfrac{1}{(n-1)!}\cdot\dfrac{d^{n-1}}{dt^{n-1}}[x^{n}]\dfrac{G(x)^{2}}{1-R(x)}\Bigg|_{t=0}\nonumber \\
0 & =\dfrac{1}{k!}\cdot\dfrac{d^{k}}{dt^{k}}[x^{n}]\dfrac{G(x)^{2}}{1-R(x)}\Bigg|_{t=0}\label{eq:gen-func-chet-1}\\
0 & =[x^{n}]\dfrac{G(x)^{2}}{1-R(x)}\Bigg|_{t=0}-b_{n}r^{n-1}\nonumber 
\end{align}
where the second last equation holds for $1\leq k\leq n-2$. The last
equation is redundant, since it represents $0=\det(t\cdot I-C_{n})|_{t=0}$
which we know to be true based on the way we have set values of the
$c_{i}$'s and $b_{i}$'s in $C_{n}$. To understand other equations,
note that there is only one term/coefficient in $G(x)$ that is dependent
on $t$, namely $q_{0}x=(t-b_{0})x$ and similarly for $R(x)$ the
term $a_{0}x=(t-c_{0})x$, and thus 
\begin{align*}
\dfrac{d}{dt}G(x)^{2} & =2G(x)\cdot x\\
\dfrac{d^{2}}{dt^{2}}G(x)^{2} & =2x^{2}\\
\dfrac{d^{i}}{dt^{i}}G(x)^{2} & =0
\end{align*}
for $i\geq3$, and similarly we can write,
\begin{align*}
\dfrac{d}{dt}(1-R(x))^{-1} & =(1-R(x))^{-2}x\\
\dfrac{d^{i}}{dt^{i}}(1-R(x))^{-1} & =i!\cdot x^{i}\cdot(1-R(x))^{-(i+1)}
\end{align*}
for all $i$. Thus we get 
\begin{align*}
\dfrac{d^{k}}{dt^{k}}\dfrac{G(x)^{2}}{1-R(x)} & =\sum_{i=0}^{k}{k \choose i}\dfrac{d^{i}}{dt^{i}}(1-R)^{-1}\dfrac{d^{k-i}}{dt^{k-i}}G^{2}\\
 & =k!x^{k}(1-R)^{-k+1}G^{2}\ +\\
 & \ \ k(k-1)!x^{k-1}(1-R)^{-k}\cdot2\cdot G\cdot x\ +\\
 & \ \ {k \choose 2}(k-2)!x^{k-2}(1-R)^{-(k-1)}2x^{2}\\
 & =k!\cdot x^{k}\cdot\left(G^{2}(1-R)^{-(k+1)}+2G(1-R)^{-k}+(1-R)^{-(k-1)}\right)\\
 & =k!\cdot x^{k}\cdot(1-R)^{-(k+1)}\cdot(G+1-R)^{2}.
\end{align*}
Thus, from equation \ref{eq:gen-func-chet-1}, we have
\begin{align}
\dfrac{1}{k!}\cdot\dfrac{d^{k}}{dt^{k}}[x^{n}]\dfrac{G(x)^{2}}{1-R(x)} & =\dfrac{1}{k!}\cdot[x^{n}]\dfrac{d^{k}}{dt^{k}}\dfrac{G(x)^{2}}{1-R(x)}\nonumber \\
 & =[x^{n-k}](1-R)^{-(k+1)}\cdot(G+1-R)^{2}\nonumber \\
 & =\sum_{j=0}^{n-k}[x^{j}](1-R)^{-(k+1)}\cdot[x^{n-k-j}]\left(G+1-R\right)^{2}.\label{eq:gen-func-chet-2}
\end{align}
To compute each of the terms in the above equation, we can write 
\begin{align*}
[x^{j}](1-R)^{-(k+1)} & =[x^{j}]\left(\sum_{i=0}{k+i \choose k}R^{i}\right)\\
 & =\sum_{i=0}{k+i \choose k}[x^{j}]R^{i}\\
 & =\sum_{i=0}^{j}{k+i \choose i}[x^{j}]R^{i}\\
 & \ \ \ [\text{since de\ensuremath{g_{x}}(R)\ensuremath{\geq1}}]\\
 & =[x^{j}]\sum_{i=0}^{j}{k+i \choose i}R^{i}
\end{align*}
Note that the first few terms look like the following:
\begin{align*}
[x^{0}](1-R)^{-(k+1)} & =1\\{}
[x^{1}](1-R)^{-(k+1)} & =(k+1)\cdot a_{0}\\{}
[x^{2}](1-R)^{-(k+1)} & ={k+1 \choose 1}\cdot ra_{1}+{k+2 \choose 2}\cdot a_{0}^{2}\\{}
[x^{3}](1-R)^{-(k+1)} & ={k+1 \choose 1}\cdot r^{2}a_{2}+{k+2 \choose 2}\cdot r\cdot2a_{0}a_{1}+{k+3 \choose 3}\cdot a_{0}^{3}\\{}
[x^{4}](1-R)^{-(k+1)} & ={k+1 \choose 1}\cdot r^{3}a_{3}+{k+2 \choose 2}\cdot r^{2}\cdot\left(2a_{0}a_{2}+a_{1}^{2}\right)\\
 & \ \ +{k+3 \choose 3}\cdot r\cdot\left(3a_{1}a_{0}^{2}\right)+{k+4 \choose 4}\cdot a_{0}^{4}
\end{align*}
Similarly, for the second term in \ref{eq:gen-func-chet-2}, the coefficients
are much easier to compute.
\begin{align*}
[x^{0}]\left(G+1-R\right)^{2} & =1\\{}
[x^{n-k-j}]\left(G+1-R\right)^{2} & =[x^{n-k-j}]\left(\dfrac{G'}{rx}+1\right)^{2}\\
 & \ \ \ [\text{from \ref{eq:gen-func-chet-0}}]\\
 & =r^{n-k-j-2}\left(2\cdot q_{n-k-j}\cdot r+\sum_{i=1}^{n-k-j-1}q_{i}q_{n-k-j-i}\right)
\end{align*}
As a consequence, our final expression for equation \ref{eq:gen-func-chet-1}
becomes the following
\begin{align*}
0 & =\dfrac{1}{k!}\cdot\dfrac{d^{k}}{dt^{k}}[x^{n}]\dfrac{G(x)^{2}}{1-R(x)}\Bigg|_{t=0}\\
 & =\sum_{j=0}^{n-k}[x^{j}](1-R)^{-(k+1)}\cdot[x^{n-k-j}]\left(G+1-R\right)^{2}\\
 & =\sum_{j=0}^{n-k}r^{n-k-j-2}\left(2\cdot q_{n-k-j}\cdot r+\sum_{i_{1}=1}^{n-k-j-1}q_{i_{1}}q_{n-k-j-i_{1}}\right)\cdot[x^{j}]\sum_{i_{2}=0}^{j}{k+i_{2} \choose i_{2}}R^{i_{2}}
\end{align*}
This gives us a closed form expression for the coefficients of the
characteristic polynomial \ref{eq:gen-func-chet-charpoly}, and can
be used to infer the values of the $c_{i}$'s of Chet Matrices in
Construction \ref{constr-chet}, but at present we
do not know how to use these equations to show the nonnegativity of
the variables using these equations. We now formally state some of
our conjectures for which these equations might be useful.

\section{\label{subsec:Conjectures-chet}Conjectures for proving nonnegativity
of Chet Matrices}

We state the main conjecture of this thesis in this section, and a
sequence of related conjectures, each interesting in their own right.
\begin{conjecture}
(Chet Conjecture) \label{conj:(Chet-Conjecture)}Let $C_{n}$ denote
the $n\times n$ Chet matrix, and let $C=\{n:C_{n}\text{ is entry-wise nonnegative}\}$.
Then the following is true:
\[
|C|=\infty.
\]
More strongly,
\[
C=\mathbb{N}.
\]
\end{conjecture}

The stronger conjecture might be easier to prove or disprove.

\noindent \textbf{Numerical Simulation.} We analytically simulated
the Chet matrices $C_{n}$ with exact precision up till $n=21$ in
Maple, and all the matrices are nonnegative, and we numerically simulated
the matrices up till $n=500$ with 100 digits of precision in Matlab using the pseudocode in Appendix \ref{chet-pseudocode},
and not only are the matrices nonnegative, but the entries $c_{i}$
are substantially far from 0, and they decay gracefully as $i$ increases.
For $n=500$, the first nonzero digit appears $6$ places after the
decimal in the smallest matrix entry, and the associated plots are given in Appendix \ref{subsec:Chet-Matrix-for-500}. The plots indicate that the numbers $c_i$ and $b_i$ behave
very smoothly, and there are no sudden jumps or discontinuities in
their values as $i$ increases, strongly
supporting the (human) intuition that the Chet Conjecture is true.
\medskip

We state some ideas, progress, and other conjectures that will help
prove the Chet Conjecture.
\medskip

\noindent \textbf{Trace conjectures.} We formulate some interesting
conjectures which if true, would be helpful (although not by direct
implication but by general understanding) in establishing the nonnegativity
of Chet Matrices $C_{n}$. The conjectures are as follows.
\begin{conjecture}
\noindent (Trace Conjecture)\label{conj:(Trace-Conjecture)} Let $A$
be a nonnegative matrix that is substochastic, that is, $\sum_{i}A_{i,j}\leq1$
and $\sum_{j}A_{i,j}\leq1$ for all $j$ and $i$. Assume the following:
Above the diagonal, $A$ has nonzero entries only for entries that
are at a distance of 1 from the diagonal, and below the diagonal,
$A$ has nonzero entries only for entries that are at a distance at
most $k$ from the diagonal, where the diagonal has distance zero
from the diagonal. Assume $\text{Tr}A^{l}\leq1$ for $l\leq k+1$,
then $\text{Tr}A^{l}\leq1$ for all $l$.
\end{conjecture}

It is easy to see that the Trace Conjecture \ref{conj:(Trace-Conjecture)}
holds for symmetric matrices provided $k\geq1$, since for a symmetric
matrix, all eigenvalues are real, and 
\[
\text{Tr}A^{2}=\sum\lambda_{i}^{2}\leq1
\]
implies that $|\lambda_{i}|\leq1$, and for any $l$, 
\[
\text{Tr}A^{l}=\sum\lambda_{i}^{l}\leq\sum_{i}|\lambda_{i}|^{l}\leq\sum_{i}|\lambda_{i}|^{2}\leq1.
\]

Note that if conjecture is true, it would mean that it is always possible
to fill in nonnegative values in the entries at distance $k+1$ from
the diagonal and ensure that $\text{Tr}A^{k+2}=1$, and inductively,
fill the matrix with nonnegative entries. The Trace conjecture does
not directly imply the Chet conjecture due to the presence of the
first row and last column in Chet matrices which are required to ensure
that the matrices $C_{n}$ remain doubly stochastic, and thus the
entries $C_{n}[k+1,1]$ and $C_{n}[n,n-(k+1)]$ are always $\geq0$
for any fixed $k$, while they are zero in the trace conjecture. In
spite of this, we feel that attempts to prove the Trace conjecture
will be helpful in proving the Chet conjecture.

The Trace Conjecture \ref{conj:(Trace-Conjecture)} is interesting
in its own right, and since we do not know its proof, the following
conjecture should be relatively easier to prove.
\begin{conjecture}
\noindent (Toeplitz Trace Conjecture) \label{conj:(Toeplitz-Trace-Conjecture)}Let
$A$ be a nonnegative matrix that is substochastic, that is, $\sum_{i}A_{i,j}\leq1$
and $\sum_{j}A_{i,j}\leq1$ for all $j$ and $i$. Assume the following:
Above the diagonal, $A$ has nonzero entries only for entries that
are at a distance of 1 from the diagonal, and below the diagonal,
$A$ has nonzero entries only for entries that are at a distance at
most $k$ from the diagonal, with entry $c_{l}$ at distance $l$
below the diagonal, where the diagonal has distance zero from the
diagonal. Assume $\text{Tr}A^{l}\leq1$ for $1\leq l\leq k+1$, then
$\text{Tr}A^{l}\leq1$ for all $l$.
\end{conjecture}

The infinite version of the Toeplitz Trace Conjecture is also interesting,
which also remains open as of this writing.
\begin{conjecture}
\noindent (Infinite Toeplitz Trace Conjecture)\label{conj:(Infinite-Toeplitz-Trace}
Consider the following Markov Chain on $\mathbb{N}$. From every vertex
$i$, there is an edge of weight $r$ to $i-1$, and edges of weight
$c_{l}$ to vertices $i+l$ for $0\leq l\leq k$, and assume that
$r+\sum_{i=0}^{k}c_{i}\leq1.$ Let $p_{l}$ be the probability of
starting at $i$ and returning to $i$ in exactly $l$ steps, and
let $p_{l}\leq\alpha/k$ for $1\leq l\leq k+1$ for some constant
$\alpha$. Then $p_{l}\leq\alpha/k$ for all $l$.
\end{conjecture}

It is also interesting to understand the behavior of eigenvalues of
Hessenberg or Hessenberg-Toeplitz matrices, which would provide further
insights into the Trace Conjectures \ref{conj:(Trace-Conjecture)},
\ref{conj:(Toeplitz-Trace-Conjecture)} and \ref{conj:(Infinite-Toeplitz-Trace},
and it leads to the following conjecture.
\begin{conjecture}
(Eigenvalues of Toeplitz-Hessenberg Matrices) Let $A$ be a nonnegative
matrix that is substochastic, that is, $\sum_{i}A_{i,j}\leq1$ and
$\sum_{j}A_{i,j}\leq1$ for all $j$ and $i$. Assume the following:
Above the diagonal, $A$ has nonzero entries only for entries that
are at a distance of 1 from the diagonal, and below the diagonal,
$A$ has nonzero entries only for entries that are at a distance at
most $k$ from the diagonal, where the diagonal has distance zero
from the diagonal. Then the eigenvalues of $A$ are contained in the
convex hull of the $k+1$ roots of unity in the complex plane.
\end{conjecture}

\noindent It is possible to construct many similar conjectures, and
insights into any would possibly be helpful to prove our original
Chet Conjecture \ref{conj:(Chet-Conjecture)}.

\chapter{Connections and Extensions}

\vspace{1.5cm}
\begin{quote}
But as more arts were invented, and some were directed to the necessities
of life, others to recreation, the inventors of the latter were naturally
always regarded as wiser than the inventors of the former, because
their branches of knowledge did not aim at utility. Hence when all
such inventions were already established, the sciences which do not
aim at giving pleasure or at the necessities of life were discovered,
and first in the places where men first began to have leisure. This
is why the mathematical arts were founded in Egypt; for there the
priestly caste was allowed to be at leisure.
\begin{flushright}
\textasciitilde{} Aristotle,  \emph{Metaphysics}
\par\end{flushright}

\end{quote}
\vspace{3cm}

In this chapter, we study some other combinatorial quantities and
relate them to edge expansion and spectral gap. Many of these connections
are also not well-studied for nonreversible chains (with some exceptions,
notably Mihail's result connecting the mixing time and edge expansion),
and our aim will be to derive new results or provide elementary proofs
of known results. Further, another aim of this section is to provide
a cohesive treatment of some topics that have appeared in different
places in different communities, and fill-in the missing lemmas. We
start by studying the most well-studied (implying we will have little
to contribute) combinatorial quantity, the Mixing time.

\section{\label{subsec:Mixing-Time}Mixing time  -- Definitions and preliminary lemmas}

To motivate the definition of mixing time for general nonnegative
matrices, we first consider the mixing time of doubly stochastic matrices.
The mixing time of a doubly stochastic matrix $A$ (i.e., of the underlying
Markov chain) is the worst-case number of steps required for a random
walk starting at any vertex to reach a distribution approximately
uniform over the vertices. To avoid complications of periodic chains,
we assume that $A$ is $\frac{1}{2}$-lazy, meaning that for every
$i$, $A_{i,i}\geq\frac{1}{2}$. Given any doubly stochastic matrix
$A$, it can be easily converted to the lazy random walk $\frac{1}{2}I+\frac{1}{2}A$.
This is still doubly stochastic and in the conversion both $\phi(A)$
and the spectral gap are halved. The mixing time will be finite provided
only that the chain is connected. Consider the indicator vector $\mathbf{1}_{\{i\}}$
for any vertex $i$. We want to find the smallest $\tau$ such that
$A^{\tau}\mathbf{1}_{\{i\}}\approx\frac{1}{n}\mathbf{1}$ or $A^{\tau}\mathbf{1}_{\{i\}}-\frac{1}{n}\mathbf{1}\approx0$,
which can further be written as $\left(A^{\tau}-\frac{1}{n}\mathbf{1}\cdot\mathbf{1}^{T}\right)\mathbf{1}_{\{i\}}\approx0$.
Concretely, for any $\epsilon$, we want to find $\tau=\tau_{\epsilon}(A)$
such that for any $i$, 
\[
\left\Vert \left(A^{\tau}-\frac{1}{n}\mathbf{1}\cdot\mathbf{1}^{T}\right)\mathbf{1}_{\{i\}}\right\Vert _{1}\leq\epsilon.
\]
Given such a value of $\tau$, for any vector $x$ such that $\|x\|_{1}=1$,
we get 
\begin{align*}
\left\Vert \left(A^{\tau}-\frac{1}{n}\mathbf{1}\cdot\mathbf{1}^{T}\right)x\right\Vert _{1} & =\left\Vert \sum\limits _{i}\left(A^{\tau}-\frac{1}{n}\mathbf{1}\cdot\mathbf{1}^{T}\right)x_{i}\mathbf{1}_{\{i\}}\right\Vert _{1}\\
 & \leq\sum\limits _{i}|x_{i}|\left\Vert \left(A^{\tau}-\frac{1}{n}\mathbf{1}\cdot\mathbf{1}^{T}\right)\mathbf{1}_{\{i\}}\right\Vert _{1}\\
 & \leq\sum\limits _{i}|x_{i}|\cdot\epsilon\\
 & =\epsilon.
\end{align*}
Thus, the mixing time $\tau_{\epsilon}(A)$ is the number $\tau$
for which $\left\Vert (A^{\tau}-J)\cdot x\right\Vert _{1}\leq\epsilon$
for any $x$ such that $\|x\|_{1}=1$.

We want to extend this definition to any nonnegative matrix $R$ with
PF eigenvalue 1 and corresponding positive left and right eigenvectors
$u$ and $v$. Fixing the largest eigenvalue to 1 was irrelevant while
studying edge expansion and spectral gap since both the quantities are
unchanged by scalar multiplications, but it will be relevant for studying
mixing time to avoid the involvement of the largest eigenvalue in
every equation. Note that if $R$ is reducible, then the mixing time
is infinite. Further, if $R$ is periodic, then mixing time is again
ill-defined. Thus, we again assume that $R$ is irreducible and $\frac{1}{2}$-lazy,
i.e. $R_{i,i}\geq\frac{1}{2}$ for every $i$. Let $x$ be any nonnegative
vector for the sake of exposition, although our final definition will
not require nonnegativity and will hold for any $x$. We want to find
$\tau$ such that $R^{\tau}x$ about the same as the component of
$x$ along the direction of $v$. Further, since we are right-multiplying
and want convergence to the right eigenvector $v$, we will define
the $\ell_{1}$-norm using the left eigenvector $u$. Thus, for the
starting vector $x$, instead of requiring $\|x\|_{1}=1$ as in the
doubly stochastic case, we will require $\|D_{u}x\|_{1}=1$. Since
$x$ is nonnegative, $\|D_{u}x\|_{1}=\langle u,x\rangle=1$. Thus,
we want to find $\tau$ such that $R^{\tau}x\approx v$, or $\left(R^{\tau}-v\cdot u^{T}\right)x\approx0$.
Since we measured the norm of the starting vector $x$ with respect
to $u$, we will also measure the norm of the final vector $\left(R^{\tau}-v\cdot u^{T}\right)x$
with respect to $u$. Thus we arrive at the following definition. 
\medskip

\begin{definition}
\label{def:gen-Mixing-time}(Mixing time of general nonnegative matrices
$R$) Let $R$ be a $\frac{1}{2}$-lazy, irreducible nonnegative matrix
with PF eigenvalue 1 with $u$ and $v$ as the corresponding positive
left and right eigenvectors, where $u$ and $v$ are normalized so
that $\langle u,v\rangle=\|D_{u}v\|_{1}=1$. Then the mixing time
$\tau_{\epsilon}(R)$ is the smallest number $\tau$ such that $\left\Vert D_{u}\left(R^{\tau}-v\cdot u^{T}\right)x\right\Vert _{1}\leq\epsilon$
for every vector $x$ with $\|D_{u}x\|_{1}=1$. 
\end{definition}

We remark that similar to the doubly stochastic case, using the triangle
inequality, it is sufficient to find mixing time of standard basis
vectors $\mathbf{1}_{\{i\}}$. We state it as an elementary lemma.

\begin{lemma}
\label{lem:mixt-std-basis}Let $R$ be as stated in Definition \ref{def:gen-Mixing-time}.
Then the mixing time $\tau_{\epsilon}(R)$ is the smallest number
$\tau$ such that $\left\Vert D_{u}\left(R^{\tau}-v\cdot u^{T}\right)y\right\Vert _{1}\leq\epsilon$
for every vector $y=\frac{\mathbf{1}_{\{i\}}}{\|D_{u}\mathbf{1}_{\{i\}}\|_{1}}$
where $\mathbf{1}_{\{i\}}$ represents the standard basis vectors.
\end{lemma}

\begin{proof}
Let $y_{i}=\frac{\mathbf{1}_{\{i\}}}{\|D_{u}\mathbf{1}_{\{i\}}\|_{1}}$,
then $y_{i}$ is nonnegative and  $\|D_{u}y_{i}\|_{1}=\langle u,y_{i}\rangle=1$. 
Then for any $x$, such that $\|D_{u}x\|_{1}=1$, we can write 
\[
x=\sum\limits _{i}c_{i}\mathbf{1}_{\{i\}}=\sum_{i}c_{i}\|D_{u}\mathbf{1}_{\{i\}}\|_{1}y_{i}
\]
with 
\[
\|D_{u}x\|_{1}=\left\Vert D_{u}\sum\limits _{i}c_{i}\mathbf{1}_{\{i\}}\right\Vert _{1}=\sum\limits _{i}|c_{i}|\|D_{u}\mathbf{1}_{\{i\}}\|_{1}=1.
\]
Thus, if for every $i$, $\left\Vert D_{u}\left(R^{\tau}-v\cdot u^{T}\right)y_{i}\right\Vert _{1}\leq\epsilon$,
then 
\begin{align*}
\left\Vert D_{u}\left(R^{\tau}-v\cdot u^{T}\right)x\right\Vert _{1} & =\left\Vert D_{u}\left(R^{\tau}-v\cdot u^{T}\right)\sum\limits _{i}c_{i}\|D_{u}\mathbf{1}_{\{i\}}\|_{1}y_{i}\right\Vert _{1}\\
 & \leq\sum\limits _{i}|c_{i}|\|D_{u}\mathbf{1}_{\{i\}}\|_{1}\left\Vert D_{u}\left(R^{\tau}-v\cdot u^{T}\right)y_{i}\right\Vert _{1}\\
 & \leq\epsilon.
\end{align*}
Thus, it is sufficient to find mixing time for every nonnegative $x$
with $\|D_{u}x\|_{1}=\langle u,x\rangle=1$, and it will hold for
all $x$. \\
\end{proof}
The next lemma is similar to \ref{lem:transRtoA}, but connects the
mixing times of $R$ and $A=D_{u}^{\frac{1}{2}}D_{v}^{-\frac{1}{2}}RD_{u}^{-\frac{1}{2}}D_{v}^{\frac{1}{2}}$
showing they are the same.

\begin{lemma}
\label{lem:The-mixing-time-A-R}The mixing time of $R$ and $A=D_{u}^{\frac{1}{2}}D_{v}^{-\frac{1}{2}}RD_{u}^{-\frac{1}{2}}D_{v}^{\frac{1}{2}}$
are the same, where $R$ is an irreducible $\frac{1}{2}$-lazy nonnegative
matrix with largest eigenvalue 1 and corresponding positive eigenvectors
$u$ and $v$ normalized so that $\langle u,v\rangle=1$.
\end{lemma}

\begin{proof}
We first show that the mixing time of $R$ and $A=D_{u}^{\frac{1}{2}}D_{v}^{-\frac{1}{2}}RD_{u}^{-\frac{1}{2}}D_{v}^{\frac{1}{2}}$
are the same. Note that if $R$ is $\frac{1}{2}$-lazy, then $A$
is also $\frac{1}{2}$-lazy, since if $R=\frac{1}{2}I+\frac{1}{2}C$
where $C$ is nonnegative, then 
\[
A=\frac{1}{2}D_{u}^{\frac{1}{2}}D_{v}^{-\frac{1}{2}}ID_{u}^{-\frac{1}{2}}D_{v}^{\frac{1}{2}}+D_{u}^{\frac{1}{2}}D_{v}^{-\frac{1}{2}}CD_{u}^{-\frac{1}{2}}D_{v}^{\frac{1}{2}}=\frac{1}{2}I+D_{u}^{\frac{1}{2}}D_{v}^{-\frac{1}{2}}CD_{u}^{-\frac{1}{2}}D_{v}^{\frac{1}{2}}.
\]
Let $w=D_{u}^{\frac{1}{2}}D_{v}^{\frac{1}{2}}\mathbf{1}$, the left
and right eigenvector of $A$ for eigenvalue 1. We will show that
for every $x$ for which $\|D_{u}x\|_{1}=1$, there exists some $y$
with $\|D_{w}y\|_{1}=1$ such that 
\[
\|D_{u}(R^{\tau}-v\cdot u^{T})x\|_{1}=\|D_{w}(A^{\tau}-w\cdot w^{T})y\|_{1}
\]
which would imply that $\tau_{\epsilon}(R)=\tau_{\epsilon}(A)$ by
definition.

Let $x$ be a vector with $\|D_{u}x\|_{1}=1$ and let $y=D_{u}^{\frac{1}{2}}D_{v}^{-\frac{1}{2}}x$.
Then since $D_{w}=D_{u}^{\frac{1}{2}}D_{v}^{\frac{1}{2}},$ we get
$\|D_{w}y\|_{1}=\|D_{u}x\|_{1}=1$. Let $R=v\cdot u^{T}+B_{R}$ and
$A=w\cdot w^{T}+B_{A}$ where $B_{A}=D_{u}^{\frac{1}{2}}D_{v}^{-\frac{1}{2}}B_{R}D_{u}^{-\frac{1}{2}}D_{v}^{\frac{1}{2}}$
and $B_{A}^{\tau}=D_{u}^{\frac{1}{2}}D_{v}^{-\frac{1}{2}}B_{R}^{\tau}D_{u}^{-\frac{1}{2}}D_{v}^{\frac{1}{2}}$.
Then 
\begin{align*}
D_{u}(R^{\tau}-v\cdot u^{T})x & =D_{u}B_{R}^{\tau}x\\
 & =D_{u}(D_{u}^{-\frac{1}{2}}D_{v}^{\frac{1}{2}}B_{A}^{\tau}D_{u}^{\frac{1}{2}}D_{v}^{-\frac{1}{2}})x\\
 & =D_{w}B_{A}^{\tau}y\\
 & =D_{w}(A^{\tau}-w\cdot w^{T})y
\end{align*}
as required.
\end{proof}
We state a final simple lemma here, that shows the inequality between
the eigenvalues of $AA^{T}$ and $\tilde{A}$ whenever $A$ is lazy.
This was shown by Fill \cite{fill1991eigenvalue} (albeit for row-stochastic
matrices), and can be obtained in multiple ways.

\begin{lemma}
\label{lem:fill-lemma} (Fill \cite{fill1991eigenvalue}) Let $A$
be a $\frac{1}{2}$-lazy irreducible nonnegative matrix with largest
eigenvalue 1 and the corresponding left and right eigenvector $w$.
Then
\[
\lambda_{2}(AA^{T})\leq\lambda_{2}\left(\frac{A+A^{T}}{2}\right).
\]
\end{lemma}

\begin{proof}
Since $A$ is $\frac{1}{2}$-lazy, we have that $2A-I$ is nonnegative,
has PF eigenvalue 1, and has the same left and right eigenvector $w$
for eigenvalue 1 implying that the PF eigenvalue is 1 by Perron-Frobenius
(Theorem \ref{thm:(Perron-Frobenius)}, part 2), and also that its
largest singular value is 1 from Lemma \ref{lem:transRtoA}. Further,
$\frac{1}{2}(A+A^{T})$ also has the same properties. Thus, for any
$x$, 
\begin{align*}
AA^{T} & =\frac{A+A^{T}}{2}+\frac{(2A-I)(2A^{T}-I)}{4}-\frac{I}{4}\\
\Rightarrow\langle x,AA^{T}x\rangle & \leq\langle x,\frac{A+A^{T}}{2}x\rangle+\|x\|_{2}\frac{\|(2A-I)\|_{2}\|(2A^{T}-I)\|_{2}}{4}\|x\|_{2}-\frac{\|x\|_{2}^{2}}{4}\\
& \leq\langle x,\frac{A+A^{T}}{2}x\rangle\\
\Rightarrow\max_{x\perp w}\langle x,AA^{T}x\rangle & \leq\max_{x\perp w}\langle x,\frac{A+A^{T}}{2}x\rangle.
\end{align*}
where the last implication followed by the variational characterization
of eigenvalues since $AA^{T}$ and $\frac{1}{2}(A+A^{T})$ are symmetric,
and thus 
\[
\lambda_{2}(AA^{T})\leq\lambda_{2}\left(\frac{A+A^{T}}{2}\right)
\]
\end{proof}
We start by giving general bounds for nonnegative matrices, and will
infer bounds for the reversible case from those equations.

\section{Mixing time and singular values}

We first show a simple lemma relating the mixing time of nonnegative
matrices to the second singular value. This lemma is powerful enough
to recover the bounds obtained by Fill \cite{fill1991eigenvalue}
and Mihail \cite{mihail1989conductance} in an elementary way, and
also give all the bounds for the mixing time of reversible chains
that depend on the second eigenvalue or the edge expansion. Since the largest
singular value of any general nonnegative matrix $R$ with PF eigenvalue
1 could be much larger than 1, the relation between mixing time and
second singular value makes sense only for nonnegative matrices with
the same left and right eigenvector for eigenvalue 1, which have largest
singular value 1 by Lemma \ref{lem:transRtoA}. 

\begin{lemma}
\label{lem:mixt-singval} (Mixing time and second singular value)
Let $A$ be a nonnegative matrix (not necessarily lazy) with PF eigenvalue
1, such that $Aw=w$ and $A^{T}w=w$ for some $w$ with $\langle w,w\rangle=1$,
and let $\kappa=\min_{i}w_{i}^{2}$. Then for every $c>0$, 
\[
\tau_{\epsilon}(A)\ \leq\ \dfrac{\ln\left(\sqrt{\dfrac{n}{\kappa}}\cdot\dfrac{1}{\epsilon}\right)}{\ln\left(\frac{1}{\sigma_{2}(A)}\right)}\ \leq\ \dfrac{c\cdot\ln\left(\frac{n}{\kappa\cdot\epsilon}\right)}{1-\sigma_{2}^{c}(A)}.
\]
\end{lemma}

\begin{proof}
Writing $\tau$ as shorthand for $\tau_{\epsilon}(A)$ and since $A=w\cdot w^{T}+B$
with $Bw=0$ and $B^{T}w=0$, we have that $A^{\tau}=w\cdot w^{T}+B^{\tau}$.
Let $x$ be a nonnegative vector such that $\langle w,x\rangle=\|D_{w}x\|_{1}=1$.
As discussed after Definition \ref{def:gen-Mixing-time}, this is
sufficient for bounding mixing time for all $x$. Then we have 
\begin{align*}
\|D_{w}(A^{\tau}-w\cdot w^{T})x\|_{1} & =\|D_{w}B^{\tau}x\|_{1}=\|D_{w}B^{\tau}D_{w}^{-1}y\|_{1}\leq\|D_{w}\|_{1}\|B^{\tau}\|_{1}\|D_{w}^{-1}\|_{1}\|y\|_{1}
\end{align*}
where $y=D_{w}x$ and $\|y\|_{1}=1$. Further, since $\|w\|_{2}=1$
and $\kappa=\min_{i}w_{i}^{2}$, we have $\|D_{w}\|_{1}\leq1$ and
$\|D_{w}^{-1}\|_{1}\leq\frac{1}{\sqrt{\kappa}}$, and using these
bounds to continue the inequalities above, we get 
\[
\|D_{w}(A^{\tau}-w\cdot w^{T})x\|_{1}\leq\frac{1}{\sqrt{\kappa}}\|B^{\tau}\|_{1}\leq\frac{\sqrt{n}}{\sqrt{\kappa}}\|B^{\tau}\|_{2}\leq\frac{\sqrt{n}}{\sqrt{\kappa}}\|B\|_{2}^{\tau}\leq\frac{\sqrt{n}}{\sqrt{\kappa}}\left(\sigma_{2}(A)\right)^{\tau}\leq\epsilon
\]
where the second inequality is Cauchy-Schwarz, the fourth inequality
used $\|B\|_{2}\leq\sigma_{2}(A)$ since $A$ has identical left and
right PF eigenvector, and the last inequality was obtained by setting
$\tau=\tau_{\epsilon}(A)=\dfrac{\ln\left(\frac{\sqrt{n}}{\sqrt{\kappa}\cdot\epsilon}\right)}{\ln\left(\frac{1}{\sigma_{2}(A)}\right)}$,
and the subsequent inequality follows from $1-x\leq e^{-x}$.
\end{proof}
For the case of $c=2$, Lemma \ref{lem:mixt-singval} was obtained
by Fill \cite{fill1991eigenvalue} in a different manner, but we find
our more general proof much simpler.

\section{Mixing time and edge expansion}

We now relate the mixing time of general nonnegative matrices $R$
to its edge expansion $\phi(R)$. The upper bound for row stochastic
matrices $R$ in terms of $\phi(R)$ were obtained by Mihail \cite{mihail1989conductance}
and simplified by Fill \cite{fill1991eigenvalue} using Lemma \ref{lem:mixt-singval}
for $c=2$. Thus, the following lemma is not new, but we prove it
here since our proof again is elementary and holds for any nonnegative
matrix $R$. 

\begin{lemma}
\label{lem:mixt-phi-rel}(Mixing time and edge expansion, upper bound by Mihail \cite{mihail1989conductance}) Let $\tau_{\epsilon}(R)$
be the mixing time of a $\frac{1}{2}$-lazy nonnegative matrix $R$
with PF eigenvalue 1 and corresponding positive left and right eigenvectors
$u$ and $v$, and let $\kappa=\min_{i}u_{i}\cdot v_{i}$. Then 
\begin{align*}
\frac{\frac{1}{2}-\epsilon}{\phi(R)}\ \leq\ \tau_{\epsilon}(R)\  & \leq\ \dfrac{4\cdot\ln\left(\frac{n}{\kappa\cdot\epsilon}\right)}{\phi^{2}(R)}.
\end{align*}
\end{lemma}

\begin{proof}
From Lemma \ref{lem:The-mixing-time-A-R}, we have that the mixing
times of $R$ and $A=D_{u}^{\frac{1}{2}}D_{v}^{-\frac{1}{2}}RD_{u}^{-\frac{1}{2}}D_{v}^{\frac{1}{2}}$
are the same. Further, from Lemma \ref{lem:transRtoA}, we have that
$\phi(R)=\phi(A)$. Thus, we show the bound for $\tau_{\epsilon}(A)$
and $\phi(A)$, and the bound for $R$ will follow.

We first upper bound $\tau_{\epsilon}(A)$ in terms of the edge expansion $\phi$, giving an essentially one-line proof of Mihail's result.
From
Lemma \ref{lem:mixt-singval} for $c=2$, we have that 
\begin{equation}
\tau_{\epsilon}(A)\leq\frac{2\cdot\ln\left(\frac{n}{\kappa\cdot\epsilon}\right)}{1-\sigma_{2}^{2}(A)}=\frac{2\cdot\ln\left(\frac{n}{\kappa\cdot\epsilon}\right)}{1-\lambda_{2}(AA^{T})},\label{eq:proof-mixtphi-2}
\end{equation}
and replacing Lemma \ref{lem:fill-lemma} in equation \ref{eq:proof-mixtphi-2},
gives
\[
\tau_{\epsilon}(A)\leq\frac{2\cdot\ln\left(\frac{n}{\kappa\cdot\epsilon}\right)}{1-\sigma_{2}^{2}(A)}=\frac{2\cdot\ln\left(\frac{n}{\kappa\cdot\epsilon}\right)}{1-\lambda_{2}(AA^{T})}\leq\frac{2\cdot\ln\left(\frac{n}{\kappa\cdot\epsilon}\right)}{1-\lambda_{2}\left(\frac{A+A^{T}}{2}\right)}\leq\frac{4\cdot\ln\left(\frac{n}{\kappa\cdot\epsilon}\right)}{\phi^{2}\left(\frac{A+A^{T}}{2}\right)}=\frac{4\cdot\ln\left(\frac{n}{\kappa\cdot\epsilon}\right)}{\phi^{2}(A)}
\]
where the second last inequality follows from Cheeger's inequality
\ref{thm:cheeger-buser} for the symmetric matrix $(A+A^{T})/2$.

We now lower bound $\tau_{\epsilon}(A)$ in terms of $\phi(A)$. We
will show the bound for nonnegative vectors $x$, and by Definition
\ref{def:gen-Mixing-time} and the discussion after, it will hold
for all $x$. By definition of mixing time, we have that for any nonnegative
$x$ such that $\|D_{w}x\|_{1}=\langle w,x\rangle=1$, since $A^{\tau}=w\cdot w^{T}+B^{\tau}$,
\[
\|D_{w}(A^{\tau}-w\cdot w^{T})x\|_{1}=\|D_{w}B^{\tau}x\|_{1}\leq\epsilon
\]
and letting $y=D_{w}x$, we get that for any nonnegative $y$ with
$\|y\|_{1}=1$, we have 
\[
\|D_{w}B^{\tau}D_{w}^{-1}y\|_{1}\leq\epsilon.
\]
Plugging the standard basis vectors for $i$, we get that for every
$i$, 
\[
\sum_{j}\left|\frac{1}{w(i)}\cdot B^{\tau}(j,i)\cdot w(j)\right|=\frac{1}{w(i)}\cdot\sum_{j}\left|B^{\tau}(j,i)\right|\cdot w(j)\leq\epsilon.
\]
Thus, for any set $S$, 
\begin{equation}
\sum_{i\in S}w(i)^{2}\cdot\frac{1}{w(i)}\cdot\sum_{j}\left|B^{\tau}(j,i)\right|\cdot w(j)=\sum_{i\in S}\sum_{j}w(i)\cdot\left|B^{\tau}(j,i)\right|\cdot w(j)\leq\sum_{i\in S}w(i)^{2}\cdot\epsilon.\label{eq:proof-mixtphi-1}
\end{equation}
Moreover, for any set $S$ for which $\sum_{i\in S}w_{i}^{2}\leq\frac{1}{2}$,
\begin{align*}
\phi_{S}(A^{\tau}) & =\frac{\langle\mathbf{1}_{S},D_{w}A^{\tau}D_{w}\mathbf{1}_{\overline{S}}\rangle}{\langle\mathbf{1}_{S},D_{w}^{2}\mathbf{1}\rangle}=\frac{\langle\mathbf{1}_{\overline{S}},D_{w}A^{\tau}D_{w}\mathbf{1}_{S}\rangle}{\langle\mathbf{1}_{S},D_{w}^{2}\mathbf{1}\rangle}\\
 & \ \ \ \ \text{[since \ensuremath{D_{w}A^{\tau}D_{w}} is Eulerian, i.e. \ensuremath{D_{w}A^{\tau}D_{w}\mathbf{1}=\ensuremath{D_{w}(A^{\tau})^{T}D_{w}}\mathbf{1}}]}\\
 & =\frac{\sum_{i\in S}\sum_{j\in\overline{S}}A^{\tau}(j,i)\cdot w(i)\cdot w(j)}{\sum_{i\in S}w_{i}^{2}}\\
 & =\frac{\sum_{i\in S}\sum_{j\in\overline{S}}w(j)\cdot w(i)\cdot w(i)\cdot w(j)+\sum_{i\in S}\sum_{j\in\overline{S}}B^{\tau}(j,i)\cdot w(i)\cdot w(j)}{\sum_{i\in S}w_{i}^{2}}\\
 & \geq\sum_{j\in\overline{S}}w_{j}^{2}-\frac{\sum_{i\in S}\sum_{j\in\overline{S}}|B^{\tau}(j,i)|\cdot w(i)\cdot w(j)}{\sum_{i\in S}w_{i}^{2}}\\
 & \geq\sum_{j\in\overline{S}}w_{j}^{2}-\frac{\sum_{i\in S}\sum_{j}|B^{\tau}(j,i)|\cdot w(i)\cdot w(j)}{\sum_{i\in S}w_{i}^{2}}\\
 & \geq\frac{1}{2}-\epsilon\\
 & \ \ \ \ \text{[since \ensuremath{\sum_{i\in S}w_{i}^{2}\leq\frac{1}{2}} and \ensuremath{\sum_{i}w_{i}^{2}=1}, and the second term follows from equation \ref{eq:proof-mixtphi-1}]}
\end{align*}
and thus 
\[
\phi(A^{\tau})\geq\frac{1}{2}-\epsilon,
\]
and using Lemma \ref{lem:phi_pow_bound}, we obtain 
\[
\phi(A^{\tau})\leq\tau\cdot\phi(A),
\]
or 
\[
\frac{\frac{1}{2}-\epsilon}{\phi(A)}\leq\tau_{\epsilon}(A).
\]

\end{proof}

\section{Mixing time and spectral gap of Reversible Matrices}

We derive the bounds between the mixing time and spectral gap of reversible
nonnegative matrices. These bounds are well-known, we state them here
to contrast them with the bounds for general matrices. Further, they
are simple to derive given Lemma \ref{lem:mixt-singval}.

\begin{lemma}
\label{lem:mixt-specgap-reversible}(Mixing time and spectral gap
of reversible matrices) Let $\tau_{\epsilon}(R)$ be the mixing time
of a $\frac{1}{2}$-lazy nonnegative reversible matrix $R$ with PF
eigenvalue 1 and corresponding positive left and right eigenvectors
$u$ and $v$, normalized so that $\langle u,v\rangle=1$, and let
$\kappa=\min_{i}u_{i}v_{i}$. Then for $\epsilon<1$, 
\begin{align*}
\frac{1-\epsilon}{2\cdot\Delta(R)}\ \leq\ \dfrac{\ln\left(\dfrac{1}{\epsilon}\right)}{\ln\left(\dfrac{1}{\lambda_{2}(R)}\right)}\ \leq\ \tau_{\epsilon}(R)\  & \leq\ \dfrac{\ln\left(\dfrac{n}{\kappa\cdot\epsilon}\right)}{\Delta(R)}.
\end{align*}
\end{lemma}

\begin{proof}
Since $R$ is reversible, all eigenvalues are real, and since it is
lazy, $\lambda_{i}(R)\geq0$ for all $i$, and thus $\sigma_{2}(R)=\lambda_{2}(R)$.
The upper bound thus follows from Lemma \ref{lem:mixt-singval}. We
show the lower bound for $A=D_{u}^{\frac{1}{2}}D_{v}^{-\frac{1}{2}}RD_{u}^{-\frac{1}{2}}D_{v}^{\frac{1}{2}}$
and it follows for $R$ from \ref{lem:The-mixing-time-A-R}. Let $x$
be the eigenvector for eigenvalue $\lambda_{2}(A)$, then $\langle x,w\rangle=0$,
and let $x$ be normalized so that $\|D_{w}x\|_{1}=1$ as required
in Definition \ref{def:gen-Mixing-time}, where $w=D_{u}^{\frac{1}{2}}D_{v}^{\frac{1}{2}}\mathbf{1}$.
Thus, 
\begin{align*}
\|D_{w}(A^{\tau}-w\cdot w^{T})x\|_{1} & =\|D_{w}(A^{\tau}x-w\cdot w^{T}x)\|_{1}\\
 & =\|D_{w}(\lambda_{2}^{\tau}(A)x-w\cdot\mathbf{0})\|_{1}\\
 & =\lambda_{2}^{\tau}(A)\|D_{w}x\|_{1}\\
 & \ \text{[since \ensuremath{A} is \ensuremath{\frac{1}{2}}-lazy and thus \ensuremath{\lambda_{2}(A)\geq0}]}\\
 & =\lambda_{2}^{\tau}(A)\\
 & =(1-\Delta(A))^{\tau}
\end{align*}
which is greater than $\epsilon$ unless 
\[
\tau\geq\dfrac{\ln\left(\dfrac{1}{\epsilon}\right)}{\ln\left(\dfrac{1}{1-\Delta(A)}\right)}.
\]

If $\Delta(A)\leq\frac{1}{2}$, then 
\begin{align*}
1-\Delta(A) &\geq\exp(-2\cdot\Delta(A)),\\
2\cdot\Delta(A) &\geq\ln\left(\dfrac{1}{1-\Delta(A)}\right)
\end{align*}
and since $\ln(\frac{1}{\epsilon})\ge1-\epsilon$, this would give
\[
\tau\geq\dfrac{1-\epsilon}{2\cdot\Delta(A)}.
\]
Now assume $\Delta(A)>\frac{1}{2}$ or $\lambda_{2}(A)<\frac{1}{2}$,
then we know that the maximum mixing time will be achieved for $\lambda_{2}(A)=\frac{1}{2}$,
and will thus be at least 
\[
\dfrac{\ln\left(\dfrac{1}{\epsilon}\right)}{\ln2}\geq\dfrac{1-\epsilon}{\ln(2)}\geq\dfrac{1-\epsilon}{2}\geq1
\]
which is always true since the mixing time at least 1, and we again
get that 
\[
\tau\geq\dfrac{1-\epsilon}{2\cdot\Delta(A)}
\]
as required.
\end{proof}

\section{Mixing time and spectral gap of Nonreversible Matrices}

We obtain bounds for the mixing time of nonnegative matrices in terms
of the spectral gap, using methods similar to the ones used to obtain
the upper bound on $\phi$ in Theorem \ref{thm:intro-phi-delta}. 

\begin{lemma}
\label{lem:mixt-lambda2} (Mixing time and spectral gap) Let $\tau_{\epsilon}(R)$
be the mixing time of a $\frac{1}{2}$-lazy nonnegative matrix $R$
with PF eigenvalue 1 and corresponding positive left and right eigenvectors
$u$ and $v$, and let $\kappa=\min_{i}u_{i}\cdot v_{i}$. Then 
\begin{align*}
(1-\Delta)\dfrac{1-\epsilon}{\Delta}\ \leq\ \tau_{\epsilon}(R)\  & \leq\ 20\cdot\frac{n+\ln\left(\dfrac{1}{\kappa\cdot\epsilon}\right)}{\Delta(R)}.
\end{align*}
\end{lemma}

\begin{proof}
Since the mixing time (as shown in the proof of Lemma \ref{lem:The-mixing-time-A-R}),
eigenvalues, edge expansion, and value of $\kappa$ for $R$ and $A=D_{u}^{\frac{1}{2}}D_{v}^{-\frac{1}{2}}RD_{u}^{-\frac{1}{2}}D_{v}^{\frac{1}{2}}$
are the same, we provide the bounds for $A$ and the bounds for $R$
follow.

To obtain the lower bound, let $y$ be eigenvector for the second
eigenvalue $\lambda_{2}$ of $A$ normalized so that $\|D_{w}y\|=1$.
Note that since $w$ is both the left and right eigenvector for eigenvalue
1 of $A$, we have that $\langle w,y\rangle=0$. To consider the mixing
time of $y$, note that 
\begin{align*}
\|D_{w}(A^{\tau}-w\cdot w^{T})y\|_{1} & =\|D_{w}(A^{\tau}-w\cdot w^{T})y\|_{1}=\|D_{w}\lambda_{2}^{\tau}y\|_{1}=|\lambda_{2}|^{\tau}
\end{align*}
which is less than $\epsilon$ for 
\[
\tau_{\epsilon}(A)\geq\dfrac{\ln(\frac{1}{\epsilon})}{\ln(\frac{1}{|\lambda_{2}|})}\geq\dfrac{1-\epsilon}{1-|\lambda_{2}|}|\lambda_{2}|\geq(1-\Delta)\dfrac{1-\epsilon}{\Delta}
\]
where the last inequality follows since $1-x\leq\ln(1/x)\leq1/x-1$
and $\text{Re}\lambda_{2}\leq|\lambda_{2}|.$

For the upper bound on $\tau_{\epsilon}(A)$, we also restrict to
nonnegative vectors, the bound for general vectors follows by the
triangle inequality as discussed in Lemma \ref{lem:mixt-std-basis}.
Similar to the proof of Lemma \ref{lem:mixt-singval}, we have for
any nonnegative vector $x$ with $\langle w,x\rangle=1$, for $A=w\cdot w^{T}+B$,
\[
\|D_{w}(A^{\tau}-w\cdot w^{T})x\|_{1}\leq\sqrt{\frac{n}{\kappa}}\cdot\|B^{\tau}\|_{2}
\]
and having 
\[
\|B^{\tau}\|_{2}\leq\frac{\epsilon\sqrt{\kappa}}{\sqrt{n}}
\]
is sufficient. Let $T$ be the triangular matrix in the Schur form
of $B$, from the proof of Lemma \ref{lem:T_norm_1_bound}, we have
that for 
\[
k\geq\dfrac{3.51n+1.385\ln\left(\frac{n}{\delta}\right)}{1-|\lambda_{m}(A)|},
\]
the norm 
\[
\|B^{k}\|_{2}\leq\delta,
\]
and thus setting $\delta=\frac{\epsilon\sqrt{\kappa}}{\sqrt{n}}$,
we get that 
\[
\tau_{\epsilon}(A)\leq\dfrac{3.51n+1.385\ln\left(n\cdot\frac{\sqrt{n}}{\sqrt{\kappa}\cdot\epsilon}\right)}{1-|\lambda_{m}(A)|}.
\]
Further, since $A$ is $\frac{1}{2}$-lazy, $2A-I$ is also nonnegative,
with the same positive left and right eigenvector $w$ for PF eigenvalue
1, thus having largest singular value 1, and thus every eigenvalue
of $2A-I$ has magnitude at most 1. Similar to the proof of Lemma
\ref{lem:real-lambda-bound}, if $\lambda_{r}=a+i\cdot b$ is any
eigenvalue of $A$, the corresponding eigenvalue in $2A-I$ is $2a-1+i\cdot2b$,
whose magnitude is at most 1, giving $(2a-1)^{2}+4b^{2}\leq1$ or
$a^{2}+b^{2}\leq a$. It further gives that 
\[
1-|\lambda_{r}|=1-\sqrt{a^{2}+b^{2}}\geq1-\sqrt{1-(1-a)}\geq1-e^{\frac{1}{2}(1-a)}\geq\frac{1}{4}(1-a)
\]
or 
\[
1-|\lambda_{m}|\geq\frac{1}{4}(1-\text{Re}\lambda_{m})\geq\frac{1}{4}(1-\text{Re}\lambda_{2}),
\]
which gives 
\begin{align*}
\tau_{\epsilon}(A) & \leq4\cdot\dfrac{3.51n+1.385\ln\left(n\cdot\frac{\sqrt{n}}{\sqrt{\kappa}\cdot\epsilon}\right)}{1-\text{Re}\lambda_{2}(A)}\\
 & \leq20\cdot\frac{n+\ln\left(\dfrac{1}{\kappa\cdot\epsilon}\right)}{1-\text{Re}\lambda_{2}(A)}
\end{align*}
completing the proof.
\end{proof}
We remark that there is only \emph{additive} and not multiplicative
dependence on $\ln\left(\frac{n}{\kappa\cdot\epsilon}\right)$. Further,
our constructions \ref{constr:Rootn-matrices-=002013}
and \ref{constr-chet} also show that the upper bound
on $\tau$ using $\text{Re}\lambda_{2}$ in Lemma \ref{lem:mixt-lambda2}
is also (almost) tight. For the construction of $A_{n}$ in Theorem
\ref{thm:rootn}, letting the columns of $U_{n}$ be $u_{1},\ldots,u_{n}$,
for $x=u_{2}$, $(A_{n}^{k}-J)u_{2}=(1-(2+\sqrt{n})^{-1})^{k}u_{3}$,
and so for $k=O(\sqrt{n})$, the triangular block of $A^{O(\sqrt{n})}$
has norm about $1/e$, which further becomes less than $\epsilon$
after about $\ln\left(\frac{n}{\epsilon}\right)$ powers. Thus for
the matrices $A_{n}$, $\tau_{\epsilon}(A_{n})\in O\left(\sqrt{n}\cdot\ln\left(\frac{n}{\epsilon}\right)\right)$.
This shows Lemma \ref{lem:mixt-lambda2} is also (almost) tight since
$\lambda_{2}(A_{n})=0$.

\section{Mixing time of a nonnegative matrix and its additive symmetrization}

We can also bound the mixing time of a nonnegative matrix $A$ with
the same left and right eigenvector $w$ for PF eigenvalue 1, with
the mixing time of its \emph{additive symmetrization} $M=\frac{1}{2}(A+A^{T})$.
Since $\phi(A)=\phi(M)$, we can bound $\tau_{\epsilon}(A)$ and $\tau_{\epsilon}(M)$
using the two sided bounds between edge expansion and mixing time
in Lemma \ref{lem:mixt-phi-rel}. For the lower bound, we get $\gamma_{1}\cdot\sqrt{\tau_{\epsilon}(M)}\leq\tau_{\epsilon}(A)$,
and for the upper bound, we get 
\[
\tau_{\epsilon}(A)\leq\gamma_{2}\cdot\tau_{\epsilon}^{2}(M),
\]
where $\gamma_{1}$ and $\gamma_{2}$ are some functions polylogarithmic
in $n,\kappa,\frac{1}{\epsilon}$. However, by bounding the appropriate
operator, we can show a tighter upper bound on $\tau_{\epsilon}(A)$,
with only a \emph{linear} instead of quadratic dependence on $\tau_{\epsilon}(M)$. 

\begin{lemma}
\label{lem:tau-A-AAT}Let $A$ be a $\frac{1}{2}$-lazy nonnegative
matrix with positive left and right eigenvector $w$ for PF eigenvalue
1, let $M=\frac{1}{2}(A+A^{T})$, and $\kappa=\min_{i}w_{i}^{2}$.
Then 
\[
\frac{1-2\epsilon}{4\cdot\ln^{\frac{1}{2}}\left(\frac{n}{\kappa\cdot\epsilon}\right)}\cdot\tau_{\epsilon}^{\frac{1}{2}}\left(M\right)\ \leq\ \tau_{\epsilon}(A)\ \leq\ \frac{2\cdot\ln\left(\frac{n}{\kappa\cdot\epsilon}\right)}{\ln\left(\frac{1}{\epsilon}\right)}\cdot\tau_{\epsilon}(M).
\]
\end{lemma}

\begin{proof}
Since $\phi(A)=\phi(M)$ and each have the same left and right eigenvector
$w$ for PF eigenvalue 1, the lower bound on $\tau_{\epsilon}(A)$
in terms of $\tau_{\epsilon}(M)$ follows immediately from Lemma \ref{lem:mixt-phi-rel},
since 
\[
\sqrt{\tau_{\epsilon}(M)}\leq\ln^{\frac{1}{2}}\left(\frac{n}{\kappa\cdot\epsilon}\right)\cdot\frac{1}{\phi(M)}=\ln^{\frac{1}{2}}\left(\frac{n}{\kappa\cdot\epsilon}\right)\cdot\frac{1}{\phi(A)}\leq\ln^{\frac{1}{2}}\left(\frac{n}{\kappa\cdot\epsilon}\right)\cdot\frac{2\tau_{\epsilon}(A)}{1-2\epsilon}.
\]

For the upper bound, we first define a new quantity for \emph{positive
semidefinite }nonnegative matrices $M$ with PF eigenvalue 1 and $w$
as the corresponding left and right eigenvector. Let $T_{\epsilon}(M)$
be defined as the smallest number $k$ for which 
\[
\|M^{k}-w\cdot w^{T}\|_{2}=\epsilon.
\]

Since $M$ is symmetric, we can write $M=w\cdot w^{T}+UDU^{*}$ where
the first column of the unitary $U$ is $w$, and the diagonal matrix
$D$ contains all eigenvalues of $M$ except the eigenvalue 1 which
is replaced by 0. Further, since $M$ is positive semidefinite, $\lambda_{2}(M)$
is the second largest eigenvalue of $M$, then for every $i>2$, $0\leq\lambda_{i}(M)\leq\lambda_{2}(M)$.
Thus we have 
\[
\|M^{k}-w\cdot w^{T}\|_{2}=\|UD^{k}U^{*}\|_{2}=\|D^{k}\|_{2}=\lambda_{2}^{k}(M)
\]
and 
\[
T_{\epsilon}(M)=\frac{\ln\left(\frac{1}{\epsilon}\right)}{\ln\left(\frac{1}{\lambda_{2}(M)}\right)}.
\]
Further, for the eigenvector $y$ of $M$ corresponding to $\lambda_{2}$,
we have for $k=T_{\epsilon}(M)$, 
\[
(M^{k}-w\cdot w^{T})y=\lambda_{2}^{k}\cdot y=\epsilon\cdot y.
\]
Setting $x=\frac{y}{\|D_{w}y\|_{1}}$, we have $\|D_{w}x\|_{1}=1$,
and we get 
\[
\|D_{w}(M^{k}-w\cdot w^{T})x\|_{1}=\left\Vert D_{w}\frac{\epsilon\cdot y}{\|D_{w}y\|_{1}}\right\Vert _{1}=\epsilon,
\]
which implies that 
\begin{equation}
\tau_{\epsilon}(M)\geq k=T_{\epsilon}(M),\label{eq:proof-AAT-1}
\end{equation}
since there is some vector $x$ with $\|D_{w}x\|=1$ that has $\|D_{w}(M^{k}-w\cdot w^{T})x\|_{1}=\epsilon$
and for every $t<k$, it is also the case that $\|D_{w}(M^{t}-w\cdot w^{T})x\|_{1}>\epsilon$.

Now we observe that since $A$ is $\frac{1}{2}$-lazy with PF eigenvalue
1, and the same left and right eigenvector $w$ for eigenvalue 1,
we have that $M=\frac{1}{2}(A+A^{T})$ is positive semidefinite. From
Lemma \ref{lem:mixt-singval}, we have 
\[
\tau_{\epsilon}(A)\leq\frac{2\cdot\ln\left(\frac{\sqrt{n}}{\sqrt{\kappa}\cdot\epsilon}\right)}{\ln\left(\frac{1}{\sigma_{2}^{2}(A)}\right)}
\]
and since $A$ is $\frac{1}{2}$-lazy, as shown in the proof of Lemma
\ref{lem:mixt-phi-rel}, 
\[
\sigma_{2}^{2}(A)=\lambda_{2}(AA^{T})\leq\lambda_{2}(M),
\]
giving 
\[
\tau_{\epsilon}(A)\leq\frac{2\cdot\ln\left(\frac{\sqrt{n}}{\sqrt{\kappa}\cdot\epsilon}\right)}{\ln\left(\frac{1}{\lambda_{2}(M)}\right)}=\frac{2\cdot\ln\left(\frac{\sqrt{n}}{\sqrt{\kappa}\cdot\epsilon}\right)}{\ln\left(\frac{1}{\epsilon}\right)}\cdot\frac{\ln\left(\frac{1}{\epsilon}\right)}{\ln\left(\frac{1}{\lambda_{2}(M)}\right)}=\frac{2\cdot\ln\left(\frac{\sqrt{n}}{\sqrt{\kappa}\cdot\epsilon}\right)}{\ln\left(\frac{1}{\epsilon}\right)}\cdot T_{\epsilon}(M)\leq\frac{2\cdot\ln\left(\frac{\sqrt{n}}{\sqrt{\kappa}\cdot\epsilon}\right)}{\ln\left(\frac{1}{\epsilon}\right)}\cdot\tau_{\epsilon}(M)
\]
where the last inequality followed from equation \ref{eq:proof-AAT-1}.
\end{proof}
One example application of Lemma \ref{lem:tau-A-AAT} is the following:
given any undirected graph $G$ such that each vertex has degree $d$,
\emph{any} manner of orienting the edges of $G$ to obtain a graph
in which every vertex has in-degree and out-degree $d/2$ cannot increase
the mixing time of a random walk (up to a factor of $\ln\left(\frac{n}{\kappa\cdot\epsilon}\right)$).

\section{Mixing time of the continuous operator}

Let $R$ be a nonnegative matrix with PF eigenvalue 1 and associated
positive left and right eigenvectors $u$ and $v$. The continuous
time operator associated with $R$ is defined as $\exp\left(t\cdot\left(R-I\right)\right)$,
where for any matrix $M$, we formally define $\exp(M)=\sum_{i=0}\frac{1}{i!}M^{i}$.
The reason this operator is considered continuous, is that starting
with any vector $x_{0}$, the vector $x_{t}$ at time $t\in\mathbb{R}_{\geq0}$
is defined as $x_{t}=\exp\left(t\cdot\left(R-I\right)\right)x_{0}$.
Since 
\[
\exp\left(t\cdot\left(R-I\right)\right)=\exp(t\cdot R)\cdot\exp(-t\cdot I)=e^{-t}\sum_{i=0}^{\infty}\frac{1}{i!}t^{i}R^{i}
\]
where we split the operator into two terms since $R$ and $I$ commute,
it follows that $\exp\left(t\cdot\left(R-I\right)\right)$ is nonnegative,
and if $\lambda$ is any eigenvalue of $R$ for eigenvector $y$,
then $e^{t(\lambda-1)}$ is an eigenvalue of $\exp\left(t\cdot\left(R-I\right)\right)$
for the same eigenvector $y$. Thus, it further follows that $u$
and $v$ are the left and right eigenvectors for $\exp\left(t\cdot\left(R-I\right)\right)$
with PF eigenvalue 1. The mixing time of $\exp\left(t\cdot\left(R-I\right)\right)$,
is the value of $t$ for which 
\[
\left\Vert D_{u}\left(\exp\left(t\cdot\left(R-I\right)\right)-v\cdot u^{T}\right)v_{0}\right\Vert _{1}\leq\epsilon
\]
for every $v_{0}$ such that $\|D_{u}v_{0}\|_{1}=1$, and thus, it
is exactly same as considering the mixing time of $\exp(R-I)$ in
the sense of Definition \ref{def:gen-Mixing-time}. 

\begin{lemma}
\label{lem:mixt-contchain}Let $R$ be a nonnegative matrix (not necessarily
lazy) with positive left and right eigenvectors $u$ and $v$ for
PF eigenvalue 1, normalized so that $\langle u,v\rangle=1$ and let
$\kappa=\min_{i}u_{i}\cdot v_{i}$. Then the mixing time of $\exp(t\cdot(R-I))$,
or $\tau_{\epsilon}\left(\exp(R-I)\right)$ is bounded as 
\[
\frac{\frac{1}{2}-\epsilon}{\phi(R)}\ \leq\ \tau_{\epsilon}\left(\exp(R-I)\right)\ \leq\ \frac{100\cdot\ln\left(\frac{n}{\kappa\cdot\epsilon}\right)}{\phi^{2}(R)}.
\]
\end{lemma}

\begin{proof}
We will first find the mixing time of the operator $\exp\left(\frac{R-I}{2}\right)$,
and the mixing time of the operator $\exp(R-I)$ will simply be twice
this number. By the series expansion of the $\exp$ function, it follows that
$\exp\left(\frac{R-I}{2}\right)$ has PF eigenvalue 1, and $u$ and
$v$ as the corresponding left and right eigenvectors. Further,

\[
\exp\left(\frac{R-I}{2}\right)=e^{-\frac{1}{2}}\left(I+\sum_{i\geq1}\frac{1}{i!}\frac{R^{i}}{2^{i}}\right)
\]
which is $\frac{1}{2}$-lazy due to the first term since $e^{-\frac{1}{2}}\geq\frac{1}{2}$
and all the other terms are nonnegative. Further, for any set $S$
for which $\sum_{i\in S}u_{i}v_{i}\leq\frac{1}{2}$, let $\delta=\frac{1}{2}$,
then 
\begin{align*}
\phi_{S}\left(\exp\left(\frac{R-I}{2}\right)\right) & =e^{-\delta}\cdot\frac{\langle\mathbf{1}_{S},D_{u}\exp\left(\frac{R}{2}\right)D_{v}\mathbf{1}_{\overline{S}}\rangle}{\langle\mathbf{1}_{S},D_{u}D_{v}\mathbf{1}\rangle}\\
 & =e^{-\delta}\cdot\sum_{i\geq1}\frac{\delta^{i}}{i!}\cdot\frac{\langle\mathbf{1}_{S},D_{u}R^{i}D_{v}\mathbf{1}_{\overline{S}}\rangle}{\langle\mathbf{1}_{S},D_{u}D_{v}\mathbf{1}\rangle}\\
 & \ \ \ \ \text{[since \ensuremath{\langle\mathbf{1}_{S},D_{u}ID_{v}\mathbf{1}_{\overline{S}}\rangle=0}]}\\
 & =e^{-\delta}\cdot\sum_{i\geq1}\frac{\delta^{i}}{i!}\cdot\phi_{S}(R^{i})\\
 & \ \ \ \ \text{[since \ensuremath{R^{i}} also has \ensuremath{u} and \ensuremath{v} as the left and right eigenvectors for eigenvalue 1]}\\
 & \leq e^{-\delta}\cdot\sum_{i\geq1}\frac{\delta^{i}}{i!}\cdot i\cdot\phi_{S}(R)\\
 & \ \ \ \ \text{[using Lemma \ref{lem:phi_pow_bound}]}\\
 & =e^{-\delta}\cdot\delta\cdot\phi_{S}(R)\sum_{i\geq1}\frac{\delta^{i-1}}{(i-1)!}\\
 & =e^{-\delta}\cdot\delta\cdot\phi_{S}(R)\cdot e^{\delta}\\
 & =\delta\cdot\phi_{S}(R)
\end{align*}
and thus, 
\begin{equation}
\phi\left(\exp\left(\frac{R-I}{2}\right)\right)\leq\frac{1}{2}\phi(R).\label{eq:mixt-contchain-1}
\end{equation}
Moreover, considering the first term in the series expansion, we get 
\begin{align*}
\phi_{S}\left(\exp\left(\frac{R-I}{2}\right)\right) & =e^{-\delta}\sum_{i\geq1}\frac{\delta^{i}}{i!}\cdot\phi_{S}(R^{i})\geq e^{-\delta}\cdot\delta\cdot\phi_{S}(R)
\end{align*}
or 
\begin{equation}
\phi\left(\exp\left(\frac{R-I}{2}\right)\right)\geq\frac{3}{10}\cdot\phi(R).\label{eq:mixt-contchain-2}
\end{equation}
Since $\exp\left(\frac{R-I}{2}\right)$ has left and right eigenvectors
$u$ and $v$, and is $\frac{1}{2}$-lazy, we get from Lemma \ref{lem:mixt-phi-rel}
and equations \ref{eq:mixt-contchain-1} and \ref{eq:mixt-contchain-2}
that 
\[
\frac{1}{2}\cdot\frac{\frac{1}{2}-\epsilon}{\frac{1}{2}\phi(R)}\leq\tau_{\epsilon}\left(\exp(R-I)\right)\leq2\cdot\frac{4\cdot\ln\left(\frac{n}{\kappa\cdot\epsilon}\right)}{\left(\frac{3}{10}\right)^{2}\phi^{2}(R)}\leq\frac{100\cdot\ln\left(\frac{n}{\kappa\cdot\epsilon}\right)}{\phi^{2}(R)}
\]
giving the result. 
\end{proof}

\section{Comparisons with the canonical paths method}

For the case of symmetric nonnegative matrices $M$ with PF eigenvalue
1, as shown in Lemma \ref{lem:mixt-specgap-reversible}, since $\tau$
varies inversely with $1-\lambda_{2}$ (up to a loss of a factor of
$\ln(\frac{n}{\kappa\cdot\epsilon})$), it follows that any lower
bound on the spectral gap can be used to upper bound $\tau_{\epsilon}(M)$.
Further, since $1-\lambda_{2}$ can be written as a minimization problem
for symmetric matrices (see Section \ref{sec:Preliminaries}), any
relaxation of the optimization problem can be used to obtain a lower
bound on $1-\lambda_{2}$, and inequalities obtained thus are referred
to as \emph{Poincare inequalities}. One such method is to use \emph{canonical
paths }\cite{sinclair1992improved} in the underlying weighted graph,
which helps to bound mixing time in certain cases in which computing
$\lambda_{2}$ or $\phi$ is infeasible. However, since it is possible
to define canonical paths in many different ways, it leads to multiple
relaxations to bound $1-\lambda_{2}$, each useful in a different
context. We remark one particular definition and lemma here, since
it is relevant to our construction in Theorem \ref{thm:rootn}, after
suitably modifying it for the doubly stochastic case. 

\begin{lemma}
\cite{sinclair1992improved}\label{lem:canpath-sinc} Let $M$ represent
a symmetric doubly stochastic matrix. Let $W$ be a set of paths in
$M$, one between every pair of vertices. For any path $\gamma_{u,v}\in S$
between vertices $(u,v)$ where $\gamma_{u,v}$ is simply a set of
edges between $u$ and $v$, let the number of edges or the (unweighted)
length of the path be $|\gamma_{u,v}|$. Let 
\[
\rho_{W}(M)=\max_{e=(x,y)}\dfrac{\sum\limits _{(u,v):e\in\gamma_{u,v}}|\gamma_{u,v}|}{n\cdot M_{x,y}}.
\]
Then for any $W$, 
\[
1-\lambda_{2}(M)\geq\frac{1}{\rho_{W}(M)}
\]
and thus, 
\[
\tau_{\epsilon}(M)\leq\rho_{W}(M)\cdot\ln\left(\frac{n}{\epsilon}\right).
\]
\end{lemma}

\begin{corollary}
\label{cor:canpath-sinc}Combining Lemma \ref{lem:phi_sing} and Lemma
\ref{lem:canpath-sinc}, it follows that for any doubly stochastic
matrix $A$, and any set $W$ of paths in the underlying graph of
$AA^{T}$, 
\[
\tau_{\epsilon}(A)\leq\dfrac{2\cdot\ln\left(\frac{n}{\epsilon}\right)}{1-\sigma_{2}^{2}(A)}=\dfrac{2\cdot\ln\left(\frac{n}{\epsilon}\right)}{1-\lambda_{2}(AA^{T})}\leq2\cdot\rho_{W}(AA^{T})\cdot\ln\left(\frac{n}{\epsilon}\right).
\]
\end{corollary}

Consider the example $A_{n}$ in Theorem \ref{thm:rootn}. It is not
difficult to see that 
\begin{equation}
\tau_{\epsilon}(A_{n})\in O\left(\sqrt{n}\cdot\ln\left(\frac{n}{\epsilon}\right)\right).\label{eq:para_mixt}
\end{equation}
This follows since the factor of $\sqrt{n}$ ensures that the only
non zero entries in the triangular matrix $T_{n}$ in the Schur form
of $A^{\lceil\sqrt{n}\rceil}$ are about $e^{-1}$, and the factor
of $\ln\left(\frac{n}{\epsilon}\right)$ further converts these entries
to have magnitude at most $\frac{\epsilon}{n}$ in $A^{\tau}$. Thus,
the operator norm becomes about $\frac{\epsilon}{n}$, and the $\ell_{1}$
norm gets upper bounded by $\epsilon$. However, from Theorem \ref{thm:rootn},
since $\phi(A_{n})\geq\frac{1}{6\sqrt{n}},$ it follows from Lemma
\ref{lem:mixt-phi-rel} that $\tau_{\epsilon}(A_{n})\in O\left(n\cdot\ln\left(\frac{n}{\epsilon}\right)\right)$,
about a quadratic factor off from the actual upper bound in equation
\ref{eq:para_mixt}. Further, from Theorem \ref{thm:rootn}, the second
eigenvalue of $A_{n}$ is $0$, and even employing Lemma \ref{lem:mixt-lambda2}
leads to a quadratic factor loss from the actual bound. However, Lemma
\ref{lem:mixt-singval} and Corollary \ref{cor:canpath-sinc} do give
correct bounds. Since $\sigma_{2}(A_{n})=1-\frac{1}{\sqrt{n}+2}$
from Theorem \ref{thm:rootn}, it follows from Lemma \ref{lem:mixt-singval}
for $c=1$ that $\tau_{\epsilon}(A_{n})\in O\left(\sqrt{n}\cdot\ln\left(\frac{n}{\epsilon}\right)\right)$,
matching the bound in equation \ref{eq:para_mixt}. Now to see the
bound given by canonical paths and corollary \ref{cor:canpath-sinc},
consider the matrix $M=A_{n}A_{n}^{T}$. Every entry of $M$ turns
out to be positive, and the set $W$ is thus chosen so that the path
between any pair of vertices is simply the edge between the vertices.
Further for $r_{n},\alpha_{n},\beta_{n}$ defined in the proof of
Theorem \ref{thm:rootn}, $M=J+r_{n}^{2}B$, where 
\[
B_{1,1}=\frac{n-2}{n},\ \ B_{n,n}=(n-2)\cdot\beta_{n}^{2},\ \ B_{i,i}=\alpha_{n}^{2}+(n-3)\cdot\beta_{n}^{2},\ \ B_{1,n}=B_{n,1}=\frac{n-2}{\sqrt{n}}\cdot\beta_{n},
\]
\[
B_{n,j}=B_{j,n}=\alpha_{n}\cdot\beta_{n}+(n-3)\cdot\beta_{n}^{2},\ \ B_{1,j}=B_{j,1}=\frac{1}{\sqrt{n}}\cdot(\alpha_{n}+(n-3)\cdot\beta_{n}),\ \ B_{i,j}=2\cdot\alpha_{n}\cdot\beta_{n}+(n-4)\cdot\beta_{n}^{2},
\]
and $2\leq i,j\leq n-1$. It follows that any entry of the matrix
$M$ is at least $c\cdot n^{-\frac{3}{2}}$ (for some constant $c$),
and from Corollary \ref{cor:canpath-sinc}, we get that $\tau_{\epsilon}(A_{n})\in O\left(\sqrt{n}\cdot\ln\left(\frac{n}{\epsilon}\right)\right)$,
matching the bound in equation \ref{eq:para_mixt}.

\section{Schur Complements}

This is the first time we talk about Laplacians in this thesis. The
reason for avoiding it so far was that it did not seem necessary to
obtain any of our required bounds for nonnegative matrices, although
it is implicitly present in all of them. Our assumptions will essentially
be the same that have been so far.

Let $R$ be an irreducible nonnegative matrix. We require the notion
of irreducibility since we will deal with inverses of submatrices
of R, which will not exist without this property. Further, $R$ has
positive $u$ and $v$ as left and right eigenvectors for eigenvalue
1, normalized so that $\sum_{i}u_{i}v_{i}=1$.
\medskip

\begin{definition}
(Laplacian of a nonnegative matrix) Let $R$ be an irreducible nonnegative
matrix with PF eigenvalue 1. Define the Laplacian of $R$ as $L=I-R$.
\end{definition}
\medskip

\begin{remark}
There are many different definitions for laplacians of nonnegative
matrices $R$ in literature. For simplicity, assume we want to understand
the Laplacian of irreducible nonnegative matrices $A$ with largest
eigenvalue 1 and left and right positive eigenvector $w$, then two
of the most prominent definitions of the laplacian of $A$ are $I-\frac{1}{2}(A+A^{T})$
and $I-A$. The first definition has been widely used to derive results
about irreducible matrices mostly since symmetric matrices are easier
to deal with and the results are cleaner. However, the first definition
loses all information about the eigenvalues of $A$, and the primary
definition for us is the second one. This definition is much more
difficult to deal with due to a lack of the variational characterization
of eigenvalues as in the symmetric case and requires entirely new
ideas and tools. This definition has also been used in many other
contexts such as statistical physics (see \cite{statmech-tolman1979principles}).
\end{remark}

\noindent \textbf{Notation.} Let $L_{S,T}$ denote the submatrix of
$L$ indexed by rows in $S$ and columns in $T$, where we will let
$L_{S}=L_{S,S}$ for brevity. Similarly, for some vector $x\in\mathbb{R}^{n}$,
we will denote that $x_{S}$ as the vector in $R^{|S|}$ containing
the entries corresponding to indices in $S$. We remark that we will
always write $x_{\{i\}}$ for singletons.

The first simple lemma is the following.

\begin{lemma}
\label{lem:L-invertible}Given a laplacian $L=I-R$ of irreducible
nonnegative $R$ with positive (left and right) eigenvectors $u$
and $v$ for largest eigenvalue 1, $L_{U,U}$ is invertible except
when $U=[n]$.
\end{lemma}

\begin{proof}
Assume not, then there is an eigenvalue $0$ for $L_{U,U}$ corresponding
to eigenvector $x_{U}\not=0$, and extend it to vector $x$ which
has 0 values outside $U$, implying $Lx=0$ or $Rx=x$ which is a
contradiction to the Perron-Frobenius Theorem \ref{thm:(Perron-Frobenius)}
since $v$ is the unique \emph{positive} vector for eigenvalue 1 of
$R$, and $x$ has zeros outside $U$ is thus not positive.
\end{proof}

\begin{definition}
\label{def:(Inverse-of-Laplacians)}(Inverse of Laplacians) A basic
problem of linear algebra, is to solve for $x$ in the equation $Lx=b$
given $b$ (where $L=I-R$ and $R$ is irreducible nonnegative with
largest eigenvalue 1). Since the kernel of $L$ has dimension 1, it
implies that for $b\not=0$, there is a unique solution for $x$ with
$\langle x,v\rangle=0$ where $x=L^{+}b$.
\end{definition}

To understand the solution $x=L^{+}b$ more clearly, for some fixed
set $U\subset[n]$, let matrix $R=\begin{bmatrix}A & B\\
C & D
\end{bmatrix}$ giving $L=\begin{bmatrix}I-A & -B\\
-C & I-D
\end{bmatrix}=\begin{bmatrix}A' & B'\\
C' & D'
\end{bmatrix}$ and assume that the vertices that $A$ is on is the set $U$, and
write $x_{U}=y$, $x_{\overline{U}}=z$ for brevity and similarly
for $b$, then from $Lx=b$, we get 
\begin{equation}
A'y+B'z=b_{1}\label{eq:schurc-1}
\end{equation}
and 
\begin{equation}
C'y+D'z=b_{2}\label{eq:schurc-2}
\end{equation}
and if $D'$ is invertible which is the case if $L$ is a laplacian
from Lemma \ref{lem:L-invertible}, then 
\begin{equation}
z=b_{2}-D'^{-1}C'y\label{eq:schurc-3}
\end{equation}
and 
\begin{equation}
(A'-B'D'^{-1}C')y=b_{1}-B'b_{2}.\label{eq:schurc-4}
\end{equation}
 This brings us to the definition of Schur Complements. \\

\begin{definition}
\label{def:(Schur-Complement)}(Schur Complement) Let the Laplacian $L=I-R$ where 
\[
R=\begin{bmatrix}A & B\\
C & D
\end{bmatrix}
\]
is an irreducible nonnegative with largest eigenvalue 1.
The Schur complement of $L$ with respect
to some set $U\subset[n]$, written as $L|_{U}$ is defined as
\[
L|_{U}=A'-B'D'^{-1}C'=(I-A)-B(I-D)^{-1}C=L_{U}-L_{U,\overline{U}}L_{\overline{U}}^{-1}L_{\overline{U},U}.
\]
We allow the operator $|$ to operate on nonnegative matrices $R$
too, and define it as 
\[
R|_{U}=I_{U}-(I-R)|_{U}=I_{U}-L|_{U}=A+B(I-D)^{-1}C.
\]
\end{definition}

\section{Properties of Schur Complements\label{text:(Schur-Complement)-PROPS}}

We now list some properties of Schur complements, or specifically
the $|$ operator as defined in  \ref{def:(Schur-Complement)}.
\begin{enumerate}
\item \textbf{Well-defined.} Schur complement is well-defined. This follows
from Lemma \ref{lem:L-invertible} since $I-D$ is invertible. 
\item \textbf{Transitivity.} Let $U$ and $V$ be two subsets of indices
with $U\subseteq V$. Then 
\[
(L|_{V})|_{U}=L|_{U}
\]
This can be easily seen from the equations \ref{eq:schurc-1}, \ref{eq:schurc-2},
\ref{eq:schurc-3}, \ref{eq:schurc-4}, since solving a subset of
equations cannot change the final solution. It can also be seen by
direct calculation from the definition but it is long and not insightful.

\item \textbf{Commutativity with the inverse.} $(L^{-1})_{U}=(L|_{U})^{-1}$.
Follows again from equations \ref{eq:schurc-1}, \ref{eq:schurc-2},
\ref{lem:L-invertible}, \ref{eq:schurc-4} or by direct calculation
from the definition. 
\item \textbf{Commutativity with the transpose.} $L^{T}|_{U}=(L|_{U})^{T}$,
again follows immediately by transposing the definition. 
\item \textbf{Closure for Laplacians.} $L|_{U}$ is the Laplacian of an
irreducible nonnegative matrix $R|_{U}$ with left and right eigenvectors
$u_{U}$ and $v_{U}$ with eigenvalue 1. To see this, first note that
\[
R|_{U}=A+B(I-D)^{-1}C=A+B\sum_{i=0}^{\infty}D^{i}C
\]
is nonnegative, where the series expansion is well defined since $I-D$ is
invertible (Lemma \ref{lem:L-invertible}) and $\|D\|<1$ since $R$
is irreducible, and each of $A$, $B$, $C$, $D$ are nonnegative
since $R$ was nonnegative. Further, note that 
\[
R|_{U}v_{U}=v_{U},\ \ R^{T}|_{U}u_{U}=u_{U}.
\]
We show for $R|_{U}$ and the other equation follows directly from
(4) above. To see this equation, note that $Lv=0$, which gives $A'v_{1}+B'v_{2}=0$
and $C'v_{1}+D'v_{2}=0$, or $v_{2}=-D'^{-1}C'v_{1}$ (since $D'$
is invertible) and substituting in the first equation gives $L|_{U}v_{U}=0$
which gives $(I-R|_{U})v_{U}=0$ as required. 
\item \textbf{Closure for PSD Laplacians.} If $L$ is positive semidefinite,
then so is $L|_{U}$. To see this, note that since $L$ is psd, $\langle x_{1},A'x_{1}\rangle+\langle x,B'y\rangle+\langle y,B^{*}x\rangle+\langle y,Dy\rangle\geq0$
for all $x,y$, and now set $y=-D^{-1}B^{*}x$, then since $L$ is
hermitian, so is $D$ and $BD^{-1}B^{*}$ (in fact also psd, considering
vectors only supported on $D$ shows $D$ is psd, and so is $D^{-1}$,
and so is $BD^{-1}B^{*}$), thus, the above change in $y$ gives $\langle x,Ax\rangle-\langle x,BD^{-1}B^{*}x\rangle\geq0$,
showing that the schur complement is psd.
\item \textbf{Eigenvalue-Eigenvector pairs.} Let $Rx=\lambda x$, then $(\lambda,x|_{U})$
is an eigenpair for $(A+B(\lambda I-D)^{-1}C)$. This easily follows
again by using equations \ref{eq:schurc-1}, \ref{eq:schurc-2}, \ref{eq:schurc-3},
\ref{eq:schurc-4} and part (5) above. 
\item \textbf{Combinatorial interpretation.} The Schur complement of the
Laplacian $L$ on some set of vertices $U$, is \emph{exactly} the
operation of replacing the weight of every edge $(i,j)$ in $R$ with
the sum of the weights of all paths of all length that go from $i$
to $\overline{U}$ to $j$, i.e. all intermediate vertices in the
path from $i$ to $j$ come from $\overline{U}$ (including the empty
set). This follows directly by looking at $R|_{U}:$
\[
I-L|_{U}=R|_{U}=A+B(I-D)^{-1}C=A+B\sum_{i=0}^{\infty}D^{i}C.
\]
This provides a robust understanding of Schur complement of the Laplacian
as a combinatorially meaningful operation. Essentially, it is arrived
at by removing vertices from the graph, and adding the ``effect''
those vertices had on the rest of the graph with regards to different
paths.
\end{enumerate}
Given these properties, we move forward to the main problem of this
section, and we will explore some more properties of schur complements
as we uncover more lemmas.

\section{\label{subsec:Capacity-and-Schur}Capacity and the Dirichlet Lemmas}

In this section, we begin exploring a different combinatorial notion
of expansion that has been studied at different points in time within
different communities under different names. Our first aim is to provide an 
exposition of a unified framework that will comprehensively cover
all concepts and definitions, hopefully from multiple perspectives.
Our next aim will be to derive many lemmas for the nonreversible versions
of these quantities, that to the best of our knowledge, are not known
as we state them. 

We start by another problem, very similar to the basic problem \ref{def:(Inverse-of-Laplacians)},  called the Dirichlet problem. We will
assume $R$ to be an irreducible nonnegative matrix with positive
$u$ and $v$ as left and right eigenvectors for eigenvalue 1, normalized
so that $\sum_{i}u_{i}v_{i}=1$. Further, since $R$ is irreducible,
we have from Lemma \ref{lem:L-invertible} that $L=I-R$ contains
only the space spanned by $v$ in its right kernel, and the space
spanned by $u$ in its left kernel. 

Note that in problem \ref{def:(Inverse-of-Laplacians)}, we dealt
with the case in which we were given $b$, and we wanted to find $x$
such that $Lx=b$. Our problem now is modified to the case in which
we have a \emph{subset} of values in $b$ and $x$, and we want to
find the remaining values. This is called the Dirichlet problem.
\medskip

\begin{definition}
\label{def:(Dirichlet-Problem)}(Dirichlet Problem) Let $U$ be a
subset of the vertices, then the Dirichlet problem given the boundary
$U$, a vector $a\in R^{|U|}$ and $b\in R^{|\overline{U}|}$ is to
find a vector $q\in R^{n}$ (called the Dirichlet vector) such that
$q_{U}=a$ and $(LD_{v}q){}_{\overline{U}}=b$. Typically, $b=0$.
\end{definition}

Note that we ``normalize'' our vector $q$ by multiplying (or ``re-weighing'')
it with $D_{v}$ before considering the effect of $L$ on it. This
is similar to what we did for the definition of edge expansion (see definition \ref{def:gen-edge-expansion}).
The direct solution of the Dirichlet problem follows directly by using
equations similar to \ref{eq:schurc-1}, \ref{eq:schurc-2}, \ref{eq:schurc-3},
\ref{eq:schurc-4}. This brings us to the most important definition
of this section.
\medskip

\begin{definition}
\label{def:(Capacity)-}(Capacity) Consider the setting of the Dirichlet
problem above with $b=0$, i.e. let $q$ be some vector with $q_{U}=a$
and $(LD_{v}q)_{\overline{U}}=0$. Then we define the capacity of
the set $U$ with values $a$ on the boundary as 
\[
\text{cap}_{U,a}(R)=\langle D_{u}q,L_{R}D_{v}q\rangle.
\]
We will write $\text{cap}_{U,a}(R)=\text{cap}(R)$ when $U$ and $a$
are clear from context.
\end{definition}

The immediate consequence of the definition is the following, similar
to Lemma \ref{lem:transRtoA} and \ref{lem:The-mixing-time-A-R}.

\begin{lemma}
\label{lem:capA-R}Let $A=D^{1/2}RD^{-1/2}$ where $D=D_{u}D_{v}^{-1}$,
then for any $U$ and any $a$, 
\[
\text{cap}_{U,a}(A)=\text{cap}_{R}(U,a)
\]
\end{lemma}

\begin{proof}
The proof is immediate from the definition, noting that $u,v,w$ are
positive vectors.
\end{proof}
Given Lemma \ref{lem:capA-R}, we will exclusively deal with $A=D_{u}^{1/2}D_{v}^{-1/2}RD_{v}^{1/2}D_{u}^{-1/2}$
onwards, since it makes calculations simpler and intelligible due
to it having the same principal left and right eigenvector.
\medskip

\begin{remark}
(Alternate definition of Capacity) Another way to define capacity
is as follows. Consider the setting of the Dirichlet problem above
with $b=0$, i.e. let $q$ be some vector with $q_{U}=a$ and $(LD_{v}q)_{\overline{U}}=0$.
Then we have defined the capacity of the set $U$ with values $a$
on the boundary as 
\[
\text{cap}_{U,a}(R)=\langle D_{u}q,L_{R}D_{v}q\rangle.
\]
We note two modifications that will not change the inherent meaning
of capacity.

First, the definition can be equivalently written with treating $y=D_{v}q$
as a vector, and imposing the condition that $y_{S}=D_{v}^{-1}a$,
but since the string $a$ was arbitrary, we can even write $y_{S}=a$,
where the $\text{cap}_{U,a}(R)$ in Definition \ref{def:(Capacity)-}
will be cap$_{U,D_{v}^{-1}a}(R)$ in the new definition. Thus we have
\[
\text{cap}_{U,a}(R)=\langle D_{u}D_{v}^{-1}y,L_{R}y\rangle.
\]
Second, the diagonal matrix $D_{u}D_{v}^{-1}$ only redefines the
inner product, since we can write the capacity as 
\[
\text{cap}_{U,a}(R)=\langle D_{u}^{1/2}D_{v}^{-1/2}y,D_{u}^{1/2}D_{v}^{-1/2}L_{R}y\rangle
\]
which is only a rescaling of capacity, and thus even the diagonal
multiplication is unnecessary. As a consequence, we can also simply
define capacity as 
\[
\text{cap}_{U,a}(R)=\langle q,L_{R}q\rangle
\]
with $q_{U}=a$ and $(L_{R}q)_{\overline{U}}=0$. However, we prefer
Definition \ref{def:(Capacity)-} to the above one, since it will
become equivalent to a new type of expansion that we will define later.
\end{remark}

Based on the remark above, we can indeed define the capacity of a
nonnegative matrix $A$ with identical left and right eigenvector $w$ as
follows.

\begin{lemma}
\label{lem:(Capacity-simpler-definition)}(Simpler definition of Capacity)
Let $L=I-A$ where $A$ is an irreducible nonnegative matrix with
largest eigenvalue 1 with $w$ as the corresponding left and right
eigenvector. Let $U\subset[n]$ and $a\in\mathbb{R}^{|U|}$, and let
$q$ be such that $q_{U}=(D_{w})_{U}a$, and $(Lq)_{\overline{U}}=0$.
Then
\[
\text{cap}_{U,a}(A)=\text{cap}(A)=\langle q,Lq\rangle.
\]
Alternately, let $q$ be a vector such that $q_{U}=a$ and $(Lq)_{\overline{U}}=0$,
then for $\overline{a}=(D_{w}^{-1})_{U}a$,
\[
\text{cap}_{U,\overline{a}}(A)=\text{cap}(A)=\langle q,Lq\rangle
\]
\end{lemma}

\begin{proof}
The proofs are immediate from the definition.
\end{proof}
These are the two definitions we will use henceforth in this section.

\begin{lemma}
\label{lem:(Equivalent-definitions-of-cap}(Equivalent definitions
of capacity) All the following are equivalent definitions of capacity.
Let $r=L_{A}D_{w}q$ or $r=L_{R}D_{v}q$ (as will be clear from context)
be the vector with $r_{\overline{U}}=0$, and $q_{U}=a$. Then,
\begin{align*}
\text{cap}_{U,a}(R) & =\langle D_{u}q,L_{R}D_{v}q\rangle\\
 & =\langle(D_{u})_{U}a,L_{R}|_{U}(D_{v})_{U}a\rangle\\
 & =\langle D_{u}D_{v}^{-1}L_{R}^{+}r,r\rangle\\
 & =\langle r_{U},(L_{R}^{+})_{U}r_{U}\rangle
\end{align*}
and
\begin{align*}
\text{cap}_{U,a}(A) & =\langle D_{w}q,L_{A}D_{w}q\rangle\\
 & =\langle D_{w}a,L_{A}|_{U}D_{w}a\rangle\\
 & =\langle r,L_{A}^{+}r\rangle\\
 & =\langle r{}_{U},(L_{A}^{+})_{U}r_{U}\rangle.
\end{align*}
\end{lemma}

\begin{proof}
The first equality is the definition, the second equality is obtained
by noting that $(L_{A}q)_{\overline{U}}=0$ and using equations similar
to \ref{eq:schurc-1}, \ref{eq:schurc-2}, \ref{eq:schurc-3}, \ref{eq:schurc-4},
the third equality follows by writing $q=D_{v}^{-1}L_{R}^{+}r$ or
$q=D_{w}^{-1}L_{A}^{+}r$ in the definition, and the last equality
follows from property 3 of Schur complements by noting that $r_{\overline{U}}=0$.
\end{proof}
We now reach the first important lemma, albeit for reversible $R$
or symmetric $A$ (see Lemma \ref{lem:transRtoA}).

\begin{lemma}
\label{lem:(Dirichlet-Lemma-for-sym}(Dirichlet Lemma for symmetric
$A$) Let $L=I-A$ be psd, and for some fixed $a$, let $q$ be such
that $q_{U}=a$ and $(Lq)_{\overline{U}}=0$. Let $x$ be any vector
such that $x_{U}=a$. Then for $\overline{a}=(D_{w}^{-1})_{U}a$,
\[
\text{cap}_{U,\overline{a}}(A)=\langle q,Lq\rangle\leq\langle x,Lx\rangle.
\]
\end{lemma}

\begin{proof}
Since $L$ is psd, we immediate get by Cauchy-Schwarz that 
\[\langle y,Lz\rangle{}^{2}\leq\langle y,Ly\rangle\langle z,Lz\rangle\]
since $L$ is psd and its square root always exists, and thus 
\[\langle y,Lz\rangle^{2}=\langle L^{1/2}y,L^{1/2}z\rangle^{2}\leq\langle L^{1/2}y,L^{1/2}y\rangle\langle L^{1/2}z,L^{1/2}z\rangle^{2}.\]
Let $y=D_{w}x$ and $z=D_{w}q$, then since $(Lz)_{\overline{U}}=0$,
we get that $\langle y,Lz\rangle=\langle z,Lz\rangle$ since $y$
and $z$ are equal in the indices  that are in $U$. Combining the
two gives the lemma, $\langle z,Lz\rangle^{2}=\langle y,Lz\rangle^{2}\leq\langle y,Ly\rangle\langle z,Lz\rangle$.
\end{proof}
This lemma is important, since it is a basic step in a large number
of other lemmas. It shows that $q$ is in fact the minimizer of $\langle x,Lx\rangle$
for all $x$ with $x_{U}=a$. This gives us our first lemma for the
capacity of nonreversible matrices.

\begin{lemma}
\label{lem:cap-A-tilddeA}(Lower bound on the capacity of nonreversible
matrices) For any $U$ and $a$, 
\[
\text{cap}_{U,a}(\tilde{A})\leq\text{cap}_{U,a}(A),
\]
where $\tilde{A}=\frac{1}{2}(A+A^{T})$.
\end{lemma}

\begin{proof}
Note that $L$ and $\tilde{L}=\frac{1}{2}(L+L^{T})$ have the same
left and right eigenvector $w$ for PF eigenvalue 1. Let $q$ be the
Dirichlet vector for $L$ and $g$ for $\tilde{L}$ where $q_{U}=a$
and $g_{U}=a$. Since $\langle x,Lx\rangle=\langle x,\tilde{L}x\rangle$
for any real $x$, setting $x=q$ in Lemma \ref{lem:(Dirichlet-Lemma-for-sym}
gives 
\[
\langle g,\tilde{L}g\rangle\leq\langle q,\tilde{L}q\rangle=\langle q,Lq\rangle
\]
which gives
\[
\text{cap}_{U,\overline{a}}(\tilde{A})\leq\text{cap}_{U,\overline{a}}(A)
\]
 as required, which holds for all $a$ since $w$ is positive and
any vector can be obtained as $\overline{a}$.
\end{proof}
In fact, at this point, we can say something much stronger than Lemma
\ref{lem:cap-A-tilddeA}. Before we proceed further, we need to understand
the idea of clumping of vertices. These are essentially linear-algebraic
modifications to ensure that the solutions to equations remain valid.
Special cases of this idea have been used at many places without sufficient
justification, and we provide complete clarification for its versatile
use.

\noindent \textbf{Vertex Clumping. }The main idea is as follows. Let
$Q\in\mathbb{C}^{n\times n}$ be any matrix, and assume $Q=\begin{bmatrix}A & B\\
C & D
\end{bmatrix}$ without loss of generality, where neither $A$ nor $D$ have size
0. Let the vertices on which $A$ is supported be $U=\{1,...,r\}$.
For some vector $x\in\mathbb{C}^{n}$, if $x\in\text{Ker}(Q)$, then
$\langle x,Qx\rangle=0$. Assume that $x\not\in\text{Ker}(Q)$. The
primary aim of vertex clumping is to write 
\[
\langle x,Qx\rangle=\dfrac{1}{\langle y,Ty\rangle}
\]
for some $y$ that is independent of $x$ and $Q$, and $T$ that
can be expressed as a function of $Q$ and $x$. Our main result here
is the following lemma.

\begin{lemma}
\textbf{\label{lem:(Vertex-Clumping)}(Vertex Clumping)} Let $Q=\begin{bmatrix}A & B\\
C & D
\end{bmatrix}\in\mathbb{C}^{n\times n}$ with $A$ supported on vertices $U\subset[n]$, a given string $a\in\mathbb{C}^{|U|}$
such that $a\not=w_{U}$ for any $w\in\text{Ker}(Q)$, and the Dirichlet
vector $x$ such that $x_{U}=a$ and $(Qx)_{\overline{U}}=0$, then
let $x_{\overline{U}}=z$ and let $y\in\mathbb{C}^{n-r+1}$ and $T\in\mathbb{C}^{(n-r+1)\times(n-r+1)}$,
with 
\[
y=\begin{bmatrix}1\\
z
\end{bmatrix}\text{\ and\ }T=\begin{bmatrix}a^{*}Aa & a^{*}B\\
Ca & D
\end{bmatrix}.
\]
Then
\[
\langle x,Qx\rangle=\langle y,Ty\rangle=\dfrac{1}{\langle e_{1},T^{+}e_{1}\rangle}
\]
where $e_{1}$ is the standard basis vector.
\end{lemma}

\begin{proof}
Note that
\begin{align*}
\langle x,Qx\rangle & =\begin{bmatrix}a^{*} & z^{*}\end{bmatrix}\begin{bmatrix}A & B\\
C & D
\end{bmatrix}\begin{bmatrix}a\\
z
\end{bmatrix}\\
 & =a^{*}Aa+a^{*}Bz+z^{*}Ca+z^{*}Dz\\
 & =\begin{bmatrix}1 & z^{*}\end{bmatrix}\begin{bmatrix}a^{*}Aa & a^{*}B\\
Ca & D
\end{bmatrix}\begin{bmatrix}1\\
z
\end{bmatrix}\\
 & =\langle y,Ty\rangle.
\end{align*}
Thus, if \textbf{$b=Qx$} and $c=Ty$ or $y=T^{+}c$, letting $c_{\{1\}}=r\not=0$
(since $x$ is not in the kernel of $Q$), and since $c_{\{2,...,n-r+1\}}=0$,
we get
\begin{align*}
\langle a,Ba\rangle+\langle a,Bz\rangle & =c_{\{1\}}=r\\
b_{\{r+1,...,n\}} & =c_{\{2,...,n-r+1\}}=0\\
\langle x,Qx\rangle=\langle y,Ty\rangle=\langle Q^{+}b,b\rangle=\langle T^{+}c,c\rangle & =r
\end{align*}
What is essentially done here, is to clump the vertices $U$ in $Q$
with respect to $x$, by taking a weighted sum of all the edges in
the cluster corresponding to $U$. Since $Qx|_{\overline{U}}=0$,
then $\langle x,Qx\rangle=\langle y,Ty\rangle=\langle y,c\rangle=r$,
and also $c=r\cdot e_{1}$ (where $e_{1}$ is the standard basis vector),
which gives
\[
\langle T^{+}c,c\rangle=r^{2}\langle T^{+}e_{1},e_{1}\rangle
\]
which gives 
\[
\langle T^{+}e_{1},e_{1}\rangle=\dfrac{1}{r}=\dfrac{1}{\langle x,Qx\rangle}.
\]
\end{proof}
As a consequence of vertex clumping, we get the following crucial
lemma. Previous lemmas of this form have existed in folklore, but
they had the following constraints.
\begin{enumerate}
\item They were always for symmetric laplacians
\item They always had the string $a$ being a bit-string corresponding to
vertices with potential 0 and 1
\item There were two clusters of vertices corresponding to where the current
was entering and leaving the graph, and instead of $e_{1}$ in lemma
\ref{lem:(Vertex-Clumping)}, they had vectors of the form $e_{i}-e_{j}$
corresponding to currents.
\end{enumerate}
In comparison, Lemma \ref{lem:(Vertex-Clumping)} holds for any laplacian
(not necessarily symmetric), requires only one cluster of vertices,
and does away altogether with the electrical point of view.

\begin{lemma}
\label{lem:cap-inverse-vertex-clumped}Let $A$ be any nonnegative
matrix with largest eigenvalue 1 for left and right eigenvector $w$,
let $L=I-A$, then for any set $U\subset[n]$ and $\overline{a}=(D_{w}^{-1})_{U}a$,
with 
\[
L=\begin{bmatrix}A & B\\
C & D
\end{bmatrix},\ \text{and}\ T=\begin{bmatrix}a^{*}Aa & a^{*}B\\
Ca & D
\end{bmatrix},
\]
we get 
\[
\text{cap}_{U,\overline{a}}(A)=\dfrac{1}{\langle e_{1},T^{+}e_{1}\rangle}
\]
where $e_{1}$ is the standard basis vector.
\end{lemma}

\begin{proof}
The proof is immediate from the vertex clumping in Lemma \ref{lem:(Vertex-Clumping)}.
\end{proof}
As a consequence of this, we can now present one of our main lemmas.

\begin{lemma}
\label{lem:main-inverseR}Let $R\in\mathbb{C}^{n\times n}$ be any
matrix such that $\tilde{R}=\frac{1}{2}(R+R^{*})$ is positive definite.
Then the following hold:
\begin{enumerate}
\item $R$ is invertible
\item $\frac{1}{2}($$R^{-1}+(R^{*})^{-1})$ is positive definite.
\item 
\[
\left(\frac{R+R^{*}}{2}\right)^{-1}\succcurlyeq\frac{R^{-1}+(R^{*})^{-1}}{2}
\]
\item 
\[
\frac{R+R^{*}}{2}\succcurlyeq\frac{R^{*}R^{-1}R^{*}+R(R^{*})^{-1}R}{2}
\]
\item Let $\tilde{R}=\frac{1}{2}(R+R^{*})$ and $\overline{R}=\frac{1}{2}(R-R^{*})$.
Let $W_{\alpha}=\tilde{R}+\alpha\overline{R}$ for $-1\leq\alpha\leq1$.
Then for any $\alpha,\beta$ such that $|\alpha|\leq|\beta|$,
\[
\frac{W_{\alpha}^{-1}+(W_{\alpha}^{*})^{-1}}{2}\succcurlyeq\frac{W_{\beta}^{-1}+(W_{\beta}^{*})^{-1}}{2}.
\]
\end{enumerate}
\end{lemma}

\begin{proof}
For (1), let $(\lambda,v)$ be some eigenvalue-eigenvector pair for
$R$ where $v$ is nonzero, then 
\[
\text{Re}\lambda=\text{Re}\langle v,Rv\rangle=\frac{1}{2}(\langle v,Rv\rangle+\overline{\langle v,Rv\rangle})=\frac{1}{2}(\langle v,Rv\rangle+\langle Rv,v\rangle)=\langle v,\frac{R+R^{*}}{2}v\rangle>0
\]
and thus $R$ has no 0 eigenvalues and is invertible.

For (2), note that for any $u\in\mathbb{C}^{n}$ such that $u$ is
nonzero,
\begin{align*}
 \left<u,\frac{R^{-1}+(R^{*})^{-1}}{2}u\right>
 & =\left<u,R^{-1}\cdot\frac{R+R^{*}}{2}\cdot(R^{*})^{-1}u\right> \\
 & =\left<(R^{*})^{-1}u,\frac{R+R^{*}}{2}\cdot(R^{*})^{-1}u\right> \\
 & =\left<v,\frac{R+R^{*}}{2}\cdot v\right> \\
 & > 0
\end{align*}

where the last inequality follows since $u\not=0$, and $(R^{*})^{-1}u$
is also nonzero since $R^{*}$ (from (1)) and thus $(R^{*})^{-1}$
is invertible, and thus full rank, implying its kernel is trivial.

For (3), let $L=\tilde{R}=\frac{1}{2}(R+R^{*})$ and $Q=\overline{R}=\frac{1}{2}(R-R^{*})$
(we define these new symbols for ease of readability of the proof).
Note that $L=L^{*}$, and since $L$ is positive definite, $L^{-1}$
exists and is well-defined. Further, $Q$ is skew-hermitian, i.e.
$Q^{*}=-Q$. We can also rewrite $R=L+Q$ and $R^{*}=L-Q$. We then
have that 

\begin{align*}
\left<u,\left(\frac{R+R^{*}}{2}\right)^{-1}u\right>-\left<u,\frac{R^{-1}+(R^{*})^{-1}}{2}u\right> & =\left<u,L^{-1}u\right>-\left<u,R^{-1}\cdot L\cdot(R^{*})^{-1}u\right>\\
 & =\left<u,L^{-1}u\right>-\left<(R^{*})^{-1}u,L\cdot(R^{*})^{-1}u\right>\\
 & =\left<R^{*}v,L^{-1}R^{*}v\right>-\left<v,Lv\right>\\
 & \ \ \ \ \text{[setting \ensuremath{u=R^{*}v}]}\\
 & =\left<(L-Q)^{*}v,L^{-1}(L-Q)^{*}v\right>-\left<v,Lv\right>\\
 & =\left<v,(L-Q)L^{-1}(L-Q)^{*}v\right>-\left<v,Lv\right>\\
 & =\left<v,(L-Q)L^{-1}(L+Q)v\right>-\left<v,Lv\right>\\
 & \ \ \ \ \text{[since \ensuremath{(L-Q)^{*}=L^{*}-Q^{*}=L+Q}]}\\
 & =\left<v,(L-Q+Q-QL^{-1}Q)v\right>-\left<v,Lv\right>\\
 & =\left<v,-QL^{-1}Qv\right>\\
 & =\left<v,Q^{*}L^{-1}Qv\right>\\
 & \ \ \ \ \text{[since \ensuremath{Q^{*}=-Q}]}\\
 & =\langle Qv,L^{-1}Qv\rangle\\
 & \geq0\\
 & \ \ \ \ \text{[since \ensuremath{L^{-1}} is positive definite]}
\end{align*}

For (4), starting with (3), we have 
\begin{align*}
 & \frac{1}{2}\left(\left(\frac{R+R^{*}}{2}\right)^{-1}+\left(\frac{R+R^{*}}{2}\right)^{-1}\right)\succcurlyeq\frac{R^{-1}+(R^{*})^{-1}}{2}\\
\Leftrightarrow & \left(\frac{R+R^{*}}{2}\right)^{-1}L\left(\frac{R+R^{*}}{2}\right)^{-1}\succcurlyeq R^{-1}L(R^{*})^{-1}\\
\Leftrightarrow & L\succcurlyeq LR^{-1}L(R^{*})^{-1}L\\
\Leftrightarrow & L\succcurlyeq\left(\frac{R+R^{*}}{2}\right)R^{-1}\left(\frac{R+R^{*}}{2}\right)(R^{*})^{-1}\left(\frac{R+R^{*}}{2}\right)\\
\Leftrightarrow & L\succcurlyeq\frac{1}{8}\left(I+R^{*}R^{-1}\right)\left(R+R^{*}\right)\left((R^{*})^{-1}R+I\right)\\
\Leftrightarrow & L\succcurlyeq\frac{1}{8}\left(R+R^{*}+R^{*}+R^{*}R^{-1}R^{*}\right)\left((R^{*})^{-1}R+I\right)\\
\Leftrightarrow & L\succcurlyeq\frac{1}{8}\left(R+R^{*}+R^{*}+R^{*}R^{-1}R^{*}+R(R^{*})^{-1}R+R+R+R^{*}\right)\\
\Leftrightarrow & L\succcurlyeq\frac{3}{4}L+\frac{1}{4}\cdot\frac{R^{*}R^{-1}R^{*}+R(R^{*})^{-1}R}{2}\\
\Leftrightarrow & \frac{R+R^{*}}{2}\succcurlyeq\frac{R^{*}R^{-1}R^{*}+R(R^{*})^{-1}R}{2}
\end{align*}

as required.

For (5), note that $W_{\alpha}^{*}=W_{-\alpha}$. In (3) above, we
showed that 
\[
W_{0}^{-1}\succcurlyeq\frac{W_{1}^{-1}+W_{-1}^{-1}}{2}.
\]
We now want to extend this inequality. First note that for every $\alpha$,
$\frac{1}{2}(W_{\alpha}+W_{\alpha}^{*})=L$ and thus from (1), $W_{\alpha}$
is invertible for every $\alpha$. To show the inequality, let $K=W_{\beta}$,
then 
\begin{align*}
W_{\alpha} & =\frac{K+K^{*}}{2}+\frac{\alpha}{\beta}\frac{K-K^{*}}{2}=\frac{1}{2}\left(1+\frac{\alpha}{\beta}\right)K+\frac{1}{2}\left(1-\frac{\alpha}{\beta}\right)K^{*}=cK+(1-c)K^{*}
\end{align*}
for 
\[
0\leq c=\frac{1}{2}\left(1+\frac{\alpha}{\beta}\right)\leq1.
\]
Then we have that
\begin{align*}
 & \frac{W_{\alpha}^{-1}+(W_{\alpha}^{*})^{-1}}{2}\succcurlyeq\frac{W_{\beta}^{-1}+(W_{\beta}^{*})^{-1}}{2}\\
\Leftrightarrow & W_{\alpha}^{-1}\frac{W_{\alpha}+W_{\alpha}^{*}}{2}(W_{\alpha}^{*})^{-1}\succcurlyeq K^{-1}\frac{K+K^{*}}{2}(K^{*})^{-1}\\
\Leftrightarrow & L\succcurlyeq W_{\alpha}K^{-1}L(K^{*})^{-1}W_{\alpha}^{*}\\
\Leftrightarrow & L\succcurlyeq\left(cK+(1-c)K^{*}\right)K^{-1}L(K^{*})^{-1}\left(cK^{*}+(1-c)K\right)\\
\Leftrightarrow & L\succcurlyeq\left(cI+(1-c)K^{*}K^{-1}\right)L\left(cI+(1-c)(K^{*})^{-1}K\right)\\
\Leftrightarrow & L\succcurlyeq\left(cI+(1-c)K^{*}K^{-1}\right)L\left(cI+(1-c)(K^{*})^{-1}K\right)\\
\Leftrightarrow & L\succcurlyeq\left(c^{2}L+(1-c)^{2}K^{*}K^{-1}L(K^{*})^{-1}K+c(1-c)L(K^{*})^{-1}K+c(1-c)K^{*}K^{-1}L\right)\\
\Leftrightarrow & L\succcurlyeq\left(c^{2}+(1-c)^{2}\right)L+c(1-c)\left(L(K^{*})^{-1}K+K^{*}K^{-1}L\right)\\
 & \ \ \ \ \ \ \text{[since \ensuremath{K^{*}K^{-1}L(K^{*})^{-1}K=L}]}\\
\Leftrightarrow & L\succcurlyeq\left(c^{2}+(1-c)^{2}\right)L+c(1-c)L+c(1-c)\left(K(K^{*})^{-1}K+K^{*}K^{-1}K^{*}\right)\\
 & \ \ \ \ \ \ \text{[on expanding \ensuremath{L=\frac{1}{2}(K+K^{*})}]}\\
\Leftrightarrow & \left(c+(1-c)\right)^{2}L\succcurlyeq\left(c^{2}+(1-c)^{2}\right)L+c(1-c)L+c(1-c)\left(K(K^{*})^{-1}K+K^{*}K^{-1}K^{*}\right)\\
\Leftrightarrow & L\succcurlyeq K(K^{*})^{-1}K+K^{*}K^{-1}K^{*}\\
 & \ \ \ \ \ \ \text{[since \ensuremath{0\leq c\leq1}]}\\
\Leftrightarrow & \frac{K+K^{*}}{2}\succcurlyeq K(K^{*})^{-1}K+K^{*}K^{-1}K^{*}
\end{align*}
which is true by (4) since $\frac{1}{2}(K+K^{*})=L$ which is positive
definite.
\end{proof}
As a consequence of these lemmas, we get the second primary theorem
of this thesis.

\begin{theorem}
\label{thm:(Monotonicity-of-Capacity)}(Monotonicity of Capacity)
Let $A$ be any nonnegative matrix with largest eigenvalue 1 for left
and right eigenvector $w$, let $L=I-A$, and let $L_{\alpha}=\tilde{L}+\alpha\overline{L}$
for some $-1\leq\alpha\leq1$ where $\tilde{L}=\frac{1}{2}(L+L^{T})$
and $\overline{L}=\frac{1}{2}(L-L^{T})$, and similarly define $A_{\alpha}$,
$\tilde{A}$ and $\overline{A}$. Then for $|\alpha|\leq|\beta|$,
we have 
\[
\text{cap}_{U,a}(A_{\alpha})\leq\text{cap}_{U,a}(A_{\beta}).
\]
\end{theorem}

\begin{proof}
Let $\alpha$ and $\beta$ be as stated. Let 
\[
L=\begin{bmatrix}A & B\\
C & D
\end{bmatrix},\ \text{and}\ H=\begin{bmatrix}a^{*}Aa & a^{*}B\\
Ca & D
\end{bmatrix},
\]
with $\tilde{H}=\frac{1}{2}(H+H^{T})$ and $\overline{H}=\frac{1}{2}(H-H^{T})$
defined similar to $\tilde{L}$ and $\overline{L}$. Note 
\[
L_{\alpha}=\begin{bmatrix}\frac{1+\alpha}{2}A+\frac{1-\alpha}{2}A^{T} & \frac{1+\alpha}{2}B+\frac{1-\alpha}{2}C^{T}\\
\frac{1+\alpha}{2}C+\frac{1-\alpha}{2}B^{T} & \frac{1+\alpha}{2}D+\frac{1-\alpha}{2}D^{T}
\end{bmatrix},
\]
and let 
\begin{align*}
H_{\alpha} & =\begin{bmatrix}a^{*}(\frac{1+\alpha}{2}A+\frac{1-\alpha}{2}A^{T})a & a^{*}(\frac{1+\alpha}{2}B+\frac{1-\alpha}{2}C^{T})\\
(\frac{1+\alpha}{2}C+\frac{1-\alpha}{2}B^{T})a & \frac{1+\alpha}{2}D+\frac{1-\alpha}{2}D^{T}
\end{bmatrix}\\
 & =\frac{1+\alpha}{2}\begin{bmatrix}a^{*}Aa & a^{*}B\\
Ca & D
\end{bmatrix}+\frac{1-\alpha}{2}\begin{bmatrix}a^{*}A^{T}a & a^{*}C^{T}\\
B^{T}a & D^{T}
\end{bmatrix}\\
 & =\tilde{H}+\alpha\overline{H}.
\end{align*}
Since $a\not=w_{U}$ where $w$ is the unique vector in the kernel
of $L_{A}$, for any vector $x$ that is not in the span of $w$,
$H_{\alpha}$ is invertible for any $\alpha$. Thus, from part (5)
of Lemma \ref{lem:main-inverseR}, we get that 
\[
\frac{H_{\alpha}^{-1}+(H_{\alpha}^{*})^{-1}}{2}\succcurlyeq\frac{H_{\beta}^{-1}+(H_{\beta}^{*})^{-1}}{2},
\]
and since $H_{\alpha}$ is real, we have 
\[
\langle e_{1},H_{\alpha}^{-1}e_{1}\rangle=\langle e_{1},(H_{\alpha}^{*})^{-1}e_{1}\rangle
\]
and thus 
\[
\langle e_{1},H_{\alpha}^{-1}e_{1}\rangle\geq\langle e_{1},H_{\beta}^{-1}e_{1}\rangle
\]
and since 
\[
\text{cap}_{U,\overline{a}}(A_{\alpha})=\frac{1}{\langle e_{1},H_{\alpha}^{-1}e_{1}\rangle}
\]
from Lemma \ref{lem:cap-inverse-vertex-clumped} (note that we are
operating in the space outside the kernel, so $H_{\alpha}^{-1}=H_{\alpha}^{+}$),
thus we get the theorem since the string $\overline{a}$ is arbitrary
since $w$ is positive.
\end{proof}
This provides a sufficiently fine-tuned understanding of the capacity
of non-symmetric nonnegative matrices, since it tells us that the
capacity strictly increases as we move away from symmetry, and even
for two non-symmetric matrices, the one farther from symmetry (in
the sense of Theorem \ref{thm:(Monotonicity-of-Capacity)}) strictly
has higher capacity.

We would also like to remark that since Theorem \ref{thm:(Monotonicity-of-Capacity)}
compares two inner products both for non-symmetric matrices, it is
highly nontrivial, as all the tools for symmetric matrices are no
longer usable.

Our next aim is to show a Dirichlet lemma similar to \ref{lem:(Dirichlet-Lemma-for-sym},
but for non-symmetric nonnegative matrices. This for shown in part
by Slowik \cite{slowik2012note} relatively recently, however, the
results there express the capacity as the infimums and supremums of
different expressions, and the explicit solutions are not presented.
We do that here, by showing the explicit expressions for the lemmas.\footnote{There is a remark regarding explicit expressions in Slowik \cite{slowik2012note}
but it is incorrect.}

\begin{lemma}
\textbf{\emph{\label{lem:(Dirichlet-Lemma-for-nonsym}(Dirichlet Lemma
for non-symmetric matrices)}} Let $L=I-A$ for some irreducible nonnegative
$A$ with largest eigenvalue 1 and the corresponding left and right
eigenvector $w$. Let $q$ be such that $q_{U}=a\not=w_{U}$ and $(Lq)_{\overline{U}}=0$.
Let $x$ be a vector such that $x_{U}=a$. Then 
\[
\text{cap}_{U,\overline{a}}(A)=\langle q,Lq\rangle\leq\langle x,L\tilde{L}^{+}L^{T}x\rangle,
\]
and the inequality is tight when $U=\{s,t\}$ with $a_{s}=1$ and
$a_{t}=0$.
\end{lemma}

\begin{proof}
Note $\langle x,Lq\rangle=\langle q,Lq\rangle$ since $(Lq)_{\overline{U}}=0$,
and $\langle q,Lq\rangle=\langle q,\tilde{L}q\rangle$ since everything
is real. By Cauchy-Schwarz inequality, for $z=\tilde{L}^{+}L^{T}x$,
we have 
\[
\langle q,\tilde{L}z\rangle^{2}\leq\langle q,\tilde{L}q\rangle\langle z,\tilde{L}z\rangle
\]
or 
\[
\langle q,L^{T}x\rangle^{2}\leq\langle q,\tilde{L}q\rangle\langle x,L\tilde{L}^{+}L^{T}x\rangle
\]
which gives the lemma.

Let $x=(L^{T})^{+}\tilde{L}q$, then note that 
\[
\langle(L^{T})^{+}\tilde{L}q,L\tilde{L}^{+}L^{T}(L^{T})^{+}\tilde{L}q\rangle=\langle q,\tilde{L}L^{+}L\tilde{L}^{+}L^{T}(L^{T})^{+}\tilde{L}q\rangle=\langle q,\tilde{L}q\rangle=\langle q,Lq\rangle=\text{cap}_{U,\overline{a}}(A).
\]
Note that the only difficulty is in verifying that $x_{U}=a$. To
do this, we first assume that the set $U=\{s,t\}$ and $a_{s}=1$
and $a_{t}=0$. Then note that if $q$ is the Dirichlet vector for
this particular $U$ and $a$, then since $(Lq)_{\overline{U}}=0$
and $L^{T}w=0$, letting $Lq=b$, we have that $w_{s}b_{s}+w_{t}b_{t}=0$.
Similarly, letting $p$ be the Dirichlet vector for $L^{T}$ and since
$(L^{T}p)_{\overline{U}}=0$ and $Lw=0$, letting $L^{T}p=c$, we
would similarly have $w_{s}c_{s}+w_{t}c_{t}=0$. But since $\langle q,Lq\rangle=\langle p,L^{T}p\rangle$
since the capacity for $L$ and $L^{T}$ are the same (use \ref{lem:(Equivalent-definitions-of-cap}
or property (4)), we have that $b_{s}=c_{s}$, giving $b_{t}=c_{t}$,
and thus $b=c$. With this observation, note that our $x$ becomes
the following.
\begin{align*}
x & =(L^{T})^{+}\tilde{L}q\\
 & =(L^{T})^{+}\cdot\frac{1}{2}(L+L^{T})q\\
 & =\frac{1}{2}q+\frac{1}{2}(L^{T})^{+}Lq\\
 & =\frac{1}{2}q+\frac{1}{2}(L^{T})^{+}b\\
 & =\frac{1}{2}q+\frac{1}{2}(L^{T})^{+}c\\
 & =\frac{1}{2}q+\frac{1}{2}p
\end{align*}
and since $p$ and $q$ agree on $U$, so does $x$, and we get that
the lemma is tight.
\end{proof}
We remark that Lemma \ref{lem:(Dirichlet-Lemma-for-nonsym} is a generalization
of Lemma \ref{lem:(Dirichlet-Lemma-for-sym}, and the latter follows
as a corollary of the former if $L$ was symmetric. A remarkable thing
about Lemma \ref{lem:(Dirichlet-Lemma-for-nonsym} is that we can
write the capacity of any matrix $L$, by using the matrix $H=L\tilde{L}^{+}L^{T}$.
We believe this matrix is interesting enough to deserve study on its
own right. We show one lemma to demonstrate its usefulness. 

\begin{lemma}
\label{lem:(Strengthening-of-H-bound}(Strengthening of Lemma \ref{lem:cap-A-tilddeA})
Let $L=I-A$ for some irreducible nonnegative $A$ with largest eigenvalue
1 and the corresponding left and right eigenvector $w$. Let $U=\{s,t\}$
with $a_{s}=1$ and $a_{t}=0$, and let $q$ be such that $q_{U}=a\not=w_{U}$
and $(Lq)_{\overline{U}}=0$, and similarly let $p$ be the Dirichlet
vector for $L^{T}$. Let $\overline{a}=(D_{w}^{-1})_{U}a$ and $H=L\tilde{L}^{+}L^{T}$.
Then 
\[
\langle p,Hp\rangle=\dfrac{\left(\text{cap}_{U,\overline{a}}(A)\right)^{2}}{\text{cap}_{U,\overline{a}}(\tilde{A})}.
\]
\end{lemma}

\begin{proof}
We first note the following. First, let $\text{cap}_{U,\overline{a}}(A)=c$
and $\text{cap}_{U,\overline{a}}(\tilde{A})=\tilde{c}$. Similar to
the proof of tightness in Lemma \ref{lem:(Dirichlet-Lemma-for-nonsym},
let $p$ be the Dirichlet vector for $L^{T}$, then we know that $Lq=L^{T}p=b$,
and we know that $b_{s}=\text{cap}(A)=c$ and $b_{t}=-c\cdot w_{s}/w_{t}$.
Thus, similar to the vertex clumping lemma, let $\mu$ be the vector
with $\mu_{s}=1$, $\mu_{t}=-w_{s}/w_{t}$, and $\mu_{i}=0$ otherwise.
Then we have 
\begin{align*}
c & =\langle q,Lq\rangle=\langle b,L^{+}b\rangle=c^{2}\langle\mu,L^{+}\mu\rangle\\
\frac{1}{c} & =\langle\mu,L^{+}\mu\rangle
\end{align*}
and similarly, we have 
\[
\langle\mu,\tilde{L}^{+}\mu\rangle=\dfrac{1}{\tilde{c}}.
\]
Now consider 
\begin{align*}
\langle p,L\tilde{L}^{+}L^{T}p\rangle & =\langle b,\tilde{L}^{+}b\rangle\\
 & =c^{2}\langle\mu,\tilde{L}^{+}\mu\rangle\\
 & =\frac{c^{2}}{\tilde{c}}
\end{align*}
as claimed.
\end{proof}

\section{Normalized Capacity and Effective Conductance}

We shall now consider the special case where $U=S\cup T$ with disjoint
$S$ and $T$, and $a_{S}=1$ and $a_{T}=0$. In this case, we refer
$\text{cap}_{U,\overline{a}}(A)=\text{cap}_{S,T}(A)$, which we can
call the \emph{effective conductance} between the sets $S$ and $T$,
and $1/\text{cap}(S,T)$ is generally referred to as the \emph{effective
resistance} between $S$ and $T$. In fact, we can normalize the capacity,
and we get the definition of normalized capacity.
\medskip

\begin{definition}
\textbf{\label{def:(Normalized-capacity)}(Normalized capacity}) Let
$A$ be an irreducible nonnegative matrix with largest eigenvalue
1 and corresponding left and right eigenvector $w$, and $L=I-A$
the corresponding Laplacian. The following quantity is referred to
as the normalized capacity. 
\[
\sigma_{A}=\min_{\langle1_{S},D_{w}^{2}1_{S}\rangle\leq\langle1_{T},D_{w}^{2}1_{T}\rangle}\dfrac{\text{cap}_{S,T}(A)}{\langle1_{S},D_{w}^{2}1_{S}\rangle}.
\]
For exposition, note that for doubly stochastic $A$,

\[
\sigma_{A}=\min_{|S|\leq|T|}\dfrac{\text{cap}_{S,T}(A)}{|S|}
\]
\end{definition}

Another completely different manner of understanding this definition
is through the combinatorial lens discussed in property (8) of Schur
complements. This will provide a beautiful perspective on understanding
Normalized Capacity.

\begin{lemma}
\label{lem:normcap_is_phi}Let $A$ be an irreducible nonnegative
matrix with largest eigenvalue 1 and corresponding left and right
eigenvector $w$ , and $L=I-A$ the corresponding Laplacian. Then
\[
\sigma_{A}=\min_{\langle1_{S},D_{w}^{2}1_{S}\rangle\leq\langle1_{T},D_{w}^{2}1_{T}\rangle}\dfrac{\text{cap}_{S,T}(A)}{\langle1_{S},D_{w}^{2}1_{S}\rangle}=\min_{U}\phi(A|_{U})
\]
\end{lemma}

\begin{proof}
Note from Lemma \ref{lem:(Equivalent-definitions-of-cap}, 
\begin{align*}
\text{cap}_{S,T}(A) & =\langle a,(D_{w})_{U}L|_{U}(D_{w})_{U}a\rangle\\
 & =\langle1_{S},(D_{w})_{U}(I-A|_{U})(D_{w})_{U}1_{S}\rangle\\
 & =\langle1_{S},(D_{w}^{2})_{U}1_{S}\rangle-\langle1_{S},(D_{w})_{U}A|_{U}(D_{w})_{U}1_{S}\rangle\\
 & =\langle1,(D_{w}^{2})_{U}1_{S}\rangle-\langle1_{S},(D_{w})_{U}A|_{U}(D_{w})_{U}1_{S}\rangle\\
 & =\langle1,(D_{w})_{U}A|_{U}(D_{w})_{U}1_{S}\rangle-\langle1_{S},(D_{w})_{U}A|_{U}(D_{w})_{U}1_{S}\rangle\\
 & \ \ \ \text{[since \ensuremath{A_{u}^{T}w_{U}=w_{U}} from properties (5) and (7) of Schur Complements]}\\
 & =\langle1_{\overline{S}},(D_{w})_{U}A|_{U}(D_{w})_{U}1_{S}\rangle
\end{align*}
and thus 
\[
\dfrac{\text{cap}_{S,T}(A)}{\langle1_{S},D_{w}^{2}1_{S}\rangle}=\frac{\langle1_{\overline{S}},(D_{w})_{U}A|_{U}(D_{w})_{U}1_{S}\rangle}{\langle1_{S},D_{w}^{2}1_{S}\rangle}=\phi_{S}(A|_{U})
\]
and first minimizing over $S$ and then $U$ gives the result.
\end{proof}
The Lemma \ref{lem:normcap_is_phi} now helps us obtain a combinatorial
understanding of the notion of (normalized) capacity similar to property
(8) of Schur complements, since it is \emph{exactly} the edge expansion
after some vertices of the graph are removed. We will first show a
simple lemma and then give 2 different interpretations. 

\begin{lemma}
\label{lem:(Maximum-principle)-Let}(Maximum principle) Let $L=I-A$
where $A$ is an irreducible nonnegative matrix with largest eigenvalue
1 and corresponding left and right eigenvector $w$. Let $a\in\{0,1\}^{|U|}$
with $U=S\cup T$ and $a_{S}=1$ and $a_{T}=0$. Let $q$ be a vector
such that $q_{U}=(D_{w})_{U}a$ and $(Lq)_{\overline{U}}=0$. Then
for every entry of $q$, 
\[
0\leq q_{i}\leq w_{i}.
\]
Further, for all $i\in S$ and $j\in T$, 
\[
(Lq)_{i}\geq0,\ (Lq)_{j}\leq0.
\]
\end{lemma}

\begin{proof}
The statement is true for entries $i\in S\cup T$, we need to show
it for entries $i\in\overline{U}$. Letting $A=\begin{bmatrix}P & Q\\
R & W
\end{bmatrix}$ be irreducible nonnegative, and let $L=I-A$, then from $(Lq)_{\overline{U}}=0$,
we have 
\begin{align*}
-Rq_{\overline{U}}+(I-W)q_{U} & =0\\
q_{\overline{U}} & =(I-W)^{-1}Rq_{U}\\
 & =\sum_{i\geq0}W^{i}Rq_{U}\\
 & =\sum_{i\geq0}W^{i}R(D_{w})_{U}a
\end{align*}
where $I-S$ is invertible due to Lemma \ref{lem:L-invertible}. Thus
from the equation above, since $W$, $R$, $w$ and $a$ are all nonnegative,
we get $q_{i}\geq0$ for all $i\in\overline{U}$. To note the upper
bound, first note that the entries of $q_{\overline{U}}$ cannot decrease
with increasing the number of ones in $a$, again because all the
matrices involved are nonnegative. Thus, the maximum entries in $q_{\overline{U}}$
are achieved when $a$ is the all ones vector. But in that case, note
that $q_{U}=w_{U}$, and the vector $q$ for which $(Lq)_{\overline{U}}=0$
is exactly $w$, which would imply $q_{\overline{U}}=w_{\overline{U}}$,
thus showing that for any $i\in\overline{U}$, $q_{i}\leq w_{i}$.

To see $(Lq)_{i}\geq0$ for $i\in S$, first let $b=(Lq)_{U}$, then
\[
b=((I-P)-Q(I-W)^{-1}R)(D_{w})_{U}a.
\]
However, from property (5) of Schur complements, we know that taking
the Schur complement creates the Laplacian of another nonnegative
matrix, specifically $A'=P+Q(I-W)^{-1}R$ with $w_{U}$ as the left
and right PF eigenvector for eigenvalue 1, and thus we have 
\[
b=(I-A')(D_{w})_{U}a.
\]
Now we know that of $a$ is the all ones vector, $b=0$. Thus, assume
$a$ has ones in the first $|S|$ entries, then for $i\in S$, 
\[
b_{i}=(1-A'_{i,i})w_{i}-\sum_{j\in S,j\not=i}A'_{i,j}w_{j},
\]
and since $A'$ is nonnegative $b_{i}=0$ exactly when the entire
row of $A'_{i,j}$ is summed, we get that $b_{i}\geq0$. Similarly,
for $i\in T$, we can write 
\[
b_{i}=-\sum_{j\in T}A_{i,j}w_{j}\leq0
\]
as required.
\end{proof}
Similar to the combinatorial interpretation in Lemma \ref{lem:normcap_is_phi}
and property (8) of Schur complements, we can now give a similar explanation
for $\text{cap}(S,T)$. However, to truly understand capacity intuitively,
we will write the interpretation for doubly stochastic matrices $A$,
for which there is no rescaling by $D_{w}$ since it has a uniform
PF eigenvector and the diagonal matrix becomes identity. It is possible
to write the interpretation for a completely general matrix, but it
will be uninformative. For doubly stochastic matrices, the interpretation
is immediately intuitive.

\begin{lemma}
\label{lem:(Interpretation-of-capacity-prob}(Interpretation of capacity
in terms of probabilities) Let $A$ be a doubly stochastic matrix
with Laplacian $I-A$. Let $S$ and $T$ be two sets with $U=S\cup T$,
and let $q$ be the (Dirichlet) vector such that $q_{i}=1$, $q_{j}=0$
for $i\in S,j\in T$, and $(Lq)_{\overline{S\cup T}}=0$. Let $\Pr_{S,T}(i)$
denote the probability that a random walk on $A$ starting at $i$
hits the set $S$ before hitting the set $T$. Let $v$ be a vector
such that $v_{i}=\Pr_{S,T}(i)$. Then $v=q$.
\end{lemma}

\begin{proof}
The proof is immediate from the following. For $i\in S$, $v_{i}=1$
and for $j\in T$, $v_{j}=0$, thus we need to compute $v_{i}$ for
$i\in\overline{U}$. Let $A=\begin{bmatrix}P & Q\\
R & W
\end{bmatrix}$ where we consider entries according to right multiplication by a
vector, that is $A_{j,i}$ represents the probability of going from
vertex $i$ to $j$. Then note that the probability of hitting the
set $S$ in exactly $t$ steps before hitting the set $T$ from vertex
$i\in U$ is exactly 
\[
(W^{t}Rq_{U})_{i}.
\]
Thus the total probability of hitting the set $S$ before $T$ is
\[
\sum_{t\geq0}(W^{t}Rq_{U})_{i}
\]
which is exactly $q_{\overline{U}}.$
\end{proof}
Before giving the second interpretation we note that we can again
clump the vertices in $S$ and $T$ to single vertices $s'$ and $t'$
and adding all the corresponding edges similar to Lemma \ref{lem:(Vertex-Clumping)}.
Note however that Lemma \ref{lem:(Vertex-Clumping)} is cleaner and
more general since it clumps everything into one super-vertex, and
does not require values 1 and 0. 

\begin{lemma}
Consider the setting of Lemma \ref{lem:(Vertex-Clumping)}. Let $L'$
be the matrix in which all the rows and columns corresponding to $S$
have been summed to a single row/column, and similarly for $T$, creating
new vertices $s'$ and $t'$ in the underlying graph (which is not
doubly stochastic now). Then 
\[
\text{cap}_{S,T}(A)=\text{cap}_{s',t'}(A)
\]
\end{lemma}

\begin{proof}
The rows can be summed up since it is just summing linear equations,
and the columns can be summed up since $v$ has the same value 1 at
every vertex in $S$, and similarly has the same value 0 for every
vertex in $T$. Since the linear equations do not change, the solution
does not change.
\end{proof}
With this, we can now give an interpretation in terms of expectation,
for the specific case where $|S|=|T|$.

\begin{lemma}
\label{lem:(Interpretation-of-capacity-exp}(Interpretation of capacity
in terms of expectations) Let $L=I-A$ with irreducible doubly stochastic
$A$ be such that $|S|=|T|=1$, and we call the vertices $s$ and
$t$. Let $q$ be the unique vector such that $(Lq)_{s}=1$, $(Lq)_{t}=-1$
and $Lq$ is zero everywhere else, and $q_{t}=0$. Then 
\[
q_{s}=\frac{1}{\text{cap}_{s,t}(A)},
\]
and 
\[
q_{i}=\text{E}_{s,t}(i),
\]
where $\text{E}_{s,t}(i)$ is the expected number of times vertex
$i$ is visited in a random walk starting from $s$ before it hits
$t$.
\end{lemma}

\begin{proof}
We provide a detailed proof for completeness. Let $s$ and $t$ be
fixed. Let $X_{u}$ the the number of times $u$ is visited in a random
walk that starts at $s$ and stops on reaching $t$, and similarly
$Y_{i,j}$ for any edge $i\rightarrow j$ be the number of times the
edge is visited in a random walk that starts at $s$ and ends at $t$.
Then $X_{u}=\sum_{i}Y_{i,u}$. Now fix some vertex $k$ and consider
$Y_{k,u}$. Let the distribution of $X_{k}$ be given by $r_{j}$,
for any $j\in\mathbb{N}.$ Thus $\sum_{j}r_{j}=1$. Now fix some $j$,
and let $X_{k}$ be visited exactly $j$ times. Let $Z_{1},...,Z_{j}$
be random variables, such that the $Z_{t}=1$ iff on the $t$'th step,
the walk went from vertex $k$ to vertex $u$, and thus $\mathbb{E}(Z_{t})=p_{k,u}$.
Then given that $X_{k}=j$, we have $Y_{k,u}=\sum_{t=1}^{t-j}Z_{t}$.
Thus, $\mathbb{E}(Y_{k,u})=j\cdot p_{k,u}$. Thus we can write 
\begin{align*}
\mathbb{E}(Y_{k,u}) & =\sum_{j}\mathbb{E}(Y_{k,u}|X_{k}=j)\cdot\Pr[X_{k}=j]\\
 & =\sum_{j=1}^{\infty}\mathbb{E}\left(\sum_{t=1}^{j}Z_{t}^{(j)}\right)\cdot r_{j}\\
 & =\sum_{j=1}^{\infty}\sum_{t=1}^{j}\mathbb{E}Z_{t}^{(j)}\cdot r_{j}\\
 & =\sum_{j=1}^{\infty}j\cdot p_{k,u}\cdot r_{j}\\
 & =p_{k,u}\cdot\mathbb{E}(X_{k})
\end{align*}
and thus
\begin{align*}
\mathbb{E}(X_{u}) & =\sum_{i}\mathbb{E}(Y_{i,u})\\
 & =\sum_{i}p_{i,u}\mathbb{E}(X_{i})
\end{align*}
for every $u$ except $s$ and $t$. Since we never visit $t$, and
alway visit $s$ at least once, we have 
\begin{align*}
\mathbb{E}(X_{t}) & =0\\
\mathbb{E}(X_{s}) & =1+\sum_{i}\mathbb{E}(X_{i})\cdot p_{i,s}
\end{align*}
Note that the equations for all vertices except $t$ completely determine
the equation for $t$, which turns out to be 
\[
\mathbb{E}(X_{t})=-1+\sum_{i}\mathbb{E}(X_{i})\cdot p_{i,s}.
\]
Note that these are exactly the equations $(Lq)_{s}=1$, $(Lq)_{t}=-1$,
$(Lq)_{i}=0$ for $i\not=s,t$, and $q_{t}=0$. Since $A$ is irreducible
and the kernel has dimension 1, setting $q_{t}=0$, we get a unique
vector satisfying the equation. Thus,
\[
q_{u}=\mathbb{E}(X_{u})
\]
 as required. We now proceed to solve the equations. With $L=I-A$,
with $\chi_{s,t}=1_{s}-1_{t}$, we get 
\[
Lq=\chi_{s,t},
\]
 with the solution 
\[
q=L^{+}\chi_{s,t}+c\cdot\mathbf{1}
\]
and with the constraint $v_{t}=0$, it gives 
\[
q=L^{+}\chi_{s,t}-\langle1_{t},L^{+}\chi_{s,t}\rangle\cdot\mathbf{1}.
\]
Note that 
\[
q_{s}=q_{s}-q_{t}=\langle\chi_{s,t},q\rangle=\langle\chi_{s,t},L^{+}\chi_{s,t}\rangle=\frac{1}{\text{cap}_{s,t}(A)}.
\]
where the last equality follows since taking the vector $v=q/q_{s}$,
we have $v_{s}=1$, $v_{t}=0$, and $(Lv)_{i}=0$ for $i\not=s,t$,
and thus 
\begin{align*}
\text{cap}_{s,t}(A) & =\langle v,Lv\rangle=\frac{1}{q_{s}^{2}}\langle q,Lq\rangle=\frac{1}{q_{s}^{2}}\langle\chi_{s,t},L\chi_{s,t}\rangle=\frac{1}{q_{s}^{2}}q_{s}=\frac{1}{q_{s}}.
\end{align*}
\end{proof}
In fact, for symmetric matrices $A$, the probability interpretation
in Lemma \ref{lem:(Interpretation-of-capacity-prob} can be written
in terms of voltages, and the expectation interpretation in Lemma
\ref{lem:(Interpretation-of-capacity-exp} can be written in terms
of currents, but we prefer the mathematical interpretations to avoid
making assumptions about reality.

Our aim now is to further explore the notion of normalized capacity
or $\sigma_{A}$. Towards this, note that we immediately get the following
lemma as a corollary of Theorem \ref{thm:(Monotonicity-of-Capacity)}.

\begin{lemma}
(Monotonicity of Normalized Capacity) Let $A$ be any nonnegative
matrix with largest eigenvalue 1 for left and right eigenvector $w$,
let $L=I-A$, and let $A_{\alpha}=\tilde{A}+\alpha\overline{A}$ for
some $-1\leq\alpha\leq1$ where $\tilde{A}=\frac{1}{2}(A+A^{T})$
and $\overline{A}=\frac{1}{2}(A-A^{T})$. Then for $|\alpha|\leq|\beta|$,
we have 
\[
\sigma_{A_{\alpha}}\leq\sigma_{A_{\beta}}.
\]
\end{lemma}

\begin{proof}
Identical to the proof of Lemma \ref{thm:(Monotonicity-of-Capacity)}.
\end{proof}

With this understanding of $\sigma_{A}$, an interesting question
is to connect $\sigma_{A}$ to the spectral gap of $A$. It turns
out, that normalized capacity is equivalent to the spectral gap, up
to constants, as shown by Schild \cite{schild2018schur} and Miller
et. al. \cite{miller2018hardy}. We present the proof for completeness
mostly following \cite{miller2018hardy}, with some changes to make
it shorter and more streamlined, and then discuss possible versions
for the non-symmetric case.

Before we proceed, we show a sequence of lemmas.

\begin{lemma}
\label{lem:miller-schild-1}$\forall$ $x\in\mathbb{R}^{n}$, $(x_{i}-x_{j})^{2}\leq\sum_{k=i}^{j-1}\dfrac{1}{p_{k}}(x_{k}-x_{k+1})^{2}$
for $\sum p_{k}\leq1$.
\end{lemma}

\begin{proof}
Using the Cauchy-Schwarz inequality, 
\[
(x_{i}-x_{j})^{2}=\left(\sum_{k}\sqrt{p_{k}}\cdot\dfrac{x_{k}-x_{k+1}}{\sqrt{p_{k}}}\right)^{2}\leq\sum_{k=i}^{j-1}\dfrac{1}{p_{k}}(x_{k}-x_{k+1})^{2}.
\]
\end{proof}

\begin{lemma}
\label{lem:miller-schild-2}Let $w_{1},...,w_{n-1}$ be any nonnegative
weights, and $f_{w}(x)=\sum_{k=1}^{n-1}w_{k}(x_{k}-x_{k+1})^{2}$.
Then for any $(i,j)$ with $i<j,$ there exists an $x$ with $x_{k}=1$
for $k\leq i$, $x_{k}=0$ with $k\geq j$, and $f_{w}(x)=\left(\sum_{k=i}^{j-1}\dfrac{1}{w_{k}}\right)^{-1}.$
\end{lemma}

\begin{proof}
Let $x$ be 1 on entries $\leq i$ and 0 on entries $\geq j$, and
for $i\leq k\leq j$, let $x_{k}=\dfrac{\sum_{r=k}^{j-1}\dfrac{1}{w_{r}}}{\sum_{r=i}^{j-1}\dfrac{1}{w_{r}}}$.
Then 
\begin{align*}
f_{w}(x) & =\sum w_{k}(x_{k}-x_{k+1})^{2}=\sum_{k=i}^{j-1}w_{k}(x_{k}-x_{k+1})^{2}=\sum_{k=i}^{j-1}w_{k}\dfrac{\dfrac{1}{w_{k}^{2}}}{\left(\sum_{r=i}^{j-1}\dfrac{1}{w_{r}}\right)^{2}}=\left(\sum_{k=i}^{j-1}\dfrac{1}{w_{k}}\right)^{-1}
\end{align*}
as required.
\end{proof}
We can now present the main lemma in \cite{miller2018hardy,schild2018schur}. We show it for doubly stochastic
matrices for simplicity.

\begin{lemma}
(Normalized Capacity and Spectral Gap \cite{miller2018hardy,schild2018schur})
\label{lem:(Normalized-Capacity-Spec-Gap}Let $A$ be a symmetric
irreducible doubly stochastic matrix with largest eigenvalue 1 and
corresponding left and right eigenvector $\mathbf{1}$. Then the normalized
capacity of $A$ is equivalent to the spectral gap of $A$, up to
constants. Quantitatively,
\[
\frac{1}{2}\cdot\Delta(A)\leq\sigma_{A}\leq4\cdot\Delta(A).
\]
\end{lemma}

\begin{proof}
We start by showing the lower bound, which is straightforward. It
is sufficient to show that $\lambda_{2}(L)\leq\lambda_{2}(L|_{U})$
for any $U$, and the inequality will follow from the Cheeger-Buser
inequality and the equivalence \ref{lem:normcap_is_phi}.

Let $z=[x,y]^{T}$, where $(Lz)(j)=0$ for $j\in\overline{M}$. Let
$x\perp1$ be the vector such that $\lambda_{2}(L|_{M})=\dfrac{\langle x,L|_{M}x\rangle}{\langle x,x\rangle}$.
Thus, there is a $y$ such that for $z=[x,y]^{T}$, $\langle z,Lz\rangle=\langle x,L|_{M}x\rangle$
where $Lz$ at any vertex in $\overline{M}$ is 0. Let $c_{0}=\langle z,1\rangle/n$,
and $z'=z-c_{0}1$, then $z'\perp1$ and $\langle z,Lz\rangle=\langle z',Lz'\rangle$,
and 
\begin{align*}
\langle z',z'\rangle & =\langle z,z\rangle-nc_{0}^{2}\\
 & =\langle x,x\rangle+\langle y,y\rangle-\dfrac{\langle y,1\rangle^{2}}{n}\\
 & \ \ \ [\text{since }x\perp1]\\
 & \geq\langle x,x\rangle
\end{align*}
since letting $m=|\overline{M}|\leq n$, and writing $y=d_{0}1+\mu$
with $\mu\perp1$ and $d_{0}=\langle y,1\rangle/m$, we get 
\[
\langle y,y\rangle-\dfrac{\langle y,1\rangle^{2}}{n}=md_{0}^{2}+\langle\mu,\mu\rangle-m\cdot\dfrac{m}{n}d_{0}^{2}\geq0.
\]
We now proceed to show the upper bound. Let $x$ be a nonnegative
vector such that $x_{i}\geq x_{i+1}$ and let $r\leq n/2$ be such
that $x_{k}=0$ for $k>r$, which exists from Lemma \ref{lem:cheeger}.

Define weights $w_{1}^{i,j},...,w_{r}^{i,j}$, with $w_{k}^{i,j}=\dfrac{x_{i}-x_{j}}{x_{k}-x_{k+1}}$
where it is understood that $w_{k}^{i,j}$ holds only for $i\leq k<j$,
and $\sum_{k=i}^{j-i}\dfrac{1}{w_{k}^{i,j}}=1$. Then note that $(x_{i}-x_{j})^{2}=\sum_{k=i}^{j-1}w_{k}^{i,j}(x_{k}-x_{k+1})^{2}$.
Define $w_{k}=\sum_{i,j}w_{k}^{i,j}$. Then we get that $\langle x,Lx\rangle=f_{w}(x)$
(defined in Lemma \ref{lem:miller-schild-2}), and for all $y$, we
get $\langle y,Ly\rangle\leq f_{w}(y)$ by applying Lemma \ref{lem:miller-schild-1}
to every edge.

Let $R_{i}=\sum_{j=i}^{r}\dfrac{1}{w_{j}},$ and let the vectors $v_{i}$
be such that they have first $i$ entries 1, last $n/2$ entries 0,
remaining entries such $f_{w}(v_{i})=R_{i}^{-1}$ (these exist by
Lemma \ref{lem:miller-schild-2}), then combining the definition of
$\sigma$ \ref{def:(Normalized-capacity)}, the Dirichlet Lemma \ref{lem:(Dirichlet-Lemma-for-sym},
and the observation above, we get

$\sigma\leq\dfrac{\langle v_{i},Lv_{i}\rangle}{i}\leq\dfrac{f_{w}(v_{i})}{i}\leq\dfrac{1}{i\cdot R_{i}}$.

Note that $\dfrac{1}{w_{j}}=R_{j}-R_{j+1}$, and using $\dfrac{a^{2}-b^{2}}{a}\leq2(a-b)$,

\begin{align*}
\sum_{i}x_{i}^{2} & =\sum_{i}^{r}(\sum_{j=i}^{r}x_{j}-x_{j+1})^{2}=\sum_{i}^{r}(\sum_{j=i}^{r}(x_{j}-x_{j+1})w_{j}^{1/2}R_{j}^{1/4}\dfrac{1}{w_{j}^{1/2}R_{j}^{1/4}})^{2}\\
 & \leq\sum_{i}^{r}\sum_{j=i}^{r}(x_{j}-x_{j+1})^{2}w_{j}R_{j}^{1/2}\sum_{j=i}^{r}\dfrac{1}{w_{j}R_{j}^{1/2}}\\
 & =\sum_{i}^{r}\sum_{j=i}^{r}(x_{j}-x_{j+1})^{2}w_{j}R_{j}^{1/2}\sum_{j=i}^{r}\dfrac{R_{j}-R_{j+1}}{R_{j}^{1/2}}\\
 & \leq2\sum_{i}^{r}\sum_{j=i}^{r}(x_{j}-x_{j+1})^{2}w_{j}R_{j}^{1/2}R_{i}^{1/2}\\
 & \leq\dfrac{2}{\sigma}\sum_{i=1}^{r}\sum_{j\geq i}^{r}(x_{j}-x_{j+1})^{2}w_{j}\dfrac{1}{\sqrt{j}}\dfrac{1}{\sqrt{i}}\\
 & =\dfrac{2}{\sigma}\sum_{j=1}^{r}(x_{j}-x_{j+1})^{2}w_{j}\dfrac{1}{\sqrt{j}}\sum_{i=1}^{j}\dfrac{1}{\sqrt{i}}\\
 & \leq\dfrac{4}{\sigma}\sum_{j=1}^{r}(x_{j}-x_{j+1})^{2}w_{j}\\
 & =\dfrac{4}{\sigma}f_{w}(x)\\
 & =\dfrac{4}{\sigma}\langle x,Lx\rangle\\
\dfrac{\sigma}{4} & \leq\dfrac{\langle x,Lx\rangle}{\langle x,x\rangle}
\end{align*}
as required.
\end{proof}
This raises the first interesting question. Since we know that $\rho_{\tilde{L}}\leq\rho_{L}$
and $\tau(\tilde{A})\leq\tau(A)$, does a relationship similar to
Lemma \ref{lem:(Normalized-Capacity-Spec-Gap} hold for nonreversible
chains, i.e., is it possible that $\sigma(A)\leq c\cdot\Delta(A)$
for some constant $c$? This would be appealing as an equivalence,
and also would imply better bounds for the mixing time since there
is a linear dependence on $\sigma$ instead of quadratic as was the
case with $\phi$ in Theorem \ref{thm:intro-phi-delta}. Unfortunately
however, we show that this relation cannot hold.

\begin{lemma}
\label{lem:normcap-cycle-counterex}It is not true that $\sigma_{A}\in O(\Delta(A)).$
Specifically, let $A$ be the directed cycle. Then
\[
\sigma_{A}\approx\sqrt{\Delta(A)}.
\]
Thus, the gap between $\sigma_{A}$ and $\sigma_{\tilde{A}}$ could
be quadratic, since if $A$ is the directed cycle, then $\tilde{A}$
is the undirected cycle, and 
\[
\sigma_{\tilde{A}}\leq\Delta(\tilde{A})\leq\Delta(A).
\]
\end{lemma}

\begin{proof}
Let $A$ be the directed cycle. Then note that for any set $S$ and
$T$ where $|S|=|T|=n/4$, and the vertices are on the opposite sides
of the square, we get that the Dirichlet vector is 1 on half the vertices,
and 0 on the other half, and $A|_{S\cup T}$ is again a cycle, leaving
the spectrum unchanged, except that the spectral gap now increases
since the second eigenvalue will become the $n/2$'th root of unity
with real part closest to 1, instead of the $n$'th root of unity
with real part closest to 1. Thus $\sigma(A)=\phi(A)$, and the lemma
follows.
\end{proof}
Thus, the gap between $\sigma_{A}$ and $\sigma_{\tilde{A}}$ could
be quadratic. In spite of this, it would be interesting to use the
matrix $H$ obtained in Lemma \ref{lem:(Dirichlet-Lemma-for-nonsym}
to obtain a lemma similar to \ref{lem:(Normalized-Capacity-Spec-Gap}
for nonreversible chains.

\section{Relation with log-Sobolev constants}

We state the relation between log-Sobolev constant and the spectral
gap, from \cite{diaconis1996logarithmic}, for the sake of comparison
with other quantities. See \cite{diaconis1996logarithmic} for related definitions Let $A$ be a symmetric matrix with largest
eigenvalue 1 and corresponding left and write eigenvector $w$. Let
$w$ be normalized so that $\sum w_{i}^{2}=1$, and let $w_{0}=\min_{i}w_{i}.$
Define the log-Sobolev constant for $A$ as follows. Let $\alpha$
be the log-Sobolev constant (for symmetric $A$, $\lambda=1-\lambda_{2}(A)$),
then
\[
\dfrac{\lambda}{1+\dfrac{1}{2}\log(1/\pi_{0})}\leq\alpha\leq\lambda
\]
where the lower bound is achieved for the clique, and the upper bound
for the hypercube. Thus,
\[
\dfrac{c_{1}.\rho}{2+\log(1/\pi_{0})}\leq\alpha\leq c_{2}\cdot\rho.
\]

\section{\label{subsec:A-different-notion}A different notion of expansion}

So far, we have seen two primary definitions of expansion. The first
definition was the standard Definition \ref{def:gen-edge-expansion}
expressed as the edge expansion $\phi$, and the second definition arose from normalized
capacity \ref{def:(Normalized-capacity)} expressed as $\sigma$.
In this subsection, we propose a new definition that will be expressed
as $\mu$. We arrived at this definition by trying to understand the
expansion of high dimensional expanders in \cite{bafna2020high}.
As we have done throughout, to conceptually grasp the definition,
we restrict to doubly stochastic matrices.
\medskip

\begin{definition}
\label{def:exp-new}Let $A$ be an irreducible doubly stochastic matrix.
Let $1_{S}$ be the vector with ones on the set $S$ and zeros elsewhere,
and let $1'_{S}=1_{S}/|S|$. Define
\[
\mu_{S}(A)=\frac{1}{2}\|A1'_{S}-A1'_{\overline{S}}\|_{1}
\]
and let 
\[
\mu(A)=\max_{S:|S|\leq n/2}\mu_{S}(A).
\]
\end{definition}

To understand the definition, observe that $\mu$ should behave approximately
like $1-\phi$. Consider $A=J=\frac{1}{n}\mathbf{1}\cdot\mathbf{1}^{T}$.
Then it is clear that $\mu(J)=0$. Consider a graph $A$ with two
almost disconnected components corresponding to $S$ and $\overline{S}$.
Then $\mu(A)\approx1$. However, $\mu$ also has behavior different
from $1-\phi$ in certain graphs. Consider a bipartite graph $A$
with two components. Then for the set $S$ supported on one component,
we get $\mu(A)=1$, and we also have that $\phi(A)\approx\text{\ensuremath{\frac{1}{2}}},$
showing that both $\mu$ and $\phi$ are constants.

\begin{lemma}
\label{lem:new-exp-and-phi}Let $A$ be an irreducible doubly stochastic
matrix. Then
\[
\mu(A)\geq1-2\cdot\phi(A).
\]
\end{lemma}

\begin{proof}
For any set $|S|=k\leq n/2$, let $k'=n-k$, and let $v$ be the vector
of sums of columns of $A$ in $\overline{S}$. By the definition of
$\phi_{S}(A)$, we have that 
\[
\phi_{S}(A)=\dfrac{1}{k}\sum_{i=1}^{k}v_{i}.
\]
Further, since $A$ is $\frac{1}{2}$-lazy, we have that $v_{i}\leq1/2$
for $i\leq k$. Thus we have

\begin{align*}
\mu_{S}(A) & =\dfrac{1}{2}\|A1'_{S}-A1'_{\overline{S}}\|_{1}\\
 & =\dfrac{1}{2}\sum_{i}\left|\dfrac{1-v_{i}}{k}-\dfrac{v_{i}}{k'}\right|\\
 & =\dfrac{1}{2}\sum_{i}\left|\dfrac{1}{k}-\dfrac{n\cdot v_{i}}{k\cdot k'}\right|\\
 & =\dfrac{1}{2}\dfrac{n}{k\cdot k'}\sum_{i}\left|\dfrac{k'}{n}-v_{i}\right|\\
 & \geq\dfrac{1}{2}\dfrac{n}{k\cdot k'}\sum_{i=1}^{k}\left(\dfrac{k'}{n}-v_{i}\right)\\
 & \text{\ \ \ [since \ensuremath{A} is \ensuremath{1/2}-lazy]}\\
 & =\dfrac{1}{2}\left(1-\dfrac{n}{k'}\phi_{S}(A)\right)\\
 & \geq\dfrac{1}{2}-\phi_{S}(A)
\end{align*}
as required. For general $A$, note that t $\sum_{i}v_{i}=k'$, and
thus, 
\[
\sum_{i>k}v_{i}=k'-k\cdot\phi_{S}(A).
\]
Thus we get,
\begin{align*}
\mu_{S}(A) & =\dfrac{1}{2}\|A1'_{S}-A1'_{\overline{S}}\|_{1}\\
 & =\dfrac{1}{2}\dfrac{n}{k\cdot k'}\sum_{i}\left|\dfrac{k'}{n}-v_{i}\right|\\
 & =\dfrac{1}{2}\dfrac{n}{k\cdot k'}\sum_{i}\left(\dfrac{k'}{n}-v_{i}\right)+2\cdot\dfrac{1}{2}\dfrac{n}{k\cdot k'}\sum_{i:v_{i}\geq k'/n}\left(v_{i}-\frac{k'}{n}\right)\\
 & =0+\dfrac{n}{k\cdot k'}\sum_{i:v_{i}\geq k'/n}\left(v_{i}-\frac{k'}{n}\right)\\
 & \geq\dfrac{n}{k\cdot k'}\sum_{i>k:v_{i}\geq k'/n}\left(v_{i}-\frac{k'}{n}\right)
\end{align*}
To understand the minimum value that the sum can take, let there be
$r$ indices $i>k$ such that $v_{i}=k'/n$, and the rest are greater
than $k'/n$. Thus we have 
\[
\sum_{i>k+r}v_{i}=k'\left(1-\frac{r}{n}\right)-k\cdot\phi_{S}(A),
\]
and continuing the summation above, we have 
\begin{align*}
\mu_{S}(A) & \geq\dfrac{n}{k\cdot k'}\sum_{i>k:v_{i}\geq k'/n}\left(v_{i}-\frac{k'}{n}\right)\\
 & =\dfrac{n}{k\cdot k'}\sum_{i>k+r}\left(v_{i}-\frac{k'}{n}\right)\\
 & =\dfrac{n}{k\cdot k'}\left(k'\left(1-\frac{r}{n}\right)-k\cdot\phi_{S}(A)-\frac{k'}{n}(k'-r)\right)\\
 & =\dfrac{n}{k\cdot k'}\left(k'\left(1-\frac{k'}{n}\right)-k\cdot\phi_{S}(A)\right)\\
 & =\frac{n}{k}\left(1-\frac{k'}{n}\right)-\frac{n}{k'}\cdot\phi_{S}(A)\\
 & =1-\frac{n}{k'}\cdot\phi_{S}(A)\\
 & \geq1-2\cdot\phi_{S}(A)
\end{align*}
as required.
\end{proof}
The second lemma we show is nicer, since it is similar to Lemma \ref{lem:(Submultiplicativity-of-phiS)},
except the binary function is minimum instead of addition. From Lemma
\ref{lem:(Submultiplicativity-of-phiS)}, we have that for doubly
stochastic $A$ and $B$, 
\[
\phi_{S}(AB)\leq\phi_{S}(A)+\phi_{S}(B).
\]
The equivalent lemma for $\mu$ is as follows.

\begin{lemma}
\label{lem:new-exp-in-bound}Let $\mu$ be as defined in \ref{lem:new-exp-and-phi},
and let $A$ and $B$ be two doubly stochastic matrices. Then 
\[
\mu_{S}(AB)\leq\min\{\mu_{S}(A),\mu_{S}(B)\}
\]
and thus 
\[
\mu(AB)\leq\min\{\mu(A),\mu(B)\}.
\]
\end{lemma}

\begin{proof}
Fix any set $S$. As the first step, note that 
\[
\mu_{S}(AB)=\dfrac{1}{2}\|AB1'_{S}-AB1'_{\overline{S}}\|_{1}\leq\|A\|_{1}\dfrac{1}{2}\|B1'_{S}-B1'_{\overline{S}}\|_{1}\leq\mu_{S}(B).
\]
Our aim is to now show that $\mu_{S}(AB)\leq\mu_{S}(A)$. Let the
columns of $B$ be $B=(b_{1},..,b_{n})$, let $S$ contain the first
$k$ vertices with $k'=n-k$, let $u=b_{1}+...+b_{k}$, then $B1'_{S}=\dfrac{1}{k}u$,
and $B1'_{\overline{S}}=\dfrac{1}{k'}(1-u)$, and

\[
B1'_{S}-B1'_{\overline{S}}=\dfrac{1}{k}u-\dfrac{1}{k'}(1-u)=\dfrac{n}{k\cdot k'}u-\dfrac{1}{k'}1=\dfrac{n}{k'}\left(\dfrac{u}{k}-\dfrac{1}{n}\mathbf{1}\right)
\]
where $0\leq u_{i}\leq1$ and $\sum_{i}u_{i}=k$ since $B$ is doubly
stochastic.

Let $v$ be a vector with 1 in the first $k$ entries and 0's otherwise,
then $u\preceq v$ (or $u$ is majorized by $v$), and thus $u=Rv$
where $R$ is doubly stochastic. By the Birkhoff-von Neumann theorem
(see Chapter \ref{sec:Preliminaries}), $R$ can be written as $R=\sum\alpha_{i}P_{i}$
with $P_{i}$ being permutations, and $\alpha_{i}\geq0,\sum\alpha_{i}=1$.
Thus 
\[
u=\sum_{i}\alpha_{i}W_{i}
\]
where $W_{i}$ is a vector which has $k$ entries that are 1 and 0's
everywhere else. Now we know that for any set $S$, 
\begin{align*}
\mu_{S}(A) & \geq\dfrac{1}{2}\|A1'_{S}-A1'_{\overline{S}}\|_{1}\\
 & =\dfrac{1}{2}\|\dfrac{1}{k}A1{}_{S}-\dfrac{1}{k'}A1{}_{\overline{S}}\|_{1}\\
 & =\dfrac{1}{2}\|\dfrac{1}{k}A1{}_{S}-\dfrac{1}{k'}(1-A1_{S})\|_{1}\\
 & =\dfrac{1}{2}\dfrac{n}{k'}\|\dfrac{1}{k}A1{}_{S}-\dfrac{1}{n}1\|_{1}
\end{align*}
 and particularly for all $|S|=k$. Thus for any set $S$ with $|S|=k$,
we have 
\begin{align*}
\mu_{S}(AB) & =\dfrac{1}{2}\|AB1'_{S}-AB1'_{\overline{S}}\|_{1}\\
 & =\dfrac{1}{2}\dfrac{n}{k'}\|\dfrac{1}{k}Au-\dfrac{1}{n}1\|_{1}\\
 & =\dfrac{1}{2}\dfrac{n}{k'}\|\sum\alpha_{i}\left(\dfrac{1}{k}AW_{i}-\dfrac{1}{n}1\right)\|_{1}\\
 & \leq\dfrac{1}{2}\dfrac{n}{k'}\sum\alpha_{i}\|\dfrac{1}{k}AW_{i}-\dfrac{1}{n}1\|_{1}\\
 & \leq\mu_{S}(A)
\end{align*}
as required.
\end{proof}

\section{\label{subsec:Tensors-and-Beyond}Tensors and Beyond}

Finally, we try to understand the expansion and mixing of systems
that are no longer linear. Our operator will now be a $k$-tensor
on $n$ dimensions, written as $T$, such that $T_{i_{1},i_{2},...,i_{k}}$
are some constant entries, and the dimension will be $n$ in each
index. We will throughout use the $3$-tensor for exposition. There
has been an upsurge of results related to tensors due to its application
in quantum information and machine learning, but relatively few
and very specific results are known with regards to random walks and
their convergence. Surprisingly, there are many versions of the Perron-Frobenius
theorem which can be instantiated and interpreted in multiple ways
for tensors (see \cite{chang2008perron,friedland2013perron}). In
fact, for the most natural definition of eigenvalue of a tensors (see
\cite{chang2013survey}), it is possible for a $k$-tensor to have
$n(k-1)^{n-1}$ eigenvalues, which is exponential when $k>2$. Our
focus in this section will be on some very specific things that we
think are very basic and important but unresolved.

The first question is to define a walk based on $T$. This is relatively
simple, and there are two equivalent and mathematically pleasing ways
to define the walk. Assume we have a distribution $p\in\mathbb{R}^{n}$,
and for simplicity assume $T$ is a nonnegative $3$-tensor. Then
we can define one step of the walk as follows:
\[
p_{t}(i)=\sum_{j=1,k=1}^{n}T_{i,j,k}p_{t-1}(j)\cdot p_{t-1}(k).
\]
Thus, our output distribution is quadratic in the input distribution.
Note that we are using the first index as the output, although any
index can be used, which is akin to taking the transpose as in the
matrix case. Another way to define the walk is to consider the last
two states, and evolve to the new state from the last two states.
We can thus write 
\[
p_{t}(i)=\sum_{j=1,k=1}^{n}T_{i,j,k}p_{t-1}(j)\cdot p_{t-2}(k).
\]
Both are interesting in their own right, and for the questions that
we are interested in, it will not matter which definition is used.
As the first step, to ensure that we are restricted to the simplex
of probability distributions, we require that $T$ is 1-line stochastic,
i.e. for every fixed $j,k$, 
\[
\sum_{i=1}^{n}T_{i,j,k}=1.
\]
As a consequence, note that if $p_{t-1}$ and $p_{t-2}$ are probability
distributions, so is $p_{t}$ since $p_{t}(i)\geq0$ since $T$ is
nonnegative, and 
\[
\sum_{i=1}^{n}p_{t}(i)=\sum_{j=1,k=1}^{n}p_{t-1}(j)\cdot p_{t-2}(k)\sum_{i=1}^{n}T_{i,j,k}=\sum_{j=1}^{n}p_{t-1}(j)\cdot\sum_{k=1}^{n}p_{t-2}(k)=1.
\]
Formally thus, we can define the tensor walk as follows.
\medskip

\begin{definition}
\label{def:tensor-walk}Let $T$ be a nonnegative $k$-tensor with
dimension $n$. Let $T$ be $1$-line stochastic, i.e., for every
fixed $j_{1},j_{2},...,j_{k-1}$, 
\[
\sum_{i=1}^{n}T_{i,j_{1},j_{2},...,j_{k-1}}=1.
\]
Then $p_{t}$ is a sequence of probability distributions due to the
tensor walk defined as follows. 
\[
p_{t}(i)=\sum_{j_{1},j_{2},...,j_{k-1}=1}^{n}T_{i,j_{1},j_{2},...,j_{k-1}}p_{t-1}(j_{1})p_{t-2}(j_{2})...p_{t-k+1}(j_{k-1})
\]
Succinctly, we write 
\[
p_{t}=T(p_{t-1},p_{t-2},...,p_{t-k}).
\]
\end{definition}

We remark that it is possible to consider any number of inputs and
outputs, where we have considered only 1 output and $k-1$ inputs
in Definition \ref{def:tensor-walk}. For evolution of probability
distributions, this is indeed the most meaningful way. We then define
the fixed point of the evolution $T$ as a probability distribution
that is preserved by $T$.
\medskip

\begin{definition}
\label{def:(Fixed-point-of-tensor-evolve}(Fixed point of tensor evolution)
Let $T$ be a $k$-tensor in $n$ dimensions that is 1 line stochastic,
with the evolution due to $T$ defined as in \ref{def:tensor-walk}.
We say that a distribution $p$ is a fixed point of $T$, if 
\begin{align*}
p & =T(p,p,...,p)\\
p(i) & =\sum_{j_{1},j_{2},...,j_{k-1}=1}^{n}T_{i,j_{1},j_{2},...,j_{k-1}}p(j_{1})p(j_{2})...p(j_{k-1}).
\end{align*}
\end{definition}

Given this definition, the first question is the following. 

\begin{center}
\emph{If $T$ is entry-wise positive, does $T$ have a unique fixed
point?}\\
\par\end{center}

Note that this is the first question to understand since positivity
implied a unique fixed point in the case of matrices due to the Perron-Frobenius
theorem. It intuitively does seem so, however, it turns out that $T$
could have multiple fixed points. The following is an example constructed
in \cite{changzhang2013counterexample}.

\begin{lemma}
\label{lem:pos-ten-counterex}\cite{changzhang2013counterexample}
Let $T$ be a positive $4$-tensor in 2 dimensions, defined as 
\[
T_{1,1,1,1}=0.872,\ T_{1,1,1,2}=2.416/3,\ T_{1,1,2,1}=2.416/3,\ T_{1,1,2,2}=0.616/3,
\]
\[
T_{1,2,1,1}=2.416/3,\ T_{1,2,1,2}=0.616/3,\ T_{1,2,2,1}=0.616/3,\ T_{1,2,2,2}=0.072,
\]
\[
T_{2,1,1,1}=0.128,\ T_{2,1,1,2}=0.584/3,\ T_{2,1,2,1}=0.584/3,\ T_{2,1,2,2}=2.384/3,
\]
\[
T_{2,2,1,1}=0.584/3,\ T_{2,2,1,2}=2.384/3,\ T_{2,2,2,1}=2.384/3,\ T_{2,2,2,2}=0.928,
\]
Note that $T$ is $1$-line stochastic. Then the following are the
2 distinct positive fixed points of $T$:
\[
p_{1}=0.2,\ p_{2}=0.8
\]
and 
\[
q_{1}=0.6,\ q_{2}=0.4.
\]
\end{lemma}

The proof can be checked by direct computation. In light of this,
it would seem that even the most basic requirement is not satisfied
by tensors. Almost every simple definition for matrices has multiple
avenues for explication when considered for tensors. In fact, even
the most intuitive way of defining the eigenvalues leads to $2^{n}$
eigenvalues even for a $3$-tensor. In spite of all this, we can actually
show the existence of a unique fixed point in the case of $2$-line
stochastic tensors.
\medskip

\begin{definition}
Let $T$ be a nonnegative $k$-tensor in $n$ dimensions. We say that
$T$ is $2$-line stochastic, if apart from the output index, there
is also an input index along which the tensor is stochastic. Formally,
we say that $T$ is 2-line stochastic, if there are some fixed input
index $j$ and the input index $i$, such that for every fixed $l_{1},...,l_{k-2}$,
we have that 
\[
\sum_{i=1}^{n}T_{i,j,l_{1},...,l_{k-2}}=1
\]
and
\[
\sum_{j=1}^{n}T_{i,j,l_{1},...,l_{k-2}}=1.
\]
\end{definition}

Given a $2$-line stochastic tensor, we can indeed show that the uniform
distribution is a unique fixed point! This is our main theorem in
this section.

\begin{theorem}
\label{thm:2-stoc-tens-uniq}Let $T$ be a positive $2$-line stochastic
tensor. Then for the tensor walk as defined in \ref{def:tensor-walk},
the uniform distribution is the unique fixed point.
\end{theorem}

\begin{proof}
We first note two crucial things before delving into the proof.

1. The fact that $T$ is line stochastic in the output index $i$
will ensure that the walk remains within the probability simplex.

2. The fact that $T$ is line stochastic in the input index $j$ will
ensure that there is a unique fixed point.

First we note that the uniform distribution is a fixed point for $T$,
by noting that 
\begin{align*}
p_{t}(i) & =\sum_{j,l_{1},...,l_{k-2}=1}^{n}T_{i,j,l_{1},...,l_{k-2}}\left(\frac{1}{n}\right)^{k-1}\\
 & =\left(\frac{1}{n}\right)^{k-1}\sum_{l_{1},...,l_{k-2}=1}^{n}\sum_{j=1}^{n}T_{i,j,l_{1},...,l_{k-2}}\\
 & =\left(\frac{1}{n}\right)^{k-1}n^{k-2}\\
 & =\frac{1}{n}.
\end{align*}
The next thing we show is that each $p_{t}$ for $t\geq k$ is a positive
vector. To see this, let $p_{1}$ to $p_{k-1}$ be the $k-1$ non-zero
initial distributions provided to us. Note $p_{i}$'s for $1\leq i\leq k-1$
could have zero entries, but do not have all zero entries. Let $x_{1},...,x_{k-1}\in[n]$
be the indices corresponding to the non-zero entries in the $p_{i}$'s
respectively. Thus we have 
\begin{align*}
p_{k}(i) & =\sum_{j_{1},j_{2},...,j_{k-1}=1}^{n}T_{i,j_{1},j_{2},...,j_{k-1}}p_{k-1}(j_{1})p_{k-2}(j_{2})...p_{1}(j_{k-1})\\
 & \geq T_{i,x_{1},x_{2},...,x_{k-1}}p_{k-1}(x_{1})p_{k-2}(x_{2})...p_{1}(x_{k-1})\\
 & >0
\end{align*}
since $T$ is positive, and further, the result holds for all $p_{t}$
inductively. Thus, it is not possible to have some starting probability
vectors $p_{1}=p_{2}=...=p_{k-1}=q$ that are not positive but the
fixed point, since $p=T(q,...,q)$ will be positive, and thus $p\not=q$.
This establishes that the fixed point must be positive. Now we show
that it must be unique.

Our next step is to define a doubly stochastic matrix $A_{t}$ using
$T$ and $p_{t}$ as follows.
\[
A_{t}(i,j)=\sum_{l_{1},...,l_{k-2}=1}^{n}T_{i,j,l_{1},...,l_{k-2}}p_{t-2}(l_{1})...p_{t-k+1}(l_{k-2}).
\]
We note that $A_{t}$ is doubly stochastic due to the 2-line stochastic
condition that can be checked by direct calculation. Further, note
that 
\[
p_{t+1}=A_{t}p_{t}.
\]
Now assume there is a positive fixed point $q$ for $T$ that is not
uniform. Let 
\[
A_{q}(i,j)=\sum_{l_{1},...,l_{k-2}=1}^{n}T_{i,j,l_{1},...,l_{k-2}}q_{l_{1}}...q_{l_{k-2}}.
\]
Note that $A_{q}$ is doubly stochastic. Further, since $q$ is positive,
$A_{q}$ is also positive and thus irreducible. Since $q$ is a fixed
point of $T$, we have that 
\[
q=A_{q}q.
\]
However, this is a contradiction since $A_{q}$ is irreducible and
doubly stochastic and only has the uniform distribution as stationary
distribution. Thus the tensor has a unique fixed point.
\end{proof}
The next interesting question is to understand how fast the tensor
converges to the fixed point, starting from any fixed distribution.
One way to achieve this is to use Mihail's proof showing that the
convergence of a nonnegative matrix will be inversely proportional
to the square of the expansion \ref{lem:mixt-phi-rel}. The nice thing
about Mihail's proof, unlike our simpler proof in Lemma \ref{lem:mixt-phi-rel}
is that it is combinatorial, and does not depend on the spectra of
the underlying object. However, it is not possible to extend it directly,
since every matrix $A_{t}$ in Lemma \ref{thm:2-stoc-tens-uniq},
although doubly stochastic, is a \emph{different} doubly stochastic
matrix, and we thus leave the problem of the speed of convergence
to the fixed point as the most important next step to Theorem \ref{thm:2-stoc-tens-uniq}
for future work.

\pagebreak{}

\chapter*{\label{sec:finth}Summary of Results and Questions}
\addcontentsline{toc}{chapter}{Summary of Results and Questions} 
We collect the main results of this thesis and some important questions
here. We restate some of our main theorems first.

\begin{theorem}
(Spectral Gap and Edge Expansion)\label{thm:intro-phi-delta-1} Let
$R$ be an $n\times n$ nonnegative matrix, with edge expansion $\phi(R)$
defined as \ref{def:gen-edge-expansion} and \ref{def:(gen-Spectral-gap-and-edgeexp},
and the spectral $\Delta(R)$ defined as \ref{def:(Spectral-gap-irred}
and \ref{def:(gen-Spectral-gap-and-edgeexp}. Then 
\[
\dfrac{1}{15}\cdot\dfrac{\Delta(R)}{n}\leq\phi(R)\leq\sqrt{2\cdot\Delta(R)}.
\]

\end{theorem}

\begin{theorem}

(Mixing Time, Singular Values, Edge Expansion, and Spectral Gap).
Let $R$ and $A$ be irreducible and $\frac{1}{2}$-lazy nonnegative
matrices with largest eigenvalue 1, where $A$ has the same corresponding
left and right eigenvector. Let the mixing time be as defined in \ref{def:gen-Mixing-time},
then 
\[
\tau_{\epsilon}(A)\ \leq\ \dfrac{\ln\left(\sqrt{\dfrac{n}{\kappa}}\cdot\dfrac{1}{\epsilon}\right)}{\ln\left(\frac{1}{\sigma_{2}(A)}\right)},
\]
\begin{align*}
\frac{\frac{1}{2}-\epsilon}{\phi(R)}\ \leq\ \tau_{\epsilon}(R)\  & \leq\ \dfrac{4\cdot\ln\left(\frac{n}{\kappa\cdot\epsilon}\right)}{\phi^{2}(R)},
\end{align*}
\[
(1-\Delta)\dfrac{1-\epsilon}{\Delta}\ \leq\ \tau_{\epsilon}(R)\ \leq\ 20\cdot\frac{n+\ln\left(\dfrac{1}{\kappa\cdot\epsilon}\right)}{\Delta(R)}.
\]

\end{theorem}

\begin{theorem}

(Rootn Matrices \ref{constr:Rootn-matrices-=002013}) There is a doubly
stochastic matrix $A_{n}$ for every $n$, such that 
\[
\phi(A_{n})\leq\dfrac{\Delta(A_{n})}{\sqrt{n}}.
\]

\end{theorem}

\begin{theorem}

(Chet Matrices \ref{constr-chet}) There is a matrix
$C_{n}$ for every $n$, such that $\sum_{i}C_{i,j}=1$ and $\sum_{j}C_{i,j}=1$
for all $i,j$, and with $\phi(C_{n})$ defined similar to that for
doubly stochastic matrices in \ref{def:(Edge-expansion-doubstoc},
it holds that 
\[
\phi(C_{n})\leq2\cdot\dfrac{\Delta(C_{n})}{n}.
\]

\end{theorem}

\begin{theorem}

(Monotonicity of Capacity) Let $A$ be an irreducible nonnegative
matrix with largest eigenvalue 1 and corresponding left and right
eigenvector $w$. Let $A_{\alpha}=\tilde{A}+\alpha\overline{A}$ where
$\tilde{A}=\frac{1}{2}(A+A^{T})$, $\overline{A}=\frac{1}{2}(A-A^{T})$,
and $-1\leq\alpha\leq1$. Then for $|\alpha|\leq|\beta|$, and every
$U\subseteq[n]$ and $a\in\mathbb{R}^{|U|}$, we have that 
\[
\text{cap}_{U,a}(A_{\alpha})\leq\text{cap}_{U,a}(A_{\beta}).
\]

\end{theorem}

\begin{theorem}

(Unique fixed point for tensors) Let $T$ be a $k$-tensor in $n$
dimensions which is entry-wise positive. Let $T$ be $2$-line stochastic,
such that for the output index $i$ and some input index $j$, it
holds that for any $i,j,l_{1},...,l_{k-2}$
\[
\sum_{i}T_{i,j,l_{1},...,l_{k-2}}=1
\]
and
\[
\sum_{j}T_{i,j,l_{1},...,l_{k-2}}=1.
\]
Let $T$ act on every probability distribution as stated in \ref{def:tensor-walk}.
Then the uniform distribution is the unique fixed point for $T$.
\end{theorem}

We collect some of our main conjectures and open problems here.

\begin{conjecture}

(Chet Conjecture) \label{conj:(Chet-Conjecture)-1}Let $C_{n}$ denote
the $n\times n$ Chet matrix defined in Construction \ref{constr-chet},
and let $C=\{n:C_{n}\text{ is entry-wise nonnegative}\}$. Then the
following is true:
\[
|C|=\infty.
\]
More strongly,
\[
C=\mathbb{N}.
\]

\end{conjecture}

\begin{conjecture}

\noindent (Trace Conjectures)\label{conj:(Trace-Conjecture)-1} Let
$A$ be a nonnegative matrix that is substochastic, that is, $\sum_{i}A_{i,j}\leq1$
and $\sum_{j}A_{i,j}\leq1$ for all $j$ and $i$. Assume the following:
Above the diagonal, $A$ has nonzero entries only for entries that
are at a distance of 1 from the diagonal, and below the diagonal,
$A$ has nonzero entries only for entries that are at a distance at
most $k$ from the diagonal, where the diagonal has distance zero
from the diagonal. Assume $\text{Tr}A^{l}\leq1$ for $l\leq k+1$.
Then for all $l$, 
\[
\text{Tr}A^{l}\leq1.
\]

The first relaxation is to show it for the Toeplitz case where every
diagonal has the same entry, and the second relaxation is to show
it for the infinite case where the condition $\text{Tr}A^{l}\leq1$
is replaced by 
\[
p_{l}\leq \alpha\cdot\frac{1}{k}
\]
for some constant $\alpha$, where $p_{l}$ is the probability of returning
back to the starting vertex after exactly $l$ steps.

\end{conjecture}

It is indeed possible to construct a large number of open problems from every chapter, for instance: Is it possible to use the matrix $H$ in Lemma \ref{lem:(Dirichlet-Lemma-for-nonsym}
to obtain a more exact bound between capacities than Theorem \ref{thm:(Monotonicity-of-Capacity)},
similar to the manner in which Lemma \ref{lem:(Strengthening-of-H-bound}
is an exact version of Lemma \ref{lem:(Dirichlet-Lemma-for-nonsym}? Is there a modified notion of normalized capacity,
such that it is equivalent to the spectral gap of a general nonnegative
$R$ up to constants, similar to Lemma \ref{lem:(Normalized-Capacity-Spec-Gap}? For the higher order notion of expansion as defined
in \cite{lee2014multiway}, what is its relation with the spectral
gap of a general nonnegative $R$? And so on. However, the most interesting and important open problem from our perspective is the following, which we believe would only be the second step into  a large and undiscovered area.  

\begin{openproblem}Is there a notion of expansion $\phi(T)$ for 2-line
stochastic tensors $T$, such that the tensor walk converges to the
uniform distribution in about $1/\phi^{2}(T)$ steps similar to lemma
\ref{lem:mixt-phi-rel}?
\end{openproblem}

We hope that these theorems, conjectures and questions get used, resolved or explored in
the future.

\printbibliography[heading=bibintoc]

\appendix

\chapter{Chet Matrices -- Exact and Numerical Computations}

We present some Chet Matrices computed exactly and numerically. The
exact computations were done in Maple, and the numerical ones in Matlab.

\section{\label{subsec:Chet-Matrix-for-16}Chet Matrix for $n=16$ computed
exactly}

The equations for $C_{16}$ where the matrices are defined in construction
\ref{constr-chet}, are as follows.
\[ r =\ \left(\dfrac{1}{16}\right)^{\frac{1}{15}} \]

$14\cdot c_{0} =\ 
2\cdot r-1$ \\

$13\cdot r\cdot c_{1} =\ 
\dfrac{8}{7}\cdot r^{2}
-\dfrac{15}{28}$ \\

$12\cdot r^{2}\cdot c_{2} =\ 
\dfrac{528}{637}\cdot r^{3}
+\dfrac{45}{637}\cdot r
-\dfrac{20}{49}$ \\

$11\cdot r^{3}\cdot c_{3} =\ 
\dfrac{38980}{57967}\cdot r^{4}
+\dfrac{7095}{231868}\cdot r^{2}
+\dfrac{110}{1029}\cdot r
-\dfrac{334605}{927472}$ \\

$10\cdot r^{4}\cdot c_{4} =\ 
\dfrac{2622960}{4463459}\cdot r^{5} 
+ \dfrac{15675}{811538}\cdot r^{3} 
+ \dfrac{3650}{93639}\cdot r^{2} 
+ \dfrac{5019075}{35707672}\cdot r 
- \dfrac{21653}{62426} $ \\
 
$9\cdot r^{5}\cdot c_{5} =\ 
\dfrac{221814104}{406174769}\cdot r^{6}
+ \dfrac{23279535}{1624699076}\cdot r^{4}
+ \dfrac{4980}{218491}\cdot r^{3}
+ \dfrac{135515025}{3249398152}\cdot r^{2} $ \\

$\ \ \ \ \ \ \ \ \ \ \ \ \ \ \ \ \ \ \ \ 
+ \dfrac{194877}{1092455}\cdot r
- \dfrac{4578373505}{12997592608}$ \\

$8\cdot r^{6}\cdot c_{6} =\ 
\dfrac{96435386896}{179123073129}\cdot r^{7}
+\dfrac{33490290}{2843223383}\cdot r^{5}
+\dfrac{10590760}{656128473}\cdot r^{4}
+\dfrac{63909555}{2843223383}\cdot r^{3}$ \\

$\ \ \ \ \ \ \ \ \ \ \ \ \ \ \ \ \ \ \ \
+ \dfrac{4070764}{99413405}\cdot r^{2}
+ \dfrac{22891867525}{102356041788}\cdot r
-\dfrac{569210253}{1530966437}$ \\
 
$7\cdot r^{7}\cdot c_{7} =\ 
\dfrac{14468197830722}{25614599457447}\cdot r^{8}
+\dfrac{4659561625}{443542847748}\cdot r^{6}
+\dfrac{8450450}{656128473}\cdot r^{5}$ \\

$\ \ \ \ \ \ \ \ \ \ \ \ \ \ \ \ \ \ \ \
+\dfrac{98852689755}{6505295100304}\cdot r^{4}
+\dfrac{2013729}{99413405}\cdot r^{3}
+\dfrac{389161747925}{10645028345952}\cdot r^{2}$ \\

$\ \ \ \ \ \ \ \ \ \ \ \ \ \ \ \ \ \ \ \
+ \dfrac{1707630759}{6123865748}\cdot r
- \dfrac{42029290110657}{104084721604864}$ \\
 
$6\cdot r^{8}\cdot c_{8} =\ 
+\dfrac{12646467661976}{19922466244681}\cdot r^{9}
+\dfrac{36526797305}{3622266589942}\cdot r^{7}
+\dfrac{2018448010}{179123073129}\cdot r^{6}$ \\
   
$\ \ \ \ \ \ \ \ \ \ \ \ \ \ \ \ \ \ \ \
+\dfrac{268437195855}{22768532851064}\cdot r^{5}
+\dfrac{74507973}{5686446766}\cdot r^{4}
+\dfrac{68675602575}{4139733245648}\cdot r^{3}$ \\
 
$\ \ \ \ \ \ \ \ \ \ \ \ \ \ \ \ \ \ \ \
+ \dfrac{15368676831}{557271783068}\cdot r^{2}
 + \dfrac{126087870331971}{364296525617024}\cdot r 
 -\dfrac{53944771607543}{120370705142688}$ \\
 
$5\cdot r^{9}\cdot c_{9}=
 \dfrac{4192421772289592}{5438833284797913}\cdot r^{10}
  +\dfrac{225785202013625}{21755333139191652}\cdot r^{8}
  +\dfrac{280853043700}{26331091749963}\cdot r^{7}$\\
  
  $\ \ \ \ \ \ \ \ \ \ \ \ \ \ \ \ \ \ \ \
+\dfrac{10494908689675}{1035968244723412}\cdot r^{6}
   +\dfrac{197518666}{19902563681}\cdot r^{5}
   +\dfrac{42922251609375}{4143872978893648}\cdot r^{4}$\\

 $\ \ \ \ \ \ \ \ \ \ \ \ \ \ \ \ \ \ \ \
+\dfrac{11384205060}{975225620369}\cdot r^{3}
 +\dfrac{210146450553285}{16575491915574592}\cdot r^{2}
 + \dfrac{269723858037715}{631946201999112}\cdot r$\\

 $\ \ \ \ \ \ \ \ \ \ \ \ \ \ \ \ \ \ \ \
-\dfrac{293344384450753313}{580142217045110720}$\\

$4\cdot r^{10}\cdot c_{10}= + 
 \dfrac{583988144046138064}{571077494903780865}\cdot r^{11}
 +\dfrac{433842562719595}{38071832993585391}\cdot r^{9}$ \\
 
 $\ \ \ \ \ \ \ \ \ \ \ \ \ \ \ \ \ \ \ \
 +\dfrac{284474332622540}{26357422841712963}\cdot r^{8}
+\dfrac{239093054413515}{25381221995723594}\cdot r^{7}$ \\

$\ \ \ \ \ \ \ \ \ \ \ \ \ \ \ \ \ \ \ \
  +\dfrac{136555337396}{16300199654739}\cdot r^{6}
 +\dfrac{504083501274005}{65265999417574956}\cdot r^{5}$ \\
 
$ \ \ \ \ \ \ \ \ \ \ \ \ \ \ \ \ \ \ \ \
+\dfrac{7074714234537}{976200845989369}\cdot r^{4}
 +\dfrac{42029290110657}{7251777713063884}\cdot r^{3}$ \\
 
 $\ \ \ \ \ \ \ \ \ \ \ \ \ \ \ \ \ \ \ \
 -\dfrac{269723858037715}{28753552190959596}\cdot r^{2}  
+ \dfrac{2640099460056779817}{5076244399144718800}\cdot r$ \\

$ \ \ \ \ \ \ \ \ \ \ \ \ \ \ \ \ \ \ \ \
 -\dfrac{487051361453894455}{843437530934814816} $ \\

$3\cdot r^{11}\cdot c_{11}= 
  \dfrac{2582373218800886646164}{1714945717196053937595}\cdot r^{12}
 +\dfrac{60910874936948457}{4619382403221694108}\cdot r^{10}$
 
 $\ \ \ \ \ \ \ \ \ \ \ \ \ \ \ \ \ \ \ \
 +\dfrac{704066732737790}{61500653297330247}\cdot r^{9}
+\dfrac{1873806347700790155}{203252825741754540752}\cdot r^{8}$\\

  $\ \ \ \ \ \ \ \ \ \ \ \ \ \ \ \ \ \ \ \
  +\dfrac{10063658054264}{1331182971803685}\cdot r^{7}
   +\dfrac{454789092312238525}{71270471363991851952}\cdot r^{6}$\\
   
    $\ \ \ \ \ \ \ \ \ \ \ \ \ \ \ \ \ \ \ \
+\dfrac{74092391002251}{13666811843851166}\cdot r^{5}
 +\dfrac{948306872766753891}{232288943704862332288}\cdot r^{4}$\\
 
 $\ \ \ \ \ \ \ \ \ \ \ \ \ \ \ \ \ \ \ \
 -\dfrac{53944771607543}{134183243557811448}\cdot r^{3}
-\dfrac{71282685421533055059}{1847752961288677643200}\cdot r^{2}$\\

$\ \ \ \ \ \ \ \ \ \ \ \ \ \ \ \ \ \ \ \
+ \dfrac{2435256807269472275}{3936041811029135808}\cdot r
  -\dfrac{250148541000967840970441}{376308088801876978306560}$ \\

$2\cdot r^{12} \cdot c_{12}=  
 \dfrac{8292429130967703269984}{3273987278283375699045}\cdot r^{13}
 +\dfrac{251191067356462921}{16167838411275929378}\cdot r^{11}$ \\

 $\ \ \ \ \ \ \ \ \ \ \ \ \ \ \ \ \ \ \ \
+\dfrac{204277892726612218}{16789678350171157431}\cdot r^{10}
 +\dfrac{52904039025561315}{5879213967736701592}\cdot r^{9}$ \\
 
 $\ \ \ \ \ \ \ \ \ \ \ \ \ \ \ \ \ \ \ \
+\dfrac{54680493901861873}{7995084928652932110}\cdot r^{8}
 +\dfrac{3520738115297033525}{654797455656675139809}\cdot r^{7}$ \\
 
 $\ \ \ \ \ \ \ \ \ \ \ \ \ \ \ \ \ \ \ \
+\dfrac{5335710693495901}{1243679877790456106}\cdot r^{6}
 +\dfrac{237927811316429277}{73910118451547105728}\cdot r^{5}$ \\
 
 $\ \ \ \ \ \ \ \ \ \ \ \ \ \ \ \ \ \ \ \
+\dfrac{1059097700970891719}{805904560808215556688}\cdot r^{4}
 -\dfrac{36081359287442657499}{6467135364510371751200}\cdot r^{3}$ \\
 
 $\ \ \ \ \ \ \ \ \ \ \ \ \ \ \ \ \ \ \ \
-\dfrac{75492961025353640525}{1074539414410954075584}\cdot r^{2}
 + \dfrac{250148541000967840970441}{359203175674518933838080}\cdot r$ \\
 
 $\ \ \ \ \ \ \ \ \ \ \ \ \ \ \ \ \ \ \ \
-\dfrac{611028831268072088101}{795955121785891907840}$ \\

$1\cdot r^{13} \cdot c_{13}= 
\dfrac{2167687960275208704611884}{417105979253302064058333}\cdot r^{14}
 +\dfrac{14252777333017827436531}{873936337483109086598412}\cdot r^{12}$ \\
 
 $\ \ \ \ \ \ \ \ \ \ \ \ \ \ \ \ \ \ \ \
+\dfrac{4397770296526056}{395716324751508761}\cdot r^{11}
 +\dfrac{43699381559221754727}{5885093181704438293592}\cdot r^{10}$ \\

  $\ \ \ \ \ \ \ \ \ \ \ \ \ \ \ \ \ \ \ \
+\dfrac{145628182696554137}{27982797250285262385}\cdot r^{9}
  +\dfrac{80226936993830858719175}{20974472099594618078361888}\cdot r^{8}$ \\

$\ \ \ \ \ \ \ \ \ \ \ \ \ \ \ \ \ \ \ \
+\dfrac{16063263041956539}{5540028546521122654}\cdot r^{7}
 +\dfrac{43750487609031995683}{20177462337272359863744}\cdot r^{6}$ \\
 
$\ \ \ \ \ \ \ \ \ \ \ \ \ \ \ \ \ \ \ \
+\dfrac{1255133000992702981}{940221987609584816136}\cdot r^{5}
 -\dfrac{984757098601178871741}{2589440999949952849180480}\cdot r^{4}$ \\
 
 $\ \ \ \ \ \ \ \ \ \ \ \ \ \ \ \ \ \ \ \
-\dfrac{2435256807269472275}{341898904585303569504}\cdot r^{3}
 -\dfrac{250148541000967840970441}{2971589907852838452660480}\cdot r^{2}$ \\
 
$\ \ \ \ \ \ \ \ \ \ \ \ \ \ \ \ \ \ \ \
+ \dfrac{1833086493804216264303}{2785842926250621677440}\cdot r
  -\dfrac{41502574124308301104426128989}{46982817503091944495530629120}$ \\

\pagebreak{}

\section{\label{chet-pseudocode} Pseudocode to compute Chet matrix for input integer $n$}

\begin{lstlisting}
chet(n):
    C = Matrix of zeroes indexed from 1 to n
    r = (1/n)^(1/(n-1))
    
    for i = 1 to n-1:
        C[i+1,i] = r
    
    bval = 1-r
    C[1,1] = bval
    C[n,n] = bval
    
    for k = 1 to n-1:
        cval = ( 1 - trace(C^k)) ) / ( k*(n-k-1)*r^(k-1) )
        bval = bval - cval
        
        for i=2 to n-k:
            C[i, i+k-1] = cval
        C[1,1+k] = bval
        C[n-k, n] = bval
    
    C[1,n] = 1
    for k = 1 to n-1:
        C[1,n] = C[1,n] - C[1,k]
    
\end{lstlisting}

\pagebreak{}

\section{\label{subsec:Chet-Matrix-for-500} Chet Matrix for $n=500$ computed
numerically}

We plot the values of $c_i$ and $b_i$ for Chet Matrices for $n=500$. The first plot contains the values of $c_i$, the second plot contains the values of $c_i$ where the y-axis has log-scale, and the third plot contains the values of $b_i$. Notice that in the third plot, the last entry is the value of $b_{499}$, which does not necessarily have to decrease with other $b_i$, since it is just a linear function of other $b_i$'s from $i=0$ to $498$ (see Construction \ref{constr-chet} for details). The same is observed for all matrices of different sizes.

\includegraphics[scale=0.8]{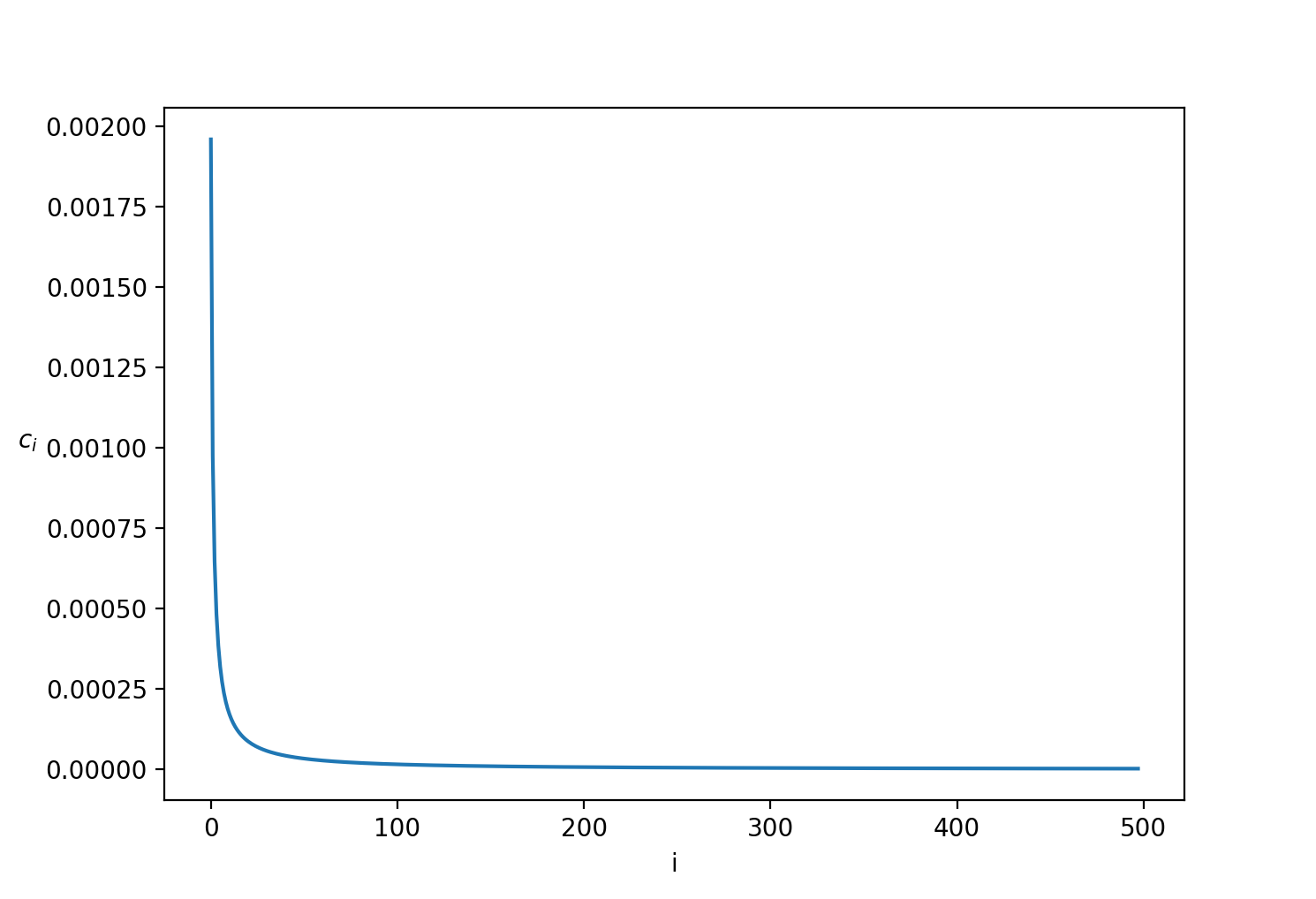}

\includegraphics[scale=0.85]{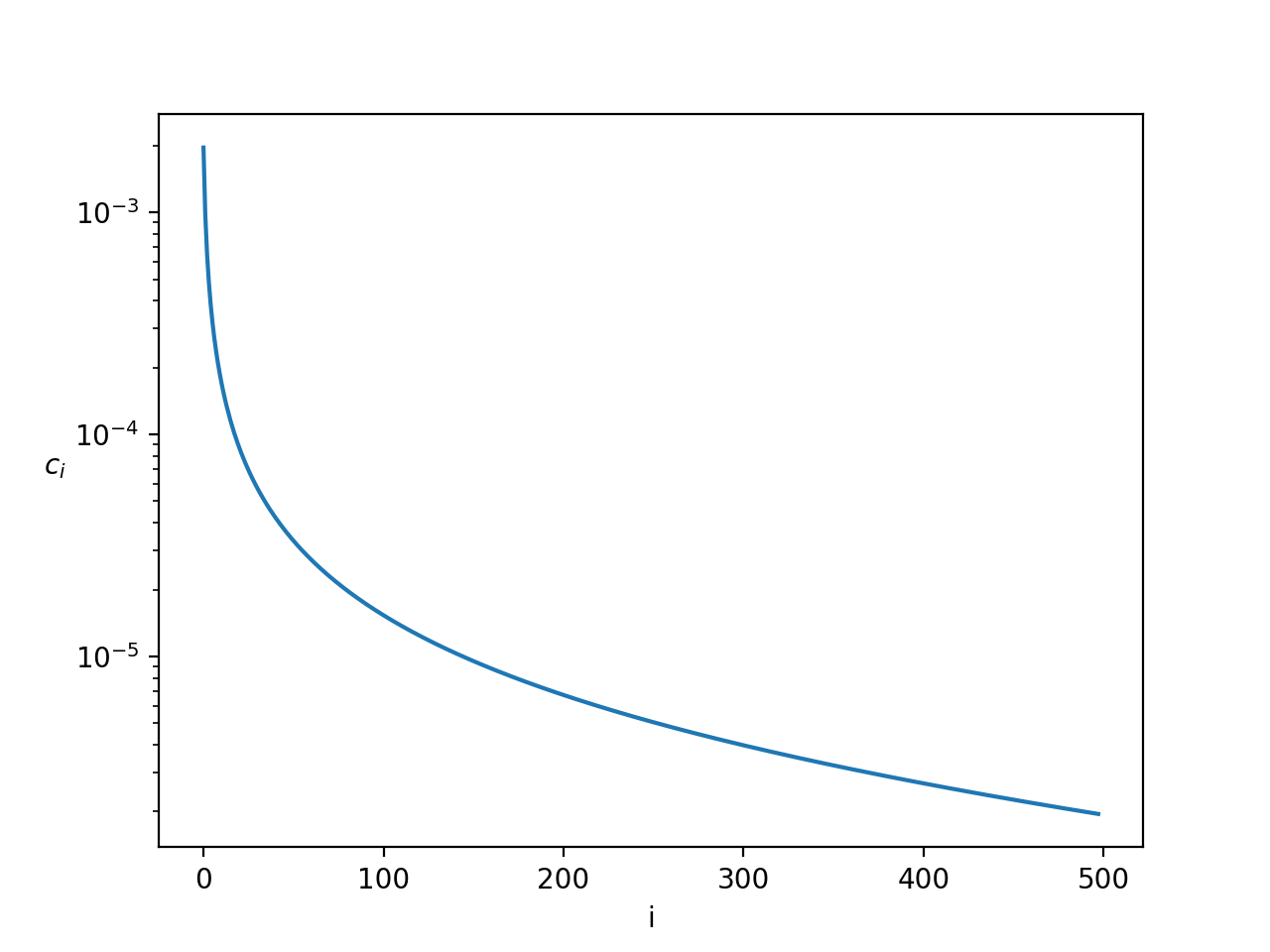}

\includegraphics[scale=0.85]{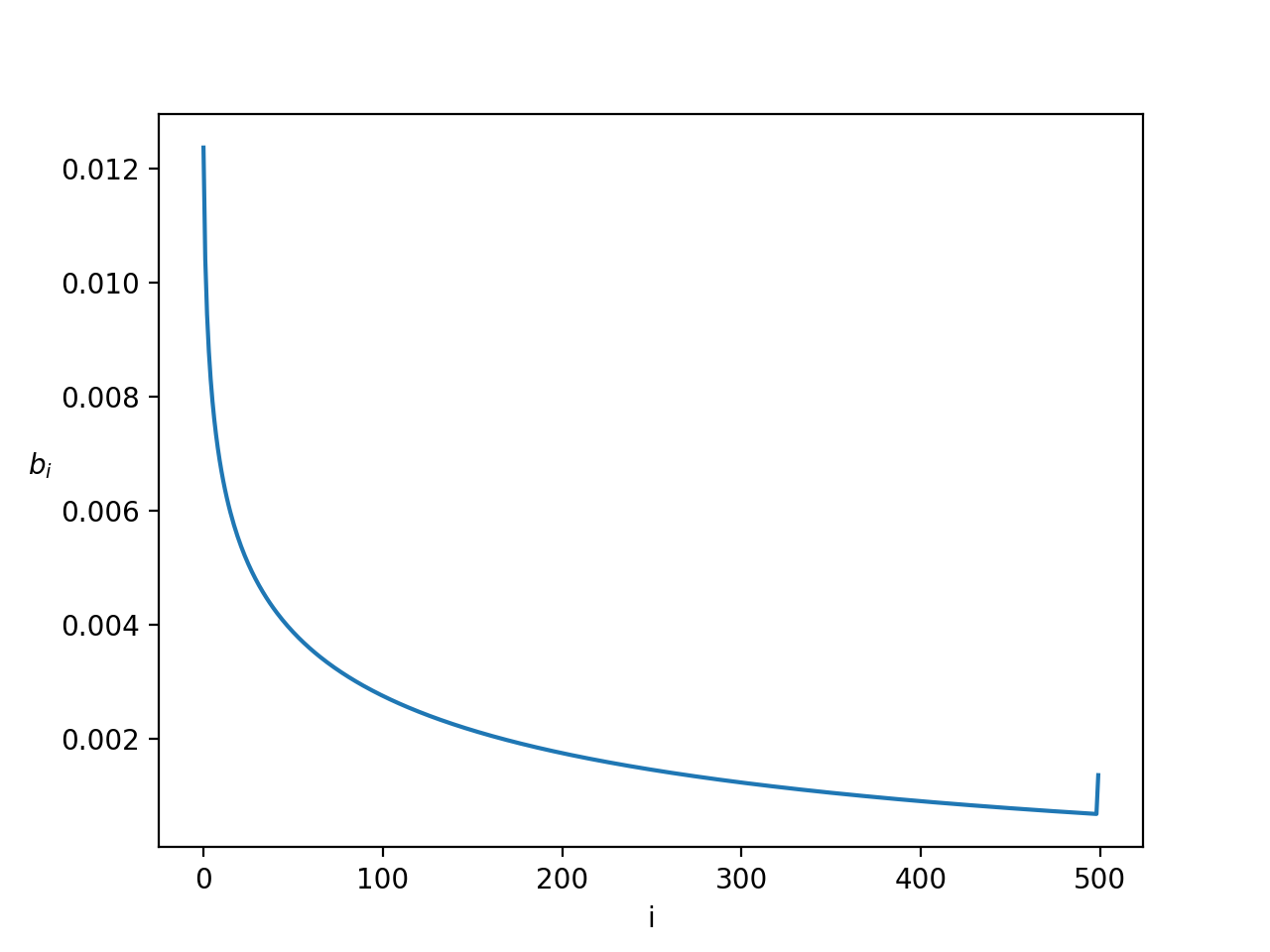}

\end{document}

%% file: thmmacros.tex
\usepackage{amsthm}
\usepackage{thmtools,thm-restate}

\numberwithin{equation}{chapter}
\declaretheoremstyle[qed=\qedsymbol]{noproofstyle}

\declaretheorem[name=Observation,numbered=no]{observation*}

\declaretheorem[numberlike=equation]{fact}

\declaretheorem[numberlike=equation]{theorem}

\declaretheorem[name=Theorem,numbered=no]{theorem*}

\declaretheorem[numberlike=equation]{lemma}
\declaretheorem[name=Lemma,numbered=no]{lemma*}

\declaretheorem[numberlike=equation]{corollary}
\declaretheorem[name=Corollary,numbered=no]{corollary*}

\declaretheorem[name=Proposition,numbered=no]{proposition*}

\declaretheorem[numberlike=equation]{claim}
\declaretheorem[name=Claim,numbered=no]{claim*}

\declaretheorem[numberlike=equation]{conjecture}
\declaretheorem[name=Conjecture,numbered=no]{conjecture*}

\declaretheorem[name=Question,numbered=no]{question*}

\declaretheorem[numberlike=equation, name=Open Problem]{openproblem}
\declaretheorem[name=Open Problem,numbered=no]{openproblem*}

\declaretheoremstyle[qed=$\lozenge$]{defstyle} 

\declaretheorem[numberlike=equation,style=defstyle]{definition}
\declaretheorem[unnumbered,name=Definition,style=defstyle]{definition*}

\declaretheorem[numberlike=equation,style=defstyle]{example}
\declaretheorem[unnumbered,name=Example,style=defstyle]{example*}

\declaretheorem[unnumbered,name=Notation=defstyle]{notation*}

\declaretheorem[numberlike=equation,style=defstyle]{construction}
\declaretheorem[unnumbered,name=Construction,style=defstyle]{construction*}

\declaretheorem[numberlike=equation,style=defstyle]{remark}
\declaretheorem[unnumbered,name=Remark,style=defstyle]{remark*}
